# The Theory of Khinchin families

Víctor J. Maciá

Departamento de Matemáticas
Universidad Autónoma de Madrid

Universidad Autónoma
de Madrid

A thesis submitted for the degree of

*Doctor in Mathematics*

Supervised by: Daniel Faraco and José L. Fernández

Madrid, 2024

*To my grandparents.*

# Contents

# Acknowledgments

First and foremost, I would like to express my deepest gratitude to my advisors, Josechu and Daniel, for their unwavering support, guidance, and encouragement throughout my doctoral studies. Working with them has been a truly rewarding experience.

A heartfelt thank you to my family for their unconditional love, patience, and encouragement. In particular I would like to thank Bea for his constant support.

To all who have contributed to this thesis, directly or indirectly, I extend my sincere appreciation. In particular I would like to thank Kari Astala and Steffen Rohde for hosting me during my visits to Helsinki and Seattle; they have taught me many things through our interactions.

# Nomenclature

| | |
|---|---|
| OGF | Ordinary generating function |
| EGF | Exponential generating function |
| $\mathbf{P}$ | Probability (referred to certain probability space) |
| $\mathbf{E}$ | Expectation (referred to certain probability space) |
| $\mathbf{V}$ | Variance (referred to certain probability space) |
| $\oplus$ | Sum of independent random variables |
| $\mathbb{D}$ | Unit disk in the complex plane $\mathbb{C}$ |
| $\mathbb{D}(z, R)$ | Disk of center $z \in \mathbb{C}$ and radius $R \geq 0$ in the complex plane $\mathbb{C}$ |
| $\stackrel{d}{=}$ | Equality in distribution of two random variables |
| $x^{\underline{k}}$ | $k$-th falling factorial of $x$ |

# Notations

- **Asymptotically equivalent:** for two functions $f, g$ we say that $f$ is asymptotically equivalent to $g$, as $t \uparrow R$, denoted $f \sim g$, as $t \uparrow R$, if
$$\lim_{t \uparrow R} \frac{f(t)}{g(t)} = 1.$$
- **Comparable:** for two functions $f, g$ we say $f$ is comparable with $g$, as $t \uparrow R$, denoted $f(t) \asymp g(t)$, as $t \uparrow R$, if there are constants $C, c > 0$ such that
$$c \leq \frac{f(t)}{g(t)} \leq C, \qquad \text{as } t \uparrow R.$$
- **Big O and Little o notation:** For two functions $f, g$ we say that $f$ is a big O of $g$, as $t \uparrow R$, denoted $f = O(g)$, as $t \uparrow R$, if there exists a constant $C > 0$ such that
$$|f(t)| \leq C|g(t)|, \qquad \text{as } t \uparrow R.$$
For the same two functions we say that $f$ is a little o of $g$, as $t \uparrow R$, denoted $f = o(g)$, as $t \uparrow R$, if
$$\lim_{t \uparrow R} \frac{f(t)}{g(t)} = 0.$$
- **Convergence in distribution of a sequence of random variables:** for a sequence of random variables $Z_n$ and a random variable $Z$, the notation $Z_n \overset{d}{\to} Z$, as $n \to \infty$, means that $Z_n$ converge in distribution towards $Z$, as $n \to \infty$.

# Introduction

The theory of Khinchin families connects Probability Theory and Complex Analysis. Along this PhD thesis, we exploit this connection to obtain, using a variety of local central limit theorems, asymptotic formulas for the coefficients of analytic functions with non-negative coefficients. We give criteria for a power series to have certain specific, and in most cases, computable, asymptotic formula.

The theory of Khinchin families was initiated and is strongly influenced by the work of Walter K. Hayman, see [48], Paul C. Rosembloom, see [83] and Luis Báez-Duarte, see [6]. We can assign to each power series $f$ in the class $\mathcal{K}$, having radius of convergence $R > 0$, a family of random variables $(X_t)_{[0,R)}$. This theory examines the behavior of this family of random variables, and also of its normalized version, as $t \uparrow R$.

It turns out that, under certain conditions, the behavior of this family of random variables encodes the asymptotic behavior of the coefficients of the power series $f$. In this thesis we develop a complete theory of Khinchin families, including operations with power series and criteria to be in some of the relevant classes of functions; we also study and organize the different connections between these families of random variables and the function theoretical properties of the functions in the class $\mathcal{K}$.

We thought it could be of interests to classify functions having certain probabilistic properties on the Khinchin family side; some of these classes of functions give sufficient conditions to have an specifics asymptotic formula for the coefficients of $f$.

A relevant result in this theory is the following:

$$f \text{ is in the Hayman class} \Rightarrow f \text{ is Strongly Gaussian } \Rightarrow f \text{ is Gaussian}$$

In this line we also have that strongly Gaussian power series conform clans. This chain of implications tells us that a characterization, or simple criteria, for being on any of these classes of functions would be really useful.

These classes of functions allow us to classify power series with non-negative coefficients and determine whether this theory applies to some of them. This question is complex, but at the same time relevant: for instance, some of the generating functions of the partitions of integers are in these classes of functions, and this allows us to find asymptotic formulas for their coefficients.

Motivated by this classifying objective we obtain, in some cases, simple criteria for a function to be in some of the classes of functions mentioned above; for instance, for one of the most relevant classes in here, the Hayman class.

This thesis aims to harmonize the various asymptotic formulas for combinatorial or probabilistic objects within a single framework. With this aim we have extended some of the classes of functions and also studied the coefficients of large powers of analytic functions. One of our main goals was to consolidate and develop some aspects of a comprehensive theory: providing a guidebook for combinatorialists or probabilists, a series of results, or basic criteria, where they can follow a step-by-step process, verify if certain straightforward conditions are met, and then derive an asymptotic formula for the object of study. However, this is not the sole primary goal; we are also profoundly interested in the functional theoretical properties of these classes of functions and their connections to the respective Khinchin families (such as the case of $f \in \mathcal{K}$ being an entire function).

The organization of this thesis is as follows:

In Chapter 1 we study the elementary properties of Khinchin families. We introduce the mean and the variance of a Khinchin family, these functions play a very relevant role in this theory and studying their growth and range is very important; for instance the mean and variance functions appear in the characteristic function of a normalized Khinchin family and also in Hayman's formula for the coefficients of a power series. Later on we introduce three critical components of this theory: the characteristic function of a Khinchin family, Hayman's formula for the coefficients, and the fulcrum of a power series in the class $\mathcal{K}$, these three objects are ubiquitous in this thesis. We conclude this first chapter by studying the moment generating function of a Khinchin family, including Chernoff bounds, and connecting the Wiman-Valiron theory with the theory of Khinchin families.

The main contributions in Chapter 1 are:

- Subsection 1.2.3 where we study the range of the variance of certain canonical products in the class $\mathcal{Q}$. This gives us a platform to illustrate the behavior of the variance.
- Lemma 1.2.1 where we characterize the range of the mean of a Khinchin family.
- A new proof of Theorem 1.2.3 which allows us to extend this result to the absolute central moments of the respective Khinchin family. See also Subsection 1.2.5 where we study the range of the variance.
- Lemma 1.3.2 which implies that the normalized Khinchin families are continuous in distribution, see also Lemma 1.3.4.
- Proposition 1.3.5 gives a zero-free region for the characteristic function of a power series $f \in \mathcal{K}$ (and also for the power series $f$).
- Lemma 1.3.9 expresses the moments of the normalized Khinchin family in terms of certain quotients of derivatives of the fulcrum (and vice versa).
- We also organize the different operations with power series and how this translates into operations in the Khinchin family side, see Subsection 1.3.6.
- In Subsection 1.4.2 we relate the moment generating function of a Khinchin family with its fulcrum. We also give Chernoff bounds for Khinchin families, see Proposition 1.4.6.

- In Proposition 1.5.7, using ideas from the Wyman-Valiron Theory, we obtain some differential inequalities for the fulcrum of an entire function $f \in \mathcal{K}$. These inequalities will be relevant in Chapter 2.

In Chapter 2 we study the moments of a Khinchin family and their asymptotic properties. Later on we introduce the concept of clan: this concept comes from [48], in that case refereed to Hayman functions. We will see later on, in Chapter 3, that strongly Gaussian power series are clans. In this chapter we study the function theoretical properties of entire functions in relation with their Khinchin families, for instance we give a formula for the order of an entire function in terms of the moments of a Khinchin family. We also give counterexamples which prove that not every entire function is a clan, though they are always weak clans. We conclude this chapter by proving that entire functions of finite order having one Borel exceptional value are clans.

The main contributions in Chapter 2 are:

- Lemma 2.3.6. This lemma allows us to obtain asymptotic formulas, in terms of the expectation, for the non-integer moments of a Khinchin family.
- In Lemma 2.3.9 we prove that certain quotients of moments of Khinchin families are increasing with respect to the respective powers. See also Corollaries 2.3.10 and 2.3.11 for comparisons of moments.
- Lemma 2.4.3 gives an equivalent definition of clan in terms of the derivatives of $f \in \mathcal{K}$. We also prove in there that transcendental entire functions $f \in \mathcal{K}$ are weak clans. See also Subsection 2.4.6 and in particular Proposition 2.4.6.
- Theorem 2.4.4 which allows us to build counterexamples proving that not every transcendental entire function $f \in \mathcal{K}$ is a clan. In this line see also Lemma **??** and Corollary 2.4.5.
- In Theorem 2.4.8 we prove that for any $p > 0$ the $p$-th moments of a Khinchin clan have an asymptotic in terms of the mean function.
- In Theorem 2.4.11 we give a characterization of clans in terms of the mean function.
- Lemma 2.4.16 expresses certain quotients of a function $f \in \mathcal{K}$ in terms of a series of factorial moments divided by powers of the mean of the Khinchin family. Later on this allows to give a characterization of clans in terms of $f$ it self, see Theorem 2.4.17.
- In Corollary 2.4.18 we obtain a characterization of clans in terms of the moment generating function (in this case $\mathbf{E}(e^{(X_t - m_f(t))\lambda})$ and $\lambda = \nu(t) = \ln(1 + 1/m_f(t)))$.
- Theorem 2.5.1 gives the order of an entire function in terms of the moments of a Khinchin family.
- Proposition 2.5.3 gives an inequality for the order of an entire function in terms of the quotient $\sigma_f^2(t)/m_f(t)$.
- Proposition 2.5.7 gives that canonical products with negative zeros in $\mathcal{K}$ conform clans.
- Theorem 2.5.9 gives sufficient conditions for an entire function of finite order $\rho > 0$ to be a clan.
- Theorem 2.5.10 proves that entire functions with one Borel exceptional value are clans. In this line Proposition 2.5.11 gives that the product of a polynomial and the exponential of an entire function of finite order in $\mathcal{K}$ is always a clan; This proves that the exponential of an entire function of finite order is always a clan.

In Chapter 3 we introduce uniform conditions on the Khinchin family side that give a local central limit theorem for the Khinchin families and therefore an asymptotic formula for the coefficients of the power series. There are, at least, three relevant classes of functions in this theory: Gaussian power series, strongly Gaussian power series and the Hayman functions. We study the function theoretical, and also the probabilistic, properties of these families. In some cases we give criteria for a power series in $\mathcal{K}$ to be in some of these classes and prove that certain specific families of power series, which are relevant from a combinatorial point of view, are in these classes.

The main contributions in Chapter 3 are:

- For non-vanishing power series: we give criteria in terms of the fulcrum to be Gaussian; see Theorems 3.1.4 and 3.1.6. See also Remark 3.1.9 for a generalization of these results.
- Under certain hypothesis, we give necessary conditions, for a non-vanishing power series $f \in \mathcal{K}$, to be a Gaussian, see Corollary 3.1.8.
- We prove that the exponential of a polynomial $e^P \in \mathcal{K}$ is a Gaussian power series, see Proposition 3.1.14.
- We prove a moment criteria for a power series in $\mathcal{K}$ to be Gaussian. This criteria applies to general power series $f \in \mathcal{K}$, see Theorem 3.1.19. We prove by means of this criteria that canonical products with negative zeros $f \in \mathcal{Q}$ such that $\lim_{t\uparrow\infty} \sigma_f^2(t) = +\infty$ are Gaussian, see Theorem 3.1.21.
- We prove that certain operations with power series preserve being Gaussian. For instance we prove that under certain simple conditions the power series $h = P \circ f$, where $P \in \mathcal{K}$ is a polynomial and $f \in \mathcal{K}$ a Gaussian power series, is Gaussian, see Theorem 3.1.32.
- In Theorem 3.2.9 we prove that strongly Gaussian power series are Gaussian.
- In Corollary 3.2.10 we prove that strongly Gaussian power series are clans.
- In Corollary 3.2.11 we prove certain asymptotic properties of strongly Gaussian power series $f \in \mathcal{K}$. As a consequence of this result we obtain that for strongly Gaussian power series we always have $M_f = +\infty$.
- We prove that certain operations with power series preserve the property of being strongly Gaussian. For instance we prove that under certain simple conditions the composition of a polynomial $P \in \mathcal{K}$ with a strongly Gaussian power series $f \in \mathcal{K}$ gives a strongly Gaussian power series $h = P \circ f$, see Theorem 3.2.14.
- In Lemma 3.3.1 we prove that the cut of a Hayman function $f \in \mathcal{K}$ should have certain asymptotic restrictions (with respect to the variance $\sigma_f^2(t)$).
- In Theorem 3.3.4 we give criteria, in terms of the fulcrum, for the major arc condition to hold. This criteria applies to non-vanishing power series in $\mathcal{K}$.
- In Sections 3.4, 3.5 and 3.6 we extend the definitions of Gaussian, strongly Gaussian and Hayman functions to uniformly Gaussian, uniformly strongly Gaussian and uniformly Hayman. Some of these extensions will be used, for instance, in Chapters 4 and 6.

In Chapter 4 we study exponential power series, that is, power series of the form $f = e^g$, with $g \in \mathcal{P}$. As explained in Section 4.1 this exponential structure is very relevant in combinatorics; an instance of this structure is given by the ordinary generating functions of the different sort of partitions of integers. Later on we give criteria for an exponential power series $f = e^g$, with $g \in \mathcal{P}$, to be Gaussian and also criteria for this power series to be in the Hayman class. This last criteria for the Hayman class gives very simple conditions in terms of the coefficients of $g$ for $f = e^g$, with $g \in \mathcal{P}$, to be in the Hayman class, see Theorems 4.3.17 and 4.3.15. We conclude this chapter by proving, using Theorem 4.3.17, that most of the generating functions of the partitions of integers are in the Hayman class.

The main contributions in Chapter 4 are:

- We adapt the criteria from Chapter 3 for a power series to be Gaussian to the case $f = e^g$ with $g \in \mathcal{P}$. These criteria are much simpler now and are given in terms of $g$, see Theorems 4.3.1 and 4.3.2. See also Remark 4.3.4 for a generalization of these results.
- Under certain hypothesis we give necessary condition, in terms of $g \in \mathcal{P}$, for $f = e^g$ to be Gaussian, see Corollary 4.3.3. See also Remark 4.3.4 for a generalization of this result.
- In Theorem 4.3.5 we give a moment criteria for an exponential power series to be Gaussian. See also Remark 4.3.6, there we give a characterization of Gaussianity for power series $f = e^g$, with $g \in \mathcal{P}$ verifying certain hypothesis.
- In Theorem 4.3.12 we give a first criteria in terms of $g \in \mathcal{P}$ for a power series $f = e^g$ to be in the Hayman class. See also Theorem 4.3.13 for further criteria.
- In Theorems 4.3.17 and 4.3.15 we give very simple criteria in terms of the coefficients of $g \in \mathcal{P}$ for $f = e^g$ to be in the Hayman class.
- Using Theorems 4.3.17 and 4.3.15 we prove that the EGF of the Bell numbers is in the Hayman class, see Theorem 4.3.18. We also prove by means of Theorem 4.3.17 that the OGF of the partitions of integers $P$, see Theorem 4.3.19, of the plane partitions $M$, see Theorem 4.3.21, and the OGF of the partitions in arithmetic progression $P_{a,b}$, see Theorem 5.1.11, are in the Hayman class.
- In Theorem 4.3.22 we give criteria for an exponential $f = e^g$, with $g \in \mathcal{K}$ to be uniformly Hayman. We also prove, using this criteria, that the exponential of a polynomial in $\mathcal{K}$ is uniformly Hayman.

In Chapter 5 we apply all the previous machinery to the exponential generating functions of the Bell numbers and also to the generating functions of the partitions of integers; partitions, partitions into distinct parts, partitions in arithmetic progression, some colored and plane partitions. All these functions are in the Hayman class, and therefore are strongly Gaussian. We also prove that these generating functions are in the hypothesis of Baéz-Duarte Theorem 3.2.8, and therefore their coefficients have an asymptotic formula given by this theorem.

The main contributions in Chapter 5 are:

- Using all the previous machinery we greatly simplify the proofs of the asymptotic formulas for the different sort of partitions of integers, see Section 5.1. For the partitions of integers see Theorem 5.1.14, this is Hardy-Ramanujan Theorem. For the partitions of integers into different parts see Theorem 5.1.15. For the partitions of integers into arithmetic progression see Theorem 5.1.16, this is Ingham's Theorem. For the plane partitions see Theorem 5.1.17, this is Wright's Theorem. For the colored partitions with coloring sequence $(j^b)_{j\geq 1}$ see Theorem 5.1.18.

- For the partitions of sets we give a very simple proof of the asymptotic formula for the Bell numbers, see Section 5.2 and in particular Theorem 5.2.1.

In Chapter 6, and searching for a unified approach, we study large powers of analytic functions. Using a variety of local central limit theorems we give asymptotic formulas for the coefficients of power series $\psi^n$, for large values of $n$, when the index of the coefficient has some asymptotic restriction. One instance of this situation is given by Lagrange's inversion formula, the Otter-Meir-Moon Theorem and its extensions. This last theorem has two interpretations: the first interpretations give an asymptotic formula for the coefficients of the generating functions of different families of rooted, weighted, combinatorial trees, this depends on $\psi \in \mathcal{K}$, and, on the Khinchin family side, for the probabilities associated to certain (Galton-Watson) random trees, this depends on the Khinchin family associated to $\psi$ (offspring distribution). We also treat the asymptotics for Lagrangian distributions, an instance of these kind of distributions being the Borel-Tanner distribution.

The main contributions in Chapter 6 are:

- We unify under the same setting different asymptotic formulas for large powers of analytic functions with non-negative coefficients. We greatly simplify the proofs using the Theory of Khinchin families.
- We simplify the proof of Theorem 6.7.1 and its extensions. The Theory of Khinchin families seems to be tailor made to deal with the extremal cases in a very simple way, see also Theorem 6.7.2 and the other cases in there.

# Chapter 1

# Elementary properties of Khinchin families

## Contents

---

In this first chapter we collect the basic definitions and properties of Khinchin families. After describing the most basic definitions, and also some examples of interest, we discuss the growth and range of the mean and variance functions. These functions play an important role in this theory; understanding their intrinsic properties will be essential.

Later on, we introduce the characteristic function and the fulcrum of a power series. We describe their properties and related results; one of these results being Hayman's formula for the coefficients of a power series: this formula represents the $n$-th coefficient of a power series in terms of the characteristic function of its Khinchin family.

We continue with the study of the moment generating function of a Khinchin family and then apply this study to the Chernoff bounds for Khinchin families. Following the lead of Rosenbloom [83] and Schumitzky [88] we conclude this chapter connecting the Wiman-Valiron's theory with the theory of Khinchin families. This connection between the Wiman-Valiron theory and the theory of Khinchin families suggests certain concentration properties of entire Khinchin families that will become apparent in subsequent chapters.

This chapter is based on the papers:

- Maciá, V.J. et al. Khinchin families and Hayman class. *Comput. Methods Funct. Theory* **21** (2021), 851–904, see [17]
- Maciá, V.J. et al. Growth of power series with nonnegative coefficients, and Moments of power series distributions. (submitted) arXiv:2401.14473v2, see [19].

## 1.1 Khinchin families

We denote by $\mathcal{K}$ the class of non-constant power series

$$f(z)=\sum_{n=0}^{\infty} a_n z^n$$

with positive radius of convergence, which have non-negative Taylor coefficients and such that $a_0>0$. Since $f\in\mathcal{K}$ is non-constant, at least one coefficient other than $a_0$ is positive.

A relevant property of the power series $f\in\mathcal{K}$ is that

$$M(f,t)=\max_{|z|\le t}|f(z)|=\max_{|z|=t}|f(z)|=f(t), \qquad \text{for any } t\in[0,R),$$

we will use this property assiduously.

The *Khinchin family* of such a power series $f\in\mathcal{K}$ with radius of convergence $R>0$ is the family of random variables $(X_t)_{t\in[0,R)}$ with values in $\{0,1,\dots\}$ and with mass functions given by

$$\mathbf{P}(X_t=n)=\frac{a_n t^n}{f(t)}\,, \quad \text{for each } n\ge 0 \text{ and } t\in(0,R)\,.$$

Notice that $f(t)>0$ for each $t\in[0,R)$. For $t=0$, we define $X_0\equiv 0$. Formally, the definition of $X_0$ is consistent with the general expression for $t\in(0,R)$, with the convention that $0^0=1$, meaning that $\mathbf{P}(X_0=0)=1$.

No hypothesis upon joint distribution of the variables $X_t$ is considered. Each $(X_t)_{t\in[0,R)}$ is a family of random variables and not a stochastic process.

### 1.1.1 Shifted Khinchin families

Sometimes it is convenient to consider power series $g(z)=\sum_{n=0}^{\infty} b_n z^n$ with radius of convergence $R>0$ and non-negative coefficients with at least two positive coefficients, but which may have $g(0)=0$. We say then that $g$ is in the shifted class $\mathcal{K}_s$. Thus $g\in\mathcal{K}_s$ if there exists an integer $l\ge 0$ such that

$$g(z)/z^l\in\mathcal{K}\,,$$

i.e., if there exists a power series $f(z)\in\mathcal{K}$ with radius of convergence $R>0$ such that $g(z)=z^l f(z)$, for $z\in\mathbb{D}(0,R)$.

To each $g \in \mathcal{K}_s$ we associate a Khinchin family $(Y_t)_{t\ \in(0,R)}$ as above:

$$\mathbf{P}(Y_t = n) = b_n t^n / g(t)\,, \quad \text{for } n \geq 0 \text{ and } t \in (0, R)\,.$$

If $(Z_t)_{t\in[0,R)}$ is the Khinchin family of $f \in \mathcal{K}$, then

$$Y_t \stackrel{d}{=} Z_t + l\,,$$

including, $Y_0 \equiv l$.

### 1.1.2 Examples of Khinchin families

Now, we will showcase some examples of Khinchin families that will be the benchmark along this thesis. In general, we will refer to these basic families, defined below, to exemplify the definitions and the applications.

#### A. The basic families

The first examples of Khinchin families are linked to the most basic and familiar discrete random variables with values in $\{0, 1, \dots\}$, that is, Bernoulli, binomial, geometric (or Pascal), negative binomial and Poisson.

- If $f(z) = 1 + z$, then for each $t > 0$, the variable $X_t$ is a Bernoulli variable with parameter $t/(1 + t)$. Observe that

$$\mathbf{P}(X_t = 0) = \frac{1}{1+t}, \quad \text{and} \quad \mathbf{P}(X_t = 1) = \frac{t}{1+t}, \quad \text{for any } t \in [0, +\infty).$$

- If $f(z) = (1 + z)^N = \sum_{n=0}^{N} \binom{N}{n} z^n$, with an integer $N \geq 1$, then for each $t > 0$, the variable $X_t$, a discrete random variable taking values on $\{0, 1, 2, \dots, N\}$, is a binomial variable with parameters $N$ and $t/(1 + t)$, that is,

$$\mathbf{P}(X_t = n) = \binom{N}{n} \frac{t^n}{(1-t)^N} = \binom{N}{n} \left(\frac{t}{1+t}\right)^n \left(1 - \frac{t}{1+t}\right)^{N-n},$$

  for any $0 \leq n \leq N$ and $t \in (0, +\infty)$.

- If $f(z) = 1/(1 - z)$, then for each $t \in [0, 1)$, the variable $X_t$ is a geometric or Pascal variable (number of failures until first success) with success probability $1 - t$, that is,

$$\mathbf{P}(X_t = n) = t^n(1 - t), \quad \text{for any } n \geq 0 \text{ and } t \in [0, 1).$$

- If $f(z) = 1/(1 - z)^N$, then for each $t \in (0, 1)$, the variable $X_t$ is a negative binomial with parameters $N$ and $1 - t$.

- If $f(z) = e^z$, then, for each $t > 0$, the random variable $X_t$ is a Poisson variable of parameter $t$.

We will refer to the Khinchin families associated to $f(z) = e^z$, $f(z) = 1 + z$ and to $f(z) = 1/(1 - z)$ as Poisson, Bernoulli and Pascal family, respectively.

- If $f(z) = a_N z^N + a_{N-1} z^{N-1} + \cdots + a_0 \in \mathcal{K}$, with an integer $N \geq 1$ and $a_N \neq 0$, is a polynomial of degree $N$. For each $t > 0$ the random variable $X_t$ takes values in the set $\{0, 1, \ldots, N\}$. The Bernoulli and binomial families are particular cases of $f \in \mathcal{K}$ a polynomial.

**B. Khinchin families and partitions**

Here we study the Khinchin families associated to the generating functions of combinatorial objects related to partitions. Our first example is the EGF of the Bell numbers, these numbers give the number of partitions of a set of certain fixed size. We continue studying the OGF of the partitions of integers (different cases: in different parts, in arithmetic progression etc.).

- The EGF $B(z)$ of $\mathcal{B}_n$, the Bell numbers, is given by the entire function

$$B(z) = e^{e^z - 1} = \sum_{n=0}^{\infty} \frac{\mathcal{B}_n}{n!} z^n, \qquad \text{for any } z \in \mathbb{C}.$$

  For any $n \geq 1$, the sequence $\mathcal{B}_n$ gives the number of partitions of a set of size $n$. We adhere to the convention that $\mathcal{B}_0 = 1$ and therefore $B \in \mathcal{K}$.

- The OGF of the number of partitions of an integer $n$, *the partition function*, given by

$$P(z) = \prod_{j=1}^{\infty} \frac{1}{1 - z^j} = \sum_{n=0}^{\infty} p(n)\, z^n\,, \quad \text{for } z \in \mathbb{D}\,,$$

  is in $\mathcal{K}$.

- The OGF $Q(z)$ of partitions into distinct parts $q(n)$ (which is also the OGF of partitions into odd parts) is given by:

$$Q(z) = \prod_{j=1}^{\infty} (1 + z^j) = \prod_{j=0}^{\infty} \frac{1}{1 - z^{2j+1}} = \sum_{n=0}^{\infty} q(n)\, z^n\,, \quad \text{for } z \in \mathbb{D}\,,$$

  is also in $\mathcal{K}$. Observe that

$$Q(z) = \frac{P(z)}{P(z^2)}\,, \quad \text{for } z \in \mathbb{D}\,. \tag{1.1.1}$$

- For integers $a, b \geq 1$, the infinite product

$$P_{a,b}(z) = \prod_{j=0}^{\infty} \frac{1}{1 - z^{aj+b}}\,, \quad \text{for } z \in \mathbb{D}\,, \tag{1.1.2}$$

the OGF of the partitions whose parts lie in the arithmetic progression $\{aj + b : j \geq 0\}$ is also in $\mathcal{K}$. Observe that $P_{1,1} \equiv P$ and that $P_{2,1} \equiv Q$.

For integers $a \geq 1, b \geq 0$, the infinite product $W_a^b$ given by

$$W_a^b(z) = \prod_{j=1}^{\infty} \Big(\frac{1}{1 - z^{j^a}}\Big)^{j^b}, \quad \text{for } z \in \mathbb{D},$$

is also in $\mathcal{K}$. We have $W_1^0 \equiv P$. Also, $W_1^1(z)$, known as the MacMahon function [63], turns out to be the OGF of *plane* partitions; see [7] for a *simple* proof. Besides, $W_a^0(z)$ is the OGF of partitions with parts which are $a$-th powers of positive integers. In general, $W_a^b(z)$ is the OGF of *colored* partitions.

### C. Other families of combinatorial interest

Here we give some examples of EGF which are exponentials of certain power series in $\mathcal{K}_s$. This structure usually appears when studying certain combinatorial objects. See Chapter 4, in particular Section 4.1, for further details.

- The EGF of the functions from $\{1, \ldots, n\}$ into $\{1, \ldots, n\}$ is given by

$$g(z) = \sum_{n=1}^{\infty} \frac{n^n}{n!} z^n.$$

  Observe that $g$ is in $\mathcal{K}_s$. The EGF of the sets of functions from $\{1, \ldots, n\}$ to itself is given by $f(z) = e^{g(z)}$, which is in $\mathcal{K}$.

- The EGF of the rooted labeled Cayley trees is given by

$$T(z) = \sum_{n=1}^{\infty} \frac{n^{n-1}}{n!} z^n$$

  and therefore $T$ is in $\mathcal{K}_s$. The EGF of the forests of rooted labeled Cayley trees is given by $F(z) = e^{T(z)}$, and is in $\mathcal{K}$.

### D. Entire functions of genus 0 with negative zeros

We consider transcendental[1] entire functions $f$ in $\mathcal{K}$ of genus 0 whose zeros are all real and negative. Due to Hadamard's factorization these functions have the form

$$f(z) = a \prod_{j=1}^{\infty} \left(1 + \frac{z}{b_j}\right), \quad \text{for any } z \in \mathbb{C},$$

[1]Transcendental entire functions are entire functions which are not polynomials.

where $a > 0$ and $(b_j)_{j\geq 1}$ is an increasing sequence of positive real numbers verifying that

$$\sum_{j=1}^{\infty} \frac{1}{b_j} < +\infty.$$

The zeros of $f$ are real and negative. We normalize $f(0) = a = 1$.

We denote by

$$N(t) = \#\{j \geq 1 : b_j \leq t\}, \qquad \text{for } t > 0,$$

the counting function of the zeros of $f$. This is a non-decreasing function verifying that $N(t) \to +\infty$, as $t \to +\infty$. *We denote by $\mathcal{Q} \subseteq \mathcal{K}$ the class of entire functions $f \in \mathcal{K}$ with genus* 0 *and (real) negative zeros which are normalized in such a way that $f(0) = 1$.*

## 1.2 Mean and variance functions

For the mean and variance of $X_t$ we reserve the notation $m_f(t) = \mathbf{E}(X_t)$ and $\sigma_f^2(t) = \mathbf{V}(X_t)$, for $t \in [0, R)$. In terms of the power series $f \in \mathcal{K}$, the mean and the variance of $X_t$ may be written as

$$\text{(1.2.1)} \qquad m_f(t) = \frac{tf'(t)}{f(t)} = t\frac{d}{dt}\ln(f(t)), \qquad \sigma_f^2(t) = tm_f'(t)\,, \quad \text{for } t \in [0, R)\,.$$

Observe that for any $t \in [0, R)$ we have

$$m_f(t) = \mathbf{E}(X_t) = \sum_{n=0}^{\infty} n\frac{a_n t^n}{f(t)} = \frac{tf'(t)}{f(t)},$$

and

$$\sigma_f^2(t) = \mathbf{E}(X_t^2) - \mathbf{E}(X_t)^2 = m_f(t)(1 - m_f(t)) + \frac{t^2 f''(t)}{f(t)} = tm_f'(t),$$

Here we have used that

$$\text{(1.2.2)} \qquad \mathbf{E}(X_t^2) = \mathbf{E}(X_t(X_t - 1)) + \mathbf{E}(X_t) = \frac{t^2 f''(t)}{f(t)} + m_f(t),$$

and also that

$$tm_f'(t) = t\left(\ln(f(t))' + t\ln(f(t))''\right) = m_f(t) + \left(\frac{t^2 f''(t)}{f(t)} - m_f(t)^2\right).$$

Therefore $\sigma_f^2(t) = tm_f'(t)$.

**A first look at the range of the mean**

For each $t \in (0, R)$, the variable $X_t$ is not a constant, and so $\sigma_f^2(t) = t m_f'(t) > 0$. Consequently, $m_f(t)$ is strictly increasing in $[0, R)$, though, in general, $\sigma_f(t)$ is not increasing. We denote

$$M_f = \lim_{t \uparrow R} m_f(t)\,. \tag{1.2.3}$$

The mean $m_f$ is an increasing diffeomorphism from $[0, R)$ to $[0, M_f)$.

For $g \in \mathcal{K}_s$, with Khinchin family $(Y_t)$, we also write $m_g(t) = \mathbf{E}(Y_t)$ and $\sigma_g^2(t) = \mathbf{V}(Y_t)$. If $g(z) = z^s h(z)$, with $s \geq 1$ and $h \in \mathcal{K}$ we have that

$$m_g(t) = s + m_h(t) \quad \text{and} \quad \sigma_g^2(t) = \sigma_h^2(t).$$

For finite radius of convergence $R > 0$ and for any integer $k \geq 1$, the function $f(t)$ and all its derivatives $f^{(k)}(t)$ are increasing functions on the interval $[0, R)$. For $k \geq 1$, we denote with $f^{(k)}(R)$ the limit

$$f^{(k)}(R) \triangleq \lim_{t \uparrow R} f^{(k)}(t)\,,$$

including

$$f(R) \triangleq \lim_{t \uparrow R} f(t)\,;$$

these limits exist, although they could be $+\infty$.

For $R = \infty$ we distinguish two cases: If $f$ is a transcendental entire function in $\mathcal{K}$, then the function $f(t)$, and all its derivatives, are increasing functions on the interval $[0, +\infty)$, besides

$$\lim_{t \to +\infty} f^{(k)}(t) = +\infty\,, \qquad \text{for any integer } k \geq 0.$$

If $f$ is a polynomial in $\mathcal{K}$ of degree $N \geq 1$, then $f^{(k)}(t)$ is an increasing function for $0 \leq k \leq N - 1$, besides

$$\lim_{t \to \infty} f^{(k)}(t) = +\infty, \qquad \text{for any integer } 0 \leq k \leq N - 1$$

and

$$\lim_{t \to +\infty} f^{(N)}(t) = N! \quad \text{and} \quad \lim_{t \to +\infty} f^{(k)}(t) = 0, \qquad \text{for any } k > N.$$

### 1.2.1 Some examples

Here we compute the mean and the variance of the basic families. We also give examples for the mean and variance of some families of combinatorial interest, for instance for the OGF of the partitions of integers $P$.

### A. Mean and variance of the basic families

Now we introduce the basic families. These families will appear recurrently along this text.

- For $f(z) = 1 + z$ the radius of convergence is $R = \infty$. The associated Khinchin family $(X_t)$ is the Bernoulli family, then, for any $t \geq 0$, the mean and variance functions are given by

$$m_f(t) = \frac{t}{1+t}, \quad \text{and} \quad \sigma_f^2(t) = \frac{t}{(1+t)^2}.$$

  In this case $M_f = 1$ and therefore $m_f$ is a diffeomorphism from $[0, +\infty)$ onto $[0, 1)$.

- For $f(z) = (1 + z)^N$, with an integer $N \geq 1$, the radius of convergence is $R = \infty$. The associated Khinchin family $(X_t)$ is the binomial family, then, for any $t \geq 0$, the mean and variance functions are given by

$$m_f(t) = \frac{Nt}{1+t}, \quad \text{and} \quad \sigma_f^2(t) = \frac{Nt}{(1+t)^2}.$$

  In this case $M_f = N$ and therefore $m_f$ is a diffeomorphism from $[0, +\infty)$ onto $[0, N)$.

- For $f(z) = 1/(1-z)$ the radius of convergence is $R = 1$. The associated Khinchin family $(X_t)$ is the Geometric (or Pascal) family, then, for any $t \in [0, 1)$, the mean and variance functions are given by

$$m_f(t) = \frac{t}{1-t}, \quad \text{and} \quad \sigma_f^2(t) = \frac{t}{(1+t)^2}.$$

  In this case $M_f = +\infty$ and therefore $m_f$ is a diffeomorphism from $[0, +\infty)$ onto $[0, +\infty)$.

- For $f(z) = 1/(1 - z)^N$, with an integer $N \geq 1$, the radius of convergence is $R = 1$. The associated Khinchin family $(X_t)$ is the Geometric (or Pascal) family, then, for any $t \in [0, 1)$, the mean and variance functions are given by

$$m_f(t) = \frac{Nt}{1-t}, \quad \text{and} \quad \sigma_f^2(t) = \frac{Nt}{(1+t)^2}.$$

  In this case $M_f = +\infty$ and therefore $m_f$ is a diffeomorphism from $[0, +\infty)$ onto $[0, +\infty)$.

- For $f(z) = e^z$ the radius of convergence is $R = \infty$. The associated family $(X_t)$ is the Poisson family, then, for any $t \geq 0$, the mean and variance functions are given by

$$m_f(t) = t, \quad \text{and} \quad \sigma_f^2(t) = t.$$

  In this case $M_f = +\infty$ and therefore $m_f$ is a diffeomorphism from $[0, +\infty)$ to $[0, +\infty)$.

- For $f(z) = a_N z^N + a_{N-1} z^{N-1} + \cdots + a_0 \in \mathcal{K}$, with an integer $N \geq 1$ and $a_N \neq 0$, a polynomial of degree $N = \deg(f)$, the associated family $(X_t)$ takes values in $\{0, 1, 2, \ldots, N\}$ and

$$m_f(t) = \frac{N a_N t^N + (N-1) a_{N-1} t^{N-1} + \cdots + a_1 t}{a_N t^N + a_{N-1} t^{N-1} + \cdots + a_0}, \tag{1.2.4}$$

therefore $M_f = \deg(f) = N$ and $m_f$ is a diffeomorphism from $[0, +\infty)$ onto $[0, N)$.

Observe that

$$\lim_{t\to+\infty} \mathbf{P}(X_t = j) = \begin{cases} 0 & \text{if } j \neq N \\ 1 & \text{if } j = N \end{cases}$$

that is, $X_t$ converges in distribution, as $t \to +\infty$, towards the constant random variable $X \equiv N = \deg f$.

We have, see equation (1.2.2), that

$$\lim_{t\to+\infty} \mathbf{E}(X_t^2) = \lim_{t\to+\infty} (\mathbf{E}(X_t(X_t-1)) + m_f(t)) = \lim_{t\to+\infty} \left(\frac{t^2 f''(t)}{f(t)} - m_f(t)\right),$$

then

$$\lim_{t\to+\infty} \mathbf{E}(X_t^2) = \lim_{t\to+\infty} \left(\frac{t^2 f''(t)}{f(t)} - m_f(t)\right) = N(N-1) - N = N^2.$$

Using the previous limits we conclude that, for a polynomial $f \in \mathcal{K}$, we always have

$$\lim_{t\to+\infty} \sigma_f^2(t) = \lim_{t\to+\infty} \left(\mathbf{E}(X_t^2) - \mathbf{E}(X_t)^2\right) = 0. \tag{1.2.5}$$

As we shall see, the reciprocal is also true, i.e. polynomials in $\mathcal{K}$ are characterized by (1.2.5), see Theorem 1.2.3 below for further details.

### B. Khinchin families and partitions

Here we study the mean and variance functions of generating functions associated to certain types of partitions. The expressions collected in here will be useful in subsequent chapters, see, for instance, Chapter 5.

- For $B(z) = e^{e^z-1}$, the EGF of the Bell numbers $\mathcal{B}_n$, the radius of convergence is $R = \infty$. For any $t \geq 0$, because of (1.2.1), the mean and the variance functions of $B$ are given by

  $$m_B(t) = te^t, \quad \text{and} \quad \sigma_B^2(t) = (t+t^2)e^t.$$

  In this case $M_P = +\infty$ and $m_P$ is a diffeomorphism from $[0, +\infty)$ onto $[0, +\infty)$.

- For $P(z)$, the OGF of the partition of integers $p(n)$, the radius of convergence is $R = 1$. Recall that

  $$P(z) = \prod_{j=1}^{\infty} \frac{1}{1-z^j}, \quad \text{for any } z \in \mathbb{D}.$$

  By using (1.2.1), and the logarithmic derivative of $P$, we find that

  $$m_P(t) = \sum_{j=1}^{\infty} \frac{jt^j}{1-t^j} \quad \text{and} \quad \sigma_P^2(t) = \sum_{j=1}^{\infty} \frac{j^2 t^j}{(1-t^j)^2}, \tag{1.2.6}$$

for any $t \in [0,1)$. Observe that

$$\frac{t}{1-t} \le m_P(t) \le \sigma_P^2(t), \qquad \text{for any } t \in [0,1),$$

therefore $M_f = +\infty$ and $m_f$ is a diffeomorphism from $[0,1)$ onto $[0,+\infty)$. We also have that $\lim_{t\uparrow 1} \sigma_P^2(t) = +\infty$.

By Euler summation, one may obtain convenient asymptotic formulas describing the behavior of $m_P(t)$ and $\sigma_P^2(t)$, as $t \uparrow 1$:

$$\text{(1.2.7)} \qquad m_P(t) \sim \frac{\zeta(2)}{(1-t)^2}\,, \qquad \sigma_P^2(t) \sim \frac{2\zeta(2)}{(1-t)^3}, \quad \text{as } t \uparrow 1.$$

See, for instance, Chapter 5.

- For $Q(z)$, *the OGF of the partitions of integers in distinct parts* $q(n)$, the radius of convergence is $R = 1$. Recall that

$$Q(z) \overset{(1)}{=} \prod_{j=1}^{\infty}(1+z^j) \overset{(2)}{=} \prod_{j=0}^{\infty} \frac{1}{1-z^{2j+1}}, \qquad \text{for any } z \in \mathbb{D}.$$

For any $t \in [0,1)$, using (1.2.1), and the logarithmic derivative of $Q$, we find that

$$m_Q(t) \overset{(1)}{=} \sum_{j=1}^{\infty} \frac{jt^j}{1+t^j} \overset{(2)}{=} \sum_{j\ge 1 \text{ odd}}^{\infty} \frac{jt^j}{1-t^j} \quad \text{and} \quad \sigma_Q^2(t) \overset{(1)}{=} \sum_{j=1}^{\infty} \frac{j^2 t^j}{(1+t^j)^2} \overset{(2)}{=} \sum_{j\ge 1 \text{ odd}}^{\infty} \frac{j^2 t^j}{(1-t^j)^2},$$

Observe that

$$\text{(1.2.8)} \qquad Q(z) = \frac{P(z)}{P(z^2)}, \qquad \text{for any } z \in \mathbb{D},$$

here we use that the partitions of integers into odd parts are in bijection with the partitions of integers into distinct parts or simply that $1 - z^{2j} = (1-z^j)(1+z^j)$:

$$Q(z) = \prod_{j=1}^{\infty}(1+z^j) = \prod_{j=1}^{\infty} \frac{(1+z^j)(1-z^j)}{1-z^j} = \frac{P(z)}{P(z^2)}.$$

Taking logarithmic derivatives in (1.2.8) we can write

$$\text{(1.2.9)} \qquad m_Q(t) = m_P(t) - 2m_P(t^2), \qquad \text{for any } t \in [0,1).$$

and therefore using the expressions for $m_Q$ we find that

$$\frac{t}{1-t} \le m_Q(t) \le m_P(t), \qquad \text{for any } t \in [0,1).$$

Here $m_P$ denotes the mean associated to $P$, the OGF of the partitions of integers.

Similarly, but now using equation (1.2.9), we find that

$$\sigma_Q^2(t) = \sigma_P^2(t) - 4\sigma_P^2(t^2), \qquad \text{for any } t \in [0,1).$$

and therefore

$$\frac{t}{(1-t)^2} \le \sigma_Q^2(t) \le \sigma_P^2(t), \qquad \text{for any } t \in [0,1).$$

In this case $M_Q = +\infty$ and $m_Q$ is a diffeomorphism from $[0,1)$ onto $[0,+\infty)$. We also have that $\lim_{t\uparrow 1} \sigma_Q^2(t) = +\infty$. Observe that

$$m_Q(t) \le \sigma_Q^2(t), \qquad \text{for any } t \in [0,1),$$

this follows from the inequality

$$\frac{jt^j}{(1-t^j)} \le \frac{j^2 t^j}{(1-t^j)^2}, \qquad \text{for any } t \ \in [0,1) \text{ and } j \ge 1.$$

- For integers $a, b \ge 1$, the power series $P_{a,b}(z)$, *the OGF of the partitions in the arithmetic progression* $\mathcal{A} = \{aj + b : j \ge 0\}$, has radius of convergence $R = 1$. Recall that $P_{a,b}$ is given by the infinite product

$$P_{a,b}(z) = \prod_{j=0}^{\infty} \frac{1}{1 - z^{aj+b}}, \qquad \text{for any } z \in \mathbb{D}\,,$$

see (1.1.2) for further details. In this case, for any $t \in [0,1)$, we have

$$m_{P_{a,b}}(t) = \sum_{j \in \mathcal{A}} \frac{jt^j}{1-t^j} \quad \text{and} \quad \sigma_{P_{a,b}}^2(t) = \sum_{j \in \mathcal{A}} \frac{j^2 t^j}{(1-t^j)^2}.$$

Taking $a = b = 1$ we have $\mathcal{A} = \{1, 2, \dots\}$ and therefore we retrieve the mean and variance of $P(z)$, the OGF of the partitions of integers.

Fixing $a = 2$ and $b = 1$ we obtain that $\mathcal{A}$ is the set of odd natural numbers, and therefore we retrieve the mean and variance of $Q(z)$, the OGF of the partitions into distinct parts (and equivalently in odd parts).

The previous example can be easily generalized by taking a subset $\mathcal{A}$ of $\{1, 2, \dots\}$. In this case we obtain

$$P_{\mathcal{A}}(z) = \prod_{j \in \mathcal{A}} \frac{1}{1 - z^j}, \qquad z \in \mathbb{D}.$$

This is the OGF of the partitions with parts in $\mathcal{A}$.

By using the same argument than above we have

$$m_{P_{\mathcal{A}}}(t)=\sum_{j\in\mathcal{A}}\frac{jt^j}{1-t^j}\quad\text{and}\quad\sigma^2_{P_{\mathcal{A}}}(t)=\sum_{j\in\mathcal{A}}\frac{j^2t^j}{(1-t^j)^2},\tag{1.2.10}$$

for any $t\in[0,1)$. We have that $M_{P_{\mathcal{A}}}=+\infty$ and therefore $m_{P_{\mathcal{A}}}$ is a diffeomorphism from $[0,1)$ onto $[0,+\infty)$. We also have that $\lim_{t\uparrow 1}\sigma^2_{P_{\mathcal{A}}}(t)=+\infty$.

Observe that, for $J=\min\mathcal{A}$, we have the inequality

$$\frac{Jt^J}{(1-t^J)}\le m_{P_{\mathcal{A}}}(t)\le\sigma^2_{P_{\mathcal{A}}}(t),\qquad\text{for any } t\in[0,1),$$

this follows combining equation (1.2.10) and the inequality

$$\frac{jt^j}{(1-t^j)}\le\frac{j^2t^j}{(1-t^j)^2},\qquad\text{for any } t\ \in[0,1)\text{ and } j\ge 1.$$

- For any pair of integers $a\ge 1$ and $b\ge 0$, the infinite product $W_{a,b}(z)$ *is the OGF of the colored partitions*, this power series has radius of convergence $R=1$. Recall that

$$W_a^b(z)=\prod_{j=1}^{\infty}\left(\frac{1}{1-z^{j^a}}\right)^{j^b},\qquad\text{for any } z\in\mathbb{D}.$$

For any $t\in[0,1)$, the mean and variance functions are given by

$$m_{W_a^b}(t)=\sum_{j=1}^{\infty}j^b\frac{j^at^{j^a}}{1-t^{j^a}}\quad\text{and}\quad\sigma^2_{W_a^b}(t)=\sum_{j=1}^{\infty}j^b\frac{j^{2a}t^{j^a}}{(1-t^{j^a})^2},$$

respectively. In this case $M_{W_a^b}=+\infty$ and therefore $m_{W_a^b}$ is a diffeomorphism from $[0,1)$ onto $[0,+\infty)$. We also have that $\lim_{t\uparrow 1}\sigma^2_{W_a^b}(t)=+\infty$.

### C. Other Khinchin families of combinatorial interests

Here we study the mean and variance functions of some exponential power series. As mentioned before this exponential structure is ubiquitous in combinatorics, see, for instance, Chapter 4.

- The EGF of the functions from $\{1,2,\dots,n\}$ into $\{1,2,\dots,n\}$ is given by

$$g(z)=\sum_{n=1}^{\infty}\frac{n^n}{n!}z^n,\qquad\text{for any } |z|<\frac{1}{e}.$$

see [32, p. 111 - II.3.2]. Consider the EGF of the sets of functions from $\{1,2,\dots,n\}$ to $\{1,2,\dots,n\}$, that is, $f(z)=e^{g(z)}$ then

$$m_f(t)=tg'(t)=\sum_{n=1}^{\infty}\frac{n^{n+1}}{n!}t^n,\quad\text{and}\quad\sigma^2_f(t)=tg'(t)+t^2g''(t),$$

for any $t \in [0, 1/e)$. Applying Stirling's formula we find that $\lim_{t\uparrow 1/e} m_f(t) = +\infty$, therefore $M_f = +\infty$ and the mean $m_f$ is a diffeomorphism from $[0, 1/e)$ onto $[0, +\infty)$.

Observe that

$$m_f(t) \le \sigma_f^2(t), \qquad \text{for any } t \in [0, 1/e),$$

then we also have that $\lim_{t\uparrow 1/e} \sigma_f^2(t) = +\infty$.

- The EGF of the bijective functions from $\{1, 2, \ldots, n\}$ onto $\{1, 2, \ldots, n\}$, that is, the EGF of the permutations is given by

$$g(z) = \sum_{n=1}^{\infty} \frac{n!}{n!} z^n = \frac{z}{1-z}, \qquad \text{for any } z \in \mathbb{D}.$$

The EGF of sets of bijective function from $\{1, 2, \ldots, n\}$ onto $\{1, 2, \ldots, n\}$ is given by $f(z) = e^{g(z)}$, therefore, for any $t \in [0, 1)$, we have

$$m_f(t) = \frac{t}{(1-t)^2} \qquad \text{and} \qquad \sigma_f^2(t) = \frac{t(1+t)}{(1-t)^3}.$$

Observe that $M_f = +\infty$ and therefore $m_f$ is diffeomorphism from $[0, 1)$ onto $[0, +\infty)$.

- The EGF of the rooted labeled Cayley trees is given by

$$T(z) = \sum_{n=1}^{\infty} \frac{n^{n-1}}{n!} z^n, \qquad \text{for any } |z| \le \frac{1}{e}.$$

This power series is in $\mathcal{K}_s$. Stirling's formula gives that $\lim_{t\uparrow 1/e} T(z) < +\infty$ and $T$ extends to be continuous on the closure of $\mathbb{D}(0, 1/e)$. Observe that

$$\frac{n^{n-1}}{n! e^n} \sim \frac{1}{\sqrt{2\pi}} \frac{1}{n^{3/2}}, \qquad \text{as } n \to \infty.$$

For the first derivative we have, using again Stirling's formula, that $\lim_{t\uparrow 1/e} T'(t) = +\infty$, therefore $M_T = +\infty$ and $m_T$ is a diffeomorphism from $[0, 1/e)$ onto $[0, +\infty)$.

### 1.2.2 Mean and variance in terms of coefficients

We may also write the mean and variance in terms of the Taylor coefficients $a_n$ of $f$. For the mean we have

$$f(t) m_f(t) = \sum_{n=0}^{\infty} n a_n t^n \,, \tag{1.2.11}$$

$$f(t)^2 m_f(t)^2 = \Big(\sum_{n=0}^{\infty} n a_n t^n\Big)^2 = \sum_{n=0}^{\infty} t^n \Big(\sum_{k=0}^{n} a_{n-k} a_k (n-k) k\Big) \,. \tag{1.2.12}$$

Besides, for the moment of order 2, we have

$$f(t)^2\mathbf{E}(X_t^2) = \big(\sum_{n=0}^{\infty} a_n t^n\big)\big(\sum_{n=0}^{\infty} n^2 a_n t^n\big) = \sum_{n=0}^{\infty} t^n \Big(\sum_{k=0}^{n} a_{n-k}a_k \frac{1}{2}\big[(n-k)^2 + k^2\big]\Big)\,, \tag{1.2.13}$$

and thus for the variance $\sigma_f^2(t)$, and any $t \in (0, R)$, we have that

$$f(t)^2\sigma_f^2(t) = f(t)^2\mathbf{E}(X_t^2) - f(t)^2 m_f(t)^2 = \sum_{n=0}^{\infty} t^n \Big(\sum_{k=0}^{n} a_{n-k}a_k \frac{1}{2}\big[((n-k) - k)^2\big]\Big)\,. \tag{1.2.14}$$

### 1.2.3 Mean and variance functions of some canonical products.

In order to explore the behavior of the variance function $\sigma_f^2(t)$, for large values of $t > 0$, we introduce now, as an illustrative example, the mean and the variance functions of some canonical products with negative zeros.

**Class $\mathcal{Q} \subseteq \mathcal{K}$ of canonical products:** We define the class $\mathcal{Q}$ of entire functions of genus 0 with negative zeros as follows: let $(b_k)_{k\geq 1}$ be a sequence of positive numbers increasing to $\infty$ in such a way that $\sum_{k=1}^{\infty} 1/b_k < +\infty$ and let the *canonical product* $f$ be

$$f(z) = \prod_{k=1}^{\infty}\Big(1 + \frac{z}{b_k}\Big), \quad \text{for } z \in \mathbb{C}. \tag{1.2.15}$$

this infinite product gives an entire function in $\mathcal{K}$, which has as zeros $\{-b_k, k \geq 1\}$, the class $\mathcal{Q} \subseteq \mathcal{K}$ is the subclass of $\mathcal{K}$ formed by these canonical products . We denote with

$$N(t) = \#\{k \geq 1 : b_k \leq t\}, \quad \text{for } t > 0\,,$$

the *counting function of zeros* of $f$. Thus $N(t)$ counts the number of zeros of $f$ in the disk $\mathbb{D}(0, t)$.

For the mean and variance functions of $f$, we have, from (1.2.1), that

$$m_f(t) = \sum_{k=1}^{\infty} \frac{t}{t + b_k} \quad \text{and} \quad \sigma_f^2(t) = \sum_{k=1}^{\infty} \frac{b_k t}{(t + b_k)^2}, \quad \text{for } t \geq 0. \tag{1.2.16}$$

Recall that

$$\frac{b_k t}{(t + b_k)^2} < \frac{t}{t + b_k}, \quad \text{for any } t > 0 \text{ and any } k \geq 1,$$

and also that

$$\frac{t}{t+b_k} < \begin{cases} 2\frac{b_k t}{(t+b_k)^2}, & \text{if } b_k > t > 0 \\ 1, & \text{if } b_k \leq t. \end{cases}$$

It follows from these two estimates that

$$\sigma_f^2(t) < m_f(t) < 2\sigma_f^2(t) + N(t), \quad \text{for any } t > 0. \tag{1.2.17}$$

***Estimates of $\sigma_f^2(t)$ for the canonical product $f$ of* (1.2.15).** We will see later that these canonical products provide interesting examples of the behavior of $\sigma_f^2(t)$, which we now estimate precisely, following the lead of Hayman in [49, Theorem 4].

We assume from the start that

$$b_{k+1} \geq 2b_k, \quad \text{for any } k \geq 1. \tag{1.2.18}$$

Consider the positive function $\varphi(x) = x/(1+x)^2$, for $x > 0$. This function $\varphi$ attains a maximum value of $1/4$ at $x = 1$, and satisfies

$$\frac{x}{4} < \varphi(x) < x, \quad \text{for } x \in (0,1), \quad \text{and} \quad \frac{1}{4x} < \varphi(x) < 1/x, \quad \text{for } x > 1\,. \tag{1.2.19}$$

We can write

$$\sigma_f^2(t) = \sum_{k=1}^{\infty} \varphi(t/b_k)\,, \quad \text{for any } t < 0\,.$$

For each $n \geq 2$, we let $I_n$ denote the interval

$$I_n := [\sqrt{b_{n-1}b_n}, \sqrt{b_n b_{n+1}}]\,.$$

For $t \in I_n, n \geq 2$, and $k > n$, we have $t/b_k < 1$ and, on account of (1.2.18), we have that $b_k \geq 2^{k-(n+1)} b_{n+1}$. Thus, using (1.2.19), we see that

$$\sum_{k>n} \varphi\Big(\frac{t}{b_k}\Big) < t \sum_{k>n} \frac{1}{b_k} \leq \frac{2t}{b_{n+1}} \leq 2\sqrt{\frac{b_n}{b_{n+1}}}\,, \quad \text{for any } t \in I_n, n \geq 2\,.$$

Besides, for $t \in I_n, n \geq 2$, and $k < n$, we have $t/b_k > 1$ and, on account of (1.2.18), we have that $b_k \leq b_{n-1}/2^{n-1-k}$. Thus, using (1.2.19), we see that

$$\sum_{k<n} \varphi\Big(\frac{t}{b_k}\Big) < \frac{1}{t} \sum_{k<n} b_k \leq \frac{2}{t} b_{n-1} \leq 2\sqrt{\frac{b_{n-1}}{b_n}}\,,$$

for any $t \in I_n, n \geq 2$.

It follows then that

$$\varphi\big(\frac{t}{b_n}\big) \le \sigma_f^2(t) \le \varphi\big(\frac{t}{b_n}\big) + 4\max\Big\{\sqrt{b_n/b_{n+1}},\ \sqrt{b_{n-1}/b_n}\Big\}, \quad \text{for any } t \in I_n\,. \tag{1.2.20}$$

and so, using that $\varphi(x)$ decreases whenever $x$ moves away from 1,

$$\begin{aligned}\sup_{t\in I_n} \sigma_f^2(t) &\le \frac{1}{4} + 4\max\Big\{\sqrt{b_n/b_{n+1}},\ \sqrt{b_{n-1}/b_n}\Big\},\\ \inf_{t\in I_n} \sigma_f^2(t) &\ge \frac{1}{4}\min\Big\{\sqrt{b_n/b_{n+1}},\ \sqrt{b_{n-1}/b_n}\Big\},\end{aligned} \tag{1.2.21}$$

for any $n \ge 2$.

### 1.2.4 Growth and range of the mean.

Since the mean function $m_f(t)$ is increasing, its range is given by $[0, M_f)$. See equation (1.2.3) for the definition of $M_f$.

The case where $M_f = \infty$ is particularly relevant. If $M_f = \infty$, then $m_f(t)$ is a diffeomorphism from $[0, R)$ onto $[0, +\infty)$ and, in particular, for each $n \ge 1$, there exists a unique $t_n \in (0, R)$ such that $m_f(t_n) = n$. These $t_n$ play an important role in Hayman's identity, see (1.3.31), and in Hayman's asymptotic formula, see equation (1.3.29).

The power series $f \in \mathcal{K}$ such that

$$M_f = \lim_{t\uparrow R} m(t) > 1 \tag{1.2.22}$$

comprise a subclass $\mathcal{K}^\star$ that will be relevant in future chapters, see Chapter 6. If $f \in \mathcal{K}^\star$, the unique value $\tau \in (0, R)$ such that $m(\tau) = 1$, is called the *apex* of $f$. This apex $\tau$ is characterized by $\tau f'(\tau) = f(\tau)$. Observe that the existence of apex for $f$ precludes $f$ from being a polynomial of degree 1, and, in particular, implies that $f''(\tau) > 0$.

The following discussion, collected later on as a Lemma, describes the quite specific cases where $M_f < +\infty$; This is Lemma 2.2 of [17].

Let $f \in \mathcal{K}$ with radius of convergence $R > 0$ and assume that $M_f < \infty$. Then

$$(\star) \quad m(t) = \frac{\sum_{n=0}^{\infty} n a_n t^n}{\sum_{n=0}^{\infty} a_n t^n} = \frac{t f'(t)}{f(t)} \le M_f\,, \quad \text{for } t \in (0, R)\,.$$

For any fixed $t_0 \in (0, R)$, we then have upon integration in $(\star)$ that

$$(\star\star) \quad \ln\Big(\frac{f(t)}{f(t_0)}\Big) \le M_f \ln\Big(\frac{t}{t_0}\Big)\,, \quad \text{for } t \in [t_0, R)\,.$$

Let us distinguish now between radius of convergence $R$ being finite or not.

• If $R < +\infty$ and $M_f < \infty$, then ($\star\star$) implies that $\lim_{t\uparrow R} f(t) < \infty$ and then $\sum_{n=0}^{\infty} a_n R^n < \infty$. And, besides, from ($\star$) we deduce also that

$$(\flat) \quad \sum_{n=0}^{\infty} n a_n R^n < \infty\,.$$

And, conversely, if ($\flat$) holds and $f$ in $\mathcal{K}$ has radius of convergence $R$, then

$$M_f = \frac{\sum_{n=0}^{\infty} n a_n R^n}{\sum_{n=0}^{\infty} a_n R^n} < \infty\,.$$

• If $R = \infty$ and $M_f < \infty$, then ($\star\star$) implies that $f$ is a polynomial. And conversely, see equation (1.2.4), and the discussion in there, for a polynomial $f$ in $\mathcal{K}$ of degree $N$, one actually has $M_f = N$.

In summary:

**Lemma 1.2.1** (Lemma 2.2 of [17]). *For $f(z) = \sum_{n=0}^{\infty} a_n z^n$ in $\mathcal{K}$ with radius of convergence $R > 0$, we have $M_f < \infty$ in just the following two cases:*

(1) *if $R < \infty$ and $\sum_{n=0}^{\infty} n a_n R^n < \infty$,*

(2) *and if $R = \infty$ and $f$ is a polynomial.*

*In the first case, we have $M_f = (\sum_{n=0}^{\infty} n a_n R^n)/(\sum_{n=0}^{\infty} a_n R^n)$. For a polynomial $f \in \mathcal{K}$, we have $M_f = \deg(f)$.*

The following Lemma gives upper and lower bounds for the quotient $f(\lambda t)/f(t)$ in terms of the mean $m_f(t)$ of the Khinchin family $(X_t)$.

**Lemma 1.2.2** (Simić, [90]). *For $\lambda > 1$ and $t > 0$, with $\lambda t < R$, we have*

$$m_f(t) \ln \lambda \le \ln\left(\frac{f(\lambda t)}{f(t)}\right) \le m_f(\lambda t) \ln \lambda\,.$$

*Proof.* Since $m_f$ is an increasing function we have that

$$m_f(t) \ln \lambda < \int_t^{\lambda t} m_f(u) \frac{du}{u} < m_f(\lambda t) \ln \lambda\,,$$

but since

$$\int_t^{\lambda t} m_f(u) \frac{du}{u} = \int_t^{\lambda t} (\ln f)'(u) du = \ln\left(\frac{f(\lambda t)}{f(t)}\right)\,,$$

the result follows. □

### 1.2.5 Growth and range of the variance $\sigma_f^2$

Some of the results below concerning $\sigma_f^2(t)$ will depend *on the gaps among the indices* of the power series. We introduce next some convenient notation.

Let $(n_k)_{k\geq 1}$ be the increasing sequence of indices so that $a_{n_k}\neq 0$, for each $k\geq 1$, and $a_n=0$, if $n\notin\{n_k : k\geq 1\}$. Thus the $(a_{n_k})$ are the nonzero Taylor coefficients of $f$, and the random variables $(X_t)$ take exactly the values $(n_k)$. Observe that $n_1=0$, and that for a polynomial $f$, the sequence $(n_k)_{k\geq 1}$ is finite. Define $\operatorname{gap}(f)$ and $\overline{\operatorname{gap}}(f)$ as

(1.2.23) $$\operatorname{gap}(f)=\sup_{k\geq 1}(n_{k+1}-n_k) \quad\text{and}\quad \overline{\operatorname{gap}}(f)=\limsup_{k\to\infty}(n_{k+1}-n_k).$$

It is always the case that $\operatorname{gap}(f)\geq\overline{\operatorname{gap}}(f)\geq 1$, for any $f\in\mathcal{K}$ which is not a polynomial. For polynomials, we still have $\operatorname{gap}(f)\geq 1$, and we can define $\overline{\operatorname{gap}}(f)=0$.

We divide the discussion on variance growth depending on whether $R$ is infinite or finite.

#### A. Entire functions, $R=\infty$

A. *Lower bounds for* $\sup_{t>0}\sigma_f^2(t)$ *and for* $\limsup_{t\to\infty}\sigma_f^2(t)$.

The (universal) lower bounds on $\sigma_f$ that we are about to discuss originate with Hayman's results in [49] quantifying Hadamard's three lines theorem for entire functions (not necessarily in $\mathcal{K}$). See also [2] and [57].

**Theorem 1.2.3** (Boĭchuk–Gol'dberg, [13]). *If $f\in\mathcal{K}$ is entire, then*

$$\sup_{t>0}\sigma_f^2(t)\geq\frac{1}{4}\operatorname{gap}(f)^2 \quad \Big(\geq\frac{1}{4}\Big).$$

*If, moreover, $f$ is transcendental, then*

$$\limsup_{t\to\infty}\sigma_f^2(t)\geq\frac{1}{4}\overline{\operatorname{gap}}(f)^2 \quad \Big(\geq\frac{1}{4}\Big).$$

Recall that an entire function $f$ is termed *transcendental* if it is not a polynomial.

This result is Theorem 2 of [13]. See also Theorem 1 in [1] and Lemma 2.5 in [74]. The probabilistic argument below is simpler than the original proof; it uses the discreteness of the random variables $X_t$.

*Proof of Theorem 1.2.3.* Let $(n_k)_{k=1}^N$ be the indices of the nonzero coefficients of $f$, with $N\leq+\infty$.

Fix $k<N$ and take $t^\star>0$ so that $m_f(t^\star)=(n_{k+1}+n_k)/2$, i.e., the midpoint of the interval $[n_k,n_{k+1}]$. Such $t^\star$ exists because $m_f(t)$ is a continuous (and increasing) function, $m_f(0)=0$, and $M_f=\infty$ or $M_f=\operatorname{degree}(f)$ if $f$ is a polynomial (recall Lemma 1.2.1).

As $X_{t^\star}$ takes the values $n_1,n_2,\ldots$, clearly $|X_{t^\star}-m_f(t^\star)|\geq\frac{1}{2}(n_{k+1}-n_k)$ with probability 1. This gives

$$\sigma_f^2(t^\star)=\mathbf{E}\big((X_{t^\star}-m_f(t^\star))^2\big)\geq\frac{1}{4}(n_{k+1}-n_k)^2.$$

The statements now follows by taking sup and limsup in the inequality above and appealing to the definitions of gap and $\overline{\operatorname{gap}}$. □

In fact, the very same argument shows, for the absolute centered moments, that if $f \in \mathcal{K}$ is entire, then

$$\sup_{t>0} \mathbf{E}(|X_t - m_f(t)|^p) \geq \frac{1}{2^p}\,\mathrm{gap}(f)^p, \quad \text{for any } p > 0,$$

and, moreover, that if $f$ is transcendental, then

$$\limsup_{t\to\infty} \mathbf{E}(|X_t - m_f(t)|^p) \geq \frac{1}{2^p}\,\overline{\mathrm{gap}}(f)^p, \quad \text{for any } p > 0.$$

As for the *sharpness of the sup part* of Theorem 1.2.3, consider the case $f(z) = a + bz$, with $a, b > 0$, for which $X_t$ is a Bernoulli variable with success probability $bt/(a + bt)$. In this case, one has $\sigma_f^2(t) = abt/(a + bt)^2$, which takes its maximum value of $1/4$ at $t = a/b$.

In fact, the converse is also true.

**Theorem 1.2.4** (Abi-Khuzzam, [1]). *For $f \in \mathcal{K}$ entire, $\sup_{t>0} \sigma_f^2(t) = 1/4$ if and only if $f(z) = a + bz$, with $a, b > 0$.*

This result is Theorem 3 of [1]. See also Lemma 2.5 in [74]. The probabilistic argument below is again simpler than the original proof.

*Proof of Theorem 1.2.4.* The 'if' part has been discussed above.

Assume that $f \in \mathcal{K}$ is entire and that $\sup_{t>0} \sigma_f^2(t) = 1/4$. Theorem 1.2.3 gives that $\mathrm{gap}(f) = 1$. Let $a_n$ and $a_{n+1}$ be any two nonzero consecutive coefficients of $f$, and let $t^\star$ be such that $m_f(t^\star) = n+1/2$ (observe that, in any case, $M_f \geq n+1$). We have that $|X_{t^\star} - m_f(t^\star)| \geq 1/2$. By hypothesis, $\mathbf{E}((X_{t^\star} - m_f(t^\star))^2) = \sigma_f^2(t^\star) \leq 1/4$, and thus $|X_{t^\star} - m_f(t^\star)| = 1/2$ with probability 1, which means that $X_{t^\star}$ only takes the values $n$ and $n + 1$, and thus $f(z) = a_n z^n + a_{n+1} z^{n+1}$. Since $f(0) > 0$, because $f \in \mathcal{K}$, it must be the case that $n = 0$ and $f(z) = a + bz$, with $a, b > 0$. □

For the *sharpness of the limsup part* of Theorem 1.2.3, consider a canonical product $h$ given by the infinite product

$$h(z) = \prod_{n=1}^{\infty} \Big(1 + \frac{z}{b_n}\Big), \tag{1.2.24}$$

where $(b_n)_{n\geq 1}$ is a sequence of positive numbers increasing to $+\infty$ so that (1.2.18) holds and, in fact, such that $\lim_{n\to\infty} b_{n+1}/b_n = +\infty$; in particular, $\sum_{k\geq 1} 1/b_k < +\infty$ holds. Obviously, $\mathrm{gap}(h) = 1$ and $\overline{\mathrm{gap}}(h) = 1$. It follows directly from the estimates (1.2.21) that $\limsup_{t\to\infty} \sigma_h^2(t) = 1/4$; this is Hayman's example in Theorem 4 of [49]. If $h$ is multiplied by a polynomial $p$ in such a way that $f = ph \in \mathcal{K}$, then it is still the case that $\limsup_{t\to\infty} \sigma_f^2(t) = 1/4$.

As it turns out, Abi-Khuzzam has characterized (see Theorem 2 in [2] and its proof) the entire functions $f \in \mathcal{K}$ with $\limsup_{t\to\infty} \sigma_f^2(t) = 1/4$ as precisely those entire functions $f \in \mathcal{K}$ which factorize as

$$f(z) = p(z) \prod_{n=1}^{\infty} \Big(1 + \frac{z}{b_n}\Big),$$

where the $b_n$ are as in Hayman's example, and where $p$ is a polynomial.

B. *Limit of* $\sigma_f^2(t)$ *as* $t \to \infty$.

Regarding the existence and possible limits of $\sigma_f^2(t)$ as $t \to \infty$, the following holds.

B.1. *Polynomials.* Polynomials $f \in \mathcal{K}$ are characterized, among the entire functions in $\mathcal{K}$, by

$$\lim_{t\to\infty} \sigma_f^2(t) = 0. \tag{1.2.25}$$

A direct calculation with the formulas (1.2.1) shows that for polynomials (1.2.25) holds; in fact $\sigma_f^2(t) = O(1/t)$, as $t \to \infty$. The converse follows, for instance, from Theorem 1.2.3.

B.2. *Transcendental functions.* As shown by Hilberdink in [51], for a transcendental entire function *it is never the case that* $\lim_{t\to\infty} \sigma_f(t)$ *exists and it is finite.*

However, there are entire functions $f \in \mathcal{K}$ for which $\limsup_{t\to\infty} \sigma_f^2(t) < +\infty$ and $\liminf_{t\to\infty} \sigma_f^2(t) > 0$. To see this, just consider a canonical product $f$ as in (1.2.15), with $b_n = 2^n$, and apply (1.2.21).

It is also possible to have $\lim_{t\to\infty} \sigma_f(t) = \infty$, as shown, for instance, by the exponential function $f(z) = e^z$, where $\sigma_f^2(t) = t$, for $t \geq 0$.

We emphasize that for $f \in \mathcal{K}$ entire, if $\lim_{t\to\infty} \sigma_f^2(t)$ exists, then that limit is 0 (just for polynomials) or $+\infty$.

C. *Boundedness of* $\sigma_f^2(t)$.

We discuss now when, if ever,

$$\sup_{t>0} \sigma_f(t) < +\infty \tag{1.2.26}$$

does hold for an entire function $f \in \mathcal{K}$.

For polynomials $f \in \mathcal{K}$, (1.2.26) holds since, in fact, in this case, $\lim_{t\to\infty} \sigma_f(t) = 0$.

In general, *if* (1.2.26) *holds, then the order* $\rho(f)$ *of* $f$ *must be zero.* (See Section 2.5 for details about the order of an entire function $f$ in $\mathcal{K}$.)

To see this, observe that, since $tm_f'(t) = \sigma_f^2(t)$, integrating, we deduce that $m_f(t) = O(\ln t)$, as $t \to \infty$. A further integration, using (1.2.1), shows that

$$\ln f(t) = O((\ln t)^2)\,, \quad \text{as } t \to \infty. \tag{1.2.27}$$

This gives that $\rho(f) = 0$, see (2.5.2).

Alternatively, Proposition 2.5.3 below shows that

$$\rho(f) \leq \big(\sup_{t>0} \sigma_f^2(t)\big) \frac{1}{M_f}\cdot$$

If $f$ is a polynomial, then $\rho(f) = 0$; and if $f$ is transcendental, $M_f = \infty$, and thus $\rho(f) = 0$.

However, $\rho(f) = 0$, or even the stronger condition (1.2.27), are not enough to ensure the boundedness of $\sigma_f^2(t)$.

Consider the canonical product

$$g(z)=\prod_{k=1}^{\infty}\Big(1+\frac{z}{b_k}\Big)^k,$$

where $(b_k)_{k\geq 1}$ is a sequence of positive numbers increasing to $\infty$ satisfying (1.2.18) and such that $\sum_{k=1}^{\infty} k/b_k<+\infty$. For each $k\geq 1$, $-b_k$ is a zero with multiplicity $k$ of the entire function $g\in\mathcal{K}$.

Additionally, we assume also that for some constant $H>0$, the $b_k$ satisfy

$$\sum_{k<n} k b_k < H b_n \quad\text{and}\quad \sum_{k>n}\frac{k}{b_k}\leq \frac{H}{b_n}\quad \text{for each } n\geq 2.$$

In this case,

$$\sigma_g^2(t)=\sum_{k=1}^{\infty} k\,\varphi\Big(\frac{t}{b_k}\Big)\,,\quad \text{for any } t>0,$$

where $\varphi(x)=x/(1+x)^2$.

Let $C$ denote a generic positive constant. With the notations of Section 1.2.3 and estimating as in there, we obtain

$$n\,\varphi\Big(\frac{t}{b_n}\Big)\leq \sigma_g^2(t)\leq n+C\,n\max\big\{\sqrt{b_{n-1}/b_n},\sqrt{b_n/b_{n+1}}\,\big\}\,,\quad \text{for } t\in I_n \text{ and } n\geq 2,$$

where $I_n=[\sqrt{b_{n-1}b_n},\sqrt{b_nb_{n+1}}]$. It follows, in particular, since $\sigma_g^2(b_n)\geq n/4$, that $\limsup_{t\to\infty}\sigma_g^2(t)=+\infty$, and also that,

$$\sup_{t\in I_n}\sigma_g^2(t)\leq Cn.$$

From (1.2.17), we see that

$$m_g(t)\leq Cn+\sum_{k\leq n} k\leq Cn^2 \quad \text{for any } t\in I_n \text{ and } n\geq 2. \tag{1.2.28}$$

Consider now the specific sequence $b_k=e^{e^k}$, $k\geq 1$, which satisfies the requirements above. In this case, (1.2.28) translates into

$$m_g(t)\leq C(\ln\ln t)^2\,,\quad \text{for } t\geq 2.$$

Integrating, this gives, for this example, that

$$\ln g(t)=O(\ln t(\ln\ln t)^2)\quad \text{as } t\to\infty,$$

which implies $\rho(g)=0$. Notice that for any function $\Phi(t)$ slowly increasing to $\infty$, the sequence $(b_k)_{k\geq 1}$ can be chosen so that $\ln g(t)=O(\Phi(t)\ln t)$ just by making $b_k$ increase fast enough. This is the best that can be expected, because if an entire function $h$ in $\mathcal{K}$ grows as $\ln h(t)=O(\ln t)$ as $t\to\infty$, then $h$ is a polynomial.

**B. Finite radius, $R < \infty$**

We now turn to functions $f \in \mathcal{K}$ with finite radius $R$ of convergence. We have the following results on the behavior of $\sigma_f^2(t)$. We divide the discussion according to whether $M_f$ is finite or not.

- *Case* $M_f = +\infty$.

In this case,

$$\sup_{t\in(0,R)} \sigma_f^2(t) = +\infty\,, \tag{1.2.29}$$

and also, of course, $\limsup_{t\uparrow R} \sigma_f^2(t) = +\infty$. To verify (1.2.29), assume that $\sup_{t\in(0,R)} \sigma_f^2(t) = S < +\infty$. Thus $tm_f'(t) \leq S$, for $t \in [0, R)$. Integrating between $R/2$ and $t \in (R/2, R)$, we would have that

$$m_f(t) \leq m_f(R/2) + S \ln\Big(\frac{2t}{R}\Big), \quad \text{for } t \in (R/2, R),$$

which implies, by letting $t \uparrow R < +\infty$, that $M_f \leq m_f(R/2) + S \ln 2 < +\infty$.

- *Case* $M_f < +\infty$.

If $M_f < +\infty$, then $\sum_{n=0}^{\infty} n a_n R^n < +\infty$, see Lemma 1.2.1, and in fact

$$\Sigma := \lim_{t\uparrow R} \sigma_f^2(t) = \frac{\sum_{n=0}^{\infty} n^2 a_n R^n}{\sum_{n=0}^{\infty} a_n R^n} - \left(\frac{\sum_{n=0}^{\infty} n a_n R^n}{\sum_{n=0}^{\infty} a_n R^n}\right)^2.$$

It is always the case that $\lim_{t\uparrow R} \sigma_f^2(t) > 0$, since $\Sigma$ is the variance of the random variable $Z$ that takes, for each integer $n \geq 0$, the value $n$ with probability $a_n R^n / (\sum_{k=0}^{\infty} a_k R^k)$, and $Z$ is a non-constant variable since $f$ is in $\mathcal{K}$.

But there is no absolute positive lower bound for $\lim_{t\uparrow R} \sigma_f^2(t)$. For $\varepsilon > 0$, the power series $f(z) = 1 + \varepsilon \sum_{n=1}^{\infty} z^n/n^4$ is in $\mathcal{K}$ and has radius of convergence $R = 1$. We have

$$M_f = \frac{\varepsilon\zeta(3)}{1+\varepsilon\zeta(4)} \quad \text{and} \quad \lim_{t\uparrow 1} \sigma_f^2(t) = \frac{\varepsilon\zeta(2)}{1+\varepsilon\zeta(4)} - \Big(\frac{\varepsilon\zeta(3)}{1+\varepsilon\zeta(4)}\Big)^2,$$

which tends to 0 as $\varepsilon \downarrow 0$.

The example $f(z) = \sum_{n=0}^{\infty} z^n/(1+n)^3$ shows that $\lim_{t\uparrow R} \sigma_f^2(t) = \infty$ may happen.

## 1.3 Basic properties of Khinchin families

Now we study the basic properties of Khinchin families; we introduce the characteristic function of a Khinchin family and its normalized version. Later on, we discuss continuous families of random variables. Normalized Khinchin families are examples of continuous families of random variables, this will be relevant when taking limits of certain differences of characteristic functions, as $t \uparrow R$. This concept will also be relevant in the proof of some results of Chapter 6.

We pass then to the study of the fulcrum of a power series $f \in \mathcal{K}$. The fulcrum is a fundamental object in the theory of Khinchin families, see, for instance, Chapter 3 and in particular Section 3.1-A and D.

Along this section we let $f(z) = \sum_{n=0}^{\infty} a_n z^n$ be a power series in $\mathcal{K}$ with radius of convergence $R > 0$ and Khinchin family $(X_t)_{t\in[0,R)}$.

### 1.3.1 Derivative power series $\mathcal{D}_f$ and its family

For a power series $f(z) = \sum_{n=0}^{\infty} a_n z^n$ in $\mathcal{K}$ with radius of convergence $R > 0$ and with at least 3 nonzero coefficients, consider the power series $\mathcal{D}_f$ given by

$$\mathcal{D}_f(z) = z f'(z)\,, \quad \text{for } |z| < R\,. \tag{1.3.1}$$

This power series $\mathcal{D}_f$ has radius of convergence $R$ and lies in $\mathcal{K}_s$. Let $(X_t)$ and $(W_t)$ be the Khinchin families associated to $f$ and $\mathcal{D}_f$, respectively; both families are defined for $0 \le t < R$.

Observe that for $t > 0$ and $n \ge 1$ we have that

$$\mathbf{P}(W_t = n) = \frac{n a_n t^n}{t f'(t)} = \frac{1}{m_f(t)} \frac{n a_n t^n}{f(t)}\,. \tag{1.3.2}$$

Thus, for $t > 0$ we have that

$$m_{\mathcal{D}_f}(t) = \mathbf{E}(W_t) = \frac{\mathbf{E}(X_t^2)}{m_f(t)} = \frac{\mathbf{E}(X_t^2)}{\mathbf{E}(X_t)}\,.$$

This gives us that the quotient $\mathbf{E}(X_t^2)/\mathbf{E}(X_t)$ is monotonically increasing in the interval $(0, R)$. Observe that if the power series $f$ has only two nonzero coefficients, say $a_N$, with $N \ge 1$, besides $a_0$, then $\mathbf{E}(X_t^2)/\mathbf{E}(X_t) = N$, for any $t \in (0, R)$.

In general, we have that

$$\mathbf{E}(W_t^p) = \frac{1}{m_f(t)} \mathbf{E}(X_t^{p+1})\,, \quad \text{for any } p > 0 \text{ and any } t \in (0, R)\,. \tag{1.3.3}$$

For the power series $\mathcal{D}_f$ a direct application of Lemma 1.2.2 gives that

**Lemma 1.3.1** (Simić [90]). *For $\lambda > 1$ and $t > 0$, with $\lambda t < R$, we have*

$$m_{\mathcal{D}_f}(t) \ln \lambda \le \ln \frac{\mathcal{D}_f(\lambda t)}{\mathcal{D}_f(t)} \le m_{\mathcal{D}_f}(\lambda t) \ln \lambda\,.$$

### 1.3.2 Normalization and characteristic functions

For each $t \in (0, R)$, we denote by $\breve{X}_t$ the normalization of $X_t$:

$$\breve{X}_t = \frac{X_t - m_f(t)}{\sigma_f(t)}, \quad \text{for } t \in (0, R)\,.$$

Observe that $\breve{X}_t$ is not defined for $t = 0$.

The characteristic function of $X_t$ may be written in terms of the power series $f$ as:

$$\mathbf{E}(e^{\imath\theta X_t}) = \sum_{n=0}^{\infty} e^{i\theta n} \frac{a_n t^n}{f(t)} = \frac{f(t e^{\imath\theta})}{f(t)}\,, \quad \text{for } t \in (0, R) \text{ and } \theta \in \mathbb{R}\,,$$

while for its normalized version $\breve{X}_t$ we have

$$\mathbf{E}(e^{\imath\theta\breve{X}_t}) = \mathbf{E}(e^{\imath\theta X_t/\sigma_f(t)})e^{-\imath\theta m_f(t)/\sigma_f(t)}\,, \quad \text{for } t \in (0,R) \text{ and } \theta \in \mathbb{R}\,,$$

and so,

$$\left|\mathbf{E}(e^{\imath\theta\breve{X}_t})\right| = \left|\mathbf{E}(e^{\imath\theta X_t/\sigma_f(t)})\right|, \quad \text{for } t \in (0,R) \text{ and } \theta \in \mathbb{R}\,.$$

### A. Some examples of characteristic function of Khinchin families

We use the basic families to illustrate the definition of characteristic function of a Khinchin family and its normalized version.

• Let $f(z) = 1+z$. In this case $R=\infty$, and the mean and variance functions are $m_f(t) = t/(1+t)$ and $\sigma_f^2(t) = t/(1+t)^2$. For each $t > 0$, the variable $X_t$ is a Bernoulli variable with parameter $t/(1+t)$ and its characteristic function is

$$\mathbf{E}(e^{\imath\theta X_t}) = \frac{1+te^{\imath\theta}}{1+t}\,, \quad \text{for } \theta \in \mathbb{R} \text{ and } t > 0\,, \tag{1.3.4}$$

and thus

$$\mathbf{E}(e^{\imath\theta\breve{X}_t}) = \frac{te^{\imath\theta/\sqrt{t}} + e^{-\imath\theta\sqrt{t}}}{1+t}\,, \quad \text{for } \theta \in \mathbb{R} \text{ and } t > 0\,. \tag{1.3.5}$$

• Let $f(z) = (1+z)^N$, with integer $N \geq 1$. In this case $R = \infty$ and the mean and variance functions are $m_f(t) = Nt/(1+t)$ and $\sigma_f^2(t) = Nt/(1+t)^2$. For each $t > 0$ the variable $X_t$ is a binomial variable with parameter $N$ and $p = t/(1+t)$ and its characteristic function is

$$\mathbf{E}(e^{i\theta X_t}) = \left(\frac{1+te^{i\theta}}{1+t}\right)^N, \quad \text{for } \theta \in \mathbb{R} \text{ and } t > 0. \tag{1.3.6}$$

This last expression follows from the fact that $X_t$ is the sum of $N$ independent Bernoulli random variables of $t/(1+t)$ or simply from (1.3.4).

As we shall see later on, combining (1.3.24) and (1.3.5) we find that

$$\mathbf{E}(e^{i\theta\breve{X}_t}) = \left(\frac{te^{\imath\theta/\sqrt{Nt}} + e^{-\imath\theta\sqrt{t/N}}}{1+t}\right)^N, \quad \text{for } \theta \in \mathbb{R} \text{ and } t > 0\,. \tag{1.3.7}$$

• Let $f(z) = 1/(1-z)$. In this case $R = 1$, and the mean and variance functions are $m_f(t) = t/(1-t)$ and $\sigma_f^2(t) = t/(1-t)^2$. For each $t \in (0,1)$, the variable $X_t$ is a geometric variable (number of failures until first success) of parameter $1-t$ and its characteristic function is

$$\mathbf{E}(e^{\imath\theta X_t}) = \frac{1-t}{1-te^{\imath\theta}}\,, \quad \text{for } \theta \in \mathbb{R} \text{ and } t \in (0,1)\,, \tag{1.3.8}$$

and thus

(1.3.9) $$\mathbf{E}(e^{\imath\theta\breve{X}_t})=\frac{1-t}{1-te^{\imath\theta(1-t)/\sqrt{t}}}e^{-\imath\theta\sqrt{t}}=\frac{1-t}{e^{\imath\theta\sqrt{t}}-te^{\imath\theta/\sqrt{t}}}\,,\quad \text{for } \theta\in\mathbb{R} \text{ and } t\in(0,1)\,.$$

• Let $f(z)=1/(1-z)^N$, with integer $N\ \geq 1$. In this case $R=1$, and the mean and the variance functions are $m_f(t)=\frac{Nt}{1-t}$ and $\sigma_f^2(t)=Nt/(1-t)^2$. For each $t\in(0,1)$, the variable $X_t$ is a negative binomial variables of parameters $N\geq 1$ and $p=1-t$, and its characteristic function is

$$\mathbf{E}(e^{i\theta X_t})=\left(\frac{1-t}{1-te^{i\theta}}\right)^N,\quad \text{for } \theta\ \in\mathbb{R} \text{ and } t\in(0,1).$$

This last expression follows since $X_t$ is the sum of $N$ independent geometric variables of parameter $t/(1+t)$ or simply from (1.3.8).

As we shall see later on, combining (1.3.24), below, with (1.3.8), we conclude that

(1.3.10) $$\mathbf{E}(e^{\imath\theta\breve{X}_t})=\left(\frac{1-t}{1-te^{\imath\theta/(1-t)/\sqrt{Nt}}}e^{-\imath\theta\sqrt{t/N}}\right)^N=\left(\frac{1-t}{e^{\imath\theta\sqrt{t/N}}-te^{\imath\theta/\sqrt{t/N}}}\right)^N,$$

for any $\theta\in\mathbb{R}$ and $t\in(0,1)$.

• Let $f(z)=e^z$. In this case $R=\infty$, and the mean and variance functions are $m_f(t)=t$ and $\sigma_f^2(t)=t$. For each $t>0$, the variable $X_t$ in its Khinchin family follows a Poisson distribution with parameter $t$ and its characteristic function is

$$\mathbf{E}(e^{\imath\theta X_t})=e^{t(e^{\imath\theta}-1)}\,,\quad \text{for } \theta\in\mathbb{R} \text{ and } t>0\,,$$

and thus

(1.3.11) $$\mathbf{E}(e^{\imath\theta\breve{X}_t})=\exp\left(t(e^{\imath(\theta/\sqrt{t})}-1-\imath\theta/\sqrt{t})\right),\quad \text{for } \theta\in\mathbb{R} \text{ and } t>0\,.$$

• Let

$$f(z)=\exp(e^z-1)=\sum_{n=0}^{\infty}\frac{\mathcal{B}_n}{n!}\,z^n\,,$$

where $\mathcal{B}_n$ is the $n$-th Bell number which counts the number of partitions of the set $\{1,\dots,n\}$. Then $R=\infty$, and the mean and variance function are $m_f(t)=te^t$ and $\sigma_f^2(t)=t(t+1)e^t$.

The characteristic function is given by

$$\mathbf{E}(e^{\imath\theta X_t})=\exp\left(e^{te^{\imath\theta}}-e^t\right),\quad \text{for } \theta\in\mathbb{R} \text{ and } t>0\,,$$

and thus

(1.3.12) $$\mathbf{E}(e^{\imath\theta\breve{X}_t})=\exp\left(e^{te^{\imath\theta\frac{e^{-t/2}}{\sqrt{t(t+1)}}}}-e^t-\imath\theta\sqrt{t/(t+1)}e^{t/2}\right),\quad \text{for } \theta\in\mathbb{R} \text{ and } t>0\,.$$

• Let $f \in \mathcal{Q}$ be a transcendental entire function of genus 0 with negative zeros, normalized in such a way that $f(0) = 1$, then

$$f(z) = \prod_{j=1}^{\infty} \left(1 + \frac{z}{b_j}\right), \qquad \text{for any } z \in \mathbb{D}.$$

The characteristic function of $(X_t)$, the Khinchin family associated to $f$, is given by

$$\mathbf{E}(e^{i\theta X_t}) = \prod_{j=1}^{\infty} \left(\frac{b_j + te^{i\theta}}{b_j + t}\right), \qquad \text{for any } t \geq 0 \text{ and } \theta \in \mathbb{R}.$$

For the normalized family $(\breve{X}_t)$ we have

$$\mathbf{E}(e^{i\theta \breve{X}_t}) = \left(\prod_{j=1}^{\infty} \left(\frac{b_j + te^{i\theta/\sigma_f(t)}}{b_j + t}\right)\right) e^{-i\theta m_f(t)/\sigma_f(t)}, \qquad \text{for any } t \geq 0 \text{ and } \theta \in \mathbb{R}.$$

**B. Continuous families**

In this section we introduce the definition of a continuous family. We want to take limits, or limits of sequences, of random variables which are part of a Khinchin family, and therefore we need some suitable concept of continuity. We will use this concept for instance to prove that certain differences of characteristic functions converge to zero, see for instance the proof of Theorem 3.2.14, or Chapter 6, where we take subsequences of Khinchin families recurrently. This concept will be used ubiquitously along this document.

A family of random variables $(Z_s)_{s\in[a,b]}$ is said to be *continuous* if it is continuous in distribution in the sense that if a sequence $(s_n)_{n\geq 1}$ with $s_n \in [a, b]$ converges to $s \in [a, b]$ then $Z_{s_n}$ converges to $Z_s$ in distribution.

The family $(Z_s)_{s\in[a,b]}$ is said to be $L^1$-*bounded* in the interval $[a, b]$ if

$$\sup_{t\in[a,b]} \mathbf{E}\left(|Z_t|\right) < +\infty.$$

**Lemma 1.3.2.** *If the family $(Z_s)_{s\in[a,b]}$ is continuous and $L^1$-bounded, then the function*

$$(s, \theta) \in [a, b] \times \mathbb{R} \longrightarrow \mathbf{E}(e^{i\theta Z_s}) \in \mathbb{D}$$

*is continuous.*

*Proof.* Take two sequences $(\theta_n)_{n\geq 1} \subseteq \mathbb{R}$ and $(s_n)_{n\geq 1} \subseteq [a, b]$ which converge, respectively, to $\theta \in \mathbb{R}$ and $s \in [a, b]$. First we write

$$\mathbf{E}(e^{i\theta_n Z_{s_n}}) = \mathbf{E}(e^{i\theta Z_{s_n}}) + \mathbf{E}(e^{i\theta_n Z_{s_n}} - e^{i\theta Z_{s_n}}). \tag{1.3.13}$$

By virtue of Lévy's convergence theorem, the first term at the right-hand side of (1.3.13) converges, as $n \to \infty$, to $\mathbf{E}(e^{\imath\theta Z_s})$. We have the bound

$$\left|\mathbf{E}(e^{i\theta_n Z_{s_n}} - e^{i\theta Z_{s_n}})\right| \le |\theta_n - \theta|\, \mathbf{E}(|Z_{s_n}|), \quad \text{for } n \ge 1\,,$$

which follows from the numerical inequality: $|e^{ix} - e^{iy}| \le |x - y|$, for all $x, y \in \mathbb{R}$.

Now using that $(Z_s)_{s\in[a,b]}$ is $L^1$-bounded, we conclude that

$$(s, \theta) \in [a, b] \times \mathbb{R} \longrightarrow \mathbf{E}(e^{\imath\theta Z_s}) \in \mathbb{D}\,,$$

is a continuous function. □

If $M_f < +\infty$ (and $R < +\infty$), and if, besides, $\sigma_f(R) < \infty$ (see Section 1.3.3), then $(\breve{X}_t)$ is continuous and $L^1$-bounded for $t \in [a, R]$, (including $t = R$) for any $a \in (0, R)$.

The (un-normalized) Khinchin family $(X_t)$, associated to any $f$ in $\mathcal{K}$ with radius of convergence $R > 0$, is continuous and $L^1$-bounded on any interval $[0, b]$ with $0 < b < R$. If $M_f < +\infty$ (and $R < +\infty$) the family $(X_t)_{t\in[0,R]}$ extended to the closed interval $[0, R]$ is continuous and bounded; observe that $\mathbf{E}(|X_t|) = \mathbf{E}(X_t) \le M_f$, for $t \in [0, R]$.

### C. Continuity of the normalized family $(\breve{X}_t)$

Let $f \in \mathcal{K}$ be a power series with radius of convergence $R > 0$ and denote $(\breve{X}_t)$ the normalized family associated to $f$.

For any $a > 0$ we have, applying Hölder's inequality, that

$$\sup_{t\in[a,R)} \mathbf{E}(|\breve{X}_t|) \le \sup_{t\in[a,R)} \mathbf{E}(\breve{X}_t^2)^{1/2} = 1, \tag{1.3.14}$$

then, for any $a > 0$, the family $(\breve{X}_t)_{t\in[a,R)}$ is $L^1$-bounded.

**Proposition 1.3.3.** *Assume that $\breve{X}_t$ converges in distribution, as $t \uparrow R$, towards certain random variable $Z$, then*

$$\lim_{t\uparrow R} \mathbf{E}(\breve{X}_t) = \mathbf{E}(Z) = 0,$$

*and*

$$\lim_{t\uparrow R} \mathbf{E}(|\breve{X}_t|^\alpha) = \mathbf{E}(|Z|^\alpha), \quad \textit{for any } \alpha \in [0, 2).$$

*Proof.* Inequality (1.3.14) gives that for any $a > 0$ and any $\alpha \in [0, 2)$, the family $(\breve{X}_t)_{t\in[a,R)}$ and the family $(|\breve{X}_t|^\alpha)_{t\in[a,R)}$ are uniformly integrable, see [10, pp. 30-32].

Combining the uniform integrability of the families with the convergence in distribution we obtain that

$$\lim_{t\uparrow R} \mathbf{E}(\breve{X}_t) = \mathbf{E}(Z) = 0,$$

and

$$\lim_{t\uparrow R} \mathbf{E}(|\breve{X}_t|^\alpha) = \mathbf{E}(|Z|^\alpha), \quad \text{for any } \alpha \in [0,2).$$

□

Now we prove that the characteristic function of a normalized Khinchin family is continous as a function of two variables. The following lemma will be useful in subsequent chapters when proving that certain, differences, or simply characteristic functions of normalized Khinchin families evaluated at $\theta = \theta(t)$, a function of $t$, for instance with $\theta(t) \to 1$ or $\theta(t) \to 0$, as $t \uparrow R$, converge, as $t \uparrow R$, to certain specific value.

**Lemma 1.3.4.** *Assume that $\breve{X}_t$ converges in distribution, as $t \uparrow R$, towards certain random variable $Z$, then, for any $a > 0$, the function*

$$(s, \theta) \in [a, R) \times \mathbb{R} \longrightarrow \mathbf{E}(e^{\imath\theta \breve{X}_s}) \in \mathbb{D}$$

*extends to a continuous function on $[a, R] \times \mathbb{R}$.*

*Proof.* Fix $a > 0$. If we take any closed subinterval $[a, b] \subseteq (0, R)$ this function is continuous in $[a, b]$, see Lemma 1.3.2.

Take a sequence $(s_n)_{n\geq 1} \subseteq [a, R)$ which converges to $R$, as $n \to \infty$, and any sequence $(\theta_n)_{n\geq 1} \subseteq \mathbb{R}$, which converges to $\theta$, as $n \to \infty$, then we have

$$\mathbf{E}(e^{i\theta_n \breve{X}_{s_n}}) = \mathbf{E}(e^{i\theta \breve{X}_{s_n}}) + \mathbf{E}(e^{i\theta_n \breve{X}_{s_n}} - e^{\imath\theta \breve{X}_{s_n}}).$$

Convergence in distribution implies that the first term at the right-hand side converges to $\mathbf{E}(e^{i\theta X})$, as $n \to +\infty$.

Now observe that

$$\left|\mathbf{E}(e^{\imath\theta_n \breve{X}_{s_n}} - e^{\imath\theta \breve{X}_{s_n}})\right| \leq |\theta_n - \theta|\, \mathbf{E}(|\breve{X}_{s_n}|) \leq |\theta_n - \theta|\,, \quad \text{for } n \geq 1\,,$$

which follows combining the inequality $|e^{ix} - e^{iy}| \leq |x - y|$, for all $x, y \in \mathbb{R}$, with the upper bound (1.3.14).

Thus

$$\lim_{n\to\infty} \mathbf{E}(e^{i\theta_n \breve{X}_{s_n}}) = \mathbf{E}(e^{i\theta Z})\,.$$

□

### 1.3.3 Extension of the Khinchin family at $t = R$ if $M_f < \infty$ (and $R < \infty$)

We assume now that $M_f < \infty$ (and $R < +\infty$). Thus we have that $\sum_{n=0}^{\infty} n a_n R^n < \infty$ and also that $f(R) = \sum_{n=0}^{\infty} a_n R^n < \infty$. The power series $\sum_{n=0}^{\infty} a_n z^n$ defines a continuous (actually, $C^1$) function on the whole closed disk $\text{cl}(\mathbb{D}(0,R))$ which extends $f$ from the open disk $\mathbb{D}(0,R)$. We let $f(R) \triangleq \sum_{n=0}^{\infty} a_n R^n = \lim_{t\uparrow R} f(t)$ and define

$$f'(R) \triangleq \sum_{n=1}^{\infty} n a_n R^{n-1} = \lim_{t\uparrow R} f'(t),$$
$$f''(R) \triangleq \sum_{n=2}^{\infty} n(n-1) a_n R^{n-2} = \lim_{t\uparrow R} f''(t).$$

The second derivative above could be $\infty$. In the present case, we may extend the Khinchin family $(X_t)_{t\in[0,R)}$ of $f$ to $t \in [0,R]$ by defining the variable $X_R$ by

$$\mathbf{P}(X_R = n) = \frac{a_n R^n}{f(R)}, \quad \text{for each } n \geq 0\,.$$

*The extended family* $(X_t)_{t\in[0,R]}$ *becomes continuous in distribution in the closed interval* $[0,R]$. Observe that $X_R$ (like any other $X_t$, with $t \in (0,R)$) is non-constant.

The variable $X_R$ has (finite) mean $\mathbf{E}(X_R) = Rf'(R)/f(R) \triangleq m_f(R)$ and variance

$$\mathbf{V}(X_R) = \frac{\sum_{n=0}^{\infty} n^2 a_n R^n}{f(R)} - \mathbf{E}(X_R)^2\,.$$

The variance $\mathbf{V}(X_R)$ is nonzero since $X_R$ is non-constant, but it could be infinite. Actually, $\mathbf{V}(X_R)$ is finite if and only if $\sum_{n=0}^{\infty} n^2 a_n R^n < +\infty$ if and only if $f''(R) < +\infty$, and in any case

$$\mathbf{V}(X_R) = \lim_{t\uparrow R} \sigma_f^2(t)\,.$$

If $\mathbf{V}(X_R)$ is finite, we write $\sigma_f^2(R) \triangleq \mathbf{V}(X_R)$.

If $M_f < \infty$ (and $R < \infty$), we may also extend the normalized family $(\breve{X}_t)$ by defining $\breve{X}_R = (X_R - m_f(R))/\sigma_f(R)$, with the understanding that $\breve{X}_R \equiv 0$, if $\sigma_f(R) = +\infty$. Recall that $\sigma_f(R) > 0$. The extended family $(\breve{X}_t)_{t\in(0,R]}$ of normalized variables is continuous (in distribution).

### 1.3.4 Fulcrum $F$ of $f$

A power series $f$ in $\mathcal{K}$ does not vanish on the real interval $[0,R)$, and so, it does not vanish in a simply connected region containing that interval, or half-line. We may consider $\ln f$, a branch of the logarithm of $f$ which is real on $[0,R)$, and the function $F$, called the *fulcrum* of $f$, which is defined and holomorphic in a region containing $(-\infty, \ln R)$ and it is given by

$$F(z) = \ln f(e^z)\,.$$

In terms of the fulcrum $F$ of $f$, the mean and variance functions of the Khinchin family associated to $f$ may be expresses as follows:

$$m_f(e^s) = F'(s) \quad \text{and} \quad \sigma_f^2(e^s) = F''(s), \quad \text{for } s < \ln(R).$$

If $f$ *does not vanish anywhere in the disk* $\mathbb{D}(0, R)$, then the fulcrum $F(z)$ of $f$ is defined everywhere in the whole half plane $\Re z < \ln R$.

### Zero-free region

The following proposition describes the specific region where $f$ is nonzero in terms of the argument and the variance of the Khinchin family associated to $f$. This region contains the interval $[0, R)$.

**Proposition 1.3.5.** *Let $f \in \mathcal{K}$ have radius of convergence $R > 0$. If for some $t \in [0, R)$ and some $\theta \in [-\pi, \pi]$, we have $f(te^{\imath\theta}) = 0$, then*

$$|\theta|\sigma_f(t) \ge \frac{\pi}{2}\,.$$

*Thus, $f \in \mathcal{K}$ does not vanish in*

$$\Omega_f = \Big\{z = te^{\imath\theta} : |\theta| < \frac{\pi}{2\sigma_f(t)}\Big\}\,.$$

*Proof.* For the proof we may use the following lemma.

**Lemma 1.3.6** (Sakovič)**.** *Let $Y$ be a random variable and $\theta \in \mathbb{R}$. If $\mathbf{E}(e^{\imath\theta Y}) = 0$, then*

$$\theta^2 \mathbf{V}(Y) \ge \frac{\pi^2}{4}\,.$$

Proposition 1.3.5 follows from Lemma 1.3.6 and from observing that if $f(te^{\imath\theta}) = 0$, then $\mathbf{E}(e^{\imath\theta X_t}) = f(te^{\imath\theta})/f(t) = 0$. □

**Remark 1.3.7** (Alternative proof of Proposition 1.3.5)**.** Alternatively, to verify Proposition 1.3.5 we may use Lemma 1 of [1], which gives that for $f \in \mathcal{K}$

$$f(t)^2 - |f(te^{\imath\theta})|^2 \le 4\sin^2(\theta/2) f(t)^2\sigma_f^2(t)\,, \quad \text{for } t \in (0, R) \text{ and } \theta \in [-\pi, \pi]\,.$$

If $f(te^{\imath\theta}) = 0$, then $1 \le 2|\sin(\theta/2)|\sigma_f(t)$ and thus $|\theta|\sigma_f(t) \ge 1$. This gives a weaker result with the constant $\pi/2$ replaced by 1. Lemma 1 of [1] is stated for entire functions, but it is valid for general $f \in \mathcal{K}$. ⊠

*A proof of Sakovič's Lemma 1.3.6.* The bound on Lemma 1.3.6 is sharp. Simply consider the random variable $Z$ which takes values $\pm 1$ with probability $1/2$; then $\mathbf{V}(Z) = 1$ and $\mathbf{E}(e^{\imath\theta Z}) = \cos\theta$, which vanishes at $\pi/2$.

The result of Sakovič appeared in [87]. As presented in the more accessible reference [84], Lemma 1.3.6 follows most ingenuously by considering the function

$$\varphi(t) = t^2 - 1 + \frac{4}{\pi}\cos\frac{\pi t}{2}\,,$$

which happens to be positive for every $t \in \mathbb{R}$, but for $t = \pm 1$. (A misprinted sign in the definition of $\varphi$ has been corrected.) Assume that $\mathbf{E}(e^{\imath\theta Y}) = 0$. Consider $W = Y - \mathbf{E}(Y)$. We have that $\mathbf{E}(e^{\imath\theta W}) = 0$, and that $\Re\mathbf{E}(e^{\imath\theta W}) = \mathbf{E}(\cos(\theta W)) = 0$ and

$$0 \le \mathbf{E}\Big(\varphi(\theta \frac{2}{\pi} W)\Big) = \theta^2 \frac{4}{\pi^2}\mathbf{E}(W^2) - 1 = \theta^2 \frac{4}{\pi^2}\mathbf{V}(Y) - 1 .$$

Actually, equality in the bound of Lemma 1.3.6 only happens for the simple symmetric variable $Z$ given above.

A more direct proof of Lemma 1.3.6, but with a weaker constant, appears in [5, Prop. 7.8] (see also [17, Lemma 2.3]).

□

**Remark 1.3.8.** *On Bŏĭchuk-Gol'berg's Theorem 1.2.3.* From Sakovič's Lemma 1.3.6 we may deduce Bŏĭchuk-Gol'berg's Theorem 1.2.3 in the weaker form that

$$\sup_{t>0} \sigma_f^2(t) \ge 1/4 ,$$

for any entire function $f \in \mathcal{K}$.

We may assume that $\sup_{t>0} \sigma_f^2(t) < +\infty$. Since $t m_f'(t) = \sigma_f^2(t)$, integrating, we deduce then that $m_f(t) = O(\ln t)$, as $t \to \infty$. A further integration shows that $\ln f(t) = O((\ln t)^2)$, as $t \to \infty$. Thus since $f \in \mathcal{K}$, we have that

$$\ln \max_{|z|\le t} |f(z)| = \ln f(t) = O((\ln t)^2) , \quad \text{as } t \to \infty .$$

The entire function $f$ is then of order 0 and Hadamard's factorization theorem gives that $f$ is an infinite canonical product or a polynomial (non-constant and not a monomial, since $f \in \mathcal{K}$). (This is the starting point of Hayman's proof of Theorem 3 of [49].) In any case, $f$ vanishes at some $z_0 \ne 0$. Write $z_0 = r_0 e^{\imath\theta_0}$, with $r_0 > 0$ and $|\theta_0| \le \pi$. Lemma 1.3.6 gives that

$$\pi \sigma_f(r_0) \ge \frac{\pi}{2} ,$$

and thus, that $\sigma_f(r_0) \ge 1/2$, and, in particular, that $\sup_{t>0} \sigma_f^2(t) \ge 1/4$.

This reasoning also gives that if $f$ is not a polynomial then

$$\limsup_{t\to\infty} \sigma_f^2(t) \ge 1/4 .$$

⊠

### 1.3.5 Characteristic functions and moments of $\breve{X}_t$

Let $f \in \mathcal{K}$ be a power series with radius of convergence $R > 0$ and Khinchin family $(X_t)$. Fix $t = e^s$, for $s < \ln(R)$. Proposition 1.3.5 gives that the characteristic function of $\breve{X}_t$ is never zero for $|\theta| < \pi/2$.

This bound also gives that the function $F(s+i\theta/\sigma_f(t))$ is well defined for any $|\theta|<\pi/2$. This follows from the fact that $f(e^s e^{i\theta/\sigma_f(t)})$ is never zero for any fixed $s<\ln(R)$ and $|\theta|<\pi/2$, see Proposition 1.3.5.

For $t=e^s$, where $s<\ln(R)$, fixed, and $|\theta|<\pi/2$, we have the relation

$$\varphi_{\breve{X}_t}(\theta)=\mathbf{E}(e^{i\theta\breve{X}_t})=e^{F(s+i\theta/\sigma_f(t))-F(s)-\frac{F'(s)}{F''(s)^{1/2}}i\theta}. \tag{1.3.15}$$

This is a $C^\infty$ function (with respect to $\theta$).

Also using that $F(s+i\theta/\sigma_f(t))-F(s)-\frac{F'(s)}{F''(s)^{1/2}}i\theta$ is well defined for any fixed $s<\ln(R)$ and $|\theta|<\pi/2$, we obtain that

$$\ln(\mathbf{E}(e^{i\theta\breve{X}_t}))=F(s+i\theta/\sigma_f(t))-F(s)-\frac{F'(s)}{F''(s)^{1/2}}i\theta. \tag{1.3.16}$$

Again we have a $C^\infty$ function with respect to $\theta$.

We want to compute the derivatives of any order of (1.3.15) and (1.3.16) at $\theta=0$. A direct application of Faà di Bruno's formula, see, for instance, [100, p. 217], to (1.3.15) and (1.3.16), respectively, gives that

**Lemma 1.3.9.** *Let $f\in\mathcal{K}$ be a power series with radius of convergence $R>0$ and denote $(X_t)$ the Khinchin family associated to $f$, then, for any $t=e^s\in(0,R)$ and any $n\geq 3$, we have*

$$\mathbf{E}(\breve{X}_t^n)=\sum_{2m_2+\cdots+nm_n=n}\frac{n!}{m_2!\dots m_n!}\prod_{j=1}^n\left(\frac{F^{(j)}(s)}{F''(s)^{j/2}}\frac{1}{j!}\right)^{m_j}, \tag{1.3.17}$$

*and*

$$\frac{F^{(n)}(s)}{F''(s)^{n/2}}=\sum_{2m_2+\dots nm_n=n}\frac{n!(l_n-1)!(-1)^{l_n+1}}{m_2!\dots m_n!}\prod_{j=1}^n\left(\frac{\mathbf{E}(\breve{X}_t^j)}{j!}\right)^{m_j}, \tag{1.3.18}$$

*here $l_n=m_2+\cdots+m_n$.*

**Remark 1.3.10** (Hidden moments)**.** As we are going to see the moments of the normal distribution are included in the sum (1.3.17). Fix $n=2k$, where $k\geq 1$ is an integer. Formula (1.3.17) can be written as

$$\mathbf{E}(\breve{X}_t^n)=\sum_{\substack{2m_2+\cdots+nm_n=n\\ m_2\neq k}}\frac{n!}{m_2!\dots m_n!}\prod_{j=1}^n\left(\frac{F^{(j)}(s)}{F''(s)^{j/2}}\frac{1}{j!}\right)^{m_j}+\frac{(2k)!}{k!\,(2!)^k}.$$

If $n$ is odd all the terms in the sum have at least one term of the form $(F^{(j)}(s)/F''(s)^{j/2})^{m_j}$ as a factor, here $1\leq j\leq n$ with $j\neq 2$ and $m_j>0$. Observe that, for $n$ odd, there always exists at least one $m_j>0$, with $1\leq j\leq n$ odd, which is positive, otherwise $n$ would be even. ⊠

**Remark 1.3.11** (The characteristic function of the standard normal)**.** For comparison and later use we apply Faà di Bruno's formula to the characteristic function of the standard normal.

Let $Z$ be a normally distributed random variable with $\mathbf{E}(Z) = 0$ and $\mathbf{V}(Z) = 1$. We have $\mathbf{E}(e^{i\theta Z}) = e^{-\theta^2/2}$, for any $\theta \in \mathbb{R}$, and therefore

$$(\star) \quad g(\theta) = \ln(\mathbf{E}(e^{i\theta Z})) = -\theta^2/2.$$

Applying Faà di Bruno's formula to $(\star)$, we obtain that, for any $n \geq 3$, we have

$$g^{(n)}(0) = 0 = \sum_{2m_2+\dots nm_n=n} \frac{n!(l_n-1)!(-1)^{l_n+1}}{m_1!\dots m_n!} \prod_{j=1}^{n} \left(\frac{\mathbf{E}(Z^j)}{j!}\right)^{m_j},$$

here $l_n = m_2 + \cdots + m_n$. Recall that $Z$ is a standard normal random variable and therefore $\mathbf{E}(Z) = 0$.

Denote $h(\theta) = -\theta^2/2$. For any $n \geq 1$, we have, applying again Faà di Bruno's formula to $\varphi_Z(\theta) = \mathbf{E}(e^{i\theta Z}) = e^{-\theta^2/2} = \exp(h(\theta))$ that

$$\mathbf{E}(Z^n)i^n = \sum_{m_1+2m_2+\cdots+nm_n=n} \frac{n!}{m_1!\dots m_n!} \prod_{j=1}^{n} \left(\frac{h^{(j)}(0)}{j!}\right)^{m_j}. \tag{1.3.19}$$

Notice that $h^{(j)}(0) = 0$, for any $j \neq 2$. In fact $h^{''}(0) = -1$, therefore, as it is well known:

$$\mathbf{E}(Z^n) = 0, \quad \text{for odd } n \geq 1,$$

and

$$\mathbf{E}(Z^{2k}) = \frac{(2k)!}{k!\,(2!)^k}, \quad \text{for any } k \geq 0.$$

⊠

**Remark 1.3.12** (Bell polynomials and moments of the standard normal)**.** There is an alternative form of Faà di Bruno's formula known as Riordan's formula, see, for instance, [100, p. 219]. This equivalent form makes use of Bell polynomials.

Applying Riordan's formula, see [100, p. 219], we obtain an equivalent expression for 1.3.17, namely:

(1.3.20)

$$\mathbf{E}(\breve{X}_t^n)i^n = \sum_{k=0}^{n} B_{n,k}\left(0, i^2, i^3\frac{F'''(s)}{F''(s)^{3/2}}, \dots, i^{n-k+1}\frac{F^{(n-k+1)}(s)}{F''(s)^{(n-k+1)/2}}\right)$$

$$= B_n\left(0, i^2, i^3\frac{F'''(s)}{F''(s)^{3/2}}, \dots, i^n\frac{F^{(n)}(s)}{F''(s)^{n/2}}\right) = i^n B_n\left(0, 1, \frac{F'''(s)}{F''(s)^{3/2}}, \dots, \frac{F^{(n)}(s)}{F''(s)^{n/2}}\right)$$

where $B_{n,k}(x_1, \ldots, x_{n-k+1})$ denotes $n$-th incomplete exponential Bell polynomial, $B_n(x_1, \ldots, x_n)$ denotes the $n$-th Bell polynomial and $n \geq 1$ is an integer.

In (1.3.20) we make use of the following property: the n-th Bell polynomial is given by

$$B_n(x_1, \ldots, x_n) = n! \sum_{j_1+2j_2+\cdots+nj_n=n} \prod_{l=1}^{n} \frac{x_l^{j_l}}{(l!)^{j_l} j_l!},$$

see, for instance, [82, p. 173] or [100, pp. 218-219], therefore

$$(\dagger) \quad B_n(ix_1, i^2x_2, \ldots, i^n x_n) = i^n B_n(x_1, \ldots, x_n).$$

Denote by $Z$ a standard normal random variable. Applying Riordan's formula, see [100, p. 219], we find that

$$\mathbf{E}(Z^n) i^n = \sum_{k=0}^{n} B_{n,k}(0, -1, 0, \ldots, 0) = B_n(0, -1, 0, \ldots, 0) = i^n B_n(0, 1, 0, \ldots, 0), \tag{1.3.21}$$

where, again, $B_{n,k}(x_1, \ldots, x_{n-k+1})$ denotes $n$-th incomplete exponential Bell polynomial and $B_n(x_1, \ldots, x_n)$ denotes the $n$-th Bell polynomial, in each of the cases for integers $n \geq 1$. For the equality at the right-hand side we use, again, the property $(\dagger)$ of Bell polynomials.

Formula (1.3.21) gives that we can express the $n$-th moment of the standard normal in terms of the $n$-th Bell polynomial. In fact we have

$$\mathbf{E}(Z^n) = B_n(0, 1, 0, \ldots, 0), \qquad \text{for any } n \geq 1,$$

and also, by virtue of equation (1.3.20), we have that

$$\mathbf{E}(\breve{X}_t^n) = B_n\left(0, 1, \frac{F'''(s)}{F''(s)^{3/2}}, \ldots, \frac{F^{(n)}(s)}{F''(s)^{n/2}}\right), \qquad \text{for any } n \geq 1.$$

⊠

### 1.3.6 Operations with power series and Khinchin families

In this section we collect results describing the Khinchin families obtained by performing operations with one, two, or more, power series in the class $\mathcal{K}$. We also study the Khinchin families of sequences of functions in $\mathcal{K}$ which converge uniformly in compact subsets to power series in $\mathcal{K}$.

#### A. Subseries and conditional distributions

Let $(X_t)$ be the Khinchin family associated to the power series $f(z) = \sum_{n=0}^{\infty} a_n z^n \in \mathcal{K}$, which has radius of convergence $R > 0$.

Let $\mathcal{L} \subseteq \{0, 1, 2, \ldots\}$ such that $0 \in \mathcal{L}$, and let

$$g(z) = \sum_{n \in \mathcal{L}} a_n z^n, \qquad \text{for any } z \in \mathbb{D}(0, R).$$

Assume that $g \in \mathcal{K}$, i.e, $a_n \neq 0$, for some $n \in \mathcal{L} \setminus \{0\}$. Let $(Y_t)$ be the Khinchin family associated to the power series $g$.

For each $t \in [0, R)$ we have that

$$\mathbf{P}(Y_t = n) = \begin{cases} 0, & \text{if } n \notin \mathcal{L}, \\ \dfrac{\mathbf{P}(X_t = n)}{\mathbf{P}(X_t \in \mathcal{L})} = \mathbf{P}(X_t = n | X_t \in \mathcal{L}), & \text{if } n \in \mathcal{L}. \end{cases}$$

In other words the law of the variable $Y_t$ is the law of $X_t$ conditioned on the event $X_t \in \mathcal{L}$.

### B. Sums of power series and mixtures

Let $f, g$ two power series in $\mathcal{K}$, both with radius of convergence at least $R > 0$. The sum $h(z) = f(z) + g(z)$ is also a power series in $\mathcal{K}$ and has radius of convergence at least $R > 0$.

Let $(X_t)$, $(Y_t)$ and $(Z_t)$ be the Khinchin families of $f, g$ and $h$, respectively. We may write

$$\mathbf{P}(Z_t = n) = \frac{f(t)}{h(t)} \mathbf{P}(X_t = n) + \frac{g(t)}{h(t)} \mathbf{P}(Y_t = n), \quad \text{for } t \in (0, R) \text{ and } n \geq 0.$$

This last expression exhibits the distribution of $Z_t$ as a mixture of the distributions of $X_t$ and $Y_t$. In this case we have the identification

$$Z_t \stackrel{d}{=} U_t X_t \oplus (1 - U_t) Y_t, \quad \text{for any } t \in [0, R),$$

where $U_t$ is a Bernoulli variable with parameter $f(t)/h(t)$ which is independent of $X_t$ and $Y_t$.

### C. Product of power series and independence of families

Let $f(z) = \sum_{n=0}^{\infty} a_n z^n$ and $g(z) = \sum_{n=0}^{\infty} b_n z^n$ be power series in $\mathcal{K}$ with radius of convergence $R > 0$ and $S > 0$. Let $(X_t)_{t \in (0,R)}$ and $(Y_t)_{t \in (0,S)}$ be the Khinchin families associated to $f$ and $g$, respectively. The product $h = fg$ is also in $\mathcal{K}$. Let $T = \min\{R, S\}$, and let $(Z_t)_{t \in (0,T)}$ be the Khinchin family associated to $h$.

The coefficients of the product $h = fg$ are given by the convolution formula

$$c_n = \text{COEFF}_n[h(z)] = \sum_{k=0}^{n} a_k b_{n-k} \, .$$

Then for each $t \in (0, T)$ the law of $Z_t$ is the same that *a sum of independent copies* of $X_t$ and $Y_t$:

$$\mathbf{P}(Z_t = k) = \sum_{k=0}^{n} \frac{a_k t^k}{f(t)} \frac{b_{n-k} t^{n-k}}{g(t)} = \sum_{j=0}^{k} \mathbf{P}(X_t = j) \mathbf{P}(Y_t = k - j) \, , \quad \text{for each } k \geq 0 \, .$$

In particular, if $(X_t)_{t \in (0,R)}$ is the Khinchin family associated with some power series $f$, and if we let $X_t^{(1)}, \ldots, X_t^{(n)}$ denote $n$ independent copies of $X_t$ then

$$\mathbf{P}\Big( \sum_{j=1}^{n} X_t^{(j)} = k \Big) = \text{COEFF}_{[k]}\big(f^n(z)\big) \, \frac{t^k}{f^n(t)} \, , \quad \text{for any } k \geq 0 \text{ and } t \in (0, R) \, .$$

We collect these claims by means of the following Lemma.

**Lemma 1.3.13.** *With the notations above, for each $t \in (0, T)$, the law of $Z_t$ is the same that the law of the sum of independent variables with the laws of $X_t$ and $Y_t$, i.e.,*

$$Z_t \stackrel{d}{=} X_t \oplus Y_t, \quad \textit{for each } t \in (0, R). \tag{1.3.22}$$

Thus, product of power series in $\mathcal{K}$, become, on the Khinchin family side, in sums of independent variables and therefore, for any $t \in [0, T)$, we have

$$m_h(t) = m_f(t) + m_g(t)\,, \quad \text{and} \quad \sigma_h^2(t) = \sigma_f^2(t) + \sigma_g^2(t)\,.$$

Using the identification (1.3.26) we have

$$\mathbf{E}(e^{i\theta \breve{Z}_t}) = \mathbf{E}(e^{i\theta X_t})\mathbf{E}(e^{i\theta \breve{Y}_t}), \quad \text{for any } t \in [0, R).$$

For the characteristic function of the powers of a power series $g \in \mathcal{K}$ we have the following: denote $(X_t)$ and $(Y_t)$ the Khinchin families of $g^N$, with integer $N \geq 1$, and $g$ respectively. We have that $X_t$ is equal, in distribution, to the sum of $N \geq 1$ independent copies of the random variable $Y_t$, that is:

$$X_t \stackrel{d}{=} Y_t^{(1)} + \cdots + Y_t^{(N)}, \quad \text{for any } t \in (0, R),$$

here $Y_t^{(1)}, \ldots, Y_t^{(N)}$ are independent identical copies of $Y_t$.

For any $t \in (0, R)$, the mean and the standard deviation of $X_t$ are given by

$$\mathbf{E}(X_t) = N m_g(t) \quad \text{and} \quad \mathbf{V}(X_t) = N\sigma_g^2(t),$$

then

$$\breve{X}_t = \frac{X_t - N m_g(t)}{\sqrt{N}\sigma_f(t)}, \quad \text{for any } t \in (0, R)\,.$$

The characteristic function of $X_t$ and $\breve{X}_t$ are given by

$$\mathbf{E}(e^{i\theta X_t}) = \mathbf{E}(e^{i\theta Y_t})^N, \quad \text{for } t \in (0, R) \text{ and } \theta \in \mathbb{R}\,, \tag{1.3.23}$$

and

$$\mathbf{E}(e^{i\theta \breve{X}_t}) = \mathbf{E}\left(e^{i\frac{\theta}{\sqrt{N}}\breve{Y}_t}\right)^N, \quad \text{for } t \in (0, R) \text{ and } \theta \in \mathbb{R}\,, \tag{1.3.24}$$

respectively.

An instance of this situation is given by $f(z) = (1+z)^N$, with integer $N \geq 1$. Denote $(X_t)$ and $(Y_t)$ the Khinchin families of $f(z) = (1+z)^N$ and $g(z) = (1+z)$, respectively, then

$$X_t \stackrel{d}{=} Y_t^{(1)} \oplus \cdots \oplus Y_t^{(N)}, \quad \text{for any } t \in (0, +\infty),$$

where $Y_t^{(1)}, \ldots, Y_t^{(N)}$ are independent identically distributed copies of $Y_t$. Here $(Y_t)$ is the Bernoulli family and $(X_t)$ is the binomial family.

Combining the items (1.3.24) and (1.3.5) we find that

$$\mathbf{E}(e^{i\theta \breve{X}_t}) = \left(\frac{t e^{\imath\theta/\sqrt{Nt}} + e^{-\imath\theta\sqrt{t/N}}}{1+t}\right)^N, \quad \text{for } \theta \in \mathbb{R} \text{ and } t > 0\,.$$

this argument proves the relation (1.3.7). Using a similar argument, but in this case with $g(z) = 1/(1-z)$ and $f(z) = 1/(1-z)^N$, with integer $N \geq 1$, we obtain equation (1.3.10).

**D. Composition of power series.**

Let $f(z) = \sum_{n=0}^{\infty} a_n z^n \in \mathcal{K}$ a power series with radius of convergence $R > 0$ and $g(z) = \sum_{n=0}^{\infty} b_n z^n \in \mathcal{K}_s$ a power series with radius of convergence $S > 0$ and such that $g(0) = 0$.

Assume that $g(\mathbb{D}(0,S)) \subset \mathbb{D}(0,R)$, then the composition $h(z) = f(g(z))$ is well defined. Denote $(Z_t)_{t\in[0,S)}$ the Khinchin family associated to $h$ and $(Y_t)_{t\in[0,R)}$ and $(X_t)_{t\in[0,S)}$ the Khinchin families associated to $f$ and $g$, respectively. For $f$ and $g$ entire functions the composition it is always well defined.

For any integer $n \geq 1$ and any $0 \leq j \leq n$, we have that

$$b_n^{(j)} = \text{COEFF}_n[g(z)^j] = \sum_{l_1+\cdots+l_j=n} b_{l_1} b_{l_2} \ldots b_{l_j},$$

then we have

$$(\star) \quad A_n = \text{COEFF}_n[h(z)] = \sum_{j=0}^{n} a_j b_n^{(j)} = \sum_{j=0}^{n} a_j \sum_{l_1+\cdots+l_j=n} b_{l_1} b_{l_2} \ldots b_{l_j}, \quad \text{for any } n > 0\,,$$

and $A_0 = a_0$, for $n = 0$.

For any $n \geq 0$, and any $t \in [0,S)$, using $(\star)$, we find that

$$\begin{aligned}\mathbf{P}(Z_t = n) = \frac{A_n t^n}{f(g(t))} &= \sum_{j=0}^{n} \frac{a_j g(t)^j}{f(g(t))} \sum_{l_1+\cdots+l_j=n} \prod_{m=1}^{j} \frac{b_{l_m} t^{l_m}}{g(t)} \\ &= \sum_{j=0}^{n} \mathbf{P}(Y_{g(t)} = j)\mathbf{P}(X_t^{(1)} \oplus \cdots \oplus X_t^{(j)} = n)\end{aligned}$$

then we can identify $Z_t$ with the random variable $X_t^{(1)} \oplus \cdots \oplus X_t^{(Y_{g(t)})}$, that is,

$$Z_t \stackrel{d}{=} X_t^{(1)} \oplus \cdots \oplus X_t^{(Y_{g(t)})}, \quad \text{for any } t \in [0,R), \tag{1.3.25}$$

where $X_t^{(1)}, X_t^{(2)} \ldots$ are a countable number of i.i.d. copies of $X_t$ which are independent of $Y_{g(t)}$.

We fix the convention that the empty sum is zero and then for $n = 0$ we have

$$\mathbf{P}(Z_t = 0) = \mathbf{P}(Y_{g(t)} = 0), \quad \text{for any } t \in [0, S).$$

We collect the previous discussion by means of the following lemma.

**Lemma 1.3.14.** *With the notations above, for each $t \in (0, T)$, the law of $Z_t$ is the same that the law of the sum of $Y_{g(t)}$ independent variables with the law of $X_t$, i.e.,*

$$Z_t \stackrel{d}{=} X_t^{(1)} \oplus \cdots \oplus X_t^{(Y_{g(t)})}, \quad \textit{for any } t \in [0, R), \tag{1.3.26}$$

For $h(z) = f(g(z))$ the fulcrum is given by the composition of the fulcrums of $g$ and $f$, respectively. In certain domain containing $[0, \ln(R))$ we have

$$H(z) = \ln(h(e^z)) = \ln(f(g(e^z))) = F(G(z)).$$

The previous composition formula gives that

$$m_h(t) = m_f(g(t)) m_g(t), \quad \text{for any } t \in [0, R), \tag{1.3.27}$$

and also that

$$\sigma_h^2(t) = \sigma_f^2(g(t)) m_g(t)^2 + \sigma_g^2(t) m_f(g(t)), \quad \text{for any } t \in [0, R). \tag{1.3.28}$$

This follows, as well, from the identification in distribution (1.3.25).

For the characteristic function of the Khinchin family associated to $h = f \circ g$ we have the following Lemma.

**Lemma 1.3.15.** *With the same notations:*

$$\mathbf{E}(e^{i\theta Z_t}) = \sum_{n=0}^{\infty} \mathbf{P}(Y_{g(t)} = n) \mathbf{E}(e^{i\theta X_t})^n, \quad \textit{for any } \theta \in \mathbb{R} \textit{ and } t \in [0, S).$$

*Proof.* Denote $f(z) = \sum_{n=0}^{\infty} a_n z^n$ the power series expansion of $f$ around $z = 0$, then

$$\mathbf{E}(e^{i\theta Z_t}) = \frac{f(g(te^{i\theta}))}{f(g(t))} = \sum_{n=0}^{\infty} \frac{g(te^{i\theta})^n}{g(t)^n} \frac{a_n g(t)^n}{f(g(t))} = \sum_{n=0}^{\infty} \mathbf{P}(Y_{g(t)} = n) \mathbf{E}(e^{i\theta X_t})^n,$$

for any $\theta \in \mathbb{R}$ and $t \in [0, S)$ □

**E. Scaling.**

For any power series $g(z)=\sum_{n=0}^{\infty} a_n z^n$ in $\mathcal{K}$ of radius $R>0$, we denote

$$Q_g \triangleq \gcd\{n\geq 1 : a_n\neq 0\} = \lim_{N\to\infty} \gcd\{1\leq n\leq N : a_n\neq 0\}\,,$$

that is, if $Q_g>1$, then only the coefficients with indices which are multiple of $Q_g$ are different from 0.

If $Q_g>1$, then we can write $g(z)=f(z^{Q_g})$ for a certain companion power series $f\in\mathcal{K}$ with radius of convergence $R^{Q_g}$. Observe that $Q_f=1$.

Let $(Y_t)_{t\in[0,R)}$ be the Khinchin family of $g$ and $(X_t)_{t\in[0,R^{Q_g})}$ be the Khinchin family associated to the power series $f$. We have that

$$Y_t \stackrel{d}{=} Q_g\cdot X_{t^{Q_g}}, \qquad \text{for any } t\in(0,R)).$$

The mean and variance functions of $g$ and $f$ are related by

$$m_g(t)=Q_g\cdot m_f(t^{Q_g})\,,$$
$$\sigma_g^2(t)=Q_g^2\cdot \sigma_f^2(t^{Q_g})\,,$$

for any $t\in[0,R)$.

**F. Identity of families**

If $\lambda>0$ and $f\in\mathcal{K}$, then $g\equiv\lambda f$ is also in $\mathcal{K}$, both power series have the same radius of convergence, say $R>0$, and the respective Khinchin families $(X_t)$ and $(Y_t)$ of $f$ and $g$ coincide in distribution: $X_t\stackrel{d}{=}Y_t$, for $t\in[0,R)$.

Conversely, assume that $f(z)=\sum_{n=0}^{\infty} a_n z^n$ and $g(z)=\sum_{n=0}^{\infty} b_n z^n$ are power series in $\mathcal{K}$ with radius of convergence $R>0$ and $S>0$, respectively. Denote $(X_t)$ and $(Y_t)$ its respective Khinchin families and suppose that there exists some $t_0\in(0,\min(R,S))$ such that $X_{t_0}\stackrel{d}{=}Y_{t_0}$, then

$$\frac{a_n t_0^n}{f(t_0)}=\frac{b_n t_0^n}{g(t_0)}, \qquad \text{for each } n\geq 0,$$

and, therefore, if we set $\lambda\triangleq g(t_0)/f(t_0)$, we have that

$$b_n=\lambda a_n, \qquad \text{for each } n\geq 0.$$

This last equality implies that $R=S$ and $g\equiv\lambda f$ in $\mathbb{D}(0,R)$.

### G. Sequences of power series

Let $(f^{[k]})_{k\ge 1}$ be a sequence of power series in $\mathcal{K}$, each of them with radius of convergence at least $R>0$ and let $(X_t^{[k]})$ denote the corresponding sequence of Khinchin families.

Assume that the sequence of power series $(f^{[k]})_{k\ge 1}$ converges uniformly on compact subsets of $\mathbb{D}(0,R)$ towards a power series $f$ in $\mathcal{K}$.

The condition $f\in\mathcal{K}$ ensures that we can associate a proper Khinchin family to $f$. The sequence of power series $f^{[k]}(z)=1/k+z^k$ is in $\mathcal{K}$, for any integer $k\ge 1$, and converges uniformly on compact subsets of $\mathbb{D}$ to the constant 0, which is not in $\mathcal{K}$.

Let $(Y_t)$ be the Khinchin family associated to the power series $f\in\mathcal{K}$. For each $n\ge 0$, the $n$-th coefficient of $f^{[k]}$ converges to the $n$-th coefficient of $f$, and, therefore,

$$\lim_{k\to+\infty}\mathbf{P}(X_t^{[k]}=n)=\mathbf{P}(Y_t=n),\qquad \text{for each } n\ge 0.$$

Consequently, for each $t\in(0,R)$ we have that $X_t^{[k]}$ converges in distribution towards the random variables $Y_t$:

$$\lim_{k\to+\infty}\mathbf{P}(X_t^{[k]}\le N)=\mathbf{P}(Y_t\le N),\qquad \text{for each } N\ge 1.$$

Also, for any moment, say, of order $q\ge 1$ we have

$$\lim_{k\to+\infty}\mathbf{E}\left(\left(X_t^{[k]}\right)^q\right)=\mathbf{E}(X_t^q),$$

since if $\mathcal{D}$ is the operator $\mathcal{D}=z(d/dz)$, then $\mathbf{E}(X_t^q)=\frac{\mathcal{D}^q(f)(t)}{f(t)}$, see equation (1.3.1) for further details.

In particular, for each $t\in[0,R)$, we have that $\lim_{k\to+\infty}m_k(t)=m_f(t)$ and $\lim_{k\to+\infty}\sigma_k^2(t)=\sigma_f^2(t)$, where $m_k$, $\sigma_k^2$ and $m_f(t)$ and $\sigma_f^2(t)$ are, respectively, the mean and variance functions of $f^{[k]}$ and $f$.

Conversely, assume that $f_k(z)=\sum_{n=0}^\infty a_n^{[k]}z^n$ and $f(z)=\sum_{n=0}^\infty a_n z^n$ are power series in $\mathcal{K}$, all of them with radius of convergence at least $R>0$.

Assume, with the notations above, that $X_t^{[k]}$ converges in distribution towards $Y_t$, for each $t\in(0,R)$, as $k\to+\infty$, and besides that for a certain $n_0\ge 0$ we have that

$$\lim_{k\to+\infty}a_{n_0}^{[k]}=a_{n_0},\qquad \text{and}\qquad a_{n_0}\neq 0.$$

Then the sequence $f^{[k]}$ converges uniformly on compact subsets of $\mathbb{D}(0,R)$ towards $g(z)$.

Let's see: from the convergence in distribution we deduce that

$$\lim_{k\to+\infty}\frac{a_{n_0}^{[k]}t^{n_0}}{f^{[k]}(t)}=\frac{a_{n_0}t^{n_0}}{f(t)},\qquad \text{for each } t\in(0,R),$$

and, therefore, that

$$(\star) \qquad \lim_{k\to+\infty} f^{[k]}(t) = f(t), \qquad \text{for each } t \in (0,R).$$

Fix $T \in (0,R)$. From $(\star)$ we deduce that

$$\sup_{k\geq 1} f^{[k]}(T) < +\infty.$$

Since the $f^{[k]}$ has non-negative Taylor coefficients, we have that

$$\sup_{k\geq 1} \sup_{|z|\leq T} |f^{[k]}(z)| < +\infty$$

Thus, by virtue of Montel's Theorem, the family $\{f^{[k]}\}_{k\geq 1}$ conforms a normal family in $\mathbb{D}(0,R)$. Since $\lim_{k\to+\infty} f^{[k]}(t) = f(t)$, for any $t \in (0,R)$, Vitali-Porter's Theorem allows us to conclude that the sequence $(f^{[k]})_{k\geq 1}$ converges uniformly on compact subsets of $\mathbb{D}(0,R)$ towards $f$.

### 1.3.7 Hayman's identity

Let $f(z) = \sum_{n=0}^{\infty} a_n z^n$ be a power series in $\mathcal{K}$. Cauchy's formula for the coefficient $a_n$ in terms of the characteristic function of its Khinchin family $(X_t)_{t\in[0,R)}$ reads

$$(1.3.29) \qquad a_n = \frac{f(t)}{2\pi t^n} \int_{|\theta|<\pi} \mathbf{E}(e^{\imath\theta X_t}) e^{-\imath\theta n}\, d\theta\,, \quad \text{for each } t \in (0,R) \text{ and } n \geq 1\,.$$

In terms of the characteristic function of the normalized variable $\breve{X}_t$, it becomes

$$(1.3.30) \qquad a_n = \frac{f(t)}{2\pi t^n \sigma_f(t)} \int_{|\theta|<\pi\sigma_f(t)} \mathbf{E}(e^{\imath\theta \breve{X}_t}) e^{-\imath\theta(n-m_f(t))/\sigma_f(t)}\, d\theta\,,$$

for each $t \in (0,R)$ and $n \geq 1$.

• If $M_f = \infty$, we may take for each $n \geq 1$ the (unique) radius $t_n \in (0,R)$ so that $m_f(t_n) = n$, to write

$$(1.3.31) \qquad a_n = \frac{f(t_n)}{2\pi t_n^n \sigma_f(t_n)} \int_{|\theta|<\pi\sigma_f(t_n)} \mathbf{E}(e^{\imath\theta \breve{X}_{t_n}})\, d\theta\,, \quad \text{for each } n \geq 1\,,$$

which we call *Hayman's identity.* This formula and its variants will be very relevant later on.

Although this identity (1.3.31) is just Cauchy's formula with an appropriate choice of radius, it encapsulates the saddle point method.

As we now check, if the Khinchin family $(X_t)$ may be extended to the closed interval $[0,R]$, then an analogous Hayman's identity is available for $t = R$.

• If $M_f < \infty$ and ($R < \infty$) the power series $f$ extends continuously to the closed disk $\text{cl}(\mathbb{D}(0, R))$, see Section 1.3.3, and the Khinchin family $(X_t)_{t\in[0,R)}$ extends to the closed interval $[0, R]$. See Section 1.2.4 and the notations therein.

We may write

$$a_n = \frac{f(R)}{2\pi R^n} \int_{|\theta|<\pi} \mathbf{E}(e^{\imath\theta X_R}) e^{-\imath\theta n}\, d\theta\,, \quad \text{for } n \geq 1\,.$$

If, moreover, $\sigma_f^2(R) < \infty$, i.e., if $\sum_{n=0}^{\infty} n^2 a_n R^n < \infty$, then $\breve{X}_R$ given by

$$\breve{X}_R = \frac{X_R - m_f(R)}{\sigma_f(R)}$$

is well defined and we also have that

$$a_n = \frac{f(R)}{2\pi R^n \sigma_f(R)} \int_{|\theta|<\pi\sigma_f(R)} \mathbf{E}(e^{\imath\theta \breve{X}_R})\, e^{-\imath\theta(n-m_f(R))/\sigma_f(R)}\, d\theta\,, \quad \text{for } n \geq 1\,. \tag{1.3.32}$$

## 1.4 Moment generating functions and Khinchin families

In this section we discuss the moment generating functions (centered or not) of Khinchin families, and also their close connection with the fulcrum of the family.

Firstly, we consider the moment generating functions of (quite) general random variables $Y$. Later on, we will apply this general concept of moment generating function to the random variables in a Khinchin family.

### 1.4.1 Moment generating functions for general random variables

Let $Y$ be a random variable that it has finite absolute moments of all orders, and we assume, in fact, that for some $\rho > 0$, we have that

$$\sum_{k=0}^{\infty} \frac{\mathbf{E}(|Y|^k)}{k!} |z|^k\,, \quad \text{for } |z| < \rho\,. \tag{1.4.1}$$

This means, in particular, that the power series

$$(\star) \quad \sum_{k=0}^{\infty} \frac{\mathbf{E}(Y^k)}{k!} z^k$$

has positive radius of convergence, at least $\rho$.

Write $m_Y = \mathbf{E}(Y)$ and $\nu_k(Y) = \mathbf{E}((Y - m_Y)^k))$, for each $k \geq 0$. Observe that $\nu_0 = 1$, $\nu_1 = 0$ and $\nu_2 = \mathbf{V}(Y)$.

The *centered moment generating function* of $Y$, defined (at least) for $\lambda \in \mathbb{D}(0, \rho)$, is given by:

$$\lambda \in \mathbb{D}(0, \rho) \mapsto \mathbf{E}(e^{\lambda(Y-m_Y)}) = \sum_{k=0}^{\infty} \frac{\nu_k(Y)}{k!} \lambda^k\,.$$

Also, the *moment generating function* of $Y$, defined (at least) for $\lambda \in \mathbb{D}(0, m\rho)$, is given by:

$$\lambda \in \mathbb{D}(0,\rho) \mapsto \mathbf{E}(e^{\lambda Y}) = \sum_{k=0}^{\infty} \frac{\mathbf{E}(Y^k)}{k!}\lambda^k .$$

We have the following identity relating the two moment generating functions:

$$e^{-\lambda m_Y}\Big(\sum_{k=0}^{\infty} \frac{\mathbf{E}(Y^k)}{k!}\lambda^k\Big) = e^{-\lambda m_Y}\mathbf{E}(e^{\lambda Y}) = \mathbf{E}(e^{\lambda(Y-m_Y)}) = \sum_{k=0}^{\infty} \frac{\nu_k(Y)}{k!}, \quad \text{for } \lambda \in \mathbb{D}(0.\rho) . \tag{1.4.2}$$

**Example 1.4.1.** As a basic example, consider *a random variable $Z$ which follows a standard normal distribution.* In this case $\rho = +\infty$ and

$$\mathbf{E}(e^{\lambda Z}) = \int_{\mathbb{R}} e^{\lambda x}\frac{e^{-x^2/2}}{\sqrt{2\pi}}dx = e^{\lambda^2/2}, \quad \text{for every } \lambda \in \mathbb{R} . \tag{1.4.3}$$

⊡

Multiplying the power series on the far left of (1.4.2) and equating coefficients with the far right we may express *the centered moments of $Y$ in terms of the moments of $Y$*:

$$\nu_k(Y) = \sum_{j=0}^{k} \binom{k}{j}\mathbf{E}(Y^j)\,(-m_Y)^{k-j}, \quad \text{for any } k \geq 0 .$$

This set of identities is a direct consequence of the binomial formula:

$$(Y-m_Y)^k = \sum_{j=0}^{k} \binom{k}{j} Y^j\,(-m_Y)^{k-j} .$$

If we move the factor $e^{-\lambda m_Y}$ to the far right of (1.4.2), multiply the two power series which are now in the right and equate coefficients, we are able to express *the moments of $Y$ in terms of its centered moments of $Y$* as follows:

$$\mathbf{E}(Y^k) = \sum_{j=0}^{k} \binom{k}{j}\nu_j m_Y^{k-j}, \quad \text{for any } k \geq 0 .$$

It is sometimes more convenient to consider the function $M_Y(\lambda)$ given by

$$M_Y(\lambda) = \ln \mathbf{E}\big(e^{\lambda(Y-m_Y)}\big) = \ln\Big(1 + \sum_{k=2}^{\infty} \frac{\lambda^k \nu_k(Y)}{k!}\Big) .$$

Since the power series $\mathbf{E}(e^{\lambda(Y-m_Y)})$ takes the value 1 at $\lambda = 0$, it does not vanish in a disk $\mathbb{D}(0,\widehat{\rho})$, for some $\widehat{\rho} \in (0,\rho]$. The function $M_Y(\lambda)$ is holomorphic in the disk $\mathbb{D}(0,\widehat{\rho})$. We take $M_Y$ to be real in the real interval $\lambda \in (-\widehat{\rho},\widehat{\rho})$.

If needed, we may reduce $\widehat{\rho}$ to guarantee further that

$$(\natural) \quad \Big|\sum_{k=2}^{\infty} \frac{\lambda^k \nu_k(Y)}{k!}\Big| < 1\,, \quad \text{for } \lambda \in \mathbb{D}(0,\widehat{\rho})\,.$$

Observe that $M_Y(0)=0$ and $M_Y'(0)=0$.

**Example 1.4.2.** *For a random variable $Z$ following a standard normal distribution,* we have $M_Z(\lambda)=\lambda^2/2$, for every $\lambda \in \mathbb{C}$. ⊡

Next we relate the Taylor coefficients of $M_Y$ with the $\nu_k(Y)$, by means of the Taylor expansion

$$\ln(1+z)=\sum_{n=1}^{\infty}(-1)^{n+1}\frac{z^n}{n}\,, \quad \text{for } |z|<1\,.$$

Appealing to ($\natural$) and to the expansion above we obtain that

$$\begin{aligned}\sum_{j=2}^{\infty}\frac{M_Y^{(j)}(0)}{j!}\lambda^j &= \ln\Big(1+\sum_{k=2}^{\infty}\frac{\nu_k(Y)}{k!}\lambda^k\Big)\\ &= \Big(\sum_{k=2}^{\infty}\frac{\nu_k(Y)}{k!}\lambda^k\Big)-\frac{1}{2}\Big(\sum_{k=2}^{\infty}\frac{\nu_k(Y)}{k!}\lambda^k\Big)^2+\frac{1}{3}\Big(\sum_{k=2}^{\infty}\frac{\nu_k(Y)}{k!}\lambda^k\Big)^3-\cdots\\ &= \frac{\nu_2(Y)}{2!}\lambda^2+\frac{\nu_3(Y)}{3!}\lambda^3+\Big(\frac{\nu_4(Y)}{4!}-\frac{\nu_2(Y)^2}{8}\Big)\lambda^4+\Big(\frac{\nu_5(Y)}{5!}-\frac{\nu_2(Y)\nu_3(Y)}{24}\Big)\lambda^5+\cdots\end{aligned}$$

Equating coefficients, we obtain a sequence of polynomials $Q_k(x_2,\ldots,x_k)$ in $k-1$ variables, starting with $k=2$, so that for any variable $Y$ satisfying (1.4.1) it holds that

$$(1.4.4) \qquad M_Y^{(k)}(0)=Q_k\big(\nu_2(Y),\ldots,\nu_k(Y)\big)\,, \quad \text{for any } k \geq 2\ .$$

For the first values of $k$ we have that

$$\begin{aligned} M_Y^{(2)} &= \nu_2(Y)\\ M_Y^{(3)} &= \nu_3(Y)\\ M_Y^{(4)} &= \nu_4(Y)-3\nu_2(Y)^2\\ M_Y^{(5)} &= \nu_5(Y)-5\nu_2(Y)\nu_3(Y)\,.\end{aligned}$$

Each polynomial $Q_k$ satisfies

$$Q_k\left(\frac{x_2}{r^2},\frac{x_3}{r^3},\ldots,\frac{x_k}{r^k}\right)=\frac{1}{r^k}Q_k(x_2,\ldots,x_k)\,, \quad \text{for any } r\neq 0 \text{ and any } x_2,\ldots,x_k\,.$$

Conversely, from

$$1+\sum_{k=2}^{\infty}\frac{\nu_k(Y)}{k!}\lambda^k=\exp\big(M_Y(\lambda)\big)\,,$$

we may deduce analogously the existence of polynomials $R_k(x_2,\dots,x_k)$ in $k-1$ variables, starting with $k=2$, so that for any variable $Y$ satisfying (1.4.1) it holds that

$$\nu_k(Y) = R_k\big(M_Y^{(2)}(0),\dots,M_Y^{(k)}(0)\big)\,,\quad \text{for any } k\geq 2\ . \tag{1.4.5}$$

Each polynomial $R_k$ satisfies the following homogeneity condition

$$R_k\left(\frac{x_2}{r^2},\frac{x_3}{r^3},\dots,\frac{x_k}{r^k}\right) = \frac{1}{r^k}R_k(x_2,\dots,x_k)\,,\quad \text{for any } r\neq 0 \text{ and any } x_2,\dots,x_k\,. \tag{1.4.6}$$

**Example 1.4.3.** *For a standard normal random variable $Z$* we have $m_Z=0$ and $\nu_2(Z)=\mathbf{V}(Z)=1$. Because of symmetry, all odd moments of $Z$ are null: $\nu_{2k+1}(Z)=0$, for $k\geq 0$, while for the even moments we have that

$$\nu_{2k}(Z) = \frac{1}{2^k}\frac{(2k)!}{k!}\,,\quad \text{for any } k\geq 1\,.$$

Since $M_Z(\lambda)=\lambda^2/2$, for any $\lambda\in\mathbb{C}$, we also have $M_Z^{(2)}(0)=1$ and $M_Z^{(k)}(0)=0$, for any $k\geq 3$.

In any case,

$$\mathbf{E}(Z^k) = \nu_k(Z) = R_k\big(1,0,\dots,0\big)\,,\quad \text{for any } k\geq 3\,, \tag{1.4.7}$$

and

$$0 = Q_k\big(1,\nu_3(Z),\dots,\nu_k(Z)\big)\,,\quad \text{for any } k\geq 3\,. \tag{1.4.8}$$

⊡

### 1.4.2 Moment generating functions of a Khinchin family and its fulcrum

Let $f$ be a power series in $\mathcal{K}$ with radius of convergence $R>0$ with Khinchin family $(X_t)_{t\in[0,R)}$. And let $F$ denote, as usual, the fulcrum of $f$, which is defined in the interval $(-\infty,\ln R)$:

$$e^{F(s)} = f(e^s)\,,\quad \text{for } s<\ln R\,.$$

We have, see Section 1.3.4, that $m_f(e^s)=F'(s)$ and $\sigma_f(s)=F''(s)$, for every $s<\ln R$.

Next we are going to obtain convenient expressions of the moment generating functions of each $X_t$ in terms of $f$ itself and more conveniently in terms of the fulcrum $F$.

Fix $t\in(0,R)$ and set $s=\ln t$. Since $X_t$ is a non-negative variable, the following change of order of summation is legitimate and we have, if $te^\lambda<R$, that

$$(\flat)\qquad \begin{aligned}\sum_{k=0}^\infty \frac{\mathbf{E}(X_t^k)}{k!}\lambda^k &= \frac{1}{f(t)}\sum_{k=0}^\infty\frac{1}{k!}\sum_{n=0}^\infty n^k a_n t^n\lambda^k = \frac{1}{f(t)}\sum_{n=0}^\infty a_n t^n\Big(\sum_{k=0}^\infty\frac{(n\lambda)^k}{k!}\Big)\\ &= \frac{1}{f(t)}\sum_{n=0}^\infty a_n (te^\lambda)^n = \frac{f(te^\lambda)}{f(t)} < +\infty\,.\end{aligned}$$

Thus with $\rho(t) = \ln R - \ln t$, condition (1.4.1) is satisfied:

$$\sum_{k=0}^{\infty} \frac{\mathbf{E}(|X_t|^k)}{k!} |z|^k \,, \quad \text{for } |z| < \rho(t)\,.$$

The identity ($\flat$) actually shows that the moment generating functions of $X_t$, centered or not, are defined for all $\lambda \in (-\infty, \rho(t))$ and may be written in terms of $f$ as

$$\text{(1.4.9)} \qquad \begin{cases} \mathbf{E}(e^{\lambda X_t}) = \dfrac{f(te^{\lambda})}{f(t)}, \\ \mathbf{E}(e^{\lambda(X_t - m_f(t))}) = \dfrac{f(te^{\lambda})}{f(t)} e^{-\lambda m_f(t)}, \end{cases} \quad \text{for each } \lambda < \ln R - \ln t\,,$$

and in terms of the fulcrum $F$ as

$$\text{(1.4.10)} \qquad \begin{cases} \mathbf{E}(e^{\lambda X_t}) = \exp\big(F(s+\lambda) - F(s)\big), \\ \mathbf{E}(e^{\lambda(X_t - m_f(t))}) = \exp\big(F(s+\lambda) - F(s) - \lambda F'(s)\big), \end{cases} \quad \text{for each } \lambda < \ln R - s\,.$$

Recall that, $t = e^s$, for $t < \ln R$, and abbreviate $\nu_k(t) = \nu_k(X_t) = \mathbf{E}((X_t - m_f(s))^k)$, for $k \geq 1$. We have,that

$$M_{X_t}(\lambda) = F(s+\lambda) - F(s) - \lambda F'(s) = \sum_{j=2}^{\infty} \frac{F^{(j)}(s)}{j!} \lambda^j \,, \quad \text{if } |\lambda| < \ln R - s\,.$$

and so

$$\text{(1.4.11)} \qquad F^{(k)}(s) = Q_k\big(\nu_2(t), \ldots, \nu_k(t)\big)\,, \quad \text{for any } k \geq 2 \;\;,$$

and

$$\text{(1.4.12)} \qquad \nu_k(t) = R_k\big(F^{(2)}(s), \ldots, F^{(k)}(s)\big)\,, \quad \text{for any } k \geq 2 \;\;.$$

### 1.4.3 Moment generating function of a Khinchin family and Chernoff bounds

We continue with the Khinchin family $(X_t)_{t\in[0,R)}$ associated to a power series $f(z)$ in $\mathcal{K}$ with radius of convergence $R$ and with fulcrum $F(s)$ defined in $(-\infty, \ln R)$.

For $t \in (0, R)$ and $s = \ln t$, and $\lambda \in \mathbb{R}$ such that $te^{\lambda} < R$ (or $s{+}\lambda < \ln R$) we have, see (1.4.10), that

$$\mathbf{E}\big(e^{\lambda(X_t - m_f(t))}\big) = \exp\big(F(s{+}\lambda) - F(s) - \lambda F'(s)\big)\,.$$

We may bound the expression $\big(F(s{+}\lambda) - F(s) - \lambda F'(s)\big)$ as follows. Fix $s < \ln R$, and let $\lambda \in \mathbb{R}$ be such that $s + \lambda < \ln R$, i.e., $\lambda \in (-\infty, (\ln R){-}s)$. We may write

$$F(s+\lambda) - F(s) - F'(s)\lambda = \int_0^{\lambda} F''(s{+}y)(\lambda - y)dy = \int_0^{\lambda} \sigma_f^2(e^s e^y)(\lambda - y)dy$$

with the standard convention that $\int_b^a = -\int_a^b$ whenever $a < b$. For $s, s+\lambda < R$, the quantity $F(s+\lambda) - F(s) - F'(s)\lambda$ is non-negative. Since $F'$ is increasing, the intermediate value theorem gives $F(s+\lambda) - F(s) < \lambda F'(s+\lambda)$ and thus

$$F(s+\lambda) - F(s) - F'(s)\lambda \le \lambda(F'(s+\lambda) - F'(s))\,.$$

Observe that if $\lambda < 0$ both factors in the bound in the right of the above inequality are negative.

Define for $s < \ln R$ and $0 \le \Lambda < (\ln R) - s$, the quantity:

$$\Sigma(s,\Lambda) = 2 \max_{|y|\le|\Lambda|} \frac{F'(s+y) - F'(s)}{y}\,.$$

For each such $s$ fixed, the expression $\Sigma(s,\Lambda)$ is positive and increasing for $\Lambda$ in the indicated range. We have, if $s, s+\lambda < \ln R$, that

$$\mathbf{E}\big(e^{\lambda(X_t - m_f(t))}\big) \le e^{(\lambda^2/2)\Sigma(s,|\lambda|)}\,.$$

The intermediate value theorem gives for $s, s+\lambda < \ln R$ that

$$\Sigma(s,\Lambda) \le 2 \max_{|y|\le|\Lambda|} F''(s+y)\,.$$

**Proposition 1.4.4.** *For the Khinchin family $(X_t)_{t\in[0,R)}$ associated to $f(z)$ in $\mathcal{K}$ with radius of convergence $R$ and fulcrum $F$, we have, for $s, s+\lambda < \ln R$ and $t = e^s$ that*

$$\mathbf{E}\big(e^{\lambda(X_t - m_f(t))}\big) \le e^{\lambda(F'(s+\lambda) - F'(s))} \le e^{(\lambda^2/2)\Sigma(s,|\lambda|)}\,. \tag{1.4.13}$$

### 1.4.4 Chernoff bounds for subexponential random variables

The most simple form of concentration for a random variable $Y$ is given by Chebyshev's inequality For random variables with finite, and small variance, this inequality ensures that the random variable $Y$ is close to its mean. Chebyshev's inequality only requires the second moment of $Y$ to be finite.

We are interested in obtaining one or two sided bounds for the tails, that is large deviation inequalities, for certain random variable $Y$; these are the Chernoff bounds. In this last case we require all the moments of $Y$ to be finite (in fact we impose that the random variable $Y$ is subexponential). We refer to Chapter 2 in [98] for further details about these concentration inequalities.

Let $U$ be a random variable. The random variable $Y$ is termed subexponential if for a certain $\Lambda > 0$ and a certain $A > 0$, it holds that

$$\mathbf{E}(e^{\lambda U}) \le e^{(\lambda^2/2)A}\,, \quad \text{for each } \lambda \text{ such that } |\lambda| \le \Lambda\,. \tag{1.4.14}$$

Proposition 1.4.4 just above tell us that each random variable of a Khinchin family is actually subexponential.

In Proposition 1.4.5 just below we are going to show that general subexponential random variables satisfy a large deviation inequality, which as a particular case will give raise in Proposition 1.4.6 to a large deviations inequality for each member of a Khinchin family.

**Proposition 1.4.5.** *Let $U$ be a random variable. Assume that for certain $\Lambda > 0$ and certain $A > 0$, it holds that*

$$\mathbf{E}(e^{\lambda U}) \le e^{(\lambda^2/2)A}\,, \quad \text{for each } \lambda \text{ such that } |\lambda| \le \Lambda\,. \tag{1.4.15}$$

*Then*

$$\mathbf{P}(|U| > y) \le 2 \cdot \begin{cases} e^{-y^2/(2A)}\,, & \text{if } 0 \le y \le \Lambda A\,, \\ e^{-\Lambda y/2}\,, & \text{if } y \ge \Lambda A\,. \end{cases} \tag{1.4.16}$$

A random variable satisfying (1.4.15) for a pair of parameters $A$ and $\Lambda$ is called *subexponential.*

*Proof.* For $y > 0$, we wish to give an upper bound for the probability $\mathbf{P}(U > y)$. To do this, take $\lambda > 0$, whose value we are going to optimize depending on $y$.

Markov's inequality gives that

$$\mathbf{P}(U > y) = \mathbf{P}\big(\lambda U > \lambda y\big) = \mathbf{P}\big(e^{\lambda U} \ge e^{\lambda y}\big) \le e^{-\lambda y}\mathbf{E}\big(e^{\lambda U}\big)\,,$$

and thus we obtain from (1.4.15) that

$$\mathbf{P}(U > y) \le e^{-\lambda y} e^{(\lambda^2/2)A}\,.$$

For any $A > 0$, the parabola $\lambda \mapsto -\lambda y + A\lambda^2/2$, decreases for $\lambda \in (-\infty, y/A)$, increases for $\lambda \in (y/A, +\infty)$ and attains a minimum value of $-y^2/(2A)$ at $\lambda = y/A$.

Thus we deduce that

$$\mathbf{P}(U > y) \le \begin{cases} e^{-y^2/(2A)}\,, & \text{if } 0 \le y \le \Lambda A\,, \\ e^{-\Lambda y/2}\,, & \text{if } y \ge \Lambda A\,. \end{cases} \tag{1.4.17}$$

For $y > 0$ and $\lambda > 0$ we may write, analogously, that

$$\mathbf{P}(U < -y) = \mathbf{P}(-\lambda U > \lambda y) \le e^{-\lambda y}\mathbf{E}(e^{-\lambda U}) \le e^{-\lambda y} e^{(\lambda^2/2)A)}\,,$$

Arguing as above, we obtain

$$\mathbf{P}(U < -y) \le \begin{cases} e^{-y^2/(2A)}\,, & \text{if } 0 \le y \le \Lambda A\,, \\ e^{-\Lambda y/2}\,, & \text{if } y \ge \Lambda A\,. \end{cases} \tag{1.4.18}$$

Combining (1.4.17) and (1.4.18) we get (1.4.16).

□

### 1.4.5 Chernoff bounds for Khinchin families

Let $(X_t)_{t\in[0,R)}$ be the Khinchin family associated to a power series $f(z)$ in $\mathcal{K}$ with radius of convergence $R$.

Let $t \in (0, R)$ and write $t = e^s$. Proposition 1.4.4 gives us that $U = X_t - m_f(t)$ satisfies the subexponentiality condition (1.4.15) for any $\Lambda$ such that $0 \le \Lambda < \ln R - s$ and any $A \ge \Sigma(s, \Lambda)$.

Thus from Proposition 1.4.5 we deduce

**Proposition 1.4.6** (Chernoff bounds)**.** *For $t \in (0, R)$ and $\Lambda \ge 0$ such that $te^\Lambda < R$, and we have for $s = \ln t$ that*

$$\mathbf{P}(|X_t - m_f(t)| > y) \le 2 \cdot \begin{cases} e^{-y^2/(2\Sigma(s,\Lambda)}, & \text{if } 0 \le y \le \Lambda\Sigma(s,\Lambda)\,, \\ e^{-\Lambda y/2}, & \text{if } y \ge \Lambda\Sigma(s,\Lambda)\,. \end{cases} \tag{1.4.19}$$

*More generally, for any $A \ge \Sigma(s, \Lambda)$, we have*

$$\mathbf{P}(|X_t - m_f(t)| > y) \le 2 \cdot \begin{cases} e^{-y^2/(2A)}, & \text{for } 0 \le y \le \Lambda A\,, \\ e^{-\Lambda y/2}, & \text{for } y \ge \Lambda A\,. \end{cases} \tag{1.4.20}$$

Observe that the restriction $te^\Lambda < R$ is vacuous in the case of entire functions $f \in \mathcal{K}$.

In terms of the power series $f$ itself, the bounds of (1.4.19) mean that

$$\sum_{|k-m_f(t)|>y} a_n t^n \le 2\, f(t) \cdot \begin{cases} e^{-y^2/(2\Sigma(s,\Lambda))}, & \text{if } 0 \le y \le \Lambda\Sigma(s,\Lambda)\,, \\ e^{-\Lambda y/2}, & \text{if } y \ge \Lambda\Sigma(s,\Lambda)\,. \end{cases}$$

## 1.5 Khinchin families and Wiman-Valiron theory

This section is devoted to show how the probabilistic point of view of the Khinchin families illuminates and provides (relatively) simple proofs of some (of the most fundamental) results of the Wiman-Valiron theory. The ideas and the essential ingredients of this approach are contained in the seminal paper [83] of Rosenbloom, and also in the paper [88] of Schumitzky. Let us mention that the paper of Hayman, [50], is a classical reference for the Wiman-Valiron method.

We are going to present proofs within the framework of Khinchin families of the classical Wiman-Valiron Theorem and of a companion theorem due to Clunie, which claim that one or a few terms of an entire power series represent quite closely the whole power series. The proofs are based respectively in the general bounds contained in Lemma 1.5.3 and Proposition 1.4.6. In those preliminary (from the point of view of the Wiman-Valiron theory) results the maximum term or a few terms of the power series are compared to the power series as a whole and the comparison is controlled by the variance function. The extra ingredient now is that this variance function $\sigma_f^2(t)$ could be bounded in terms of $f(t)$, but not for all values of $t$, if fact for all but a set of finite logarithmic measure. These comparisons of $\sigma_f$ with $f$ are based on some differential inequalities and also on a classical lemma of Borel, which we present first.

### 1.5.1 Probability bounds and variance of the variables $X_t$

The following bound is a useful observation about general random variables with values in the lattice $\mathbb{Z}$. The bound follows simply from Chebyshev's inequality.

**Lemma 1.5.1.** *For some absolute constant $H$ and any random variable $U$ taking values in $\mathbb{Z}$, we have that*

$$\max_{k\in\mathbb{Z}} \mathbf{P}(U=k) \geq \frac{H}{1+\sqrt{\mathbf{V}(U)}}\,.$$

The proof below works with $H = 1/(4\sqrt{2})$.

Lemma 1.5.1 quantifies the observation that if all the probabilities $\mathbf{P}(U=k)$ are small, then the variable $U$ is quite spread out and its variance must be large.

*Proof.* We may assume that $\mathbf{V}(U) < +\infty$. And also that $\sqrt{\mathbf{V}(U)} > 0$. Let us denote $\mu = \mathbf{E}(U)$ and $\sigma^2 = \mathbf{V}(U)$.

Chebyshev's inequality gives that

$$\mathbf{P}\big((U-\mu)^2 \geq 2\sigma^2\big) \leq \frac{1}{2}\,, \quad \text{and thus that} \quad \mathbf{P}\big((U-\mu)^2 \leq 2\sigma^2\big) \geq \frac{1}{2}\,.$$

Therefore

$$\begin{aligned}
\frac{1}{2} &\leq \mathbf{P}\big((U-\mu)^2 \leq 2\sigma^2\big) = \sum_{-\sqrt{2}\sigma \leq |k-\mu| \leq \sqrt{2}\sigma} \mathbf{P}(U=k)\\
&\leq \big(\max_{k\in\mathbb{Z}} \mathbf{P}(Z=k)\big)\cdot \#\big\{k\in\mathbb{Z} : -\sqrt{2}\sigma \leq |k-\mu| \leq \sqrt{2}\sigma\big\}\\
&\leq \big(\max_{k\in\mathbb{Z}} \mathbf{P}(Z=k)\big)\cdot(1+2\sqrt{2}\sigma) \leq \Big(\max_{k\in\mathbb{Z}} \mathbf{P}(Z=k)\Big)2\sqrt{2}\cdot(1+\sigma)\,,
\end{aligned}$$

which gives the result.

We have used above that if $I$ is a closed interval in $\mathbb{R}$, then $\#\big(I\cap\mathbb{Z}\big) \leq 1+|I|$. $\square$

In the inequality of Lemma 1.5.1, the quantity $1+\sqrt{\mathbf{V}(U)}$ cannot be replaced by $\sqrt{\mathbf{V}(U)}$. For instance, for $p\in(0,1/3)$, let $U_p$ be the random variable which takes the values $\pm1$ with probability $p$, and the value $0$ with probability $1-2p$. Observe that $\sqrt{\mathbf{V}(U_p)} = \sqrt{2p}$ and $\max(p,1-2p) = 1-2p$ (recall that $p\in(0,1/3)$). We have

$$\lim_{p\downarrow 0}\Big(\big(\max_{k\in\mathbb{Z}} P(U_p=k)\big)\cdot\sqrt{\mathbf{V}(U_p)}\Big) = \lim_{p\downarrow 0}\max(p,1-2p)\cdot\sqrt{2p} = 0\,.$$

**Remark 1.5.2.** *Continuous variables.* For continuous random variables, the following analogue of Lemma 1.5.1 holds. Let $Y$ be a continuous random variable, with density $\phi(x)$, (finite) standard deviation $\sigma$ and mean $m$. Let $\nu$ denote $\nu = \sup_{x\in\mathbb{R}}\phi(x)$. Chebyshev's inequality gives that

$$\int_{|x-m|\leq\lambda\sigma}\phi(x)dx = \mathbf{P}(|U-m|\leq\lambda\sigma) \geq (1-\frac{1}{\lambda^2}), \qquad \text{for any } \lambda>1.$$

Therefore

$$\sigma \cdot \nu \geq \frac{1}{2\lambda}\left(1 - \frac{1}{\lambda^2}\right), \quad \text{if } \lambda > 1,$$

maximizing the lower bound in $\lambda$ (taking $\lambda = 1/\sqrt{3}$) we find that

$$\sigma \cdot \nu \geq Q,$$

with $Q = 1/(3\sqrt{3})$.

For a random variable $U$ distributed uniformly in the interval $[-1, 1]$ we have $\nu = 1/2$ and $\sigma = 1/\sqrt{3}$, and thus $\sigma \cdot \nu = 1/(2\sqrt{3})$, in this case. ⊠

For the Khinchin family $(X_t)_{t\in[0,R)}$ associated to $f \in \mathcal{K}$, Lemma 1.5.1 translates into the following inequality:

$$\left(\max_{n\geq 0}\left(\frac{a_n t^n}{f(t)}\right)\right) \cdot (1 + \sigma_f(t)) \geq H\,, \quad \text{for any } t \in [0, R)\,,$$

which we write as a Lemma.

**Lemma 1.5.3.** *For a power series $f(z) = \sum_{n=0}^{\infty} a_n z^n \in \mathcal{K}$ with radius of convergence $R > 0$ we have*

$$f(t) \leq \frac{1}{H}\left(\max_{n\geq 0} a_n t^n\right)(1 + \sigma_f(t))\,, \quad \textit{for any } t \in [0, R)\,.$$

### 1.5.2 Differential inequalities

The following lemma is Proposition (the only one) in Fuchs' paper [34].

**Lemma 1.5.4** (Fuchs)**.** *Let $\mu$ be a regular Borel (positive) measure in $\mathbb{R}$. Let $\mathcal{I}$ be a* collection *of open intervals indexed by a set $\Omega$:*

$$\mathcal{I} = \{I_\omega : \omega \in \Omega\}\,,$$

*which is such that for a constant $B > 0$ and any finite subcollection of* pairwise *disjoint intervals $I_{\omega_j}, 1 \leq j \leq n$, it holds that*

$$\sum_{j=1}^{n} \mu(I_{\omega_j}) \leq B\,.$$

*Then for any Borel set $E$ such that $E \subset \bigcup_{\omega\in\Omega} I_\omega$, it holds that $\mu(E) \leq 2B$.*

In the proof we will use Besicovitch covering Theorem (valid in $\mathbb{R}^k$) in the elementary one-dimensional case. This one-dimensional covering result claims that if $(J_j)_{j=1}^m$ is a finite collection of finite open intervals, then we can find *two* disjoint sets of indices $U, V \subset \{1, \ldots, m\}$, so that both the intervals $\{J_j; j \in U\}$ and the intervals $\{J_j; j \in V\}$ are pairwise disjoint and so that

$$\bigcup_{j=1}^{m} I_j = \Big(\bigcup_{j\in U} I_j\Big) \cup \Big(\bigcup_{j\in V} I_j\Big)\,.$$

In our applications below the measures $\mu$ are integrals against non-negative continuous functions in $\mathbb{R}$, which are locally finite ($\mu(K) < +\infty$ for any compact set $K$ in $\mathbb{R}$) and thus regular, see, for instance, Theorem 2.18 and Theorem 2.17 of [86].

*Proof.* By (inner) regularity, it is enough to show that $\mu(K) \le 2B$, for any compact $K \subset E$.

Let $\omega_1, \dots, \omega_m$ be such that

$$K \subset \bigcup_{j=1}^{m} I_{\omega_j} \,.$$

Write $J_j = I_{\omega_j}$ Applying the covering Theorem with the notations above we have

$$\sum_{j\in U} \mu(J_j) \le B \quad \text{and} \quad \sum_{j\in V} \mu(J_j)$$

and thus that

$$\mu(K) \le \mu\Big(\bigcup_{j\in U} I_j\Big) + \mu\Big(\bigcup_{j\in V} I_j\Big) \le \sum_{j\in U} \mu(J_j) + \sum_{j\in V} \mu(J_j) \le 2B\,,$$

as claimed. □

A typical application of Lemma 1.5.4 follows. It is an adaptation of [34, Lemma 1].

**Lemma 1.5.5.** *Let $G$ be a $C^1$ function defined in $(a, b) \subseteq (-\infty, +\infty)$ and increasing. Let $A = \lim_{x\downarrow a} G(x)$ and $B = \lim_{x\uparrow b} G(x)$. Let $\phi$ be a positive continuous function defined in $(A, B)$ and such that $\int_A^B \phi(v)dv < +\infty$.*

*Then there exists an open set $E \subset (a, b)$ of finite Lebesgue measure so that*

$$\frac{1}{\lambda}\int_{G(x)}^{G(x+\lambda)} \phi(u)du \le 1\,, \quad \text{for every } x \in (a, b) \setminus E \text{ and any } \lambda > 0 \text{ such that } x + \lambda < b,$$

*and*

$$\frac{1}{\lambda}\int_{G(x-\lambda)}^{G(x)} \phi(u)du \le 1\,, \quad \text{for every } x \in (a, b) \setminus E \text{ and any } \lambda > 0 \text{ such that } x - \lambda > a\,.$$

In this lemma we allow $a$ and $A$ to be $-\infty$ and $b$ and $B$ to be $+\infty$. Of course if $(a, b)$ is a finite interval, there is nothing to prove.

*Proof.* We define $E_+$ by

$$E_+ = \Big\{x \in (a, b) : \text{there exists } y \in (x, b) \text{ such that } y - x < \int_x^y \phi(G(u))G'(u)du\Big\}\,.$$

The set $E_+$ is open, since if $x \in E_+$ and $y$ is as in the definition above, then by continuity and for some $\delta > 0$, we have that

$$y - z < \int_z^y \phi(G(u))G'(u)du, \quad \text{for any } z \in (x - \delta, x + \delta)\,.$$

The set $E_+$ is contained in the union of a collection $\mathcal{I}$ of open intervals (contained in $(a,b)$) which satisfy

$$|I| < \int_I \phi(G(u))G'(u)du\,.$$

Thus for any finite collection of intervals $I_1,\dots,I_m$ of pairwise disjoint intervals of $\mathcal{I}$ we have that

$$\sum_{j=1}^m |I_j| < \sum_{j=1}^m \int_{I_j} \phi(G(u))G'(u)du < \int_a^b \phi(G(u))G'(u)du = \int_A^B \phi(v)dv\,.$$

Given this, Lemma 1.5.4 shows that $|E_+| < +\infty$.

For any $x \in (a,b)\setminus E_+$ and if $x < x+\lambda < b$ we then have that

$$\lambda = (x+\lambda) - x \ge \int_x^{x+\lambda} \phi(G(u))G'(u)du = \int_{G(x)}^{G(x+\lambda)} \phi(v)dv\,.$$

Arguing analogously or applying the above with $-G(-x)$ and $\phi(-u)$, we obtain $E_- \subset (a,b)$ of finite measure so that if $x \in (a,b)\setminus E_-$ and if $x > x-\lambda > a$ then

$$\lambda = x - (x-\lambda) \ge \int_{x-\lambda}^{x} \phi(G(u))G'(u)du = \int_{G(x-\lambda)}^{G(x)} \phi(v)dv\,.$$

□

**Corollary 1.5.6.** *Let $G$ be a $C^1$ function defined in $[a,+\infty)$ and increasing so that $G(a) = A$ and $\lim_{x\to\infty} G(x) = +\infty$. Then for any $\varepsilon > 0$ there exists $c > a$ and a set $E \subset [c,+\infty)$ of finite Lebesgue measure so that for each $x \in (c,+\infty)\setminus E$ we have*

$$\frac{G(x+\lambda) - G(x)}{\lambda} \le G(x+\lambda)^{1+\varepsilon}\,, \quad \textit{for every } \lambda > 0\,,$$

$$\frac{G(x) - G(x-\lambda)}{\lambda} \le G(x)^{1+\varepsilon}\,, \ \textit{for every } \lambda > 0 \textit{ such that } x-\lambda > c\,,$$

*and besides*

$$(\flat) \qquad G'(x) \le G(x)^{1+\varepsilon}\,.$$

*Proof.* Take $c > a$ so that $G(c) > 1$ and take $\phi(u) = 1/u^{1+\varepsilon}$, for $u \ge G(c)$. Then $\int_{G(c)}^\infty \phi(v)dv < +\infty$.

Apply Lemma 1.5.5 to obtain an exceptional set $E \subset [c,+\infty)$ of finite Lebesgue measure, such that for $x \in [c,+\infty)\setminus E$ and $\lambda > 0$ we have that

$$\lambda \ge \int_{G(x)}^{G(x+\lambda)} \phi(u)du \ge \frac{1}{G(x+\lambda)^{1+\varepsilon}} \int_{G(x)}^{G(x+\lambda)} du = \frac{G(x+\lambda) - G(x)}{G(x+\lambda)^{1+\varepsilon}}\,,$$

and for $x \in [c,+\infty)\setminus E$ and $\lambda > 0$ such that $x-\lambda > c$ we have that

$$\lambda \ge \int_{G(x-\lambda)}^{G(x)} \phi(u)du \ge \frac{1}{G(x)^{1+\varepsilon}} \int_{G(x-\lambda)}^{G(x)} du = \frac{G(x) - G(x-\lambda)}{G(x)^{1+\varepsilon}}\,.$$

The inequality $(\flat)$ follows from by fixing $x \in (c,+\infty)\setminus E$ and letting $\lambda \to 0$. □

Next, we apply Corollary 1.5.6 to entire functions in $\mathcal{K}$. Let $f$ be an entire power series in $\mathcal{K}$, and let $F$ be its fulcrum defined in $(-\infty, +\infty)$. We have, see Section 1.3.4, that

$$F(s) = \ln f(e^s), \ F'(s) = m_f(e^s) \quad \text{and} \quad F''(s) = \sigma_f^2(e^s)\,, \quad \text{for any } s \in \mathbb{R}\,.$$

**Proposition 1.5.7.** *For $f$ entire in $\mathcal{K}$ and for $\varepsilon > 0$, there exists $S \in \mathbb{R}$ and a set $E \subset [S, +\infty)$ of finite measure such that for $s \in [S, +\infty) \setminus E$ we have that*

$$F'(s) \le F(s)^{1+\varepsilon} \quad \textit{and} \quad \begin{cases} \dfrac{F(s+\lambda) - F(s)}{\lambda} \le F(s+\lambda)^{1+\varepsilon}\,, & \textit{for every } \lambda > 0 \\ \dfrac{F(s) - F(s-\lambda)}{\lambda} \le F(s)^{1+\varepsilon}\,, & \textit{if } s > s-\lambda > S\,, \end{cases}$$

*and also that*

$$F''(s) \le F'(s)^{1+\varepsilon} \quad \textit{and} \quad \begin{cases} \dfrac{F'(s+\lambda) - F'(s)}{\lambda} \le F'(s+\lambda)^{1+\varepsilon}\,, & \textit{for every } \lambda > 0 \\ \dfrac{F'(s) - F'(s-\lambda)}{\lambda} \le F'(s)^{1+\varepsilon}\,, & \textit{if } s > s-\lambda > S\,, \end{cases}$$

*Besides,*

$$F''(s) \le F(s)^{1+\varepsilon}\,, \quad \textit{for every } s \in [S, +\infty) \setminus E\,. \tag{1.5.1}$$

*Proof.* This follows from applying Corollary 1.5.6 with $G = F$ and with $G = F'$. The last inequality follows from applying the two previous bounds with $\eta$ (instead of $\varepsilon$) such that $(1+\eta)^2 = 1+\varepsilon$. $\square$

### 1.5.3 Borel's Lemma

Borel's Lemma is a standard estimate on the growth of increasing functions which is of frequent use in the theory of entire functions. See, for instance, [85, Chapter 9]. We shall need the following variant of Borel's Lemma which appears in Hayman's paper [50] as Lemma 1: Borel's Lemma would correspond to the case $\varepsilon = 1$.

**Lemma 1.5.8** (Borel, Hayman)**.** *Let $G$ be an increasing function defined in $[a, +\infty)$ for some $a \in \mathbb{R}$, and such that $G(a) > 0$. For every $\varepsilon > 0$, there exists a set $E \subset [a, +\infty)$ of finite Lebesgue measure so that*

$$G\Big(x + \frac{1}{G(x)^{\varepsilon}}\Big) \le 2G(x)\,, \quad \textit{for any } x \in [a, +\infty) \setminus E\,.$$

The proof below follows the lines of the proof of Borel's lemma in [85]. In general formulations of this kind of results the functions $G$ are not assumed to be continuous; for certain applications not necessarily continuous functions are needed. But the assumption of continuity simplifies the proof a bit, and in our applications below the function $G$ are actually $C^{\infty}$.

*Proof.* Denote $y_0 = a$ and

$$E := \{x \geq y_0 : G\big(x + 1/G(x)^\varepsilon\big) \geq 2G(x)\}\,.$$

We are going to show that $E$ is contained in a countable (maybe finite) collection of closed intervals of finite total length and this that $|E| < \infty$.

Consider $E_0 = [y_0 + 1, +\infty) \cap E$. If $E_0 = \emptyset$ we are done, and we stop. If not, let $x_0 = \inf E_0$. Write $h_0 = 1/G(x_0)^\varepsilon$ and let $y_1 = x_0 + h_0$. We have

$$y_1 \geq x_0 \geq y_0 + 1 \quad \text{and} \quad G(y_1) \geq 2G(x_0)\,.$$

Consider $E_1 = [y_1 + 1/2, +\infty) \cap E$. If $E_1 = \emptyset$ we are done, and we stop. If not, let $x_1 = \inf E_1$. Write $h_1 = 1/G(x_1)^\varepsilon$ and let $y_2 = x_1 + h_1$. We have

$$y_2 \geq x_1 \geq y_1 + 1/2 \quad \text{and} \quad G(y_2) \geq 2G(x_1)\,.$$

It the process does not stop, we obtain two sequences $(x_n)_{n\geq 0}$ and $(y_n)_{n\geq 0}$ such that

$$y_{n+1} \geq x_n \geq y_n + 1/2^n \quad \text{and} \quad G(y_{n+1}) \geq 2G(x_n) \quad \text{for any } n \geq 0\,.$$

Continue with this process indefinitely or stop in for some $n \geq 0$ we have $E_n = \emptyset$, in which case the proof is finished. Since $G$ the function is non-decreasing we also have

$$G(x_{n+1}) \geq 2G(x_n)\,, \quad \text{for any } n \geq 0\,,$$

thus

$$G(x_n) \geq 2^n G(x_0)\,, \quad \text{for any } n \geq 0\,.$$

Also

$$y_{n+1} = x_n + h_n \quad \text{and} \quad h_n = \frac{1}{G(x_n)^\varepsilon} \leq \frac{1}{G(x_0)^\varepsilon}\frac{1}{2^{n\varepsilon}}\,, \quad \text{for any } n \geq 0\,.$$

The increasing sequence $(x_n)_{n\geq 0}$ tends to $\infty$, because if $\lim_{n\to\infty} x_n = u \in \mathbb{R}$ then $G(u) = +\infty$. And also $\lim_{n\to\infty} y_n = +\infty$.

We have that $[x_0, +\infty) = \bigcup_{n\geq 0}[x_n, x_{n+1}]$. And that

$$E \subset \bigcup_{n\geq 0}[x_n, y_{n+1} + 2^{-(n+1)}]$$

because $E \cap [y_{n+1} + 2^{-(n+1)}, x_{n+1}) = \emptyset$, by the very definition of $x_{n+1}$.

Finally,

$$|E| \leq \sum_{n=0}^{\infty}(y_{n+1} - x_n) + \sum_{n=0}^{\infty}\frac{1}{2^{(n+1)}} \leq 1 + \sum_{n=0}^{\infty} h_n \leq 1 + \frac{1}{G(x_0)^\varepsilon}\sum_{n=0}^{\infty}\frac{1}{2^{n\varepsilon}} < \infty\,.$$

□

### 1.5.4 Maximum term

The maximum term of an entire power series is a key ingredient of the Wiman-Valiron theory.

Let $g$ be an entire entire function $g(z) = \sum_{n=0}^{\infty} b_n z^n$, (no a priori assumptions on the coefficients other than $g$ being entire). The *maximum term* is the function $\mu(g,t)$ defined by

$$t \geq 0 \mapsto \mu(g,t) = \max_{t \geq 0} |b_n| t^n \,.$$

Since $g$ is entire, we have that $\lim_{n\to\infty} |b_n| t^n = 0$, for any $t \geq 0$, and thus that $\mu(g,t)$ is well defined: it is a (finite) maximum, not a supremum.

Recall also the standard notation: $M(g,t) = \max_{|z|\leq t} |g(t)|$. The Wiman-Valiron Theorem relates these two functions, $\mu(g,t)$ and $M(g,t)$, both defined in $(0,+\infty)$.

Cauchy's coefficient inequality gives that

$$|b_n| t^n \leq M(g,t) \,, \quad \text{for any } n \geq 0 \text{ and any } t \geq 0 \,, \tag{1.5.2}$$

so that

$$\mu(g,t) \leq M(g,t), \quad \text{for any } t \geq 0 \,.$$

The interesting comparison in the Wiman-Valiron theory goes in the other direction. As we shall see, the Wiman-Valiron Theorem claims that for every entire function $g$ the maximum term $\mu(g,t)$ is *not much less* than the maximum modulus $M(g,t)$ for *most values of* $t$, in a sense that we shall make precise later on, of course.

Observe that the function $f$ defined by the power series $\sum_{n=0}^{\infty} |b_n| z^n$, is also entire, (by the formula of Cauchy-Hadamard), and that

$$\mu(g,t) = \mu(f,t) \quad \text{and} \quad M(g,t) \leq M(f,t) = f(t) \,, \quad \text{for any } t \geq 0 \,. \tag{1.5.3}$$

This observation will mean that for the Wiman-Valiron comparison of $\mu(g,t)$ and $M(g,t)$ we could assume that $g$ has non-negative Taylor coefficients.

Observe also that if $h$ is an entire function in $\mathcal{K}$ and if $f(z) = z^k h(z)$ for some integer $k \geq 1$, then

$$\mu(f,t) = t^k \mu(h,t) \quad \text{and} \quad M(f,t) = t^k M(h,t) \,, \quad \text{for any } t \geq 0 \,. \tag{1.5.4}$$

This observation will mean that for the Wiman-Valiron comparison for an entire function $g$ we could assume that $g$ is in $\mathcal{K}$.

For an entire function $f \in \mathcal{K}$, we have that

$$\frac{\mu(f,t)}{f(t)} = \max_{k \geq 0} \mathbf{P}(X_t = k) \,, \quad \text{for any } t > 0 \,.$$

For a polinomial $g(z) = \sum_{n=0}^{N} b_n z^n$ of degree $N$, we have

$$\lim_{t\to\infty} \frac{\mu(g,t)}{t^N} = |b_N| = \lim_{t\to\infty} \frac{M(g,t))}{t^N} \,. \tag{1.5.5}$$

The maximum term of the exponential function is given (asymptotically) by the following lemma.

**Lemma 1.5.9.** *We have*

$$\sup_{k\geq 0} \frac{t^k}{k!} \asymp \frac{e^t}{\sqrt{t}}, \qquad \text{as } t \to +\infty.$$

The lemma claims that for the exponential function, $g(z) = e^z$, then

$$\mu(g,t) \asymp \frac{M(g,t)}{\ln^{1/2} M(g,t)}, \quad \text{as } t \to \infty\,.$$

*Proof.* We prove that there are constants $0 < C \leq D$ and $T > 0$ such that

$$C \sup_{k\geq 0} \frac{t^k}{k!} \leq \frac{e^t}{\sqrt{t}} \leq D \sup_{k\geq 0} \frac{t^k}{k!}, \qquad \text{for any } t \geq T.$$

Along this proof we will denote by $C$ different absolute constants.

Fix $t \geq 1$. Using Stirling's lower bound we have

$$\frac{t^k}{k!} \leq C \frac{(te)^k}{k^k \sqrt{k}}, \qquad \text{for any } k \geq 1.$$

Now we find the value $x > 0$ where the function $h(x)$, given by

$$h(x) = \frac{(te)^x}{x^x \sqrt{x}}, \qquad \text{for any } x > 0,$$

reach its maximum. This value is given by $\hat{x}$, the solution of $\ln(t) = \ln(\hat{x}) + 1/2\hat{x}$, which implies that $t = \hat{x} e^{-1/(2\hat{x})}$. This last equality gives $t < \hat{x}$ and also that $t \to \infty$ implies that $\hat{x} \to \infty$, therefore $t \sim \hat{x}$, as $t \to \infty$. In fact: $t = \hat{x} + 1/2 + O(1/\hat{x})$.

Using this last expression we find that

$$h(\hat{x}) = e^{-1/2} \frac{e^{\tilde{x}}}{\sqrt{\hat{x}}} \sim \frac{e^t}{\sqrt{t}}, \qquad \text{as } t \to +\infty,$$

and this implies that there exists certain $T > 0$ such that

$$\sup_{k\geq 0} \frac{t^k}{k!} \leq C \frac{e^t}{\sqrt{t}}, \qquad \text{for any } t > T.$$

Now we prove the lower bound. Fix $t > 1$ and take $n = \lfloor t \rfloor$. Stirling's asymptotic formula gives, for a different constant $C$, that

$$\frac{t^n}{n!} \geq C \frac{e^n t^n}{n^n \sqrt{n}} \geq C \frac{e^n}{\sqrt{n}} \geq C \frac{e^t}{\sqrt{t}},$$

for $t$ large enough, therefore there exists $T > 0$ such

$$\sup_{k\geq 0} \frac{t^k}{k!} \geq C \frac{e^t}{\sqrt{t}}, \qquad \text{for any } t \geq T.$$

□

### 1.5.5 Wiman-Valiron Theorem

The logarithmic measure of a set $A \subset (0, +\infty)$ is defined as $\int_A \frac{dx}{x}$. If $A = \{e^y : y \in B\}$, then $\int_A \frac{dx}{x} = \int_B dy = |B|$.

We restate Lemma 1.5.3 in terms of the maximum term of $f$ entire in $\mathcal{K}$ as

$$f(t) \le \frac{1}{H}\mu(f,t)(1+\sigma_f(t)), \quad \text{for any } t \ge 0\,. \tag{1.5.6}$$

where $H$ is some absolute constant.

**Theorem 1.5.10** (Wiman-Valiron)**.** *Let $g(z) = \sum_{n=0}^{\infty} b_n z^n$ be a non-constant entire power series. For each $\varepsilon > 0$, there exists a set $W = W_\varepsilon \subset (0, +\infty)$ of* finite logarithmic measure *such that*

$$\frac{M(g,t)}{\ln^{1/2+\epsilon} M(g,t)} \le \mu(g,t), \qquad \textit{for any } t \in (0,+\infty) \setminus W.$$

The following proof comes from [83].

*Proof.* Some preliminary reductions. Because of (1.5.3), we may assume that $g$ has non-negative Taylor coefficients. Also, because of (1.5.5), we may assume that $g$ is not a polynomial and thus that $\lim_{t\to\infty} \ln M(g,t)/\ln t = +\infty$. And, finally, because of (1.5.4), taking into account the limit above and by reducing $\varepsilon$ if necessary, we may assume that $g \in \mathcal{K}$.

We are thus to prove that for a given entire function $g$ in $\mathcal{K}$ and $\varepsilon > 0$, there is a set $W \subset (0, +\infty)$ of finite logarithmic measure so that

$$g(t) \le \mu(g,t) \ln^{1/2+\varepsilon} g(t)\,, \quad \text{for any } t \notin W\,.$$

We already know, from (1.5.6), that, for some absolute constant $H > 0$,

$$(\star_1) \quad g(t) \le \frac{1}{H}\mu(g,t)(1+\sigma_g(t))\,, \quad \text{for any } t \ge 0\,.$$

For the fulcrum $F$ of $g$ we have, from (1.5.1) of Proposition 1.5.7, that there exists $S \in \mathbb{R}$ and a set $E \subset [S, +\infty)$ of finite Lebesgue measure so that

$$F''(s) \le F(s)^{1+\varepsilon}, \quad \text{for } s \in [S,+\infty) \setminus E\,,$$

which translates into

$$(\star_2) \quad \sigma_g(e^s) \le \ln^{1/2+\varepsilon/2} g(e^s), \quad \text{for } s \in [S,+\infty) \setminus E\,.$$

Increasing the value $S$, if necessary, we have, since $\lim_{s\to\infty} g(e^s) = +\infty$, that

$$(\star_3) \quad \frac{1}{H}\big(1+\ln^{1/2+\varepsilon/2} g(e^s)\big) \le \ln^{1/2+\varepsilon} g(e^s)\,, \quad \text{for } s \in [S,+\infty) \setminus E\,.$$

Combining $(\star_1), (\star_2)$ and $(\star_3)$, we obtain that

$$g(e^s) \le \mu(g, e^s) \ln^{1/2+\varepsilon} g(e^s)\,, \quad \text{for } s \in [S,+\infty) \setminus E\,.$$

The set $W = \{t > e^S : \ln t \in E\}$ has finite logarithmic measure and we have that

$$g(t) \le \mu(g,t) \ln^{1/2+\varepsilon} g(t)\,, \quad \text{for } t \in [e^S,+\infty) \setminus W\,,$$

as desired. $\square$

### 1.5.6 Clunie's Theorem

**Theorem 1.5.11** (Clunie)**.** *Let $f(z)=\sum_{n=0}^{\infty} a_n z^n$ be an entire function in $\mathcal{K}$. For every $\varepsilon>0$ there exist $T>0$ and a set $C\subset[T,+\infty)$ of finite logarithmic measure such that for every $t\in[T,+\infty)\setminus C$ we have that*

$$\mathbf{P}\big(|X_t-m_f(t)|>y\big)\le 2\begin{cases} e^{-y^2/(2m_f(t)^{1+\varepsilon})}\,, & \textit{for } 0\le y\le m_f(t)\,,\\ \\ e^{-y/(2m_f(t)^{\varepsilon})}\,, & \textit{for } 0\le y\le m_f(t)\,.\end{cases}$$

In other terms, for $t\in[T,+\infty)\setminus C$, we have

$$\sum_{|k-m_f(t)|>y} a_n t^n \le 2f(t)\, e^{-y^2/(2m_f(t)^{1+\varepsilon})}\,, \quad \text{for } 0\le y\le m_f(t)\,.$$

The following proof of Theorem 1.5.11 is patterned after [88].

*Proof.* We may assume that $f$ is not a polynomial and thus that

$$M_f=\lim_{t\to\infty} m_f(t)=\lim_{s\to\infty} F'(s)=+\infty\,.$$

Here $F$ denotes the fulcrum of $f$.

From Proposition 1.5.7 with $\varepsilon/2$, we obtain $S\in\mathbb{R}$ and a set $E\subset[S,+\infty)$ of finite Lebesgue measure such that for $s\in[S,+\infty)\setminus E$ we have that

$$\Sigma(s,\Lambda)<F'(s+\Lambda)^{1+\varepsilon/2}\,, \quad \text{for } \Lambda>0 \text{ such that } s-\Lambda\ge S\,.$$

To reach this conclusion we have used that $F'(s)\le F'(s+\lambda)$, if $\lambda\ge 0$.

By increasing $S$ and augmenting $E$ if necessary, we have from the Borel-Hayman Lemma 1.5.8, that

$$F'\Big(s+\frac{1}{F'(s)^{\varepsilon}}\Big)<2F'(s)\,, \quad \text{for } s\in[S,+\infty)\setminus E\,.$$

Thus

$$\Sigma\Big(s,\frac{1}{F'(s)^{\varepsilon}}\Big)<\frac{2^{1+\varepsilon/2}}{F'(s)^{\varepsilon/2}}F'(s)^{1+\varepsilon}\,, \quad \text{for } s\notin E \text{ and such that } s-1/F'(s)\ge S\,.$$

Since $\lim_{s\to\infty} F'(s)=+\infty$, we may increase the value of $S$ again so that

$$\Sigma\Big(s,\frac{1}{F'(s)^{\varepsilon}}\Big)<F'(s)^{1+\varepsilon}\,, \quad \text{for } s\in[S,\infty)\setminus E\,.$$

Define $T=e^S$ and the set $C=\{e^s : s\in E\}$. The set $C$ has finite logarithmic measure.

Finally, we appeal to Proposition 1.4.6, as remarked in the formulation (1.4.20), to deduce that for any $t\in[T,+\infty)\setminus C$ it holds that

$$\mathbf{P}\big(|X_t-m_f(t)|>y\big)\le 2e^{-y^2/(2m_f(t)^{1+\varepsilon})}\,, \quad \text{for every } y \text{ such that } 0\le y\le m_f(t)\,,$$

and

$$\mathbf{P}\big(|X_t-m_f(t)|>y\big)\le 2e^{-y^2/(2m_f(t)^{\varepsilon})}\,, \quad \text{for every } y \text{ such that } y\ge m_f(t)\,.$$

□

Chapter 2

# Moments of a Khinchin family

## Contents

---

This chapter deals with the beautiful interplay between function theoretical properties of the holomorphic functions $f \in \mathcal{K}$ and certain probabilistic properties of their associated Khinchin families $(X_t)$ as $t \uparrow R$, in particular, the behavior as $t \uparrow R$ of the mean $m_f(t) = \mathbf{E}(X_t)$, of the variance $\sigma_f^2(t) = \mathbf{V}(X_t)$ or, in general, of the moments $\mathbf{E}(X_t^p)$, for $p > 0$.

For instance, in Section 1.2.5 we will provide simple proofs, *of probabilistic nature*, of some lower bounds and asymptotic lower bounds of $\sigma_f(t)$, as $t \uparrow R$, quantifying Hadamard's three lines theorem, due to Hayman [49], Bŏĭchuk and Gol'dberg [13], Abi-Khuzzam [2] and others. Later, see Theorem 2.5.1 and Proposition 2.5.3, we will show, for entire functions $f$ in $\mathcal{K}$, how the growth of the mean $m_f(t)$ and, in general, of the moments $\mathbf{E}(X_t^p)$, with $p > 0$, and of the quotient $\sigma_f^2(t)/m_f(t)$, relate to the order of $f$, generalizing some results of Pólya and Szegő, and of Báez-Duarte.

It turns out, Theorem 2.4.8, that for clans, $\lim_{t\uparrow R} \mathbf{E}(X_t^p)/\mathbf{E}(X_t)^p = 1$, for any $p > 0$. For instance, the ordinary generating function of the partitions of integers $P(z) = \prod_{k=1}^{\infty} 1/(1-z^k)$ is a clan, and for the family $(X_t)_{t>0}$ associated to $P$ we obtain readily, for any $p > 0$, that $\mathbf{E}(X_t^p) \sim \zeta(2)^p/(1-t)^{2p}$, as $t \uparrow 1$.

For entire functions in $\mathcal{K}$, Pólya and Szegő [81] showed that entire functions with non-negative coefficients of finite order $\rho$ such that $\lim_{t\to\infty} \ln f(t)/t^\rho$ exists and is positive are clans; we show, more generally, see Theorem 2.5.9 in Section 2.5, that entire functions in $\mathcal{K}$ of regular growth (in a precise sense) are clans.

The entire gap series with non-negative coefficients presented in Section 2.5.2 furnish the basic examples of entire functions of any order $\rho$, with $0 \leq \rho \leq +\infty$, which are not clans. A classical result of Pfluger and Pólya [80] ensures that these entire gap series have no Borel exceptional values. We show in Theorem 2.5.10 that entire functions in $\mathcal{K}$ with one Borel exceptional value are always clans.

This chapter is mainly based on the paper:

- Maciá, V.J. et al. Growth of power series with nonnegative coefficients, and Moments of power series distributions. (submitted) arXiv:2401.14473v2, see [19].

## 2.1 Moments

In this section we compute the moments of any order of a Khinchin family $(X_t)_{t\in[0,R)}$. For integers $n \geq 0$, these moments admit closed expressions in terms of $f$ itself, but these expressions are much more intricate than those for the mean $m_f(t)$ and the variance $\sigma_f^2(t)$.

Using the operator $\mathcal{D}$, recall that $\mathcal{D}f(z) = zf'(z)$, see equation (1.3.1), we have that

$$\mathbf{E}(X_t^q) = \frac{\mathcal{D}^q f(t)}{f(t)} = \sum_{n=0}^{\infty} n^q \frac{a_n t^n}{f(t)}, \quad \text{for any } t \in [0, R).$$

To express these moments in terms of $f$ and some of its derivatives we make use of falling factorials. Stirling number of the second kind allow to express the powers $x^n$ in terms of falling factorials $x(x-1)\dots(x-j+1)$. Recall that the Stirling number of the second kind $S(n,j)$ counts the number of partitions of a set of size $n$ in $j$ (non-empty) subsets. In fact we have:

$$x^n = \sum_{j=0}^{n} S(n,j) x(x-1)\dots(x-j+1), \quad \text{for any } n \geq 1.$$

Here $S(n,0) = 0$, for any $n \geq 1$. For $n = 0$ we stick to the conventions that $S(0,0) = 1$, $x^0 = 1$ and also that $x(x-1)\dots(x-j+1)$ if $j = 0$. For more details see, for instance, [REF Stirling Second Kind - Flajolet].

We will use the notation $x^{\underline{j}} = x(x-1)\dots(x-j+1)$ to denote the $j$-th falling factorial of $x$ and, analogously, $X^{\underline{j}} = X(X-1)\dots(X-j+1)$ to denote the $j$-th falling factorial of a random variable $X$. With this notation in mind, for any Khinchin family $(X_t)$, we have

$$X_t^n = \sum_{j=0}^{n} S(n,j) X_t^{\underline{j}}, \quad \text{for any } n \geq 1 \text{ and } t \in [0, R). \tag{2.1.1}$$

For any $j \geq 0$, the *$j$-th factorial moment* of the Khinchin family $(X_t)$ is given by

$$\mathbf{E}(X_t^{\underline{j}}) = \sum_{n=j}^{\infty} n^{\underline{j}} \frac{a_n t^n}{f(t)} = t^j \frac{f^{(j)}(t)}{f(t)}, \quad \text{for any } t \in (0, R). \tag{2.1.2}$$

Using (2.1.1) and (2.1.2), we conclude that

$$\mathbf{E}(X_t^n) = \sum_{j=0}^{n} S(n,j) t^j \frac{f^{(j)}(t)}{f(t)}, \quad \text{for any } n \geq 0 \text{ and } t \in [0, R).$$

We collect this expression by means of the following Lemma.

**Lemma 2.1.1.** *Let $f \in \mathcal{K}$ be power series with radius of convergence $R > 0$ and denote $(X_t)_{[0,R)}$ its Khinchin family, then*

$$\mathbf{E}(X_t^n) = \sum_{j=0}^{n} S(n,j) t^j \frac{f^{(j)}(t)}{f(t)},$$

*for any $n \geq 0$ and $t \in [0, R)$.*

**Some examples: moments of the most basic families**

- For $f(z) = 1 + z$ and its Bernoulli family we have
$$\mathbf{E}(X_t^n) = S(n,1)\frac{t}{1+t} = \frac{t}{1+t}, \qquad \text{for any } t \geq 0 \text{ and } n \geq 1,$$
Observe that in this case $X_t^n \equiv X_t$, for any $n \geq 1$.

- For $f(z) = 1/(1-z)$ and its Pascal family we have
$$\mathbf{E}(X_t^n) = \sum_{j=1}^{n} S(n,j)j!\frac{t^j}{(1-t)^j}, \text{ for any } t \in [0,1) \text{ and } n \geq 1.$$
In particular, for $t = 1/2$, we have that $n$-th moment of the geometric variable of parameter $1/2$ is $\tilde{B}_n$, the $n$-th ordered Bell's number.

- For $f(z) = e^z$ and its Poisson family we have
$$\mathbf{E}(X_t^n) = \sum_{j=1}^{n} S(n,j)t^j, \qquad \text{for any } t \in [0,1) \text{ and } n \geq 1.$$
In particular, for $t = 1$, we have that the $n$-th moment of the Poisson variable of parameter 1 is $\mathcal{B}_n$, the $n$-th Bell number.

### 2.1.1 Central and absolute central moments

The expression for the general central moments in terms of $f$ is a little bit more involved:

$$\begin{aligned}\mathbf{E}((X_t - m(t))^k) &= \sum_{l=0}^{k} \binom{k}{l} (-1)^{k-l} m(t)^{k-l} \mathbf{E}(X_t^l) \\ &= \sum_{0\leq j\leq l\leq k} \binom{k}{l} S(l,j)(-1)^{k-l} m(t)^{k-l} t^j \frac{f^{(j)}(t)}{f(t)}\end{aligned} \tag{2.1.3}$$

Using the previous equality we can find upper bounds for the absolute central moments, these moments are difficult to compute in terms of $f$.

We will use the following, special, notations for the third and fourth absolute central moments:

$$\gamma^3(t) = \mathbf{E}(|X_t - m(t)|^3) \quad \text{and} \quad \kappa^4(t) = \mathbf{E}(|X_t - m(t)|^4), \tag{2.1.4}$$

for any $t \in [0, R)$.

Jensen's inequality implies that

$$\gamma(t) \leq \kappa(t), \qquad \text{for any } t \in [0, R).$$

This inequality will prove useful in determining upper bounds for $\gamma(t)$. While the fourth central moment, $\kappa^4(t)$, can be expressed in a closed form, the third absolute central moment, $\gamma^3(t)$, remains, a priori, unmanageable (at least for $t$ close to $R$, see equation (2.2.2), for the expression of $\gamma(t)^3$, for $t$ close to 0).

### 2.1.2 The moments determine the family

The following proposition claims that Khinchin families are determined by their moments at any $t = t_0$.

**Proposition 2.1.2.** *Let $(X_t)$ and $(Y_t)$ be the Khinchin families associated to $f, g \in \mathcal{K}$, respectively. Assume that $f$ and $g$ have radius of convergence $R, S$, respectively, and also that for $t_0 \in (0, \min\{R, S\})$ we have that*

$$\mathbf{E}(X_{t_0}^k) = \mathbf{E}(Y_{t_0}^k), \quad \text{for any } k \geq 0,$$

*then $R = S$ and there exists $\lambda > 0$ such that $f(z) = \lambda g(z)$, for any $z \in \mathbb{D}(0, R)$. In fact $X_t \stackrel{d}{=} Y_t$, for any $t \in (0, R)$.*

*Proof.* We have that $\mathbf{E}(X_{t_0}^k) = \mathbf{E}(Y_{t_0}^k)$, for any $k \geq 0$, this implies that all the factorial moments of $X_{t_0}$ and $Y_{t_0}$ are equal. Equation (2.1.2) gives that

$$\frac{f^{(k)}(t_0)}{f(t_0)} = \frac{g^{(k)}(t_0)}{g(t_0)}, \quad \text{for any } k \geq 0.$$

The previous equality implies that, for $z$ in a small disk around $z = t_0$, we have

$$\sum_{k=0}^{\infty} f^{(k)}(t_0)(z - t_0)^k = \frac{f(t_0)}{g(t_0)} \sum_{k=0}^{\infty} g^{(k)}(t_0)(z - t_0),$$

then both power series are identical in a small disk around $z = t_0$ and the identity principle gives that $R = S$ and $f(z) = \lambda g(z)$, for any $z \in \mathbb{D}(0, R)$, with $\lambda = f(t_0)/g(t_0)$. □

### 2.1.3 Recurrence for the moments and the central moments

For a Khinchin family $(X_t)_{t \in [0,R)}$ and any integer $r \geq 0$ denote

$$\mu_r(t) = \mathbf{E}(X_t^r) \quad \text{and} \quad \tilde{\mu}_r(t) = \mathbf{E}((X_t - m_f(t))^r).$$

Thus $\mu_1(t) = m_f(t)$ and $\tilde{\mu}_2(t) = \sigma_f^2(t)$. We add the convention that $\mu_0(t) \equiv 1 \equiv \tilde{\mu}_0(t)$, for every $t \in (0, R)$.

For the moments and the central moments we have the following differential recurrences.

**Lemma 2.1.3.** *With the same notations:*

$$\mu_{r+1}(t) = m_f(t)\mu_r(t) + t\mu_r'(t), \quad \text{for } t \in [0, R) \text{ and } r \geq 0.$$

*and*

$$\tilde{\mu}_{r+1}(t) = t\tilde{\mu}_r'(t) + r\tilde{\mu}_{r-1}(t)\sigma_f^2(t), \quad \text{for } t \in [0, R) \text{ and any integer } r \geq 1.$$

*Proof.* Taking derivatives in

$$\mu_r(t) = \frac{1}{f(t)} \sum_{n=1}^{\infty} n^r a_n t^n$$

we deduce that

$$\mu_{r+1}(t) = m_f(t)\mu_r(t) + t\mu_r'(t), \qquad \text{for any } t \in [0, R) \text{ and } r \geq 0.$$

For the centered moments we have that

$$\tilde{\mu}_r(t) = \sum_{n=0}^{\infty} \frac{a_n t^n}{f(t)} (n - m_f(t))^r, \qquad \text{for any } t \in [0, R) \text{ and any integer } r \geq 0\,,$$

taking derivatives we find that

$$\begin{aligned} \tilde{\mu}_{r+1}(t) &= t(\tilde{\mu}_r'(t) + r\tilde{\mu}_{r-1}(t)m_f'(t)) \\ &= t\tilde{\mu}_r'(t) + r\tilde{\mu}_{r-1}(t)\sigma_f^2(t), \end{aligned}$$

for $t \in [0, R)$ and any integer $r \geq 1$. □

### 2.1.4 Sequences of power series and moments

Assume that $(f^{[k]})_{k\geq 1}$ is a sequence of functions in $\mathcal{K}$, all of them with radius of convergence at least $R > 0$ and also that $f^{[k]}$ converges uniformly on compact subsets to a function $f \in \mathcal{K}$, with radius of convergence $R > 0$.

For each $k \geq 1$ denote $(X_t^{[k]})$ the Khinchin family of $f^{[k]}$ and by $(X_t)$ the Khinchin family associated to $f$. In the same vein denote $m_k$, $\sigma_k^2$ and $m_f(t)$ and $\sigma_f^2(t)$ the respective mean and variance.

We already know, see Chapter 1, that for each $t \in (0, R)$ the sequence of random variables $(X_t^{[k]})$ converges in distribution towards a variables $X_t$.

For $\alpha > 0$, we denote $m_{\alpha,k}(t) = \mathbf{E}\left(\left(X_t^{[k]}\right)^\alpha\right)$ and $m_\alpha(t) = \mathbf{E}(X_t^\alpha)$. For moments of integer order $q$ we have

$$\mu_q(t) = \frac{1}{f(t)} \mathcal{D}^q f(t),$$

see equation (1.3.1) for the definition of $\mathcal{D}$, and thus

$$m_{q,k}(t) = \frac{1}{f^{[k]}(t)} \mathcal{D}^q f^{[k]}(t).$$

converges to $m_\alpha(t)$, for $t \in (0, R)$.

**Proposition 2.1.4.** *With the same notation than above, for each $t \in (0, R)$, and each $\alpha > 0$ we have that*

$$\lim_{k\to+\infty} \mathbf{E}\left(\left(X_t^{[k]}\right)^\alpha\right) = \mathbf{E}(X_t^\alpha).$$

*Proof.* Fix $t \in (0, R)$ and $\alpha > 0$. Let $\varepsilon > 0$ be given and fix $r \in (t, R)$. Take $K \geq 1$ such that

$$|f^{[k]}(z) - f(z)| \leq \varepsilon, \quad \text{if } |z| = r \text{ and } k \geq K.$$

Cauchy's formula gives that

$$|a_n^{[k]} - a_n| \leq \frac{\varepsilon}{r^n}, \quad \text{for each } n \geq 0 \text{ and } k \geq K.$$

Fixing a bigger $K$, if necessary, we can assume that

$$\left|\frac{1}{f^{[k]}(t)} - \frac{1}{f(t)}\right| \leq \varepsilon, \quad \text{and} \quad \frac{1}{f^{[k]}(t)} \leq \frac{2}{f(t)}, \quad \text{for any } k \geq K.$$

We write

$$\begin{aligned}\mathbf{E}\left(\left(X_t^{[k]}\right)^\alpha\right) - \mathbf{E}(X_t^\alpha) &= \sum_{n=0}^{\infty} \frac{n^\alpha a_n^{[k]} t^n}{f^{[k]}(t)} - \sum_{n=0}^{\infty} \frac{n^\alpha a_n t^n}{f(t)} \\ &= \sum_{n=0}^{\infty} \frac{(a_n^{[k]} - a_n) n^\alpha t^n}{f^{[k]}(t)} + \sum_{n=0}^{\infty} a_n n^\alpha t^n \left(\frac{1}{f^{[k]}(t)} - \frac{1}{f(t)}\right).\end{aligned}$$

Estimating, for $k \geq K$:

$$\left|\mathbf{E}\left(\left(X_t^{[k]}\right)^\alpha\right) - \mathbf{E}(X_t^\alpha)\right| \leq \varepsilon \left(\frac{2}{f(t)} \sum_{n=0}^{\infty} n^\alpha \left(\frac{t}{r}\right)^n + \sum_{n=0}^{\infty} a_n n^\alpha t^n\right).$$

□

**Corollary 2.1.5.** *With the previous notation, for any $t \in (0, R)$ we have*

$$\lim_{k \to \infty} m_k(t) = m_f(t), \quad \text{and} \quad \lim_{k \to \infty} \sigma_k^2(t) = \sigma_f^2(t).$$

Another corollary, which will be useful later on, register the convergence of the even centered moments of $\left(X_t^{[k]}\right)$ to the even centered moments of $(X_t)$.

**Corollary 2.1.6.** *With the previous notations, for $t \in (0, R)$ and any integer $j \geq 1$ we have*

$$\lim_{k \to \infty} \mathbf{E}\left(|X_t^{[k]} - m_k(t)|^{2j}\right) = \mathbf{E}(|X_t - m_f(t)|^{2j}).$$

*Proof.* It is enough to observe that $\mathbf{E}(|X_t - m_f(t)|^2 j) = \mathbf{E}((X_t - m(t))^{2j})$ is written as a linear combination of terms (a total of $2j + 1$ terms) of the form $\mathbf{E}(X_t^s) m_f(t)^{2j-s}$.

□

## 2.2 Moments of a Khinchin family near 0

Let $f(z) = \sum_{n=0}^{\infty} a_n z^n$ be a power series in $\mathcal{K}$ with radius of convergence $R > 0$. Let $(X_t)_{t \in [0,R)}$ be its Khinchin family and $m_f(t)$ and $\sigma_f^2(t)$ be, respectively, the mean and variance of $X_t$, for $t \in [0, R)$.

### 2.2.1 Analycity of the moments near 0

Let's study the moments $\mathbf{E}(X_t^k)$, for $t$ close to 0. Take $S \in (0, R)$ such that $f(z) \neq 0$ for any $z$ in the open disk $\mathbb{D}(0, S)$. For any integer $k \geq 1$ we have

$$\mathbf{E}(X_t^k) = \frac{1}{f(t)} \sum_{j=0}^{\infty} a_j j^k t^j.$$

Therefore, the moments $\mathbf{E}(X_t^k)$ are the restriction to the interval $[0, S)$ of the holomorphic function

$$z \mapsto \frac{1}{f(z)} \sum_{j=0}^{\infty} a_j j^k z^j, \qquad \text{for any } z \in \mathbb{D}(0, S).$$

This implies that the moments $\mathbf{E}(X_t^k)$ and the central moments $\mathbf{E}((X_t - m(t))^k)$, see equation (2.1.3), are analytic functions (i.e the restriction of certain analytic function to the real interval $(0, S)$).

The absolute central moments $\mathbf{E}(|X_t - m(t)|^k)$ are also analytic for $t \in (0, S')$ (that is, these moments are given by the restriction of an analytic function to the interval $(0, S')$). The argument justifying this is the following: the mean $m(t)$ is an increasing function; we distinguish two cases: there exists $\tau \in (0, R)$ such that $m(\tau) = 1$ (this is the case $M_f > 1$) or $m(t) < 1$ for any $t \in (0, R)$; in this case we choose $\tau = R$.

For any $t \in [0, \tau)$, the mean, $m(t)$, verifies the inequality $m(t) < 1$, and therefore, in this interval, we have

$$\begin{aligned} \mathbf{E}(|X_t - m(t)|^k) &= m(t)^k \frac{a_0}{f(t)} + \sum_{n=1}^{\infty} (n - m(t))^k \frac{a_n t^n}{f(t)} \\ &= \mathbf{E}((X_t - m(t))^k) + (1 + (-1)^{k+1}) m(t)^k \frac{a_0}{f(t)}, \qquad \text{for any } t \in (0, \tau). \end{aligned} \tag{2.2.1}$$

We take $S' = \min\{\tau, S\}$ and then the absolute central moment $\mathbf{E}(|X_t - m(t)|^k)$ is the restriction of certain holomorphic function in $\mathbb{D}(0, S')$ to the open interval $(0, S')$.

An immediate consequence of the previous analysis is that, for $t$ close to $s = 0$, we can find an explicit expression for $\gamma^3(t)$. Equation (2.2.1) gives that

$$\gamma^3(t) = \mathbf{E}((X_t - m(t))^3) + 2m(t)^3 \frac{a_0}{f(t)}, \qquad \text{for any } t \in (0, S'). \tag{2.2.2}$$

### 2.2.2 Asymptotic behavior of the moments near 0

We study now the asymptotic behavior of the moments for $t$ close to zero. Let $f \in \mathcal{K}$ be a power series with radius of convergence $R > 0$. Assume first that $f'(0) = a_1 > 0$. In this case, we have that

$$m_f(t) = \frac{a_1}{a_0} t + O(t^2) \quad \text{and} \quad \sigma_f^2(t) = \frac{a_1}{a_0} t + O(t^2)\,, \quad \text{as } t \downarrow 0\,,$$

and, besides, that

$$\mathbf{E}(|X_t - m_f(t)|^3) = \frac{a_1}{a_0} t + O(t^2) \quad \text{and} \quad \sqrt{t}\, \mathbf{E}(|\breve{X}_t|^3) = \sqrt{\frac{a_0}{a_1}} + O(t)\,, \quad \text{as } t \downarrow 0\,.$$

In general, if for $k \geq 1$, we have $a_k > 0$, but $a_j = 0$, for $1 \leq j < k$, then

$$m_f(t) = \frac{k a_k}{a_0} t^k + O(t^{k+1}) \quad \text{and} \quad \sigma_f^2(t) = \frac{k^2 a_k}{a_0} t^k + O(t^{k+1})\,, \quad \text{as } t \downarrow 0\,,$$

and, besides, that

$$\mathbf{E}(|X_t - m_f(t)|^3) = \frac{k^3 a_k}{a_0} t^k + O(t^{k+1}) \quad \text{and} \quad t^{k/2}\, \mathbf{E}(|\breve{X}_t|^3) = \sqrt{\frac{a_0}{a_k}} + O(t)\,, \quad \text{as } t \downarrow 0\,.$$

## 2.3 Asymptotic behavior of the moments of the most basic families

For $\beta \geq 0$ and $t \in [0, R)$, *the moment* $\mathbf{E}(X_t^\beta)$ *of exponent* $\beta$ *of* $X_t$ maybe be written in terms of the power series of $f$ as

$$\mathbf{E}(X_t^\beta) = \sum_{n=0}^{\infty} n^\beta \frac{a_n t^n}{f(t)}\,.$$

We are particularly interested in the comparison of the moment $\mathbf{E}(X_t^\beta)$ of exponent $\beta$ with the mean $\mathbf{E}(X_t)$, which acts as a sort of unit, as $t \uparrow R$.

We would like to consider also *moments of negative exponents* of the $X_t$. But since $P(X_t = 0) = a_0/f(t) > 0$, we consider instead of $X_t$ the variable $Y_t$ *which is just* $X_t$ *conditioned on being positive*: $Y_t \triangleq (X_t | X_t \geq 1)$, which takes values in $\{1, 2, \ldots\}$:

$$\mathbf{P}(Y_t = n) = \frac{a_n t^n}{f(t) - a_0} = \mathbf{P}(X_t = n) \frac{f(t)}{f(t) - a_0}\,, \quad \text{for } n \geq 1 \text{ and } t \in (0, R)\,.$$

These variables $Y_t$, conditioned versions of the $X_t$, are defined for $t \in (0, R)$.

We describe next the asymptotics of the moments of the most basic families.

### 2.3.1 Moments of geometric and negative binomial variables

For the geometric family associated to $f(z) = 1/(1 - z)$ we have

**Proposition 2.3.1.** *If* $(X_t)_{0 \leq t < 1}$ *is the Khinchin family of* $1/(1 - z)$*, then for* $\beta \geq 0$

$$\lim_{t \uparrow 1} \frac{\mathbf{E}(X_t^\beta)}{\mathbf{E}(X_t)^\beta} = \Gamma(\beta + 1).$$

Observe that for this family $\mathbf{E}(X_t) = t/(1-t)$, for $0 \le t < 1$, and thus $\mathbf{E}(X_t) \sim 1/(1-t)$ as $t \to 1$.

This estimate of moments is actually equivalent to the following growth estimate of certain power series.

**Lemma 2.3.2.** *For $\beta > 0$, we have* $\displaystyle\sum_{n=1}^{\infty} n^{\beta-1} t^n \sim \Gamma(\beta)\frac{1}{(1-t)^\beta}\,,\quad$ *as* $t \uparrow 1\,.$

*Proof.* From the binomial expansion, it follows that

$$\frac{1}{(1-z)^\beta} = \sum_{n=0}^{\infty} \frac{\Gamma(n+\beta)}{\Gamma(\beta)\, n!} z^n\,, \quad \text{for } |z| < 1\,,$$

and form Stirling's formula that

$$\frac{\Gamma(n+\beta)}{\Gamma(\beta)\, n!} \sim \frac{n^{\beta-1}}{\Gamma(\beta)}\,, \quad \text{as } n \to \infty\,, \tag{2.3.1}$$

□

For moments of negative exponents for the conditioned versions $Y_t$ of the $X_t$ of the geometric family, we have, as $t \to 1$,

$$\mathbf{E}(Y_t^{-\beta}) \sim \begin{cases} \zeta(\beta)(1-t), & \text{for } \beta > 1\,, \\ (1-t)\ln \frac{1}{1-t}, & \text{for } \beta = 1\,, \\ \Gamma(1-\beta)(1-t)^\beta, & \text{for } \beta \in [0,1)\,. \end{cases}$$

The following estimates of moments of the negative binomial of parameter $N$ associated to $1/(1-z)^N$ are a consequence again of Lemma 2.3.2.

**Proposition 2.3.3.** *Let $N$ be an integer $N \ge 1$ and let $(X_t)_{0\le t<1}$ be the Khinchin family of $1/(1-z)^N$. For $\beta > 0$,*

$$\lim_{t\uparrow 1} \frac{\mathbf{E}(X_t^\beta)}{\mathbf{E}(X_t)^\beta} = \frac{\Gamma(\beta+N)}{\Gamma(N)N^\beta}.$$

For this negative binomial family we have that $\mathbf{E}(X_t) = Nt/(1-t)$, for $t \in [0,1)$.

As for moments of negative exponents for the conditioned versions $Y_t$ of the $X_t$ of this negative binomial family of parameter $N$, we have, as $t \to 1$,

$$\mathbf{E}(Y_t^{-\beta}) \sim \begin{cases} \dfrac{\zeta(\beta-N+1)}{\Gamma(N)}\,(1-t)^N, & \text{for } \beta > N\,, \\ \dfrac{1}{\Gamma(N)}\,(1-t)^N \ln \dfrac{1}{1-t}, & \text{for } \beta = N\,, \\ \dfrac{\Gamma(N-\beta)}{\Gamma(N)}\,(1-t)^\beta, & \text{for } \beta \in [0,N)\,. \end{cases}$$

### 2.3.2 Moments of Bernoulli and binomial variables

For the Khinchin family $(X_t)_{0\le t<\infty}$ of a *polynomial* $f$ *of degree* $N$ we have that $X_t$ tends in distribution to the constant $N$ as $t\to\infty$ and also that for $\beta\ge 0$:

$$\lim_{t\to\infty}\mathbf{E}(X_t^\beta)=N^\beta=\lim_{t\to\infty}\mathbf{E}(X_t)^\beta\,.$$

Particular instances of this polynomial case are the Bernoulli family and the binomial family which are the Khinchin families, respectively, of $f(z)=1+z$ (of degree 1) and of $f(z)=(1+z)^N$ (of degree $N$).

For moments of negative exponents with the conditioned versions $Y_t=\mathbf{E}(X_t|X_t\ge 1)$ we have that

$$\lim_{t\to\infty}\mathbf{E}(Y_t^{-\beta})=N^{-\beta}=\lim_{t\to\infty}\mathbf{E}(Y_t)^{-\beta}\,.$$

### 2.3.3 Moments of Poisson variables

The Poisson family $(X_t)_{t\ge 0}$ is the Khinchin family associated to the exponential function $e^z$. For $t>0$, the variable $X_t$ is a Poisson variable with parameter $t$, and thus $\mathbf{E}(X_t)=t$.

For $\beta>0$ and $t>0$ the moment of exponent $\beta$ of $X_t$ is given by

$$\mathbf{E}(X_t^\beta)=\sum_{n=0}^\infty n^\beta\frac{t^n e^{-t}}{n!}\,.$$

The relevant power series comparison for the Poisson family corresponding to Lemma 2.3.2 for the case of the geometric family is the following proposition.

**Proposition 2.3.4.** *For $\beta\ge 0$, we have that*

$$\sum_{n=0}^\infty n^\beta\frac{t^n}{n!}\sim t^\beta e^t\,,\qquad \textit{as } t\to\infty\,,$$

*and*

$$\sum_{n=1}^\infty \frac{1}{n^\beta}\frac{t^n}{n!}\sim \frac{e^t}{t^\beta}\,,\qquad \textit{as } t\to\infty\,.$$

The Proposition 2.3.4, of analytic character, is equivalent to the combination of the propositions 2.3.7 and 2.3.8 which record the asymptotics results of the moments of the Poisson variables. Actually, we shall derive these results about moments probabilistically and not using Proposition 2.3.4 which we leave as a corollary.

**Remark 2.3.5** (Analytic proof of Proposition 2.3.4)**.** A proof of analytic character is the following: upon derivation of $e^t=\sum_{n=0}^\infty t^n/n!$, we deduce that for any integer $k\ge 1$,

$$\sum_{n=0}^\infty\frac{n(n-1)\dots(n-k+1)}{n!}t^n=t^k e^t,\qquad \text{for any } t\in\mathbb{R}.$$

From there, it follows that

$$\sum_{n=0}^{\infty} \frac{n^k}{n!} t^n \sim t^k e^t, \qquad \text{as } t \to \infty.$$

Convexity and Jensen's inequality imply that, for any $\beta > 0$,

$$\sum_{n=0}^{\infty} \frac{n^\beta}{n!} t^n \sim t^\beta e^t, \qquad \text{as } t \to \infty.$$

⊠

An application of Jensen's inequality will allow to reduce some discussions below to the case of integer exponents. This lemma will also be useful in other sections.

**Lemma 2.3.6.** *Let $(U_t)_{0 \le t < R}$ and $0 < R \le \infty$, be a family of non-negative random variables such that for some $p > 1$,*

$$\mathbf{E}(U_t^p) \sim \mathbf{E}(U_t)^p\,, \ \textit{as } t \uparrow R.$$

*Then, for any $\beta \in (0, p]$,*

$$\mathbf{E}(U_t^\beta) \sim \mathbf{E}(U_t)^\beta\,, \ \textit{as } t \uparrow R.$$

*Proof.* Let $1 < \beta \le p$. We have from Jensen's inequality that

$$\mathbf{E}(U_t) \le \mathbf{E}(U_t^\beta)^{1/\beta} \le \mathbf{E}(U_t^p)^{1/p}\,.$$

The result follows in this case, since, by hypothesis, $\mathbf{E}(U_t^p)^{1/p} \sim \mathbf{E}(U_t)$, as $t \uparrow R$.

For $\beta \in (0, 1)$, consider $u \in (0, 1)$ such that $1 = \beta u + p(1 - u)$. Jensen's inequality gives that

$$\mathbf{E}(U_t^\beta) \le \mathbf{E}(U_t)^\beta\,.$$

By Hölder's inequality

$$\mathbf{E}(U_t) \le \mathbf{E}(U_t^\beta)^u\, \mathbf{E}(U_t^p)^{1-u}\,;$$

The result follows from these two inequalities, since, by hypothesis and the definition of $u$, we have that $\mathbf{E}(U_t^p)^{1-u} \sim \mathbf{E}(U_t)^{1-\beta u}$, as $t \uparrow R$. □

We split the discussion of the moments of $X_t$ into those of positive exponents, Proposition 2.3.7, and those of negative exponents, Proposition 2.3.8, which together prove Proposition 2.3.4.

• *Moments of positive exponent.*

**Proposition 2.3.7.** *For the family $(X_t)_{0 \le t < \infty}$ associated to $f(z) = e^z$, and for any $\beta > 0$, we have*

$$\mathbf{E}(X_t^\beta) \sim t^\beta = \mathbf{E}(X_t)^\beta\,, \qquad \textit{as } t \to \infty\,.$$

This statement is the first part of Proposition 2.3.4.

*Proof.* We will first show that the result holds when $\beta = k$ is an integer and after that we will extend it to any $\beta > 0$ by means of Lemma 2.3.6.

Notice that $\mathbf{E}(X_t) = m_f(t) = t$ and $\mathbf{V}(X_t) = \sigma_f^2(t) = t$, for $t > 0$, and thus $\mathbf{E}(X_t^2) = t^2 + t$, for $t > 0$, which means that the result holds for $k = 2$. For $k = 1$, it holds trivially.

Let $k \geq 1$ be an integer. Upon differentiating $k$ times the identity $e^t = \sum_{n=0}^{\infty} t^n/n!$, we see that

$$\begin{aligned}\mathbf{E}\big(X_t(X_t-1)\dots(X_t-k+1)\big) &= \sum_{n=k}^{\infty} n(n-1)\cdots(n-k+1)\frac{t^n e^{-t}}{n!}\\ &= t^k e^{-t}\sum_{n=k}^{\infty} n(n-1)\cdots(n-k+1)\frac{t^{n-k}}{n!}\\ &= t^k\,, \quad \text{for } t \geq 0\,.\end{aligned}$$

For integer $k \geq 1$, we may write the monomials $x^k$ in terms of the falling factorials $x(x-1)\dots(x-j+1)$ and Stirling numbers of the second kind as

$$x^k = \sum_{j=0}^{k} S(k,j)x(x-1)\dots(x-j+1)$$

and thus the $k$th moment of $X_t$ is given by the polynomial

$$\mathbf{E}(X_t^k) = \sum_{j=0}^{k} S(k,j)t^j\,.$$

Since $S(k,k) = 1$, for $k \geq 1$, we see that for each integer $k \geq 1$

$$\mathbf{E}(X_t^k) \sim t^k, \quad \text{as } t \to \infty\,.$$

Lemma 2.3.6 finishes the proof. □

- *Moments of negative exponent.*

Consider next the variables $Z_t$ which are the reciprocals of the truncated Poisson variables:

$$Z_t = \begin{cases} 1/X_t\,, & \text{if } X_t \geq 1\,,\\ 0\,, & \text{if } X_t = 0\,,\end{cases} \quad \text{for } t > 0\,.$$

We consider now moments of positive exponents of the reciprocals $(Z_t)_{t>0}$. Observe that for $\beta > 0$,

$$\mathbf{E}(Z_t^\beta) = e^{-t}\sum_{n=1}^{\infty}\frac{1}{n^\beta}\frac{t^n}{n!}\,, \quad \text{for any } t > 0\,.$$

**Proposition 2.3.8.** *For $\beta > 0$,*

$$\mathbf{E}(Z_t^\beta) \sim \frac{1}{t^\beta}\,, \quad \textit{as } t \to \infty\,.$$

In other terms,

$$\sum_{n=1}^{\infty} \frac{1}{n^\beta} \frac{t^n}{n!} \sim \frac{e^t}{t^\beta}\,, \quad \text{as } t \to \infty\,,$$

which is the second half of Proposition 2.3.4. For precise asymptotic *expansions* of the negative moments of Poisson variables we refer to [55].

For the variables $Y_t = (X_t | X_t \geq 1)$ and $\beta > 0$, we have that

$$\mathbf{E}(Y_t^{-\beta}) = \mathbf{E}(Z_t^\beta) \frac{e^t}{e^t - 1}\,, \quad \text{for } t > 0\,,$$

and thus we also have

$$\mathbf{E}(Y_t^{-\beta}) \sim \frac{1}{t^\beta}\,, \quad \text{as } t \to \infty\,.$$

*Proof.* On account of Lemma 2.3.6, to prove Proposition 2.3.8 it suffices to consider integer exponents $\beta > 0$.

For $t > 0$, the first moment, $\mathbf{E}(Z_t)$, is given by

$$\mathbf{E}(Z_t) = e^{-t} \sum_{n=1}^{\infty} \frac{1}{n} \frac{t^n}{n!}\,.$$

Now taking the derivative with respect to $t$ of $\mathbf{E}(Z_t) e^t$, we see that

$$\frac{d}{dt} \sum_{n=1}^{\infty} \frac{1}{n} \frac{t^n}{n!} = \frac{e^t - 1}{t}\,,$$

for $t > 0$. Hence,

$$\mathbf{E}(Z_t) = e^{-t} \int_0^t \frac{e^s - 1}{s}\, ds = e^{-t} \alpha_1(t),$$

where $\alpha_1(t) = \int_0^t (e^s - 1)/s\, ds$.

The second moment, $\mathbf{E}(Z_t^2)$, is given by

$$\mathbf{E}(Z_t^2) = e^{-t} \sum_{n=1}^{\infty} \frac{1}{n^2} \frac{t^n}{n!}\,.$$

and similarly,

$$\frac{d}{dt} \sum_{n=1}^{\infty} \frac{1}{n^2} \frac{t^n}{n!} = \sum_{n=1}^{\infty} \frac{1}{n} \frac{t^{n-1}}{n!} = \frac{\alpha_1(t)}{t}.$$

Thus, if $\alpha_2(t) = \int_0^t \alpha_1(s)/s\, ds$ then for $t > 0$,

$$\mathbf{E}(Z_t^2) = e^{-t} \alpha_2(t)\,.$$

In general, for any integer $k \geq 1$,

$$\mathbf{E}(Z_t^k) = e^{-t} \sum_{n=1}^{\infty} \frac{1}{n^k} \frac{t^n}{n!} =: e^{-t} \alpha_k(t)\,, \quad \text{for } t > 0\,,$$

with $\alpha_k(t) = \int_0^t \alpha_{k-1}(s)/s\, ds$ and $\alpha_0(t) = e^t - 1$. To finish the proof it is enough to show that for any integer $k \geq 0$,

$$\alpha_k(t) \sim \frac{e^t}{t^k}\,, \quad \text{as } t \to \infty\,.$$

By induction, we readily see that $\lim_{t\downarrow 0} \alpha_k(t)/t = 1$, for $k \geq 0$. We prove the asymptotic formula also by induction. Let $k \geq 1$, and assume that the formula holds for $k-1$. (For $k = 0$, it holds trivially.)

For $M > 0$ we let

$$\varepsilon_M = \sup_{t \geq M} \left| \alpha_{k-1}(t) t^{k-1} e^{-t} - 1 \right|.$$

We have by the induction hypothesis that $\lim_{M\to\infty} \varepsilon_M = 0$.

Fix $M > k$. For $t \geq M$, we have that

$$(1 - \varepsilon_M) \int_M^t \frac{e^s}{s^k} ds \leq \alpha_k(t) - \int_0^M \frac{\alpha_{k-1}(s)}{s} ds \leq (1 + \varepsilon_M) \int_M^t \frac{e^s}{s^k} ds \tag{2.3.2}$$

Now, integrating by parts we obtain that

$$\int_M^t \frac{e^s}{s^k} ds = \frac{e^t}{t^k} - \frac{e^M}{M^k} + k \int_M^t \frac{e^s}{s^{k+1}} ds\,.$$

Since

$$\int_M^t \frac{e^s}{s^{k+1}} ds \leq \frac{1}{M} \int_M^t \frac{e^s}{s^k} ds\,,$$

we further obtain the bound

$$\frac{e^t}{t^k} - \frac{e^M}{M^k} \leq \int_M^t \frac{e^s}{s^k} ds \leq \frac{M}{M-k} \left( \frac{e^t}{t^k} - \frac{e^M}{M^k} \right). \tag{2.3.3}$$

From (2.3.2) and (2.3.3),

$$1 - \varepsilon_M \leq \liminf_{t\to\infty} \alpha_k(t) t^k e^{-t} \leq \limsup_{t\to\infty} \alpha_k(t) t^k e^{-t} \leq (1 + \varepsilon_M) \frac{M}{M-k}\,.$$

Letting $M \to \infty$, we finally conclude that

$$\alpha_k(t) \sim \frac{e^t}{t^k}\,, \quad \text{as } t \to \infty\,,$$

as desired. □

### 2.3.4 Relative growth of moments

We are comparing now, *in the case* $M_f = +\infty$, the growth of moments of different exponents and also the growth of the factorial moments and the moments of the $X_t$; this is covered, respectively, by Corollaries 2.3.10 and 2.3.11. The basic tool is the following lemma.

**Lemma 2.3.9.** *For* $1 < \alpha < \beta$,

$$\frac{\mathbf{E}(X_t^\alpha)}{\mathbf{E}(X_t)^\alpha} \le \frac{\mathbf{E}(X_t^\beta)}{\mathbf{E}(X_t)^\beta}, \quad \textit{for any } t \in (0, R). \tag{2.3.4}$$

*Proof.* Recall that the $X_t$ are non-negative variables. Write $\alpha = u + (1-u)\beta$, with $u = (\beta - \alpha)/(\beta - 1) \in (0,1)$. Hölder's inequality gives that

$$\mathbf{E}(X_t^\alpha) = \mathbf{E}(X_t^u\, X_t^{(1-u)\beta}) \le \mathbf{E}(X_t)^u \mathbf{E}(X_t^\beta)^{1-u},$$

and so that

$$\frac{\mathbf{E}(X_t^\alpha)}{\mathbf{E}(X_t^\beta)} \le \Big(\frac{\mathbf{E}(X_t)}{\mathbf{E}(X_t^\beta)}\Big)^u.$$

Jensen's inequality gives that $\mathbf{E}(X_t^\beta) \ge \mathbf{E}(X_t)^\beta$, and therefore,

$$\frac{\mathbf{E}(X_t^\alpha)}{\mathbf{E}(X_t^\beta)} \le \mathbf{E}(X_t)^{(1-\beta)u} = \frac{\mathbf{E}(X_t)^\alpha}{\mathbf{E}(X_t)^\beta},$$

as stated. □

From Lemma 2.3.9, we deduce the following two corollaries.

**Corollary 2.3.10.** *Assume* $M_f = +\infty$. *If* $1 < \alpha < \beta$, *then*

$$\lim_{t \uparrow R} \frac{\mathbf{E}(X_t^\alpha)}{\mathbf{E}(X_t^\beta)} = 0\,.$$

*Proof.* From Lemma 2.3.9, we have that

$$\frac{\mathbf{E}(X_t^\alpha)}{\mathbf{E}(X_t^\beta)} \le \frac{\mathbf{E}(X_t)^\alpha}{\mathbf{E}(X_t)^\beta} = m_f(t)^{\alpha-\beta}, \quad \text{for any } t \in (0, R).$$

The statement follows since $\lim_{t\uparrow R} m_f(t) = +\infty$ and $\alpha - \beta < 0$. □

**Corollary 2.3.11.** *Assume* $M_f = +\infty$. *For any integer* $k \ge 1$, *we have that*

$$\lim_{t \uparrow R} \frac{\mathbf{E}(X_t^{\underline{k}})}{\mathbf{E}(X_t^k)} = 1\,.$$

*Proof.* From Corollary 2.3.10, we deduce that

$$\lim_{t \uparrow R} \frac{\mathbf{E}(X_t^j)}{\mathbf{E}(X_t^k)} = 0\,, \quad \text{for } 0 \le j < k.$$

The statement follows by expanding $\mathbf{E}(X_t^{\underline{k}})$ as $\mathbf{E}(X_t^k)$ plus a linear combination of the moments $\mathbf{E}(X_t^j)$ with $0 \le j < k$. □

## 2.4 Khinchin Clans

In the previous section, we have compared the growth of $\mathbf{E}(X_t^\beta)$ with $\mathbf{E}(X_t)^\beta$, as $t \uparrow R$, for the basic examples of Khinchin families. Motivated by Hayman, [48], we introduce next a particular kind of Khinchin families, that we call *clans*, for which $\mathbf{E}(X_t^2) \sim \mathbf{E}(X_t)^2$ as $t \uparrow R$. Concretely,

**Definition 1.** *Let $f$ be in $\mathcal{K}$ have radius of convergence $R \leq \infty$. We say that $f$ is a* clan (*and also that the associated family $(X_t)_{t\in[0,R)}$ is a clan*) *if*

$$\lim_{t\uparrow R} \frac{\sigma_f(t)}{m_f(t)} = 0. \tag{2.4.1}$$

For a clan, the normalized variables $Y_t = X_t/\mathbf{E}(X_t)$, for $t \in (0, R)$, converge in probability to the constant 1 as $t \uparrow R$, since its variance $\mathbf{V}(Y_t) = \sigma_f^2(t)/m_f(t)^2$ converges to 0 as $t \uparrow R$.

This clan condition is equivalent to

$$\lim_{t\uparrow R} \frac{\mathbf{E}(X_t^2)}{\mathbf{E}(X_t)^2} = 1.$$

In terms of just the mean $m_f$, being a clan is equivalent to

$$\lim_{t\uparrow R} \frac{t m_f'(t)}{m_f(t)^2} = 0,$$

while in terms of the derivative power series $\mathcal{D}_f$, the condition for being clan becomes

$$\lim_{t\uparrow R} \frac{m_{\mathcal{D}_f}(t)}{m_f(t)} = 1.$$

Alternatively, if we define

$$L_f(t) := \frac{f(t) f''(t)}{f'(t)^2}, \quad \text{for } t \in (0, R), \tag{2.4.2}$$

and since

$$\frac{\mathbf{E}(X_t^2)}{\mathbf{E}(X_t)^2} = \frac{1}{m_f(t)} + L_f(t), \tag{2.4.3}$$

we have that $f$ is a clan if and only if

$$\lim_{t\uparrow R} L_f(t) = 1 - 1/M_f. \tag{2.4.4}$$

For a clan, the normalized variables $Y_t = X_t/\mathbf{E}(X_t)$, for $t \in (0, R)$, converge in probability to the constant 1, as $t \uparrow R$, since its variance $\mathbf{V}(Y_t) = \sigma_f^2(t)/m_f(t)^2$ converges to 0, as $t \uparrow R$.

**Remark 2.4.1.** For any power series $f(z) = \sum_{n=0}^\infty a_n z^n$ in $\mathcal{K}$, if $N \geq 1$ is the first index $n \geq 1$ so that $a_n \neq 0$, then

$$\lim_{t\downarrow 0} t^N \frac{\sigma_f^2(t)}{m_f^2(t)} = \frac{a_0}{a_N}, \quad \text{and thus} \quad \lim_{t\downarrow 0} \frac{\sigma_f^2(t)}{m_f^2(t)} = +\infty.$$

see Section 2.2 for further details. ⊠

### 2.4.1 First examples of clans

The Poisson family associated to $f(z) = e^z$ is a clan, since in this case $m_f(t) = t$ and $\sigma_f^2(t) = t$ and the radius of convergence is $R = \infty$.

The Bernoulli and binomial families are also clans. In fact, all polynomials $f$ are clans, since for them $\sigma_f(t) \to 0$ and $m_f(t) \to \text{degree}(f)$, as $t \to \infty$. Further, we have:

**Lemma 2.4.2.** *Let $f$ be a clan. Then $M_f < +\infty$ if and only if $f$ is a polynomial.*

*Proof.* Polynomials in $\mathcal{K}$ are clans and have $M_f$ finite, in fact $M_f$ coincides with its degree. Conversely, let $f(z) = \sum_{n=0}^{\infty} a_n z^n$ be in $\mathcal{K}$, not a polynomial and with $M_f < \infty$. Thus, because of Lemma 1.2.1, we have that $R < \infty$ and $\sum_{n=1}^{\infty} n a_n R^n < +\infty$.

Then since

$$\frac{\mathbf{E}(X_t^2)}{\mathbf{E}(X_t)^2} = \frac{(\sum_{n=0}^{\infty} n^2 a_n t^n)(\sum_{n=0}^{\infty} a_n t^n)}{(\sum_{n=0}^{\infty} n a_n t^n)^2}\,, \quad \text{for } t \in (0, R)\,,$$

taking limit as $t \uparrow R$, we obtain

$$(\star) \quad \frac{(\sum_{n=0}^{\infty} n^2 a_n R^n)(\sum_{n=0}^{\infty} a_n R^n)}{(\sum_{n=0}^{\infty} n a_n R^n)^2} = 1\,.$$

If we set $b_n = a_n R^n / (\sum_{j=0}^{\infty} a_j R^j)$, para cada $n \geq 0$, which satisfy $\sum_{n=0}^{\infty} b_n = 1$, then $(\star)$ becomes

$$\sum_{n=0}^{\infty} n^2 b_n = \Big(\sum_{n=0}^{\infty} n b_n\Big)^2.$$

This means that $b_n = 0$, for each $n \geq 0$ except for one value of $n$. This would imply that $f$ is a monomial, a contradiction. □

In Section 3.2.3 it will be shown that every $f \in \mathcal{K}$ of the ample class of strongly Gaussian functions conforms a clan. Primordial examples of strongly Gaussian functions are the generating function of the partitions and its variants (see, for example, [17] and [18]).

### 2.4.2 Clan in terms of $f$, directly

Observe that we may write

$$\text{(2.4.5)} \qquad \frac{\mathbf{E}(X_t^2)}{\mathbf{E}(X_t)^2} = \frac{1}{m_f(t)} + \frac{f(t) f''(t)}{f'(t)^2}\,.$$

Thus, we have, in terms of $f$ directly, that $f$ conforms a clan if and only

$$\text{(2.4.6)} \qquad \lim_{t \uparrow R} \Big(\frac{1}{m_f(t)} + \frac{f(t) f''(t)}{f'(t)^2}\Big) = 1\,.$$

As a consequence of (2.4.6) we have that

**Lemma 2.4.3.** *Let $f$ be a power series in $\mathcal{K}$ with radius of convergence $R>0$ and such that $M_f=+\infty$. Then $f$ is a clan if and only if*

$$\lim_{t\uparrow R} L_f(t)=\lim_{t\uparrow R}\frac{f(t)f''(t)}{f'(t)^2}=1\,. \tag{2.4.7}$$

It follows, from instance, that the exponential generating function of the Bell numbers (number of partitions of sets) $e^{e^z-1}$ is a clan.

In particular, and because of Lemma 2.4.2, if $f$ is an *entire transcendental function in* $\mathcal{K}$, then $f$ is a clan if and only if

$$\lim_{t\to\infty}\frac{f(t)f''(t)}{f'(t)^2}=1\,. \tag{2.4.8}$$

For a polynomial $f$ of degree $N$, we have that

$$\lim_{t\to\infty}\frac{f(t)f''(t)}{f'(t)^2}=1-\frac{1}{N}\,.$$

For any transcendental entire function it is always the case that

$$\liminf_{t\to\infty} L_f(t)=\liminf_{t\to\infty}\frac{f(t)f''(t)}{f'(t)^2}\geq 1\,. \tag{2.4.9}$$

This follows from Lemma 2.4.2 and the general identity (2.4.5). This fact has been pointed out in [91, p. 682]. Even further, assume that $M_f=+\infty$ and let $c$ be such that $m_f(c)=1$. From the identity

$$L_f(s)=\Big(s\Big(1-\frac{1}{m_f(s)}\Big)\Big)'\,,\quad \text{for } s\in(0,R)\,,$$

we deduce that

$$\int_c^t L_f(s)ds=t\Big(1-\frac{1}{m_f(t)}\Big)\,,\quad \text{for } t\in[c,R)$$

and so, that

$$\lim_{t\uparrow R}\frac{1}{t}\int_c^t L_f(s)ds=1\,.$$

This implies, in particular and given (2.4.9), that

$$\liminf_{t\to\infty} L_f(t)=1\,,\quad \text{for any transcendental entire function } f \text{ in } \mathcal{K}\,. \tag{2.4.10}$$

### 2.4.3 Growth of the quotient $\sigma_f/m_f$ and gaps

Let $f\in\mathcal{K}$, not a polynomial, have radius of convergence $R>0$. As in Section 1.2.5, we denote by $(n_k)_{k=1}^{\infty}$ the increasing sequence of indices of the nonzero coefficients of $f$. We define $\overline{G}(f)$ by

$$\overline{G}(f)=\limsup_{k\to\infty}\frac{n_{k+1}}{n_k}\cdot \tag{2.4.11}$$

Clearly, $\overline{G}(f) \geq 1$. If for a given $f$ we had $\liminf_{k\to\infty} n_{k+1}/n_k > 1$, then $f$ would be the sum of a polynomial and a power series with Hadamard gaps; but notice that $\overline{G}(f)$ calls for a 'lim sup'. Observe that the definition of $\overline{G}(f)$ involves the quotients $n_{k+1}/n_k$, and not the differences $n_{k+1}-n_k$ as it is the case in $\operatorname{gap}(f)$ and $\overline{\operatorname{gap}}(f)$.

**Theorem 2.4.4.** *Assume that $f \in \mathcal{K}$, with radius of convergence $R > 0$, is not a polynomial and that $M_f = +\infty$. Then*

$$\limsup_{t\uparrow R} \frac{\sigma_f(t)}{m_f(t)} \geq \frac{\overline{G}(f)-1}{\overline{G}(f)+1}\cdot$$

The proof below mimics our proof of Theorem 1.2.3.

*Proof.* Since $M_f = +\infty$, we have that $m_f(t)$ is a homeomorphism from $[0, R)$ onto $[0, +\infty)$, and thus for any integer $k$ there is $t_k \in (0, R)$ so that

$$m_f(t_k) = \frac{n_k + n_{k+1}}{2}\cdot$$

For the random variable $X_{t_k}$, we have that $|X_{t_k} - m_f(t_k)| \geq (n_{k+1} - n_k)/2$ with probability 1, and thus

$$\sigma_f^2(t_k) = \mathbf{E}\big((X_{t_k} - m_f(t_k))^2\big) \geq \frac{1}{4}\,(n_{k+1} - n_k)^2,$$

and also

$$\frac{\sigma_f^2(t_k)}{m_f^2(t_k)} \geq \frac{(n_{k+1} - n_k)^2}{(n_{k+1} + n_k)^2}\cdot$$

The result follows. □

Similarly, for the general absolute central moments of the family of functions $f \in \mathcal{K}$ as in the statement of Theorem 2.4.4, we have that

$$\limsup_{t\uparrow R} \frac{\mathbf{E}\left(|X_t - m_f(t)|^p\right)}{m_f(t)^p} \geq \Big(\frac{\overline{G}(f)-1}{\overline{G}(f)+1}\Big)^p, \quad \text{for any } p > 0.$$

### 2.4.4 Some non examples of clans

The geometric and negative binomial families *are not clans*, since for $f(z) = 1/(1-z)^N$ we have that $\lim_{t\uparrow 1} \sigma_f(t)/m_f(t) = 1/\sqrt{N}$. If fact, for each $\alpha > 0$ the function $f(z) = 1/(1-z)^\alpha$, which is in $\mathcal{K}$ is not a clan, because $\lim_{t\uparrow 1} \sigma_f(t)/m_f(t) = 1/\sqrt{\alpha}$.

Theorem 2.4.4 provides us with many examples of power series in $K$ which are not clans:

**Corollary 2.4.5.** *Let $f$ be in $\mathcal{K}$. If $M_f = +\infty$ and $\overline{G}(f) > 1$, then $f$ is* not *a clan.*

Power series, with radius of convergence $R = 1$, like $1 + \sum_{k=1}^{\infty} z^{2^k}$, or, entire, like $1 + \sum_{k=1}^{\infty} z^{2^k}/(2^k!)$ are in $\mathcal{K}$ and are not clans.

### 2.4.5 Some basic properties of clans

We register now a few properties of power series in $\mathcal{K}$ that are clans, that is, functions in $\mathcal{K}$, with radius of convergence $R > 0$, and satisfying the limit condition given in (2.4.1).

• It follows from Chebyshev's inequality that if $(X_t)_{t\in[0,R)}$ is a clan, then for any $\varepsilon > 0$,

$$\lim_{t\uparrow R} \mathbf{P}\Big(\Big|\frac{X_t}{\mathbf{E}(X_t)} - 1\Big| > \varepsilon\Big) = 0,$$

and thus that $X_t/\mathbf{E}(X_t)$ converges in probability to the constant 1 as $t \uparrow R$: the random variable $X_t$ concentrates about its mean $m_f(t)$ as $t \uparrow R$.

• By Lemma 2.3.6, we have that if $f$ is a clan, then

$$\lim_{t\uparrow R} \frac{\mathbf{E}(X_t^p)}{\mathbf{E}(X_t)^p} = 1 \quad \text{for any } p \in (0, 2].$$

Theorem 2.4.8 below will show that if $f$ is a clan, then this limit result actually holds for any $p > 0$.

• If $f$ is a clan (with at least three nonzero coefficients), then $\mathcal{D}_f = zf'(z)$ is also a clan. This will be proved right after Theorem 2.4.8. (The condition of three nonzero coefficients excludes the case in which the variables associated to $\mathcal{D}_f$ are constant.)

• Notice that if $g \in \mathcal{K}$ with radius of convergence $R \leq \infty$ conforms a clan, then for any integer $N \geq 1$, $f(z) = g(z^N)$ also conforms a clan, since for $t \in [0, R^{1/N})$, $m_f(t) = N m_g(t^N)$ and $\sigma_f(t) = N\sigma_g(t^N)$, where $m_g$, $m_f$, $\sigma_g$, and $\sigma_f$ denote the means and variances of the Khinchin families of $g$ and $f$.

If $f, g \in \mathcal{K}$ are clans, then their product $h \equiv fg$ is also a clan. For we have $m_h = m_f + m_g$ and $\sigma_h^2 = \sigma_f^2 + \sigma_g^2 = o(m_f^2 + m_g^2)$ and thud $\sigma_h = o(m_h)$. In particular, if $f \in \mathcal{K}$ and $N \geq 1$ is an integer, then $f^N$ is also a clan. Later, see Section 2.4.8, we will verify that if $f$ is a clan then $\mathcal{D}f$ is a clan.

• If $f$ is an entire function in $\mathcal{K}$ and $f$ is a clan, then $e^f$ is a clan and $f^k$ is a clan for every integer $k \geq 1$. To show this, we distinguish between $f$ polynomial or not.

If $f$ is a polynomial in $\mathcal{K}$, then $f^k$ is also polynomial in $\mathcal{K}$ and thus a clan. For $h \equiv e^f$, we have $M_h = +\infty$. Also

$$\frac{h''(t)h(t)}{h'(t)^2} = \frac{f''(t)}{f'(t)^2} + 1,$$

and since clearly $\lim_{t\to\infty} f''(t)/f'(t)^2 = 0$ for the polynomial $f$, the condition of Lemma 2.4.3 holds.

If $f$ is not a polynomial then $M_f = +\infty$. Write $h \equiv e^f$. Since

$$\frac{h''(t)h(t)}{h'(t)^2} = \frac{f''(t)}{f'(t)^2} + 1 = \Big(\frac{f''(t)f(t)}{f'(t)^2}\Big)\frac{1}{f(t)} + 1\,,$$

the result follows from Lemma 2.4.3. For powers $f^k$ the argument to show that they are clans is analogous.

• Finally, if $f$ and $g$ are entire functions which are clans, then the composition $f \circ g$ is a clan. This is clear if both $f$ and $g$ are polynomials. Otherwise, this follows from the identity

$$L_{f\circ g}(t) = L_f(g(t)) + \frac{L_g(t)}{m_f(g(t))}, \quad \text{for any } t > 0,$$

and by combining Lemma 2.4.3 and the fact that, for a polynomial $h$ in $\mathcal{K}$ of degree $N$, we have that $\lim_{t\to\infty} L_h(t) = 1 - 1/N$, while $\lim_{t\to\infty} m_h(t) = N$.

In particular, if $g$ is an entire function which is a clan, then $e^g$ is a clan.

### 2.4.6 Weak clans

We say that $f \in \mathcal{K}$ gives raise to a *weak clan* if

$$\liminf_{t\uparrow R} \frac{\mathbf{E}(X_t^2)}{\mathbf{E}(X_t)^2} = 1\,,$$

or, equivalently,

$$\liminf_{t\uparrow R} \frac{\sigma_f(t)}{m_f(t)} = 0\,.$$

Every entire function $f \in \mathcal{K}$ is a weak clan. In fact, we have

**Proposition 2.4.6.** *Let $f$ be an entire function in $\mathcal{K}$. Then for every $\varepsilon > 0$, there exists a set $G_\varepsilon \subset (0, R)$ of logarithmic measure not exceeding $\varepsilon$, such that*

$$(\flat) \qquad \lim_{\substack{t\to\infty;\\ t\notin G_\varepsilon}} \frac{\sigma_f(t)}{m_f(t)} = 0\,.$$

*Proof.* Polynomials are clans, so we may assume that $M_f = +\infty$.

For $a > 0$, consider $H_a = \{x \geq a : \sigma_f^2(x) \geq m_f(x)^{3/2}\}$. For $x \in H_a$ we have that

$$\frac{m_f'(x)}{m_f(x)^{3/2}} \geq \frac{1}{x}\,,$$

and thus, by observing that $M_f = +\infty$, we deduce that

$$\int_{H_a} \frac{dx}{x} \leq \frac{2}{m_f(a)^{1/2}}\,.$$

Let $a = a(\varepsilon)$ be such that $2/m_f(a)^{1/2} \leq \varepsilon$, then for $G_\varepsilon = H_{a(\varepsilon)}$, we have that $G_\varepsilon$ has logarithmic measure at most $\varepsilon$ and for $t \notin G_\varepsilon$ we have that

$$\frac{\sigma_f^2(t)}{m_f^2(t)} \leq \frac{1}{m_f^{1/2}(t)}\,,$$

which gives $(\flat)$

□

As a consequence of this result, we have that, in general, for entire functions $f \in \mathcal{K}$, which are not clans, the quantity $\sigma_f(t)/m_f(t)$ must oscillate and has not limit as $t \to \infty$.

The same proof gives the same restricted limit (avoiding a set of arbitrarily small logarithmic measure) with $\sigma_f^2(t)/m_f(t)^{1+\delta}$, with $\delta > 0$, instead of $\sigma_f^2(t)/m_f(t)^2$.

In contrast, power series $f \in \mathcal{K}$ with finite radius of convergence *need not be weak clans.* In fact, for $f(z) = 1/(1-z)$, we have that $\lim_{t\uparrow 1} \sigma_f(t)/m_f(t) = 1$, and for $f(z) = 1+\ln(1/(1-z))$, we have $\lim_{t\uparrow 1} \sigma_f(t)/m_f(t) = +\infty$.

In the next two examples, with radius $R = 1$, the quotient $\sigma_f^2(t)/m_f^2(t)$ has a positive (finite) limit as $t \uparrow 1$.

• Let $g(z) = 1 + \sum_{n=1}^{\infty} n^4 z^n$. In this case, we have $M_g = +\infty$, in fact,

$$m_g(t) \sim \frac{5}{1-t}, \quad \text{as } t \uparrow 1\,.$$

Also

$$g^{(j)}(t) \sim \frac{(4+j)!}{(1-t)^{5+j}}, \quad \text{and so} \quad \lim_{t\uparrow 1} \frac{g''(t)g(t)}{g'(t)^2} = 6/5\,,$$

and thus, via (2.4.5), we deduce that

$$\lim_{t\uparrow 1} \frac{\mathbf{E}(X_t^2)}{\mathbf{E}(X_t)^2} = \frac{6}{5} > 1\,.$$

Thus $g$ is not a weak clan.

• Let $h(z) = 1 + \sum_{n=1}^{\infty} z^n/z^4$. In this case, we have $M_h = \zeta(3)/(1+\zeta(4)) < +\infty$. Since $th'(t) = \sum_{n=1}^{\infty} z^n/n^3$, and $th'(t) + t^2 h''(t) = \sum_{n=1}^{\infty} z^n/n^2$, we readily obtain that

$$\lim_{t\uparrow 1} \frac{\mathbf{E}(X_t^2)}{\mathbf{E}(X_t)^2} = \frac{\zeta(2)(1+\zeta(4))}{\zeta(3)^2} > 1\,.$$

**Remark 2.4.7.** In general, if $\alpha \in \mathbb{R}$, the power series in $\mathcal{K}$, of radius $R = 1$, given by

$$f(z) = 1 + \sum_{n=1}^{\infty} n^{\alpha} z^n$$

is not a weak clan, and thus not a clan. In fact, for each $\alpha$, the limit $L := \lim_{t\uparrow 1} \dfrac{\mathbf{E}(X_t^2)}{\mathbf{E}(X_t)^2}$ exists (including the possibility of $\infty$ as limit) and it is not 1.

We have

$$M_f = \begin{cases} +\infty, & \text{for } \alpha \geq -2, \\ \dfrac{\zeta(-\alpha-1)}{1+\zeta(-\alpha)} & \text{for } \alpha < -2\,, \end{cases}$$

and

$$L = \begin{cases} \dfrac{\alpha+2}{\alpha+1}, & \text{for } \alpha > -1, \\ +\infty, & \text{for } -1 \geq \alpha \geq -3, \\ \dfrac{\zeta(-\alpha-2)(1+\zeta(-\alpha))}{\zeta(-\alpha-1)} & \text{for } -3 > \alpha\,. \end{cases}$$

⊠

### 2.4.7 Moments of clans

Our next result shows that for a clan $f$, any moment, not just the second one, is asymptotically equivalent, as $t \uparrow R$, to the corresponding power of $m_f(t)$.

**Theorem 2.4.8.** *If $f \in \mathcal{K}$ with radius of convergence $R \leq \infty$ is a clan with associated family $(X_t)_{t\in[0,R)}$, then*

$$\lim_{t\uparrow R} \frac{\mathbf{E}(X_t^p)}{\mathbf{E}(X_t)^p} = 1, \quad \textit{for every } p > 0. \tag{2.4.12}$$

As a consequence of Theorem 2.4.8, we can now prove that, as anticipated in Section 2.4.5, if $f$ is a clan (with at least three nonzero coefficients), then $\mathcal{D}_f = zf'(z)$ is also a clan. Let $(X_t)_{t\in[0,R)}$ and $(W_t)_{t\in[0,R)}$ denote, respectively, the Khinchin families of $f$ and of $\mathcal{D}_f$.

As observed in (1.3.3), the moments of the families $(W_t)$ and $(X_t)$ are related by

$$\mathbf{E}(W_t^p) = \frac{1}{m_f(t)}\,\mathbf{E}(X_t^{p+1}), \quad \text{for any } p > 0 \text{ and any } t \in (0, R).$$

Thus,

$$\frac{\mathbf{E}(W_t^2)}{\mathbf{E}(W_t)^2} = \frac{\mathbf{E}(X_t^3)}{m_f(t)}\,\frac{m_f(t)^2}{\mathbf{E}(X_t^2)^2} = \frac{\mathbf{E}(X_t^3)}{m_f(t)^3}\left(\frac{m_f(t)^2}{\mathbf{E}(X_t^2)}\right)^2.$$

Since $f$ is a clan, both fractions on the far right tend towards 1 as $t \uparrow R$ and, consequently, $\mathcal{D}_f$ is, as claimed, also a clan.

As another consequence, which we have also anticipated, observe that since the partition function $P(z) = \prod_{k=1}^{\infty} 1/(1-z^k)$ is a clan, Theorem 2.4.8 and (1.2.7) give for the moments of its associated family $(X_t)$ that, for any $p > 0$,

$$\mathbf{E}(X_t^p) \sim \mathbf{E}(X_t)^p \sim \frac{\zeta(2)^p}{(1-t)^{2p}}, \quad \text{as } t \uparrow 1.$$

**Remark 2.4.9.** Hayman, in Theorem III of [48], shows that the successive derivatives of Hayman (admissible) functions satisfy asymptotic formulas which are equivalent to the conclusion of Theorem 2.4.8, $\lim_{t\uparrow R} \mathbf{E}(X_t^k)/\mathbf{E}(X_t)^k = 1$ for $k \geq 1$ integer. Our probabilistic proof below shows that this conclusion is valid under the simple and more general notion of clan. ⊠

*Proof of Theorem* 2.4.8. Due to Lemma 2.3.6, it is enough to prove (2.4.12) for any integer $k \geq 1$.

If $f$ is a polynomial of degree $N$, we have, for any integer $k \geq 1$, that $\lim_{t\to\infty} \mathbf{E}(X_t^k) = N^k$. In particular,

$$\lim_{t\to\infty} \frac{\mathbf{E}(X_t^k)}{\mathbf{E}(X_t)^k} = \frac{N^k}{N^k} = 1.$$

We assume now that $f$ is not a polynomial and, consequently, that $M_f = +\infty$.

We first check that

$$\limsup_{t\uparrow R} \frac{\mathbf{E}(X_t^k)}{\mathbf{E}(X_t)^k} \leq e k!. \tag{2.4.13}$$

Corollary 2.3.11, using that $M_f = +\infty$, gives that the inequality (2.4.13) is equivalent to

$$\limsup_{t\uparrow R} \frac{\mathbf{E}(X_t^k)}{\mathbf{E}(X_t)^k} \leq e k!,$$

which follows from (2.4.21) and Theorem 2.4.17.

Denote $V_t := X_t/\mathbf{E}(X_t)$, for $t \in (0, R)$, so that, for any integer $k \geq 1$ and any $t \in (0, R)$

$$\frac{\mathbf{E}(X_t^k)}{\mathbf{E}(X_t)^k} = \mathbf{E}(V_t^k).$$

We aim to show that for a clan, $\lim_{t\uparrow R} \mathbf{E}(V_t^k) = 1$ for any integer $k \geq 1$. For $k = 1$, we have that $\mathbf{E}(V_t) = 1$, for any $t \in (0, R)$, and the case $k = 2$ is just the definition of clan.

By (2.4.13), for any integer $k \geq 1$, the moments of $V_t$ satisfy that

$$\limsup_{t\uparrow R} \mathbf{E}(V_t^{2k}) \leq e\,(2k)!. \tag{2.4.14}$$

Fix an integer $k \geq 3$. Consider a constant $\omega > 0$ and apply the Jensen, Cauchy–Schwarz and Chebyshev inequalities:

$$\begin{aligned} 1 = \mathbf{E}(V_t)^k \leq \mathbf{E}(V_t^k) &= \mathbf{E}\big(V_t^k\,\mathbf{1}_{\{|V_t-1|>\omega\}}\big) + \mathbf{E}\big(V_t^k\,\mathbf{1}_{\{|V_t-1|\leq\omega\}}\big) \\ &\leq \mathbf{E}(V_t^{2k})^{1/2}\,\mathbf{P}\big(|V_t-1|>\omega\big)^{1/2} + (1+\omega)^k \leq \mathbf{E}(V_t^{2k})^{1/2}\,\frac{\sigma_f(t)}{m_f(t)}\,\frac{1}{\omega} + (1+\omega)^k, \end{aligned}$$

where with $\mathbf{1}_A$ we denote the indicator function of the event $A$. Since $f$ is a clan, $\lim_{t\uparrow R} \sigma_f(t)/m_f(t) = 0$, and this and the bound in (2.4.14) combine to imply that

$$1 \leq \limsup_{t\uparrow R} \mathbf{E}(V_t^k) \leq (1+\omega)^k,$$

for any $\omega > 0$. Therefore, $\limsup_{t\uparrow R} \mathbf{E}(V_t^k) = 1$. Now since $k \geq 1$, we have that $\mathbf{E}(V_t^k) \geq \mathbf{E}(V_t)^k = 1$ for any $t \in (0, R)$, we get, as desired, that

$$\lim_{t\uparrow R} \mathbf{E}(V_t^k) = 1$$

for every integer $k \geq 1$. □

**Remark 2.4.10.** Theorem 2.4.8 does not hold for general families or sequences of random variables, i.e., if $(V_n)_{n\geq 1}$ is a sequence of non-negative random variables, then the condition $\lim_{n\to\infty} \mathbf{E}(V_n^2)/\mathbf{E}(V_n)^2 = 1$ does not imply that $\lim_{n\to\infty} \mathbf{E}(V_n^p)/\mathbf{E}(V_n)^p = 1$ for $p > 2$ (although this would be the case for $p < 2$, because of Lemma 2.3.6).

Indeed, for $n \geq 1$, define $V_n$ taking the value $\sqrt{n}/\ln(n+1)$ with probability $1/(n+1)$ and the value $\frac{1}{n}\big((n+1)-\sqrt{n}/\ln(n+1)\big)$ with probability $n/(n+1)$. For this sequence of random variables, we have that $\mathbf{E}(V_n) = 1$ for any $n \geq 1$ and $\lim_{n\to\infty} \mathbf{E}(V_n^p) = 1$, if $p \leq 2$, but $\lim_{n\to\infty} \mathbf{E}(V_n^p) = +\infty$, if $p > 2$. ⊠

### 2.4.8 Derivative $\mathcal{D}f$ of a clan is a clan

Next we show that *if $f$ is a clan (with at least 3 nonzero coefficients) then $\mathcal{D}f$ is also a clan.* Let $(X_t)_{t\in[0,R)}$ and $(W_t)_{t\in[0,R)}$ denote, respectively, the Khinchin families of $f$ and of $zf'(z) = \mathcal{D}_f$.

As observed in (1.3.3) the moments of the families $(W_t)$ and $(X_t)$ are related by

$$\mathbf{E}(W_t^p) = \frac{1}{m_f(t)}\mathbf{E}(X_t^{p+1})\,, \quad \text{for any } p > 0 \text{ and any } t \in (0,R)\,.$$

In particular,

$$\mathbf{E}(W_t) = \frac{1}{m_f(t)}\mathbf{E}(X_t^2)\,,$$

and thus,

$$\frac{\mathbf{E}(W_t^2)}{\mathbf{E}(W_t)^2} = \frac{\mathbf{E}(X_t^3)}{m_f(t)^3}\frac{m_f(t)^4}{\mathbf{E}(X_t^2)^2}\,.$$

Since $f$ is a clan, both fractions on the right tend towards 1 as $t \uparrow R$ and, consequently, $\mathcal{D}f$ is, as claimed, also a clan.

### 2.4.9 On the mean function of a clan.

Hayman showed in Lemmas 2 and 3 of [48]) that the mean $m_f(t)$ of functions $f \in \mathcal{K}$ in the Hayman class cannot grow too slowly. Recall, from Section 3.2.3, that the Hayman class is a subclass of the class of strongly Gaussian functions, and thus, by Corollary 3.2.10, functions in the Hayman class are clans. Next, building upon Hayman's approach, we present a characterization of clans in terms of $m_f$ alone.

Observe that to first order approximation we have that

$$m_f(t + t/m_f(t)) \approx m_f(t) + m_f'(t)\frac{t}{m_f(t)} = m_f(t) + \frac{\sigma_f^2(t)}{m_f(t)},$$

and thus that

$$\frac{m_f(t + t/m_f(t))}{m_f(t)} \approx 1 + \frac{\sigma_f^2(t)}{m_f(t)^2}.$$

This suggests that clans, for which $\lim_{t\uparrow R} \sigma_f(t)/m_f(t) = 0$, may be characterized in terms of the behavior of the quotient $m_f(t + t/m_f(t))/m_f(t)$ as $t \uparrow R$. This is the content of Theorem 2.4.11.

Also, to second order approximation we have, using (1.2.1), that

$$\begin{aligned} \ln f(t + t/m_f(t)) - \ln f(t) &\approx (\ln f)'(t)\frac{t}{m_f(t)} + \frac{1}{2}\,(\ln f)''(t)\,\frac{t^2}{m_f(t)^2} \\ &= 1 + \frac{1}{2}\frac{\sigma_f^2(t)}{m_f(t)^2} - \frac{1}{2}\frac{1}{m_f(t)}, \end{aligned} \tag{2.4.15}$$

which in turn suggests that clans, at least when $M_f = +\infty$, may be characterized in terms of the behavior of the difference $\ln f(t + t/m_f(t)) - \ln f(t)$, as $t \uparrow R$. This is the content of Theorem 2.4.17.

**Theorem 2.4.11.** *Let $f \in \mathcal{K}$ have radius of convergence $0 < R \leq \infty$.*

*If $R = \infty$, the power series $f$ is a clan if and only if*

$$\lim_{t\uparrow R} \frac{m_f(t + t/m_f(t))}{m_f(t)} = 1. \tag{2.4.16}$$

*If $R < \infty$, the power series $f$ is a clan if and only if (2.4.16) holds, and besides,*

$$\lim_{t\uparrow R} (R - t)\, m_f(t) = \infty. \tag{2.4.17}$$

**Remark 2.4.12.** *The classical Borel lemma, see for instance [85, Chapter 9], claims that for any function $\mu(t)$ continuous and increasing in $[t_0, \infty)$ for some $t_0$ and if $a > 1$, then there exists an exceptional set $E$ of logarithmic measure at most $a/(a-1)$ so that*

$$\mu\big(t + t/\mu(t)\big) \leq a\mu(t)\,, \quad \textit{for any } t \in [r_0, +\infty) \setminus E\,.$$

⊠

*Proof of Theorem 2.4.11.* If $f$ is a polynomial then $\lim_{t\to\infty} m_f(t) = \deg(f)$, and thus (2.4.16) holds. We may assume thus that $f$ is not a polynomial.

For the *direct part* we assume that $f$ is a clan and not a polynomial, and thus that $M_f = +\infty$. Consider $t \in (0, R)$ and denote by $\Delta(t)$ the supremum

$$\Delta(t) = \sup_{s\geq t} \frac{\sigma_f^2(s)}{m_f(s)^2}.$$

Since $f$ conforms a clan, $\lim_{t\uparrow R} \Delta(t) = 0$. Now, take $0 < t \leq r \leq s < R$. Since $\sigma_f^2(r) = r m_f'(r)$, we have that

$$\frac{m_f'(r)}{m_f(r)^2} \leq \frac{\Delta(t)}{r}\,.$$

After integration in the interval $(t,s)$ and using that $\Delta$ is decreasing, we get that

$$\frac{1}{m_f(t)}-\frac{1}{m_f(s)}\le \Delta(t)\ln\frac{s}{t}\le \frac{\Delta(t)}{t}(s-t)\,. \tag{2.4.18}$$

If $R<\infty$, by taking the limit as $s\uparrow R$, and using that $M_f=+\infty$, we get that

$$\frac{1}{m_f(t)}\le \frac{\Delta(t)}{t}(R-t)\,,\quad \text{for } t\in[0,R)$$

and since $\lim_{t\uparrow R}\Delta(t)=0$, we deduce that (2.4.17) holds.

So in either case, whether $R$ is finite or not, there exists $t_0\in(0,R)$ such that $t+t/m_f(t)<R$ for every $t\in(t_0,R)$. In (2.4.18), take $t\in(t_0,R)$ and $s=t+t/m_f(t)<R$ and multiply by $m_f(t)>0$ in both sides to get

$$0\le 1-\frac{m_f(t)}{m_f\big(t+t/m_f(t)\big)}\le \Delta(t)\,.$$

The stated result follows since $\lim_{t\uparrow R}\Delta(t)=0$.

Next, the *converse part* of the statement:

We have that $M_f=+\infty$. For $R=+\infty$, this follows since $f$ is not a polynomial. For $R<+\infty$, this follows from assumption (2.4.17).

Observe that whether $R$ is finite or not, we always have that $t+t/m_f(t)<R$ for $t\in[t_0,R)$, for appropriate $t_0$, and thus that $m_f(t+t/m_f(t))$ is well defined.

Denote $\lambda(t)=1+1/m_f(t)$. Consider the power series $\mathcal{D}f$.

From Lemma 1.3.1 we have that

$$(\natural)\quad m_{\mathcal{D}f}(t)\ln\lambda(t)<\ln\frac{\mathcal{D}f\big(\lambda(t)t\big)}{\mathcal{D}f(t)}\,,\quad \text{for } t\in[t_0,R)\,.$$

Since $\mathcal{D}f(t)=m_f(t)f(t)$, for $t\in[0,R)$, appealing to Lemma 1.2.2 we have that

$$(\natural\natural)\quad \begin{aligned}\ln\frac{\mathcal{D}f\big(\lambda(t)t\big)}{\mathcal{D}f(t)}&=\ln\frac{m_f\big(\lambda(t)t\big)}{m_f(t)}+\ln\frac{f\big(\lambda(t)t\big)}{f(t)}\\ &\le \ln\frac{m_f\big(\lambda(t)t\big)}{m_f(t)}+m_f\big(\lambda(t)t\big)\ln\lambda(t)\,,\quad \text{for } t\in[t_0,R)\,.\end{aligned}$$

Combining the inequalities $(\natural)$ and $(\natural\natural)$ and dividing by $m_f(t)\ln\lambda(t)$ we obtain that

$$1\le\frac{\mathbf{E}(X_t^2)}{\mathbf{E}(X_t)^2}=\frac{m_{\mathcal{D}f}(t)}{m_f(t)}\le\frac{\ln\dfrac{m_f(\lambda(t)t)}{m_f(t)}}{m_f(t)\ln\lambda(t)}+\frac{m_f(\lambda(t)t)}{m_f(t)}\,,\quad \text{for } t\in[t_0,R)\,.$$

Since $\lim_{t\uparrow R}m_f(\lambda(t)t)/m_f(t)=1$, by hypothesis, and since $\lim_{t\uparrow}m_f(t)\ln\lambda(t)=1$ because $\lim_{t\uparrow R}m_f(t)=\infty$, we deduce that $\lim\limits_{t\uparrow R}\mathbf{E}(X_t^2)/\mathbf{E}(X_t)^2=1$, and thus that $f$ is a clan. $\square$

### 2.4.10 On the quotient $f(\lambda t)/f(t)$

Let $f$ be in $\mathcal{K}$ and with radius of convergence $R \leq +\infty$. For $\lambda > 1$, we consider the quotient $f(\lambda t)/f(t)$. Bounds for this quotient appeared already in Lemma 1.2.2.

**Lemma 2.4.13.** *If $f$ is a power series in $\mathcal{K}$ with radius of convergence $R \leq +\infty$, then for $\lambda > 1$, the function $f(\lambda t)/f(t)$ is increasing for $t \in (0, R/\lambda)$.*

*Proof.* Fix $\lambda > 1$. Consider the function $g(t) = \ln f(\lambda t) - \ln f(t)$, for $t \in (0, R/\lambda)$. The derivative of $g$ multiplied by $t$ is $m_f(\lambda t) - m_f(t)$, which is always positive, since $m_f$ is strictly increasing in $(0, R)$. □

**Lemma 2.4.14.** *Let $f$ be an entire power series in $\mathcal{K}$ and $\lambda > 1$. If $f$ is a polynomial of degree $d$, then $f(\lambda t)/f(t)$ is bounded for $t \in (0, \infty)$, in fact, $\lim_{t\to\infty} f(\lambda t)/f(t) = \lambda^d$. If $f$ is transcendental, then $\lim_{t\to\infty} f(\lambda t)/f(t) = +\infty$.*

*Proof.* By Lemma 2.4.13, we have that $\lim_{t\to\infty} f(\lambda t)/f(t) = \sup_{t\to\infty} f(\lambda t)/f(t) := H$, which could be $+\infty$ If $H < +\infty$, then $f(\lambda^n) \leq H^n f(1)$, para cada $n \geq 1$, which, by Cauchy estimates, implies that $f$ is a polynomial of degree at most $\ln H/\ln \lambda$. Conversely, if $f \in \mathcal{K}$ is a polynomial of degree $d$ then $\lim_{t\to\infty} f(\lambda t)/f(t) = \lambda^d$. □

See [81, Item 24] and also Lemma 1.2.2.

**Lemma 2.4.15.** *Let $f$ be an entire power series in $\mathcal{K}$. Assume that for a continuous function $\lambda(t)$ defined in $(0, +\infty)$ with values in $(1, +\infty)$ we have that $f(\lambda(t)t)/f(t)$ is bounded for $t \in (0, +\infty)$. Then*

$$\lambda(t) = 1 + O(1/m_f(t)), \quad \textit{as } t \to \infty.$$

*Proof.* Assume first that $f$ is not a polynomial. Let $H > 0$ be such that $f(\lambda(t)t)/f(t) \leq H$, for every $t \in (0, \infty)$. Lemma 1.2.2 gives us that

$$(\star) \quad m_f(t) \ln \lambda(t) \leq H, \quad \text{for } t \in (0, \infty).$$

Write $\lambda(t) \triangleq 1 + \delta(t)$. It is the case that $\lim_{t\to\infty} \delta(t) = 0$; for if for a sequence $t_n$ converging to $\infty$ we have that $\delta(t_n) \geq \delta_0 > 0$, then, using Lemma 2.4.13

$$\frac{f(\lambda(t_n)t_n)}{f(t_n)} \geq \frac{f((1+\delta_0)t_n)}{f(t_n)}$$

but the left hand side is bounded by $H$ while the right hand side tends to $\infty$, because of Lemma 2.4.14.

Now, from $(\star)$ we have that

$$\delta(t) m_f(t) \leq H\delta(t)/\ln \lambda(t),$$

and thus

$$\limsup_{t\to\infty} \delta(t) m_f(t) \leq H,$$

since $\lim_{t\to\infty} \delta(t)/\ln\lambda(t) = 1$.

For a polynomial we have that if $f(\lambda(t)t)/f(t)$ is bounded for $t \in (0, +\infty)$ then $\lambda(t)$ is bounded, which is consistent with the conclusion since $\lim_{t\to\infty} m_f(t) = \text{degree}(f)$. □

The discussion above leads us to consider most naturally the case $\lambda(t) = 1 + 1/m_f(t)$.

**Lemma 2.4.16.** *Let $f \in \mathcal{K}$ with radius of convergence $R \leq \infty$. Assume that*

$$\text{there exists } T \in (0, R) \text{ such that } t + t/m_f(t) < R, \text{ for any } t \in [T, R). \tag{2.4.19}$$

*Then*

$$\frac{f(t + t/m_f(t))}{f(t)} = \sum_{k=0}^{\infty} \frac{1}{k!} \frac{\mathbf{E}(X_t^{\underline{k}})}{\mathbf{E}(X_t)^k}, \quad \text{for any } t \in (T, R) \tag{2.4.20}$$

*and, in particular,*

$$\frac{f(t + t/m_f(t))}{f(t)} \geq \frac{1}{k!} \frac{\mathbf{E}(X_t^{\underline{k}})}{\mathbf{E}(X_t)^k}, \quad \text{for any } k \geq 0. \tag{2.4.21}$$

The condition (2.4.19) in Lemma 2.4.16 is satisfied whenever $f$ is a clan. To see this, observe first that if $R = +\infty$, then (2.4.19) is obvious. Now, if $R < +\infty$, and $f$ is a clan, then (2.4.17) of Theorem 2.4.11 gives, in particular, that $(R - t)m_f(t) > R > t$, for any $t \in [T, R)$, for some $T \in (0, R)$, and thus that $t + t/m_f(t) < R$, for $t \in [T, R)$, which is (2.4.19).

*Proof of Lemma 2.4.16.* Fix $t \in [T, R)$. The radius of convergence of the Taylor expansion of $f$ around $t$ exceeds $t/m_f(t)$, and this gives that

$$f(t + t/m_f(t)) = \sum_{k=0}^{\infty} \frac{1}{k!} \frac{t^k f^{(k)}(t)}{m_f(t)^k}. \tag{2.4.22}$$

Dividing by $f(t)$ and appealing to the expression (2.1.2) of the factorial moments of the $(X_t)$ in terms of $f$, we obtain (2.4.20). The inequalities in (2.4.21) follow since all the summands in (2.4.20) are non-negative. □

In the next result, as suggested by the second order approximation (2.4.15), clans are characterized as in Theorem 2.4.11, but involving $\ln f(t)$ instead of $m_f(t)$.

**Theorem 2.4.17.** *Let $f \in \mathcal{K}$ with radius of convergence $R \leq \infty$. Assume that $M_f = +\infty$. Then $f$ is a clan if and only if condition (2.4.19) holds and*

$$\lim_{t\uparrow R} \ln\Big(\frac{f(t + t/m_f(t))}{f(t)}\Big) = 1. \tag{2.4.23}$$

Concerning the hypothesis $M_f = +\infty$ of Theorem 2.4.17, observe that for any polynomial $f \in \mathcal{K}$ of degree $N$, we have that

$$\lim_{t\to\infty} \frac{f(t + t/m_f(t))}{f(t)} = \Big(1 + \frac{1}{N}\Big)^N < e.$$

*Proof.* Assume first that $f$ is a clan.

We have observed after the statement of Lemma 2.4.16 that condition (2.4.19) is satisfied by clans. To verify (2.4.23), denote $\lambda(t) := 1 + 1/m_f(t)$, for $t \in (0, R)$. Observe that since $\lim_{t \uparrow R} m_f(t) = M_f = +\infty$, we have that

$$\lim_{t \uparrow R} m_f(t) \ln \lambda(t) = 1. \tag{2.4.24}$$

From Lemma 1.2.2, we have, for $T \in (0, R)$ as in condition (2.4.19), that

$$m_f(t) \ln \lambda(t) \le \ln \frac{f(\lambda(t)t)}{f(t)} \le \frac{m_f(\lambda(t)t)}{m_f(t)} \, (m_f(t) \ln \lambda(t)), \quad \text{for any } t \in [T, R).$$

Using (2.4.24) and (2.4.16) of Theorem 2.4.11, the limit (2.4.23) follows.

For the converse implication, assuming now that (2.4.19) (2.4.23) hold, we will verify that $f$ is a clan. It is enough to show that

$$\limsup_{t \uparrow R} \frac{\mathbf{E}(X_t^2)}{\mathbf{E}(X_t)^2} \le 1,$$

or, because of Corollary 2.3.11 and the hypothesis $M_f = +\infty$, that

$$\limsup_{t \uparrow R} \frac{\mathbf{E}(X_t^2)}{\mathbf{E}(X_t)^2} \le 1. \tag{2.4.25}$$

Using hypothesis (2.4.19) and Lemma 2.4.16 we obtain that

$$\frac{f(t + t/m_f(t))}{f(t)} = \sum_{k=0}^{\infty} \frac{1}{k!} \frac{\mathbf{E}(X_t^k)}{\mathbf{E}(X_t)^k}, \quad \text{for } t \in [T, R).$$

Fix an integer $N \ge 3$. Since the summands above are all non-negative we may bound

$$\frac{f(t + t/m_f(t))}{f(t)} \ge \sum_{k=0}^{N} \frac{1}{k!} \frac{\mathbf{E}(X_t^k)}{\mathbf{E}(X_t)^k}, \quad \text{for } t \in [T, R).$$

We split the sum on the right separating the summands corresponding to $k \le 2$ and those with $3 \le k \le N$:

$$\frac{f(t + t/m_f(t))}{f(t)} \ge 1 + 1 + \frac{1}{2} \frac{\mathbf{E}(X_t^2)}{\mathbf{E}(X_t)^2} + \sum_{k=3}^{N} \frac{1}{k!} \frac{\mathbf{E}(X_t^k)}{\mathbf{E}(X_t)^k}.$$

For $3 \le k \le N$, we have that

$$\liminf_{t \uparrow R} \frac{\mathbf{E}(X_t^k)}{\mathbf{E}(X_t)^k} = \liminf_{t \uparrow R} \frac{\mathbf{E}(X_t^k)}{\mathbf{E}(X_t^k)} \, \frac{\mathbf{E}(X_t^k)}{\mathbf{E}(X_t)^k} \ge 1, \tag{2.4.26}$$

since Jensen's inequality gives that $\mathbf{E}(X_t^k) \ge \mathbf{E}(X_t)^k$ and $\lim_{t \uparrow R} \mathbf{E}(X_t^k)/\mathbf{E}(X_t^k) = 1$, because of Corollary 2.3.11 and the hypothesis $M_f = +\infty$.

Fix now $\tau \in (0,1)$. Because of (2.4.26), there exists $S = S(N,\tau) \in [T,R)$ so that $\mathbf{E}(X_t^k)/\mathbf{E}(X_t)^k \geq \tau$, for any $t \in [S,R)$ and any $3 \leq k \leq N$. Thus, we have that

$$\frac{f(t+t/m_f(t))}{f(t)} \geq 1+1+\frac{1}{2}\frac{\mathbf{E}(X_t^2)}{\mathbf{E}(X_t)^2}+\tau\sum_{k=3}^{N}\frac{1}{k!}, \quad \text{for } t \in [S,R).$$

From this inequality and the hypothesis (2.4.23), we deduce that

$$e \geq 1+1+\frac{1}{2}\limsup_{t\uparrow R}\frac{\mathbf{E}(X_t^2)}{\mathbf{E}(X_t)^2}+\tau\sum_{k=3}^{N}\frac{1}{k!}.$$

Letting now $\tau \uparrow 1$, and then $N \to \infty$, we obtain that

$$e \geq 1+1+\frac{1}{2}\limsup_{t\uparrow R}\frac{\mathbf{E}(X_t^2)}{\mathbf{E}(X_t)^2}+\big(e-1-1-\frac{1}{2}\big),$$

and conclude that (2.4.25) holds. □

In terms of the moment generation function of $X_t$, clans are characterized by the following corollary.

**Corollary 2.4.18.** *Let $f \in \mathcal{K}$ with radius of convergence $R \leq \infty$. Assume that $M_f = +\infty$. Then $f$ is a clan if and only if condition (2.4.19) holds and*

$$\lim_{t\uparrow R}\ \mathbf{E}\Big(e^{(X_t-m_f(t))\nu(t)}\Big) = 1, \tag{2.4.27}$$

*where $\nu_f(t) := \ln(1+1/m_f(t))$ for $t>0$.*

*Proof.* The results follows from the expression

$$\mathbf{E}\Big(e^{(X_t-m_f(t))\nu(t)}\Big) = \frac{f(t+t/m_f(t))}{f(t)}\Big(1+\frac{1}{m_f(t)}\Big)^{-m_f(t)}$$

and Theorem 2.4.17. □

Observe that for polynomials in $\mathcal{K}$, (2.4.27) holds.

## 2.5 Entire functions in $\mathcal{K}$

In this section, we relate the order (of growth) of an entire function $f$ in $\mathcal{K}$ with the growth of the mean $m_f(t)$ and the variance $\sigma_f^2(t)$ of the associated Khinchin family.

We shall resort occasionally to the fact that for $f$ in $\mathcal{K}$, we have that

$$\max\{|f(z)| : |z| \leq t\} = \max\{|f(z)| : |z| = t\} = f(t), \quad \text{for every } t \in [0,R). \tag{2.5.1}$$

In particular, we shall be interested in the relation between the order (of growth) of $f$ and the growth of the mean $m_f(t) = \mathbf{E}(X_t)$, of the variance $\sigma_f^2(t) = \mathbf{V}(X_t)$, and of moments $\mathbf{E}(X_t^p)$, with $p > 0$, of the associated family of probability distributions.

Recall that the order $\rho(f)$ of an entire function $f$ in $\mathcal{K}$ is given by

$$\rho(f) := \limsup_{t\to\infty} \frac{\ln\ln \max_{|z|=t}\{|f(z)|\}}{\ln t} = \limsup_{t\to\infty} \frac{\ln\ln f(t)}{\ln t}, \tag{2.5.2}$$

where we have used (2.5.1) in the second expression. On the other hand, for any entire function $f(z) = \sum_{n=0}^\infty a_n z^n$, Hadamard's formula (see Theorem 2.2.2 in [12]) gives $\rho(f)$ in terms of the coefficients of $f$:

$$\rho(f) = \limsup_{n\to\infty} \frac{n\ln n}{\ln(1/|a_n|)}\cdot \tag{2.5.3}$$

## 2.5.1 The order of $f$ entire and the moments $\mathbf{E}(X_t^p)$

We can express the order of an entire function $f \in \mathcal{K}$ in terms of the mean $m_f(t)$ and, in fact, of any moment $\mathbf{E}(X_t^p)$, as follows.

**Theorem 2.5.1.** *Let $f \in \mathcal{K}$ be an entire function of order $\rho(f) \le +\infty$. Then*

$$\limsup_{t\to\infty} \frac{\ln[\mathbf{E}(X_t^p)^{1/p}]}{\ln t} = \rho(f)\,, \quad \textit{for any } p \ge 1. \tag{2.5.4}$$

The case $p = 1$ of Theorem 2.5.1, i.e., $\limsup_{t\to\infty} \ln m_f(t)/\ln t = \rho(f)$, appears in Pólya and Szegő, see item 53 in p. 9 of [81].

*Proof.* We abbreviate and write

$$\Lambda_p := \limsup_{t\to\infty} \frac{\ln[\mathbf{E}(X_t^p)^{1/p}]}{\ln t}\,, \quad \text{for } p > 0.$$

Observe that $\Lambda_p \le \Lambda_q$ if $0 < p \le q$, by Jensen's inequality.

(a) First we show that

$$\Lambda_p \le \rho(f)\,, \quad \text{for any } p > 0.$$

Fix $p > 0$. The inequality trivially holds if $\rho(f) = +\infty$, so we may assume $\rho(f) < +\infty$. Let $\omega > \rho(f)$ and take $\tau = (\omega + \rho(f))/2$. If $f(z) = \sum_{n=0}^\infty a_n z^n$, for $z \in \mathbb{C}$, then Hadamard's formula (2.5.3) gives $N = N_\tau > 0$ such that

$$a_n \le \frac{1}{n^{n/\tau}}, \quad \text{if } n \ge N.$$

For $t$ such that $t^\tau > N$, we have

$$\begin{aligned} f(t)\,\mathbf{E}(X_t^p) &= \sum_{n=1}^\infty n^p a_n\, t^n = \sum_{n\le t^\omega} n^p a_n\, t^n + \sum_{n>t^\omega} n^p a_n\, t^n \\ &\le t^{p\omega} f(t) + \sum_{n\ge 1} n^p \frac{1}{n^{n/\tau}}\, n^{n/\omega} = t^{p\omega} f(t) + C\,, \end{aligned}$$

where $C = C(\omega, \rho(f), p) < +\infty$. Thus,

$$\mathbf{E}(X_t^p) \le t^{p\omega} + C/f(t), \quad \text{if } t^\tau > N,$$

and so $\Lambda_p \le \omega$, for every $\omega > \rho(f)$, which implies, as desired, that $\Lambda_p \le \rho(f)$ for $p > 0$ fixed above.

(b) To finish the proof, it is enough to show that $\rho(f) \le \Lambda_1$, because $\Lambda_1 \le \Lambda_p$ for $p \ge 1$, and this, combined with part (a), would give (2.5.4).

We may assume that $\Lambda_1 < +\infty$, since otherwise there is nothing to prove. We observe first that, for any $\omega > \Lambda_1$, there exists $T = T_\omega$ such that

$$m_f(t) = \mathbf{E}(X_t) \le t^\omega, \quad \text{for } t \ge T,$$

which, recall (1.2.1), can be written in terms of $f$ as

$$\frac{f'(t)}{f(t)} < t^{\omega-1}, \quad \text{for } t \ge T.$$

Upon integration, the above inequality gives that

$$\ln f(t) - \ln f(T) \le \frac{1}{\omega}\,(t^\omega - T^\omega), \quad \text{if } t \ge T,$$

which implies, by the very definition (2.5.2) of order, that $\rho(f) \le \omega$. From this, we deduce, as desired, that $\rho(f) \le \Lambda_1$. □

**Remark 2.5.2.** We mention that, if $\rho(f)$ is finite, then (2.5.4) holds also for $p \in (0,1)$.

To see this, fix $p \in (0,1)$. We just need to show that $\rho(f) \le \Lambda_p$. We are going to interpolate $\Lambda_1$ between $\Lambda_p$ and $\Lambda_2$.

Let $u = 1/(2-p) \in (1/2, 1)$, so that

$$1 = pu + 2(1-u).$$

By Hölder's inequality, we have that

$$\mathbf{E}(X_t) = \mathbf{E}(X_t^{pu}\, X_t^{2(1-u)}) \le \mathbf{E}(X_t^p)^u\, \mathbf{E}(X_t^2)^{1-u}, \quad \text{for any } t > 0,$$

and thus,

$$\Lambda_1 \le pu\,\Lambda_p + 2(1-u)\,\Lambda_2.$$

By (2.5.4), $\Lambda_1 = \Lambda_2 = \rho(f) < +\infty$, so

$$\rho(f) \le pu\,\Lambda_p + 2(1-u)\,\rho(f).$$

This, substituting the value of $u$, gives that $\rho(f) \le \Lambda_p$. ⊠

### 2.5.2 Entire gaps series, order and clans

Recall, from Definition 1, that an entire function $f$ in $\mathcal{K}$ is a clan if

$$\lim_{t\to\infty} \frac{\sigma_f(t)}{m_f(t)} = 0.$$

We shall now exhibit examples of entire functions in $\mathcal{K}$ of *any given order* $\rho$, $0 \le \rho \le \infty$, *which are not clans*. These (counter)examples will be used in forthcoming discussions.

Fix $0 < \rho < +\infty$.

Consider the increasing sequence of integers given by $n_1 = 0$, $n_2 = 1$ and $n_{k+1} = kn_k$, for any $k \ge 2$, and let

$$f(z) = 1 + \sum_{k=2}^{\infty} \frac{1}{n_k^{n_k/\rho}} z^{n_k}.$$

The function $f$ is entire and belongs to $\mathcal{K}$. Hadamard's formula (2.5.3) gives that $\rho(f) = \rho$.

Recall the gap parameter $\overline{G}(f)$ given in (2.4.11), and observe that, in this case, $\overline{G}(f) = \limsup_{k\to\infty} k = +\infty$, so from Theorem 2.4.4, we deduce that

$$\limsup_{t\to\infty} \frac{\sigma_f(t)}{m_f(t)} \ge 1 \tag{2.5.5}$$

holds, and, in particular, that $f$ is not a clan.

With the same specification of the sequence $n_k$, the power series $g$ and $h$ in $\mathcal{K}$ given by

$$g(z) = 1 + \sum_{k=2}^{\infty} \frac{1}{n_k^{n_k^2}} z^{n_k} \quad \text{and} \quad h(z) = 1 + \sum_{k=2}^{\infty} \frac{1}{n_k^{n_k/\sqrt{\ln n_k}}} z^{n_k}$$

are entire, of respective orders $\rho(g) = 0$ and $\rho(h) = +\infty$, and are such that (2.5.5) holds, and so, in particular, they are not clans.

The examples above of entire power series which are not clans are based on the seminal examples of Borel [15], see also [101] of Whittaker, of entire functions whose lower order does not coincide with the order.

### 2.5.3 The order of $f$ entire and the quotient $\sigma_f^2(t)/m_f(t)$

The next result compares the order of the entire function with the quotient $\sigma_f^2(t)/m_f(t)$ as $t \to \infty$.

**Proposition 2.5.3.** *Let $f \in \mathcal{K}$ be an entire function of order $\rho(f) \le \infty$. Then*

$$\liminf_{t\to\infty} \frac{\sigma_f^2(t)}{m_f(t)} \le \liminf_{t\to\infty} \frac{\ln m_f(t)}{\ln t} \le \liminf_{t\to\infty} \frac{\ln \ln f(t)}{\ln t}$$
$$\le \limsup_{t\to\infty} \frac{\ln \ln f(t)}{\ln t} = \limsup_{t\to\infty} \frac{\ln m_f(t)}{\ln t} = \rho(f) \le \limsup_{t\to\infty} \frac{\sigma_f^2(t)}{m_f(t)}.$$

The equality statements in the middle of the comparisons of Proposition 2.5.3 are the very definition of order and the case $p = 1$ of Theorem 2.5.1. On the other hand, Báez-Duarte shows in Proposition 7.7 of [5] that

$$\liminf_{t\to\infty} \frac{\sigma_f^2(t)}{m_f(t)} \le \rho(f) \le \limsup_{t\to\infty} \frac{\sigma_f^2(t)}{m_f(t)}\cdot$$

*Proof of Proposition* 2.5.3. We will check first that

$$\limsup_{t\to\infty} \frac{\ln\ln f(t)}{\ln t} \le \limsup_{t\to\infty} \frac{\ln m_f(t)}{\ln t} \le \limsup_{t\to\infty} \frac{\sigma_f^2(t)}{m_f(t)}\cdot \tag{2.5.6}$$

Of course, we already know that the first two limsup coincide with the order $\rho(f)$.

For the inequality on the right of (2.5.6), let us denote

$$L = \limsup_{t\to\infty} \frac{\sigma_f^2(t)}{m_f(t)}$$

and assume that $L < +\infty$, since otherwise there is nothing to prove. Recall that

$$\frac{\sigma_f^2(t)}{m_f(t)} = \frac{t m_f'(t)}{m_f(t)}\cdot$$

Take $\omega > L$. Then there exists $T > 0$ such that

$$\frac{t m_f'(t)}{m_f(t)} \le \omega, \quad \text{for any } t \ge T,$$

and, by integration, for $t > T$,

$$\ln m_f(t) - \ln m_f(T) \le \omega(\ln t - \ln T),$$

and thus

$$\limsup_{t\to\infty} \frac{\ln m_f(t)}{\ln t} \le \omega, \quad \text{and consequently,} \quad \limsup_{t\to\infty} \frac{\ln m_f(t)}{\ln t} \le L.$$

The proof of the inequality on the left of (2.5.6),

$$\limsup_{t\to\infty} \frac{\ln\ln f(t)}{\ln t} \le \limsup_{t\to\infty} \frac{\ln m_f(t)}{\ln t},$$

follows as above using that

$$t(\ln f(t))' = m_f(t).$$

The proof of the inequalities for the liminf,

$$\liminf_{t\to\infty} \frac{\sigma_f^2(t)}{m_f(t)} \le \liminf_{t\to\infty} \frac{\ln m_f(t)}{\ln t} \le \liminf_{t\to\infty} \frac{\ln\ln f(t)}{\ln t},$$

is analogous. □

Concerning Proposition 2.5.3, a few observations are in order.

(a) The three liminf in the statement, in general, do not give the order of $f$, since the third one is the lower order of $f$, which, for instance and again, by Borel [15], does not have to coincide with the order.

(b) As for the third limsup in the second line: as pointed out by Báez-Duarte (see [5], p. 100), in general, the order $\rho(f)$ is not given by $\limsup_{t\to\infty} \sigma_f^2(t)/m_f(t)$, as proposed by Kosambi in Lemma 4 of [59]. The example of Baéz-Duarte is the canonical product $f$ given by

$$f(z) = \prod_{n=1}^{\infty} \Big(1 + \frac{z^{n^2}}{n^{4n^2}}\Big), \tag{2.5.7}$$

that is an entire function in $\mathcal{K}$. Borel's theorem (see, for instance, Theorem 2.6.5 in [12]) tells us that the order $\rho(f)$ of any canonical product coincides with the exponent of convergence of its zeros. The zeros of $f$ have exponent of convergence is $3/4$, and thus $f$ has order $\rho(f) = 3/4$. From the case $p = 1$ of Theorem 2.5.1, we see that $f$ satisfies that, say, $m_f(t) \leq C\, t^{7/8}$, for some $C > 0$ and every $t \geq 1$. As $f(n^4 e^{i\pi/n^2}) = 0$, we obtain from Lemma 1.3.6 that

$$\sigma_f^2(n^4) \geq \frac{1}{4}\, n^4, \quad \text{for any } n \geq 1,$$

and thus that

$$\frac{\sigma_f^2(n^4)}{m_f(n^4)} \geq \frac{1}{4C}\, n^{1/2}, \quad \text{for any } n \geq 1.$$

So, for this particular function $f$, we have that

$$\limsup_{t\to\infty} \sigma_f^2(t)/m_f(t) = +\infty,$$

but $\liminf_{t\to\infty} \sigma_f^2(t)/m_f(t) \leq \rho(f) = 3/4$.

Alternatively, and more generally, recall that for any $\rho \in [0, +\infty]$, we have exhibited in Section 2.5.2 an entire transcendental function $f \in \mathcal{K}$, with order $\rho(f) = \rho$, and such that $\limsup_{t\to\infty} \sigma_f(t)/m_f(t) \geq 1$, and thus such that $\limsup_{t\to\infty} \sigma_f^2(t)/m_f(t) = +\infty$. Also,

$$\liminf_{t\to\infty} \sigma_f^2(t)/m_f(t) \leq \rho(f) < +\infty,$$

because of Proposition 2.5.3.

(c) On the other hand, as a consequence of Proposition 2.5.3, for an entire function $f \in \mathcal{K}$, we have that *if the limit* $\lim_{t\to\infty} \sigma_f^2(t)/m_f(t)$ *exists* (including the possibility of being $\infty$), then the order $\rho(f)$ of $f$ is precisely

$$\rho(f) = \lim_{t\to\infty} \frac{\sigma_f^2(t)}{m_f(t)} \cdot \tag{2.5.8}$$

In Báez-Duarte's example (2.5.7), the limit of $\sigma_f^2(t)/m_f(t)$ as $t \to \infty$ does not exist, and (2.5.8) does not hold. In general, (2.5.8) does not hold for any entire function $f$ in $\mathcal{K}$ for which the lower order does not coincide with the order.

For the class of entire functions in $\mathcal{K}$ of genus zero which we are to discuss in the next Section 2.5.4, the identity (2.5.8) holds (see the comments after Proposition 2.5.4).

The identity (2.5.8) also holds for the class of nonvanishing entire functions in $\mathcal{K}$ of finite order. To see this, let $f(z) = e^{g(z)} \in \mathcal{K}$, where $g$ is entire (not necessarily in $\mathcal{K}$). We may assume that $g(t) \in \mathbb{R}$ for $t > 0$, and further that, for some $T > 0$, $g(t) > 0$ for $t \geq T$. Assume that $f$ has finite order. Hadamard's factorization theorem gives that $g$ is a polynomial, say of order $N$. Thus $\rho(f) = N$. If the leading coefficient of the polynomial $g$ is $b$, then

$$\begin{aligned} m_f(t) &= tf'(t)/f(t) = tg'(t) \sim b\,N\,t^N, \quad \text{as } t \to \infty, \\ \sigma_f^2(t) &= tm_f'(t) = tg'(t) + t^2 g''(t) \sim bN^2\,t^N \quad \text{as } t \to \infty. \end{aligned}$$

Therefore, (2.5.8) holds.

### 2.5.4 Entire functions of genus 0 with negative zeros

We consider now entire transcendental functions $f$ in $\mathcal{K}$ of genus 0 (and thus of order $\leq 1$) whose zeros are all real and negative. See Chapter 4 of [12]. All such $f$ are canonical products of the form

$$f(z) = a \prod_{j=1}^{\infty} \Big(1 + \frac{z}{b_j}\Big), \quad \text{for } z \in \mathbb{C}, \tag{2.5.9}$$

where $a > 0$ and $(b_j)_{j \geq 1}$ is an increasing sequence of positive real numbers such that

$$\sum_{j=1}^{\infty} \frac{1}{b_j} < \infty. \tag{2.5.10}$$

The zeros $-b_1, -b_2, -b_3, \ldots$ of $f$ all lie on the negative real axis.

We normalize $f(0) = a = 1$ in the discussion that follows. We have considered those canonical products in Section 1.2.3 with that normalization. We will use the notation $N(t)$ of the counting function of the zeros of $f$:

$$N(t) = \#\{j \geq 1 : b_j \leq t\}, \quad \text{for } t > 0\,,$$

which is a non-decreasing function such that $N(t) \to \infty$ as $t \to \infty$.

Recall, also, that Borel's theorem tells us that the order $\rho(f)$ of the canonical product coincides with the exponent of convergence of its zeros, which, in turn (see Theorem 2.5.8 in [12]), is given by

$$\limsup_{t \to \infty} \frac{\ln N(t)}{\ln t} = \rho(f).$$

The function $\ln f(t)$, for $t \in (0, +\infty)$, may be expressed in terms of the counting function $N(t)$; concretely, by integration by parts, one obtains

$$\ln f(t) = \sum_{j=1}^{\infty} \ln\Big(1 + \frac{t}{b_j}\Big) = \int_0^{\infty} \frac{tN(x)}{x(x+t)}\,dx = \int_0^{\infty} \frac{N(ty)}{y(y+1)}\,dy. \tag{2.5.11}$$

This representation (2.5.11) of $\ln f(t)$ in terms of the counting function $N(t)$ comes from [96].

From precise asymptotic information of the counting function $N(t)$ of $f$, one may obtain asymptotic information of the mean and variance function of the family associated to $f$.

**Proposition 2.5.4.** *Let $f$ be an entire function of genus 0 with only negative zeros, and given by (2.5.9). Assume that for $\rho \in (0,1)$ we have*

$$N(t) \sim C\, t^{\rho} \quad \text{as } t \to \infty. \tag{2.5.12}$$

*Then*

$$\text{(a)} \qquad \ln f(t) \sim \frac{C\pi}{\sin(\pi\rho)}\, t^{\rho}, \qquad \text{as } t \to \infty,$$

$$\text{(b)} \qquad m_f(t) \sim \frac{C\pi\rho}{\sin(\pi\rho)}\, t^{\rho}, \qquad \text{as } t \to \infty,$$

$$\text{(c)} \qquad \sigma_f^2(t) \sim \frac{C\pi\rho^2}{\sin(\pi\rho)}\, t^{\rho}, \qquad \text{as } t \to \infty.$$

*Conversely, if* (a), (b) *or* (c) *holds, then* (2.5.12) *holds.*

It follows that for entire functions of genus 0 given by formula (2.5.9), and whenever $N(t) \sim Ct^{\rho}$ as $t \to \infty$, with $\rho \in (0,1)$, then

$$\lim_{t\to\infty} \frac{\sigma_f^2(t)}{m_f(t)} = \rho,$$

as was pointed out (but left unproved) by Báez-Duarte in Proposition 7.9 of [5], and also that equality holds in Proposition 2.5.3. Furthermore, it follows that the entire functions of Proposition 2.5.4 are clans; but see Proposition 2.5.7 for a more general statement.

*Proof.* That (c) $\Rightarrow$ (b) $\Rightarrow$ (a) follows immediately by integration, since $m_f(t) = t(\ln f)'(t)$ and $\sigma_f^2(t) = t m_f'(t)$; recall the formulas (1.2.1).

That (a) implies (2.5.12) is a classical Tauberian theorem, see Valiron [96] and Titchmarsh [93, 94].

The proof of the direct part of Proposition 2.5.4 follows from the representation (2.5.11) and the following standard identity for the Euler Beta function:

$$\int_0^{\infty} \frac{y^{\eta}}{(1+y)^2}\, dy = \text{Beta}(1+\eta, 1-\eta) = \frac{\pi\eta}{\sin(\pi\eta)}, \quad \text{for any } \eta \in [0,1). \tag{2.5.13}$$

The representation (2.5.11) gives that

$$\ln f(t) = t^{\rho} \int_0^{\infty} \frac{y^{\rho}}{y(y+1)}\, \frac{N(ty)}{(ty)^{\rho}}\, dy.$$

Since $N(ty)/(ty)^{\rho} \to C$ as $t \to \infty$ and $y^{\rho}/(y(y+1))$ is integrable in $[0,\infty)$,

$$\lim_{t\to\infty} \frac{\ln f(t)}{t^{\rho}} = C \int_0^{\infty} \frac{y^{\rho}}{y(y+1)}\, dy.$$

Integrating by parts, using that $\rho < 1$ and (2.5.13), we obtain that

$$\int_0^\infty \frac{y^\rho}{y(y+1)}\,dy = \frac{1}{\rho}\int_0^\infty \frac{y^\rho}{(y+1)^2}\,dy = \frac{\pi}{\sin(\pi\rho)}\cdot$$

That is,

$$\ln f(t) \sim \frac{C\pi}{\sin(\pi\rho)}\,t^\rho, \quad \text{as } t \to \infty.$$

For the mean and the variance, we have the representations

$$m_f(t) = t^\rho \int_0^\infty \frac{y^\rho}{(y+1)^2}\,\frac{N(ty)}{(ty)^\rho}\,dy \quad \text{and} \quad \sigma_f^2(t) = t^\rho \int_0^\infty \frac{y^\rho(y-1)}{(y+1)^3}\,\frac{N(ty)}{(ty)^\rho}\,dy.$$

Arguing as above, (b) and (c) follow from (2.5.12). □

**Remark 2.5.5.** For general canonical products $f$ of genus $p \geq 0$ with only negative zeros, there is a representation of $\ln f(t)$ analogous to (2.5.11), which is the case $p = 0$ (see, for instance, Theorem 7.2.1 in [11]):

$$\ln f(t) = (-1)^p \int_0^\infty \frac{(t/x)^{p+1}}{1+t/x}\,N(x)\,\frac{dx}{x}\cdot$$

This expression means in particular that for $p$ odd, the canonical product $f$ is bounded by 1, for $t \in (0,\infty)$, and thus shows that $f$ is not in $\mathcal{K}$. Observe, in any case, that for a primary factor $E_p(z)$, which for $|z| < 1$ is $E_p(z) = \exp(-\sum_{j=p+1}^\infty z^j/j)$, we have that its $(2p+2)$-coefficient is $-p/(2(p+1)^2)$, and therefore $E_p(-z/a)$ with $a > 0$, which vanishes at $-a$, is not in $\mathcal{K}$, if $p \geq 1$; in general, thus, canonical products of nonzero genus with only negative zeros are not in $\mathcal{K}$. ⊠

**Remark 2.5.6.** For entire functions $f$ of genus zero and only negative zeros, if the number of zeros is *comparable* to a power $t^\rho$ with $\rho \in (0,1)$ ($N(t) \asymp t^\rho$, as $t \to \infty$), so are the mean and variance of its Khinchin family. This follows most directly from the representation (2.5.11). ⊠

Entire functions of genus 0 with only negative zeros are always clans.

**Proposition 2.5.7.** *Every entire function $f$ in $\mathcal{K}$ defined by (2.5.9) with $\sum_{j\geq 1} 1/b_j < +\infty$ is a clan.*

*Proof.* Assume that $f$ is not a polynomial. For $f$ given by (2.5.9), we have, using (1.2.17), that

$$\sigma_f^2(t) < m_f(t),$$

and thus $\sigma_f^2(t)/m_f^2(t) \leq 1/m_f(t)$. Since $m_f(t) \to \infty$ as $t \to \infty$, we obtain that $f$ is a clan. □

### 2.5.5 Entire functions, proximate orders and clans

For $\rho \geq 0$, a *proximate $\rho$-order* $\rho(t)$ is a continuously differentiable function defined in $(0,+\infty)$ and such that

$$\lim_{t\to\infty} \rho(t) = \rho \quad \text{and} \quad \lim_{t\to\infty} \rho'(t)t\ln t = 0.$$

Traditionally, proximate orders are allowed to have a discrete set of points where they are not differentiable, but have both one-sided derivatives at those points, see Section 7.4 in [11].

If for a proximate $\rho$-order $\rho(t)$ we write $V(t) = t^{\rho(t)}$, for $t > 0$, then (see, for instance, Lemma 5 in Section 12, Chapter I, of [60]), for every $\lambda > 0$ we have that

$$\lim_{t\to\infty} \frac{V(\lambda t)}{V(t)} = \lambda^{\rho}. \tag{2.5.14}$$

In other terms, the function $V(t)/t^{\rho}$ is a *slowly varying* function.

**Theorem 2.5.8** (Valiron's proximate theorem for $\mathcal{K}$)**.** *If $f$ is an entire function in $\mathcal{K}$ of finite order $\rho \geq 0$, then there is a proximate $\rho$-order $\rho(t)$ such that*

$$\limsup_{t\to\infty} \frac{\ln f(t)}{t^{\rho(t)}} = 1.$$

See Theorem 7.4.2 in [11]; the smoothness which we require in our definition of proximate $\rho$-order $\rho(t)$ is provided by Proposition 7.4.1 and Theorem 1.8.2 (the smooth variation theorem) in [11].

We have the following.

**Theorem 2.5.9.** *Let $f$ be an entire function in $\mathcal{K}$ of finite order $\rho > 0$, let $\rho(t)$ be a proximate $\rho$-order, and let $\tau > 0$. Then*

$$\lim_{t\to\infty} \frac{\ln f(t)}{t^{\rho(t)}} = \tau \tag{2.5.15}$$

*if and only if*

$$\lim_{t\to\infty} \frac{m_f(t)}{t^{\rho(t)}} = \tau\rho. \tag{2.5.16}$$

*If either (2.5.15) or (2.5.16) holds, then $f$ is a clan.*

Observe that in both (2.5.15) and (2.5.16) a 'lim' is assumed, and not a 'lim sup' as in Valiron's Theorem 2.5.8, which encompass all finite order entire functions. Condition (2.5.15) of Theorem 2.5.9 concerns entire functions that are said to have *regular growth*.

Comparing with Valiron's theorem, the limit $\tau$ instead of 1 amounts no extra generality, since replacing $\rho(t)$ by $\rho^{\star}(t) = \rho(t) + \ln\tau/\ln t$, we have that $\rho^{\star}(t)$ is also a proximate $\rho$-order and

$$\lim_{t\to\infty} \frac{\ln f(t)}{t^{\rho^{\star}(t)}} = 1.$$

That (2.5.16) implies that $f$ is a clan is due to Simić, [90]. In [90], see also [91], it is claimed, in the terminology used here, that any entire function of *finite order* in $\mathcal{K}$ is a clan, which is not the case; see, for instance, Section 2.5.2. The error in the argument originates in a misprint in the statement of Theorem 2.3.11 in p. 81 of [11]: the lim sup appearing in that statement should be a lim inf (which is what is actually proved in [11]). See also [61], p. 101, for a similar warning. The argument of [90] shows precisely that (2.5.16) implies that $f$ is a clan.

As for the implication (2.5.15) ⇒ (2.5.16), compare with Lemma 3.1 in [73].

For constant proximate $\rho$-order, ($\rho(t) = \rho$, for $t > 0$) that (2.5.15) implies (2.5.16) and then that $f$ is a clan is due to Pólya and Szegö with an argument involving some delicate estimates: combine items 70 and 71, of page 12, of theirs [81].

*Proof.* Fix $\lambda > 1$. Lemma 1.2.2 gives us that

$$m_f(t) \ln \lambda \le \ln f(\lambda t) - \ln f(t) \le m_f(\lambda t) \ln \lambda, \quad \text{for } t > 0,$$

and thus dividing by $V(t) = t^{\rho(t)}$,

$$\frac{m_f(t)}{V(t)} \le \frac{1}{\ln \lambda}\left[\frac{\ln f(\lambda t)}{V(\lambda t)}\frac{V(\lambda t)}{V(t)} - \frac{\ln f(t)}{V(t)}\right] \le \frac{m_f(\lambda t)}{V(\lambda t)}\frac{V(\lambda t)}{V(t)}, \quad \text{for } t > 0. \tag{2.5.17}$$

We first prove that (2.5.15) ⇒ (2.5.16). From (2.5.14), (2.5.15) and (2.5.17), and letting $t \to \infty$, we deduce that

$$\limsup_{t\to\infty} \frac{m_f(t)}{t^{\rho(t)}} \le \tau \frac{\lambda^\rho - 1}{\ln \lambda} \le \lambda^\rho \liminf_{t\to\infty} \frac{m_f(t)}{t^{\rho(t)}}\cdot$$

Letting $\lambda \downarrow 1$, equation (2.5.16) follows.

We now prove that (2.5.16) ⇒ (2.5.15). Assume first that

$$\limsup_{t\to\infty} \frac{\ln f(t)}{t^{\rho(t)}} < +\infty. \tag{2.5.18}$$

If (2.5.18) holds, then from the first inequality of (2.5.17), and using (2.5.16), we deduce that

$$\tau\rho \ln \lambda + \liminf_{t\to\infty} \frac{\ln f(t)}{t^{\rho(t)}} \le \lambda^\rho \liminf_{t\to\infty} \frac{\ln f(t)}{t^{\rho(t)}},$$

while the second inequality of (2.5.17) gives

$$\lambda^\rho \limsup_{t\to\infty} \frac{\ln f(t)}{t^{\rho(t)}} \le \limsup_{t\to\infty} \frac{\ln f(t)}{t^{\rho(t)}} + \tau\rho\lambda^\rho \ln \lambda.$$

Writing the two inequalities above as

$$\tau\rho \le \frac{\lambda^\rho - 1}{\ln \lambda} \liminf_{t\to\infty} \frac{\ln f(t)}{t^{\rho(t)}} \quad \text{and} \quad \frac{\lambda^\rho - 1}{\ln \lambda} \limsup_{t\to\infty} \frac{\ln f(t)}{t^{\rho(t)}} \le \tau\rho\lambda^\rho,$$

and by letting $\lambda \downarrow 1$, we get that

$$\tau \le \liminf_{t\to\infty} \frac{\ln f(t)}{t^{\rho(t)}} \quad \text{and} \quad \limsup_{t\to\infty} \frac{\ln f(t)}{t^{\rho(t)}} \le \tau,$$

so (2.5.15) follows.

To show that (2.5.18) holds, we first observe that

$$\frac{d}{ds} s^{\rho(s)} = (s \ln s \cdot \rho'(s) + \rho(s))\, s^{\rho(s)-1}, \quad \text{for } s > 0,$$

and from the defining properties of the proximate orders, we deduce that, for appropriately large $A > 0$,

$$\frac{d}{ds} s^{\rho(s)} \geq \frac{\rho}{2}\, s^{\rho(s)-1}, \quad \text{for each } s \geq A. \tag{2.5.19}$$

From (2.5.16), by incrementing $A$ if necessary, we deduce that

$$\frac{d}{ds} \ln f(s) = \frac{m_f(s)}{s} \leq 2\tau\rho\, s^{\rho(s)-1} \leq 4\tau \frac{d}{ds}\, s^{\rho(s)}, \quad \text{for } s \geq A,$$

where (2.5.19) was used in the last inequality. We thus have that

$$\ln f(t) \leq \ln f(A) + 4\tau\, t^{\rho(t)} - 4\tau A^{\rho(A)}, \quad \text{for } t \geq A,$$

from which we obtain that

$$\limsup_{t\to\infty} \frac{\ln f(t)}{t^{\rho(t)}} \leq 4\tau,$$

as wanted.

This argument shows in fact that

$$\frac{1}{\rho} \liminf_{t\to\infty} \frac{m_f(t)}{t^{\rho(t)}} \leq \liminf_{t\to\infty} \frac{\ln f(t)}{t^{\rho(t)}} \leq \limsup_{t\to\infty} \frac{\ln f(t)}{t^{\rho(t)}} \leq \frac{1}{\rho} \limsup_{t\to\infty} \frac{m_f(t)}{t^{\rho(t)}}.$$

For another proof of this last chain of inequalities, see, for instance, Theorem 4 in [95].

Finally, we prove that (2.5.16) $\Rightarrow$ $f$ is a clan. From (2.5.14) and (2.5.16), and taking into account that $\rho > 0$, we deduce that

$$\lim_{t\to\infty} \frac{m_f(\lambda t)}{m_f(t)} = \lambda^\rho, \quad \text{for any } \lambda > 0, \tag{2.5.20}$$

and thus that the mean $m_f$ is a regularly varying function, see Section 1.4 of [11].

To show that $f$ is a clan, we may assume that $f$ is not a polynomial, and thus that $\lim_{t\to\infty} m_f(t) = +\infty$.

Now, (2.5.20) implies that for any function $\lambda(t)$ such that $\lambda(t) > 1$ and such that $\lim_{t\to\infty} \lambda(t) = 1$, we have that

$$\lim_{t\to\infty} \frac{m_f(\lambda(t)\, t)}{m_f(t)} = 1. \tag{2.5.21}$$

To see this, fix $\varepsilon > 0$. Then we have that $1 < \lambda(t) \leq 1 + \varepsilon$, for $t \geq t_\varepsilon$. Therefore,

$$1 \leq \frac{m_f(\lambda(t)\, t)}{m_f(t)} \leq \frac{m_f((1+\varepsilon)\, t)}{m_f(t)}, \quad \text{for } t \geq t_\varepsilon.$$

Thus, $\limsup_{t\to\infty} m_f(\lambda(t)\, t)/m_f(t) \leq (1+\varepsilon)^\rho$, and thus (2.5.21) holds.

Applying (2.5.21) with $\lambda(t) = 1 + 1/m_f(t)$ and appealing to Theorem 2.4.11, we conclude that $f$ is a clan. $\square$

### 2.5.6 Exceptional values and clans

The entire gap series of Section 2.5.2, which are our basic examples of entire functions in $\mathcal{K}$ which are not clans, have no Borel exceptional values. This follows, for instance, from a classical result of Pfluger and Pólya [80]. We show next that, in general, entire functions which are not clans have no Borel exceptional values.

Recall that, by definition, $a$ is a Borel exceptional value of an entire function $f$ of finite order if the exponent of convergence of the $a$-values of $f$ (i.e., the zeros of $f(z) - a = 0$) is strictly smaller than the order of $f$; a theorem of Borel claims that a non-constant entire function of finite order can have at most one Borel exceptional value.

**Theorem 2.5.10.** *If the entire function $f \in \mathcal{K}$ has finite order and has one Borel exceptional value, then $f$ is a clan.*

*Proof.* Let $\rho$ be the order of $f$. Assume that $a \in \mathbb{C}$ is the Borel exceptional value for $f$. Denote with $s$ the exponent of convergence of the zeros of $f(z) - a$. Thus $s < \rho$, since $a$ is a Borel exceptional value for $f$.

Let $f(z) = a + P(z)e^{Q(z)}$ be the Hadamard factorization of $f$, where $P$ is the canonical product formed with the zeros of $f - a$, and $Q$ is a polynomial of degree $d$ and leading coefficient $c \neq 0$. The order of $P$ is $s$, and thus the order $\rho$ of $f$ must be the integer $d$.

Now,

$$|f(t) - a| = |P(t)|e^{\Re Q(t)} \quad \text{and} \quad \ln|f(t) - a| = \ln|P(t)| + \Re Q(t), \quad \text{for } t > 0.$$

Take $s' \in (s, d)$. For a certain $t'$ depending on $s'$, we have for $t \geq t'$ that $\ln|P(z)| \leq |z|^{s'}$, if $|z| = t$. Besides,

$$\frac{\Re Q(t)}{t^d} = \Re c + O\Big(\frac{1}{t}\Big), \quad \text{as } t \to \infty.$$

We conclude that

$$\lim_{t\to\infty} \frac{\ln|f(t) - a|}{t^d} = \Re c,$$

and therefore that

$$\lim_{t\to\infty} \frac{\ln f(t)}{t^d} = \Re c.$$

Observe that if $\Re c = 0$, then $\Re Q(t) = O(t^{d-1})$ as $t \to \infty$, and that would mean that

$$\limsup_{t\to\infty} \frac{\ln f(t)}{t^h} = 0,$$

for some $h$ such that $(s <)\, h < d$, and thus, in particular, that $f$ would be of order at most $h$, which is not the case.

Thus $\Re c > 0$, and condition (2.5.15) of Theorem 2.5.9 holds, so $f$ is a clan. □

For an entire function, not necessarily of finite order, a Picard exceptional value is a value that is taken just a finite number of times. For Picard exceptional values and clans, we have the following result, which came out in a conversation of one of the authors with Walter Bergweiler.

**Proposition 2.5.11.** *If $f = Pe^g$ is in $\mathcal{K}$, where $P$ is a polynomial and $g$ is an entire function in* $\mathcal{K}$ of finite order*, then $f$ is a clan.*

The value 0 is Picard exceptional for $f = Pe^g$. It is not assumed that $P$ is in $\mathcal{K}$, but it is assumed that $g$ is in $\mathcal{K}$. Observe also that the assumption is that $g$ is of finite order; if $e^g$ were of finite order, that $f$ is a clan would follow from Theorem 2.5.10.

*Proof.* The entire function $f$ is transcendental, since $g \in \mathcal{K}$ is not a constant. From Lemma 2.4.3, we have that $f$ is a clan if and only if $\lim_{t\to\infty} L_f(t) = 1$. To show this, we verify first that $g$ satisfies

$$\lim_{t\to\infty} \frac{g''(t)}{g'(t)^2} = 0. \tag{2.5.22}$$

Condition (2.5.22) clearly holds if $g$ is a polynomial.

Assume thus that $g$ is not a polynomial. From the case $p = 1$ of Theorem 2.5.1 applied to the derivative $g'$, which is also of finite order, it follows that for some finite constant $S > 0$ and radius $R_1 > 0$, we have that

$$\frac{g''(t)}{g'(t)} \leq t^S, \quad \text{for } t > R_1.$$

Besides, since $g'$ is not a polynomial, we have, for some radius $R_2$, that $g'(t) > t^{S+1}$, for $t > R_2$. And thus (2.5.22) holds.

Next, a calculation, recall (2.4.2), gives that

$$L_f(t) = \Big(\frac{P''(t)}{P(t)}\frac{1}{g'(t)^2} + 2\frac{P'(t)}{P(t)}\frac{1}{g'(t)} + \frac{g''(t)}{g'(t)^2} + 1\Big) \Big/ \Big(\frac{P'(t)}{P(t)}\frac{1}{g'(t)} + 1\Big)^2, \quad \text{for } t > 0.$$

Since $P$ is a polynomial, we have that $P'(t)/P(t)$ and $P''(t)/P(t)$ tend to 0 as $t \to \infty$. Besides, since $g \in \mathcal{K}$, we have that $\lim_{t\to\infty} g'(t) = +\infty$. Using now (2.5.22), it is deduced that $\lim_{t\to\infty} L_f(t) = 1$. □

Chapter 3

# Gaussian and strongly Gaussian power series

## Contents

---

Let $f \in \mathcal{K}$ be a power series with radius of convergence $R > 0$ and Khinchin family $(X_t)$. In this chapter we study, and give different criteria, for the convergence in distribution of the Khinchin family $(X_t)$, as $t \uparrow R$, to a standard normal. These are the so called Gaussian power series.

We also study uniform conditions that, combined with Hayman's formula (1.3.31), give information about the asymptotic behavior of the coefficients of $f$. These conditions are codified by means of two classes of functions: strongly Gaussian power series and Hayman functions. Along this section we also prove that, under certain conditions, some operations with power series preserve being in each of these classes.

This chapter is partially based on the paper:

- Maciá, V.J. et al. Khinchin families and Hayman class. *Comput. Methods Funct. Theory* **21** (2021), no. 4, 851–904, see [17].

## 3.1 Gaussian power series

For a power series $f \in \mathcal{K}$, with radius of convergence $R > 0$, we say that $f$, or its associated Khinchin family $(X_t)$, is Gaussian if

$$\lim_{t\uparrow R} \mathbf{E}(e^{i\theta \breve{X}_t}) = e^{-\theta^2/2}, \qquad \text{for any } \theta \in \mathbb{R}.$$

Here $\breve{X}_t = (X_t - m_f(t))/\sigma_f(t)$, for any $t \in (0, R)$. Recall that $e^{-\theta^2/2}$ is the characteristic function of the standard normal.

### A. Some basic tools in probability

We collect here, for later use, two basic tools in probability: Lévy's convergence theorem and the method of moments. We restrict both results to the case of the standard normal distribution.

Theorem 3.1.1 gives, in our context, that for $f \in \mathcal{K}$ a Gaussian power series, the family of random variables $\breve{X}_t$ converges in distribution to a standard normal random variable, as $t \uparrow R$.

The next result follows combining Lévy's convergence theorem with the fact that the characteristic function and the cumulative distribution function of a standard normal are continuous functions, see [102, p. 185] for the general statement.

**Theorem 3.1.1** (Lévy's convergence theorem for the standard normal)**.** *For a sequence of random variables $Z_n$ and $Z$ a standard normal random variable we have that:*

$$Z_n \xrightarrow{d} Z, \text{ as } n \to \infty, \quad \text{if and only if} \quad \lim_{n\to\infty} \mathbf{E}(e^{i\theta Z_n}) = e^{-\theta^2/2}, \quad \text{for any } \theta \in \mathbb{R}.$$

Now we introduce the method of moments for the standard normal distribution, see [9, Theorem 30.2] for the general case. In the next result we use that the standard normal distribution is determined by its moments, see, for instance, [9, p. 389].

**Theorem 3.1.2** (Method of moments: convergence to the standard normal)**.** *Let $Z_n$ be a sequence of random variables with finite moments of any order. Denote by $Z$ a standard normal random variable and assume that*

$$\lim_{n\to\infty} \mathbf{E}(Z_n^m) = \mathbf{E}(Z^m), \quad \text{for any integer } m \geq 1,$$

*then $Z_n$ converges in distribution towards $Z$, as $n \to \infty$.*

### B. Some examples

We start by checking the concept of Gaussianity on some of our basic families: the Poisson, Bernoulli and Pascal families.

- For the Poisson family, $f(z) = e^z$, we have a explicit expression for the characteristic functions of each $\breve{X}_t$:

  $$\mathbf{E}(e^{i\theta\breve{X}_t}) = \exp(t(e^{i\theta/\sqrt{t}} - 1 - i\theta/\sqrt{t})), \quad \text{for each } t > 0 \text{ and } \theta \in \mathbb{R}.$$

  see (1.3.11).

  Fixing $\theta \in \mathbb{R}$ and making $t \to +\infty$, we conclude that

  $$\lim_{t\to+\infty} \mathbf{E}(e^{i\theta\breve{X}_t}) = e^{-\theta^2/2}, \quad \text{for any } \theta \in \mathbb{R}.$$

  Therefore, as $t \to +\infty$, the Poisson random variable of parameter $t$, and normalized, that is $\breve{X}_t = (X_t - t)/\sqrt{t}$, converges in distribution towards the standard normal.

- For the Bernoulli family, $f(z) = 1 + z$, we have, see (1.3.5), that

$$\mathbf{E}(e^{\imath\theta\breve{X}_t}) = \frac{te^{\imath\theta/\sqrt{t}} + e^{-\imath\theta\sqrt{t}}}{1+t}, \quad \text{for } \theta \in \mathbb{R} \text{ and } t > 0,$$

therefore, fixing $\theta \in \mathbb{R}$, we have

$$\lim_{t\to+\infty} \mathbf{E}(e^{i\theta\breve{X}_t}) = 1, \quad \text{for any } \theta \in \mathbb{R},$$

and $\breve{X}_t$ converges in distribution towards the constant 0, as $t \to +\infty$.

- The previous example is a particular case of $R \in \mathcal{K}$ a polynomial of degree $N \geq 1$. In this case the family of random variables $(\breve{X}_t)$, also, converges in distribution towards the constant 0, as $t \to +\infty$.

  Indeed: recall that for any polynomial $R \in \mathcal{K}$ we have

$$\lim_{t\to+\infty} m_R(t) = \deg(R) = N,$$

and that the characteristic function is given by

$$\mathbf{E}(e^{i\theta X_t}) = \frac{R(te^{i\theta})}{R(t)}, \quad \text{for any } t \geq 0 \text{ and } \theta \in \mathbb{R},$$

therefore, using again that $R$ is a polynomial with $\deg(R) = N$, we find that

$$\lim_{t\to+\infty} \mathbf{E}(e^{i\theta X_t}) = e^{iN\theta}, \quad \text{for any } \theta \in \mathbb{R}. \tag{3.1.1}$$

The family of random variables $X_t$ converges in distribution towards the constant random variable $N \geq 1$. This follows also from the fact that any polynomial in $\mathcal{K}$ conform a clan and therefore the family of random variables $X_t/m_R(t)$ converges in probability to the constant 1, as $t \to +\infty$, that is, $X_t$ converges in probability, as $t \to +\infty$, to the constant $N$.

Denote $j = N - \max\{0 \leq s < N : a_s \neq 0\}$. Using this notation we have

$$R(z) = a_N z^N + a_{N-j} z^{N-j} + o(z^{N-j}), \quad \text{as } z \to +\infty,$$

and therefore

$$(m_R(t) - N) = \frac{tR'(t) - NR(t)}{R(t)} \sim \frac{-ja_{N-j}}{a_N}\frac{1}{t^j}, \quad \text{as } t \to +\infty.$$

Now observe that

$$\begin{aligned}\sigma_R^2(t) &= m_R(t) - m_R(t)^2 + \frac{t^2R''(t)}{R(t)} \\ &= (2N+1)(m_R(t) - N) - (m_R(t) - N)^2 + \left(\frac{t^2R''(t)}{R(t)} - N(N-1)\right).\end{aligned}$$

Using again that $R(z) = a_N z^N + a_{N-j} z^{N-j} + o(z^{N-j})$, as $z \to +\infty$, we conclude that

$$\sigma_R^2(t) \sim \frac{j^2 a_{N-j}}{a_N} \frac{1}{t^j}, \qquad \text{as } t \to +\infty.$$

We collect the previous discussion by means of the following lemma.

**Lemma 3.1.3.** *Let $R \in \mathcal{K}$ be a polynomial of degree $N \geq 1$, then*

$$(N - m_R(t)) \sim \frac{j a_{N-j}}{a_N} \frac{1}{t^j}, \quad \text{and} \quad \sigma_R^2(t) \sim \frac{j^2 a_{N-j}}{a_N} \frac{1}{t^j}, \quad \text{as } t \to \infty, \tag{3.1.2}$$

*and therefore*

$$\lim_{t\to\infty} \frac{N - m_R(t)}{\sigma_R(t)} = 0, \quad \text{and} \quad \lim_{t\to\infty} \frac{N - m_R(t)}{\sigma_R^2(t)} = \frac{1}{j}. \tag{3.1.3}$$

*Recall that $j = N - \max\{0 \leq s < N : a_s \neq 0\}$.*

The characteristic function of $\breve{X}_t$ is given by

$$\mathbf{E}(e^{i\theta \breve{X}_t}) = \mathbf{E}(e^{i\theta(X_t - N)/\sigma_R(t)}) e^{-i\theta(m_R(t)-N)/\sigma_R(t)}.$$

Using that $R$ is a polynomial of degree $N$:

$$(\star) \quad \lim_{t\to+\infty} \mathbf{E}(e^{i\theta(X_t-N)/\sigma_R(t)}) = \lim_{t\to+\infty} \frac{R(te^{i\theta/\sigma_R(t)})}{R(t)} e^{-i\theta N/\sigma_R(t)} = 1\,,$$

therefore $(\star)$ gives that

$$\lim_{t\to+\infty} \mathbf{E}(e^{i\theta\breve{X}_t}) = 1\,, \qquad \text{for any } \theta \in \mathbb{R}\,,$$

that is, for a polynomial $R \in \mathcal{K}$, Lévy's convergence theorem gives that the family of random variables $(\breve{X}_t)$ converges in distribution towards the constant random variable $X \equiv 0$.

We conclude that the polynomials in $\mathcal{K}$ are not Gaussian.

- For the Pascal family, $f(z) = 1/(1-z)$, we have a explicit expression for the characteristic function of the family $\breve{X}_t$:

$$\mathbf{E}(e^{\imath\theta\breve{X}_t}) = \frac{1-t}{1 - te^{\imath\theta(1-t)/\sqrt{t}}} e^{-\imath\theta\sqrt{t}} = \frac{1-t}{e^{\imath\theta\sqrt{t}} - te^{\imath\theta/\sqrt{t}}}\,, \tag{3.1.4}$$

for any $\theta \in \mathbb{R}$ and $t \in (0,1)$, see (1.3.9), therefore

$$\lim_{t\uparrow 1} \mathbf{E}(e^{i\theta\breve{X}_t}) = \frac{e^{-i\theta}}{1 - i\theta}, \qquad \text{for any } \theta \in \mathbb{R}\,.$$

Observe that

$$m_f(t) = \frac{t}{1-t} \quad \text{and} \quad \sigma_f^2(t) = \frac{t}{(1-t)^2},$$

for any $t \in [0,1)$, and then

$$\lim_{t \uparrow 1} \frac{m_f(t)}{\sigma_f(t)} = 1.$$

In summary: the family of random variables $X_t/\sigma_f(t)$ converges in distribution as $t \uparrow 1$ towards a random variable $Z$, where $Z$ is a exponentially distributed random variable with parameter $\lambda = 1$, and therefore the family of random variables $\breve{X}_t$ converges in distribution, as $t \uparrow 1$, towards $Z+1$, that is, the Pascal family is not Gaussian.

### 3.1.1 Criteria for non-vanishing power series to be Gaussian

In this section we give criteria for a non-vanishing power series to be Gaussian.

#### A. Global criteria for a non-vanishing power series to be Gaussian

Let $f \in \mathcal{K}$ be a power series with radius of convergence $R > 0$. *Assume that $f$ has no zeros in the disk $\mathbb{D}(0,R)$.*

We can take the holomorphic logarithm $\ln(f(z))$ of $f(z)$ which is real at $z=0$ and consider the auxiliary function, the fulcrum of $f$, recall that $F(z) = \ln(f(e^z))$, see Subsection 1.3.4. The fulcrum of $f$ is holomorphic in $\Omega = \{z \in \mathbb{C} : \Re(z) < \ln(R)\}$ (in practice this will be a half-plane or the entire complex plane).

For $s < \ln(R)$ and $\theta \in \mathbb{R}$ we have

$$F(s+i\theta) - F(s) - F'(s)i\theta + F''(s)\frac{\theta^2}{2} = i^3 \int_0^\theta \int_0^\alpha \int_0^\beta F'''(s+i\phi)d\phi d\beta d\alpha,$$

therefore, fixing some $A > 0$, for any $s < \ln(R)$ and $|\theta| \le A$ we have

$$|F(s+i\theta) - F(s) - F'(s)i\theta + F''(s)\frac{\theta^2}{2}| \le \sup_{|\phi|\le A} |F'''(s+i\phi)| \frac{|\theta|^3}{6}.$$

Recall that $F'(s) = m_f(e^s)$ and also that $F''(s) = \sigma_f^2(e^s)$, for any $s < \ln(R)$.

Fixing $t = e^s$, for $s < \ln(R)$, and changing $\theta$ by $\theta/\sigma_f(t)$, we obtain that

$$\left| \ln(f(te^{i\theta/\sigma_f(t)})) - \ln(f(t)) - i\frac{m_f(t)}{\sigma_f(t)}\theta + \frac{\theta^2}{2} \right| \le \frac{\sup_{|\phi|\le A} |F'''(s+i\phi)|}{\sigma_f^3(t)} \frac{|\theta|^3}{6},$$

and therefore

$$\left| \ln \mathbf{E}(e^{i\theta \breve{X}_t}) + \frac{\theta^2}{2} \right| \le \frac{\sup_{|\phi|\le A} |F'''(s+i\phi)|}{\sigma_f^3(t)} \frac{|\theta|^3}{6}. \tag{3.1.5}$$

This bound bears resemblance to the bound in Berry-Esseen's Theorem, which gives an upper bound for the rate of convergence to a standard normal of a sum of independent and identically distributed random variables in terms of a third moment.

As a consequence of inequality (3.1.5) we deduce the following theorem.

**Theorem 3.1.4.** *With the previous notations, if $f \in \mathcal{K}$ has radius of convergence $R > 0$, has no zeros in the disk $\mathbb{D}(0, R)$, and verifies*

$$(\star) \qquad \lim_{s \uparrow \ln(R)} \frac{\sup_{|\phi| \leq A} |F'''(s + i\phi)|}{F''(s)^{3/2}} = 0, \qquad \textit{for any } A > 0,$$

*then $f$ is Gaussian.*

In some favorable situations we will have $|F'''(s + i\phi)| \leq F'''(s)$ and therefore condition $(\star)$ translates into

$$\lim_{s \uparrow \ln(R)} \frac{F'''(s)}{F''(s)^{3/2}} = \lim_{s \uparrow \ln(R)} \frac{F'''(s)}{\sigma_f^3(t)} = 0.$$

This happens, for instance, for power series $f = e^g$, with $g \in \mathcal{K}_s$, see Chapter 4.

**Remark 3.1.5.** Observe that, under the hypothesis of the Theorem 3.1.4, the characteristic function $\mathbf{E}(e^{i\theta \check{X}_t})$ converges uniformly, in compact subsets of $\mathbb{R}$, as $t \uparrow R$, to the characteristic function of the standard normal $e^{-\theta^2/2}$. In fact we have

$$\lim_{t \uparrow R} \int_{-M}^{M} |\mathbf{E}(e^{i\theta \check{X}_t}) - e^{-\theta^2/2}| d\theta = 0, \qquad \text{for any } M \geq 0.$$

see the definition (3.2.1) of strongly Gaussian power series below. ⊠

We can generalize Theorem 3.1.4 in the following way: for $s < \ln(R)$ and $\theta \in \mathbb{R}$ we have

$$|F(s + i\theta) - F(s) - F'(s) i\theta + F''(s)\frac{\theta^2}{2} - F'''(s) i \frac{\theta^3}{3!}| \leq \sup_{\phi \in \mathbb{R}} |F^{(4)}(s + i\phi)| \frac{|\theta|^4}{4!}.$$

and therefore for any $s < \ln(R)$, changing $\theta$ by $\theta/\sigma_f(e^s)$, we find that

$$(\star\star) \qquad \left| \ln(\mathbf{E}(e^{i\theta \check{X}_t})) + \frac{\theta^2}{2} - \frac{F'''(s)}{\sigma_f^3(e^s)} i \frac{\theta^3}{3!} \right| \leq \frac{\sup_{\phi \in \mathbb{R}} |F^{(4)}(s + i\phi)|}{\sigma_f^4(e^s)} \frac{|\theta|^4}{4!}.$$

As a consequence of the bound $(\star\star)$ we obtain the following result.

**Theorem 3.1.6.** *Let $f \in \mathcal{K}$ be a power series with radius of convergence $R > 0$ and having no zeros in the disk $\mathbb{D}(0, R)$. Assume that*

$$(\dagger) \qquad \limsup_{s \uparrow \ln(R)} \frac{\sup_{\phi \in \mathbb{R}} |F^{(4)}(s + i\phi)|}{\sigma_f^4(e^s)} = 0,$$

*then*

$$f \textit{ is Gaussian} \qquad \textit{if and only if} \qquad \lim_{s \uparrow \ln(R)} \frac{F'''(s)}{F''(s)^{3/2}} = 0.$$

*Proof.* Denote $t = e^s$, for $s < \ln(R)$ and assume that $f$ is Gaussian. Using inequality ($\star\star$), combined with the triangle inequality, we find that

$$\limsup_{s\uparrow\ln(R)} \frac{|F'''(s)|}{\sigma_f^3(e^s)} \frac{|\theta|^3}{3!} \leq \limsup_{s\uparrow\ln(R)} \frac{\sup_{\phi\in\mathbb{R}} |F^{(4)}(s+i\phi)|}{\sigma_f^4(e^s)} \frac{|\theta|^4}{4!},$$

therefore, for any $\theta \neq 0$, we have

$$\limsup_{s\uparrow\ln(R)} \frac{|F'''(s)|}{\sigma_f^3(e^s)} \frac{1}{3!} \leq \limsup_{s\uparrow\ln(R)} \frac{\sup_{\phi\in\mathbb{R}} |F^{(4)}(s+i\phi)|}{\sigma_f^4(e^s)} \frac{|\theta|}{4!}.$$

Making $\theta \to 0$, we conclude that

$$\limsup_{s\uparrow\ln(R)} \frac{F'''(s)}{\sigma_f^3(e^s)} = 0.$$

Conversely, assume that

$$\limsup_{s\uparrow\ln(R)} \frac{F'''(s)}{\sigma_f^3(e^s)} = 0,$$

using, again, the triangle inequality we find that

$$\limsup_{s\uparrow\ln(R)} \left|\ln \mathbf{E}(e^{i\theta\breve{X}_t}) + \frac{\theta^2}{2}\right| \leq \limsup_{s\uparrow\ln(R)} \frac{\sup_{\phi\in\mathbb{R}} |F^{(4)}(s+i\phi)|}{\sigma_f^4(e^s)} \frac{|\theta|^4}{4!},$$

and therefore our hypothesis (†) implies that $f$ is Gaussian. □

**Remark 3.1.7.** In general we have the inequality

$$|F^{(4)}(s)| \leq \sup_{\phi\in\mathbb{R}} |F^{(4)}(s+i\phi)|, \quad \text{for any } s < \ln(R),$$

then, using (†), we also obtain that

$$\limsup_{s\uparrow\ln(R)} \frac{|F^{(4)}(s)|}{\sigma_f(e^s)^4} \leq \limsup_{s\uparrow\ln(R)} \frac{\sup_{\phi\in\mathbb{R}} |F^{(4)}(s+i\phi)|}{\sigma_f(e^s)^4} = 0.$$

⊠

The following corollary is implicit in the proof of Theorem 3.1.6.

**Corollary 3.1.8.** *Let $f \in \mathcal{K}$ be a power series with radius of convergence $R > 0$ and having no zeros in the disk $\mathbb{D}(0,R)$. Assume that*

$$(\dagger) \qquad \limsup_{s\uparrow\ln(R)} \frac{\sup_{\phi\in\mathbb{R}} |F^{(4)}(s+i\phi)|}{\sigma_f^4(e^s)} < +\infty,$$

*then*

$$f \text{ Gaussian} \quad \text{implies that} \quad \lim_{s\uparrow\ln(R)} \frac{F'''(s)}{F''(s)^{3/2}} = 0.$$

**Remark 3.1.9** (Generalization of the previous Theorems)**.** Any of the previous Theorems can be generalized. First assume that for some $k \geq 4$ we have

$$\limsup_{s \uparrow \ln(R)} \frac{\sup_{\phi \in \mathbb{R}} |F^{(k)}(s+i\phi)|}{F''(s)^{k/2}} = 0,$$

in this case

$$f \text{ is Gaussian} \qquad \text{if and only if} \qquad \lim_{s \uparrow \ln(R)} \frac{F^{(j)}(s)}{F''(s)^{j/2}} = 0, \text{ for any } 3 \leq j \leq k-1.$$

If we assume that for some $k \geq 4$ we have

$$\limsup_{s \uparrow \ln(R)} \frac{\sup_{\phi \in \mathbb{R}} |F^{(k)}(s+i\phi)|}{F''(s)^{k/2}} < +\infty,$$

then

$$f \text{ is Gaussian} \qquad \text{implies that} \qquad \lim_{s \uparrow R} \frac{F^{(j)}(s)}{F''(s)^{j/2}} = 0, \text{ for any } 3 \leq j \leq k-1.$$

⊠

**Remark 3.1.10** (Comparing derivatives)**.** In this remark we compare different derivatives of the fulcrum.

**Lemma 3.1.11.** *Let $f \in \mathcal{K}$ be a power series with radius of convergence $R > 0$ and* having no zeros in the disk $\mathbb{D}(0, R)$, *then*

$$|F^{(k)}(s+i\theta)| \leq |F^{(k)}(s)| + \left( \sup_{\alpha \in \mathbb{R}} |F^{(k+1)}(s+i\alpha)| \right) |\theta| \,,$$

*for any $s < \ln(R)$ and $\theta \in \mathbb{R}$.*

*Proof.* For any $s < \ln(R)$ and $\theta \in \mathbb{R}$ we have

$$F^{(k)}(s+i\theta) - F^{(k)}(s) = \int_0^\theta F^{(k+1)}(s+i\alpha) i d\alpha$$

therefore

$$|F^{(k)}(s+i\theta)| \leq |F^{(k)}(s)| + \left( \sup_{\alpha \in \mathbb{R}} |F^{(k+1)}(s+i\alpha)| \right) |\theta|.$$

□

**Lemma 3.1.12.** *Let $f \in \mathcal{K}$ be a power series with radius of convergence $R > 0$ and having no zeros in the disk $\mathbb{D}(0, R)$. Fix $k \geq 3$ and assume that*

$$\limsup_{t\uparrow R} \frac{\sup_{\alpha\in\mathbb{R}} |F^{(k+1)}(s+i\alpha)|}{F''(s)^{k/2}} = 0,$$

*then*

$$\limsup_{t\uparrow R} \frac{|F^{(k)}(s+i\theta)|}{F''(s)^{k/2}} \leq \limsup_{t\uparrow R} \frac{|F^{(k)}(s)|}{F''(s)^{k/2}}, \qquad \textit{for any } \theta \in \mathbb{R}.$$

*and, for any $A > 0$,*

$$\limsup_{s\uparrow \ln(R)} \frac{\sup_{|\theta|\leq A} |F^{(k)}(s+i\theta)|}{F''(s)^{k/2}} \leq \limsup_{s\uparrow \ln(R)} \frac{|F^{(k)}(s)|}{F''(s)^{k/2}},$$

*Proof.* Applying Lemma 3.1.11, we conclude that

$$\limsup_{s\uparrow \ln(R)} \frac{|F^{(k)}(s+i\theta)|}{F''(s)^{k/2}} \leq \limsup_{s\uparrow \ln(R)} \frac{|F^{(k)}(s)|}{F''(s)^{k/2}}.$$

and also that for any $A > 0$ we have

$$\limsup_{s\uparrow \ln(R)} \frac{\sup_{|\theta|\leq A} |F^{(k)}(s+i\theta)|}{F''(s)^{k/2}} \leq \limsup_{s\uparrow \ln(R)} \frac{|F^{(k)}(s)|}{F''(s)^{k/2}}.$$

□

**Corollary 3.1.13.** *Let $f \in \mathcal{K}$ be a power series with radius of convergence $R > 0$ and having no zeros in the disk $\mathbb{D}(0, R)$. Fix $k \geq 3$, and assume that*

$$C_{k+1} = \limsup_{s\uparrow \ln(R)} \frac{\sup_{\alpha\in\mathbb{R}} |F^{(k+1)}(s+i\alpha)|}{F''(s)^{k/2}} < +\infty,$$

*then, for any $A > 0$, we have that*

$$\limsup_{s\uparrow \ln(R)} \frac{\sup_{|\theta|\leq A} |F^{(k)}(s+i\theta)|}{F''(s)^{k/2}} \leq \limsup_{s\uparrow \ln(R)} \frac{|F^{(k)}(s)|}{F''(s)^{k/2}} + C_{k+1}A\,.$$

Fix $A > 0$. For any $s < \ln(R)$ we always have that

$$|F^{(k)}(s)| \leq \sup_{|\theta|\leq A} |F^{(k)}(s+i\theta)|\,,$$

and therefore, under the conditions of Corollary 3.1.13, we have that

$$\limsup_{s\uparrow \ln(R)} \frac{|F^{(k)}(s)|}{F''(s)^{k/2}} < +\infty \quad \text{if and only if} \quad \limsup_{s\uparrow \ln(R)} \frac{\sup_{|\theta|\leq A} |F^{(k)}(s+i\theta)|}{F''(s)^{k/2}} < +\infty.$$

### B. Direct applications of Theorem 3.1.4

Now we apply Theorems 3.1.4 to some examples of Khinchin families.

**B.1. Exponential of a polynomial not in $\mathcal{K}$** Here we prove, using the global criteria (3.1.4), that the exponential of a polynomial $p$, not necessarily in $\mathcal{K}$, but with $e^p \in \mathcal{K}$ is Gaussian.

**Proposition 3.1.14.** *Let $P$ be a polynomial and denote $f(z) = e^{P(z)}$. If the power series $f = e^P$ is in $\mathcal{K}$, then $f = e^P$ is Gaussian.*

*Proof.* Denote $P(z) = \sum_{j=0}^{N} b_j z^j$, a polynomial of degree $N \geq 1$. The fulcrum of $f = e^P$ is given by $F(z) = P(e^z)$ and the radius of convergence is $R = \infty$.

We have

$$F'''(z) = e^{3z}P'''(z) + 3e^{2z}P''(z) + e^z P'(e^z) = \sum_{j=0}^{N} j^3 b_j e^{jz}\,,$$

and therefore

$$|F'''(s+i\theta)| \leq \sum_{j=0}^{N} j^3 |b_j| e^{js}\,.$$

Moreover,

$$\sum_{j=0}^{N} j^3 |b_j| e^{js} \sim N^3 |b_N| e^{Ns}, \qquad \text{as } s \to +\infty.$$

On the other hand, we have

$$F''(s) = \sum_{j=0}^{N} j^2 b_j e^{js} \sim N^2 b_N e^{Ns}, \qquad \text{as } s \to +\infty,$$

therefore

$$\lim_{s\to+\infty} \frac{\sup_{\theta\in\mathbb{R}} |F'''(s+i\theta)|}{F''(s)^{3/2}} = 0.$$

The previous limit gives, by virtue of Theorem 3.1.4, that $f = e^P$ is Gaussian. Observe that $F''(s) = \sigma_f^2(e^s) > 0$, and this implies that $b_N$ is a positive real number. □

**B.2. Gaussianity of the Poisson family** The Poisson family is the Khinchin family associated to $f(z) = e^z$. The fulcrum of $f$ is given by $F(z) = e^z$, for any $z \in \mathbb{C}$, then we have

$$\lim_{s\to+\infty} \frac{\sup_{\theta\in\mathbb{R}} |F'''(s+i\theta)|}{F''(s)^{3/2}} = \lim_{s\to+\infty} \frac{F'''(s)}{F''(s)^{3/2}} = \lim_{s\to+\infty} e^{-s/2} = 0,$$

and therefore, applying Theorem 3.1.4, we conclude that $f(z) = e^z$ is a Gaussian power series.

**B.3. Gaussianity of the OGF of the involutions** The generating function of the involutions is given by $I(z) = \exp(z + z^2/2)$. The fulcrum of $I(z)$ is given by $F(z) = e^z + e^{2z}/2$, therefore

$$\lim_{s\to+\infty} \frac{\sup_{\theta\in\mathbb{R}} |F'''(s+i\theta)|}{F''(s)^{3/2}} = \lim_{s\to+\infty} \frac{F'''(s)}{F''(s)^{3/2}} = \lim_{s\to+\infty} \frac{e^s + 4e^{2s}}{(e^s + 2e^{2s})^{3/2}} = 0\,.$$

Applying Theorem 3.1.4, we conclude that $I$ is a Gaussian power series.

**B.4. Gaussianity of the Euler family** Recall that the Euler family is the Khinchin family associated to the generating function of the partitions of integers $P(z)$. Also recall that

$$P(z) = \prod_{j=1}^{\infty} \frac{1}{1-z^j}, \qquad \text{for any } |z| < 1.$$

The fulcrum of $F$ is given by

$$F(z) = \ln(P(e^z)) = -\sum_{j=1}^{\infty} \log(1 - e^{jz}) = \sum_{j,k\geq 1} \frac{e^{kjz}}{k}, \qquad \text{for any } \Re(z) < 0, \tag{3.1.6}$$

and therefore

$$F'(z) = \sum_{j,k\geq 1} j e^{kjz}\,,$$
$$F''(z) = \sum_{j,k\geq 1} kj^2 e^{jkz}\,,$$
$$F'''(z) = \sum_{j,k\geq 1} k^2 j^3 e^{kjz}\,.$$

For any $q \geq 1$, the Fulcrum of $f$, and all its derivatives, verify the inequality

$$|F^{(q)}(s+i\theta)| \leq F^{(q)}(s)\,, \qquad \text{for any } s < 0 \text{ and } \theta \in \mathbb{R}.$$

We want to check that

$$\lim_{s\downarrow 0} \frac{F'''(-s)}{F''(-s)^{3/2}} = 0.$$

Let's see: for any $s > 0$ we have

$$s^2 F'(-s) = \sum_{k\geq 1} \sum_{j\geq 1} (js) e^{-k(js)}\, s,$$

therefore

$$\lim_{s\downarrow 0} s^2 F'(-s) = \sum_{k\geq 1} \int_0^\infty x e^{-kx} dx = \sum_{k\geq 1} \frac{1}{k^2} \int_0^\infty y e^{-y} dy = \zeta(2)\Gamma(2). \tag{3.1.7}$$

In general, by using the same argument, we obtain that for any $q \geq 1$, we have

$$\lim_{s \downarrow 0} s^{q+1} F^{(q)}(-s) = \zeta(2)\Gamma(q+1),$$

then

$$F^{(q)}(-s) \sim \zeta(2)\Gamma(q+1)s^{-(q+1)}, \quad \text{as } s \downarrow 0. \tag{3.1.8}$$

We conclude then that

$$\frac{F'''(-s)}{F''(-s)^{3/2}} \sim \frac{3}{\sqrt{2\zeta(2)}} s^{1/2}, \quad \text{as } s \downarrow 0,$$

therefore $P$ is a Gaussian power series.

**Remark 3.1.15.** More generally, for the generating function of the partitions of integers $P(z)$, and for any $n \geq 3$, using an analogous argument to that used in (3.1.7), we have

$$\frac{F^{(n)}(-s)}{F''(-s)^{n/2}} \sim \frac{\zeta(2)\Gamma(n+1)}{(\zeta(2)\Gamma(3))^{n/2}} s^{n/2-1}, \quad \text{as } s \downarrow 0,$$

therefore

$$\lim_{s \downarrow 0} \frac{\sup_{\theta \in \mathbb{R}} |F^{(n)}(-s+i\theta)|}{F''(-s)^{n/2}} \leq \lim_{s \downarrow 0} \frac{F^{(n)}(-s)}{F''(-s)^{n/2}} = 0.$$

Here we use that $|F^{(n)}(s+i\theta)| \leq F^{(n)}(s)$, for any $s < 0$ and $\theta \in \mathbb{R}$, see equation (3.1.6). ⊠

### C. A local criteria for Gaussianity for non-vanishing power series

Let $f \in \mathcal{K}$ be a power series with radius of convergence $R > 0$. The function $f$ has no zeros in the interval $[0, R)$ and, therefore, $f$ has no zeros in a simply connected domain containing $[0, R)$. For instance, $f$ has no zeros, in the region

$$\Omega_f = \left\{ z = te^{\imath\theta} : |\theta| < \frac{\pi}{2\sigma_f(t)} \right\}.$$

see Subsection 1.3.4 for further details.

For each $t \in (0, R)$ we give a disk $\Gamma_t = \mathbb{D}(t, 2t\mu(t))$, of radius $2t\mu(t)$, where $0 < \mu(t) \leq 1$ and $2t\mu(t) < R - t$, that is, $\Gamma_t \subset \mathbb{D}(0, R)$, and such that $f$ has no zeros in $\Gamma_t$.

We insists: *the function $f$ has no zeros in the disk $\Gamma_t$, for any $t \in (0, R)$.*

The fulcrum of $f$, that is, the function $F(z) = \log(f(e^z))$, is well defined and is holomorphic in the disk $\mathbb{D}(\ln(t), \mu(t))$, for any $t > 0$. Observe that for any $z \in \mathbb{D}(\ln(t), \mu(t))$, we have $e^z \in \Gamma_t$, because if $|z - \ln(t)| < \mu(t)$, then

$$|e^z - t| = t|e^{z-\ln(t)} - 1| \leq t(e^{\mu(t)} - 1) \leq 2t\mu(t),$$

here we use that $e^x - 1 \leq 2x$, for any $x \in (0,1]$.

We can extend the variance $\sigma_f^2$ to a holomorphic function in a domain by defining $\sigma_f^2(z) = F''(\ln(z))$, for any $z \in \Gamma_t$, and any $t \in (0,R)$ or, equivalently, by defining $\sigma_f^2(e^z) = F''(z)$, for any $z \in \mathbb{D}(\ln(t), \mu(t))$ and any $t \in (0,R)$.

We denote by

$$A(t) = \sup_{z\in\Gamma_t} \frac{|\sigma_f^2(z)|}{\sigma_f^2(t)} = \sup_{z\in\Gamma_t} \frac{|F''(\ln(z))|}{F''(\ln(t))} < +\infty, \qquad \text{for any } t \in (0,R).$$

The next lemma gives an approximation of the holomorphic function $h$ by it's Taylor polynomial of second order and an upper bound for $h''$ instead of an upper bound for $h'''$. See [48, p. 78].

**Lemma 3.1.16.** *Let $h(z) = \sum_{n=0}^{\infty} a_n z^n$ be a power series with radius of convergence at least 1. Suppose that $|h''(z)| \leq M_2$, for any $z \in \mathbb{D}$, then*

$$\left| h(z) - h(0) - h'(0)z - \frac{h''(0)}{2} z^2 \right| \leq M_2 |z|^3, \qquad \textit{for any } z \in \mathbb{D}. \tag{3.1.9}$$

*Proof.* Cauchy's estimates for the coefficients of $h''$ give that, for any $n \geq 3$ and any $0 < r < 1$, we have

$$n(n-1)|a_n| \leq \frac{M_2}{r^{n-2}},$$

therefore making $r \uparrow 1$ we obtain that

$$|a_n| \leq \frac{M_2}{n(n-1)}, \qquad \text{for any } n \geq 3,$$

then the left-hand side of (3.1.9) verifies that

$$\left| \sum_{n=3}^{\infty} a_n z^{n-3} \right| |z|^3 \leq M_2 \left( \sum_{n=3}^{\infty} \frac{1}{n(n-1)} \right) |z|^3 = \frac{M_2}{2} |z|^3 \leq M_2 |z|^3, \qquad \text{for any } z \in \mathbb{D}.$$

□

For any holomorphic function in a disk $\mathbb{D}(a,s)$ and verifying that $|h''(s)| \leq M_2$, for any $z \in \mathbb{D}(a,s)$, Lemma 3.1.16 gives, rescaling, that

$$\left| h(z) - h(a) - h'(a)(z-a) - \frac{h''(s)}{2}(z-a)^2 \right| \leq \frac{M_2}{s} |z-a|^3, \qquad \text{for any } z \in \mathbb{D}(a,s).$$

Our function $h$ and the disk $\mathbb{D}(a,s)$ will be the fulcrum $F(z)$ and the disk $\mathbb{D}(\ln(t), \mu(t))$. We obtain, writing in terms of $m_f(t)$ and $\sigma_f^2(t)$, that

$$\left| F(z) - \ln(f(t)) - m_f(t)(z - \ln(t)) - \frac{\sigma_f^2(t)}{2}(z - \ln(t))^2 \right| \leq \frac{A(t)\sigma_f^2(t)}{\mu(t)} |z - \ln(t)|,$$

for any $z \in \mathbb{D}(\ln(t), \mu(t))$.

In particular, for any $z = \ln(t) + i\theta$ with $|\theta| < \mu(t)$, we have that

$$\left| \ln(f(te^{i\theta})) - \ln(f(t)) - m_f(t) i\theta + \sigma_f^2(t)\frac{\theta^2}{2} \right| \leq \frac{A(t)\sigma_f^2(t)}{\mu(t)} |\theta|^3.$$

If we take $\theta = \phi/\sigma_f(t)$, and impose that $|\theta| \leq \sigma_f(t)\mu(t)$, we obtain, passing to characteristic functions, that

$$\left| \ln(\mathbf{E}(e^{i\phi \breve{X}_t})) + \frac{\phi^2}{2} \right| \leq \frac{A(t)}{\mu(t)\sigma_f(t)} |\phi|^3 , \qquad \text{for any } |\phi| \leq \mu(t)\sigma_f(t) .$$

**Theorem 3.1.17.** *With the previous notations, suppose that*

$$\limsup_{t \uparrow R} A(t) < +\infty, \qquad \textit{and} \qquad \lim_{t \uparrow R} \mu(t)\sigma_f(t) = +\infty,$$

*then $f$ is Gaussian.*

*Proof.* Fix $\phi \in \mathbb{R}$. There exists $t_0 > 0$ such that $|\phi| \leq \mu(t)\sigma_f(t)$, for any $t \geq t_0$, therefore,

$$\lim_{t \uparrow R} \left( \ln(\mathbf{E}(e^{i\phi \breve{X}_t})) + \frac{\phi^2}{2} \right) = 0 .$$

□

### 3.1.2 Moment criteria for general $f \in \mathcal{K}$ to be Gaussian

Let $f$ be a power series in $\mathcal{K}$ with radius of convergence $R > 0$ with Khinchin family $(X_t)_{t \in [0,R)}$. *Notice that in this case we don't assume that $f$ is non-vanishing.* Let $F$ denote, as usual, the fulcrum of $f$, which is defined in the interval $(-\infty, \ln R)$:

$$e^{F(s)} = f(e^s) , \quad \text{for } s < \ln R .$$

To check gaussianity of the family $X_t$ we will use the so-called method of moments of convergence of probability distributions. See, for instance, Billingsley, [9, Section 30]. We let $Z$ denote a random variable with standard normal distribution. Since the standard normal is determined by its moments and since the normalized variables $\breve{X}_t$ have moments of all orders, we have, see [9, Theorem 30.2] that

**Lemma 3.1.18.**

$$\textit{if } \lim_{t \uparrow R} \frac{\nu_k(t)}{\sigma_f(t)^k} = \mathbf{E}(Z^k) , \textit{ for every } k \geq 3, \quad \textit{then the family } X_t \textit{ is gaussian} .$$

Recall that $\nu_k(t)$ denotes the $k$-th central moment of $\breve{X}_t$.

Using the previous lemma we obtain that

**Theorem 3.1.19.** *With the notations above, if*

$$\lim_{s\uparrow \ln R} \frac{F^{(k)}(s)}{F''(s)^{k/2}} = 0\,, \quad \textit{for every } k \geq 3\,,$$

*then $f$ is gaussian.*

*Proof.* For $k \geq 3$ we have, appealing to (1.4.6) and to (1.4.12), we have that

$$\begin{aligned} R_k\Big(\frac{F''(s)}{F''(s)^{2/2}}, \frac{F^{(3)}(s)}{F''(s)^{3/2}}, \ldots, \frac{F^{(k)}(s)}{F''(s)^{k/2}}\Big) &= \frac{1}{F''(s)^{k/2}} R_k\Big(F''(s), F^{(3)}(s), \ldots, F^{(k)}(s)\Big) \\ &= \frac{\nu_k(s)}{F''(s)^{k/2}}\,. \end{aligned}$$

From the hypothesis of the theorem and because of the identities (1.4.7), we deduce, as desired, that

$$\lim_{s\uparrow \ln R} \frac{\nu_k(s)}{F''(s)^{k/2}} = R_k(1, 0, \ldots, 0) = \mathbf{E}(Z^k) \quad \text{for every } k \geq 2\,.$$

then Lemma 3.1.18 gives that $f$ is Gaussian. □

**Canonical products in the Class $\mathcal{Q}$** As an example of the use of the gaussianity criterium of Theorem 3.1.19 we next consider the case of an entire function $f$ in $\mathcal{Q} \subseteq \mathcal{K}$, that is, entire functions of genus 0 whose only zeros are negative and normalized, i.e., $f(0) = 1$. Because of Hadamard's factorization Theorem, such an $f$ is a canonical product and may be expressed as

$$(\star) \quad f(z) = \prod_{j=1}^{\infty} \left(1 + \frac{z}{b_j}\right), \quad \text{for any } z \in \mathbb{C}.$$

Recall, see Section 2.5.4, that

$$m_f(t) = \sum_{j=1}^{\infty} \frac{t}{(b_j + t)} \quad \text{and} \quad \sigma_f(t)^2 = \sum_{j=1}^{\infty} \frac{t b_j}{(b_j + t)^2}\,, \quad \text{for any } t \geq 0\,,$$

and also that $\sigma_f^2(t) < m_f(t)$, for any $t > 0$.

If in the expression for $\sigma_f^2(t)$ and for each $t > 0$ we split the summands into those with $b_j < t$ and those with $b_j \geq t$ we obtain that

$$\frac{1}{4}\Big(\frac{1}{t}\sum_{b_j < t} b_j + t\sum_{b_j \geq t} \frac{1}{b_j}\Big) \leq \sigma_f^2(t) \leq \Big(\frac{1}{t}\sum_{b_j < t} b_j + t\sum_{b_j \geq t} \frac{1}{b_j}\Big)\,, \quad \text{for } t > 0\,.$$

The fulcrum $F$ of $f$ is given by

$$F(s) = \ln(f(e^s)) = \sum_{j=1}^{\infty} \ln\left(1 + \frac{e^s}{b_j}\right)\,, \quad \text{for any } s \in \mathbb{R}\,.$$

**Lemma 3.1.20.** *For each integer $k \geq 2$, there exists a constant $C_k > 0$ such that for any $f$ as in $(\star)$ and any $s \in \mathbb{R}$*

$$\left|F^{(k)}(s)\right| \leq C_k F''(s)\,.$$

*Proof.* The first two derivatives of $F$ are given by

$$F'(s) = \sum_{j=1}^{\infty} \frac{e^s}{(b_j + e^s)}\,, \quad \text{and} \quad F''(s) = \sum_{j=1}^{\infty} \frac{b_j e^s}{(b_j + e^s)^2}\,, \quad \text{for any } s \in \mathbb{R}\,.$$

For any integers $n \geq 1$ and $b > 0$, and any $s \in \mathbb{R}$, denote

$$A_n(b, s) = \frac{b e^{ns}}{(b + e^s)^{n+1}}\,.$$

For any $n \geq 1, b > 0$ and $s \in \mathbb{R}$ we have the inequality

$$A_n(b, s) = \frac{b e^{ns}}{(b + e^s)^{n+1}} = \left(\frac{e^s}{b + e^s}\right)^{n-1} \frac{b e^s}{(b + e^s)^2} \leq A_1(b, s)\,.$$

and also the differential identity

$$\frac{d}{ds} A_n(b, s) = n A_n(b, s) - (n + 1) A_{n+1}(b, s)\,.$$

By iteration, we obtain, from the identity above, and for any $k \geq 2$, real coefficients $\delta_{k,l}$, for $k \geq l \geq 0$, so that

$$\frac{d^k}{ds^k} A_1(b, s) = \sum_{l=0}^{k} \delta_{k,l} A_{1+l}(b, s)\,, \quad \text{for any } s \in \mathbb{R} \text{ and any } b > 0\,,$$

and thus that

$$\left|\frac{d^k}{ds^k} A_1(b, s)\right| \leq \sum_{l=0}^{k} |\delta_{k,l}| A_{1+l}(b, s) \leq \left(\sum_{l=0}^{k} |\delta_{k,l}|\right) A_1(b, s) > 0, \quad \text{for any } s \in \mathbb{R} \text{ and any } b > 0\,.$$

For $k \geq 2$, since $F^{(k)}(s) = \sum_{j=1}^{\infty} A_1^{(k-2)}(b_j, s)$, we conclude that

$$|F^{(k)}(s)| \leq \sum_{j=1}^{\infty} \left|\frac{d^{(k-2)}}{ds} A_1(b_j, s)\right| \leq \underbrace{\left(\sum_{l=0}^{k-2} |\delta_{k-2,l}|\right)}_{C_k :=} \sum_{j=1}^{\infty} A_1(b_j, s) = C_k F''(s)\,, \quad \text{for any } s \in \mathbb{R}\,.$$

□

We use the lemma and Theorem 3.1.19 to prove

**Theorem 3.1.21.** *Let $f$ be a canonical product in $\mathcal{Q}$ given by $(\star)$. If $\lim_{t\to+\infty}\sigma_f^2(t)=+\infty$, then $f$ is Gaussian.*

*Proof.* We use Theorem 3.1.19. We are going to show that for any $k\geq 3$

$$\lim_{s\to\infty}\frac{|F^{(k)}(s)|}{F''(s)^{k/2}}=0\,.$$

But from the lemma above

$$\frac{|F^{(k)}(s)|}{F''(s)^{k/2}}=\frac{|F^{(k)}(s)|}{F''(s)}\cdot F''(s)^{1-k/2}\leq C_k\sigma_f(e^s)^{2-k},$$

which tends to 0, as $s\to\infty$ since $\lim_{s\to\infty}\sigma_f(e^s)=+\infty$ and $k\geq 3$. □

**Example 3.1.22.** Let $f\in\mathcal{Q}$ be a canonical product as in $(\star)$. Assume that

$$N(t)\sim t^\rho,\qquad \text{as } t\to+\infty,$$

then Proposition 2.5.4 gives that $\lim_{t\to+\infty}\sigma_f^2(t)=+\infty$ and therefore Theorem 3.1.21 implies that $f$ is a Gaussian power series. ⊡

### 3.1.3 Operations with Gaussian power series

Below these lines we study the Gaussianity of power series obtained by performing certain operations with one, two, or more, Gaussian power series. Among these operations we find the product of Gaussian power series or the compositions of a polynomial with a Gaussian power series.

#### A. Subordination of Gaussian power series

Let $f\in\mathcal{K}$ be a power series with radius of convergence $R>0$.

**Proposition 3.1.23.** *Fix an integer $N\geq 1$, then $f$ is Gaussian if and only if $h(z)=f(z^N)\in\mathcal{K}$, and $h$ is Gaussian.*

*Proof.* Assume that $f\in\mathcal{K}$ is Gaussian, and denote $g(z)=f(z^N)$, for any $|z|<R^{1/N}$. If $X_t$ is the Khinchin family associated to $f$ and $Y_t$ is the Khinchin family associated to $g$ we have

$$Y_t\stackrel{d}{=}NX_{t^N},\qquad \text{for any } t\in[0,R^{1/N}),$$

then

$$m_g(t)=Nm_f(t^N),\qquad \text{and}\qquad \sigma_g^2(t)=N^2\sigma_f^2(t^N),$$

for any $t\in[0,R^{1/N})$. Observe that

$$\breve{Y}_t=\frac{Y_t-m_g(t)}{\sigma_g(t)}\stackrel{d}{=}\frac{N\,X_{t^N}-N\,m_f(t^N)}{N\sigma_f(t^N)}\stackrel{d}{=}\breve{X}_{t^N},\qquad \text{for any } t\in(0,R^{1/N}),\tag{3.1.10}$$

therefore

$$\lim_{t\uparrow R^{1/N}} \mathbf{E}(e^{i\theta \check{Y}_t}) = \lim_{t\uparrow R^{1/N}} \mathbf{E}(e^{i\theta \check{X}_{t^N}}) = e^{-\theta^2/2}\,,$$

then $g(z) = f(z^N)$ is a Gaussian power series. The other implication follows by using again the equality in distribution (3.1.10). □

### B. Products of Gaussian power series

Here we study the product of two, or more, Gaussian power series. Under certain conditions, the product of Gaussian power series is also a Gaussian power series. We distinguish several cases: products of power series with the same radius of convergence and products of power series with different radius of convergence.

**B.1. Power series with the same radius of convergence** Assume that $f, g \in \mathcal{K}$ are Gaussian power series with radius of convergence $R > 0$. Denote $h = fg$, the product power series.

Let $(X_t), (Y_t)$ and $(Z_t)$ the Khinchin families of $f, g$ and $h$ respectively, then

$$Z_t \stackrel{d}{=} X_t \oplus Y_t, \qquad \text{for any } t \in [0, R),$$

that is, $Z_t$ is equal, in distribution, to the sum of two independent random variables. In fact we have

$$m_h(t) = m_f(t) + m_g(t), \qquad \text{and} \qquad \sigma_h^2(t) = \sigma_f^2(t) + \sigma_g^2(t),$$

for any $t \in [0, R)$.

Assume that

$$\lim_{t\uparrow R} \frac{\sigma_f^2(t)}{\sigma_f^2(t) + \sigma_g^2(t)} = \lambda^2, \qquad \text{with } \lambda \in [0, 1].$$

Therefore the limit

$$\lim_{t\uparrow R} \frac{\sigma_g^2(t)}{\sigma_f^2(t) + \sigma_g^2(t)} = \mu^2, \qquad \text{with } \mu \in [0, 1].$$

also exists; notice that $\lambda^2 + \mu^2 = 1$.

We have

$$\check{Z}_t \stackrel{d}{=} \frac{\sigma_f(t)}{\sigma_h(t)} \check{X}_t \oplus \frac{\sigma_g(t)}{\sigma_h(t)} \check{Y}_t, \qquad \text{for any } t \in (0, R),$$

then

$$\mathbf{E}(e^{i\theta \check{Z}_t}) = \mathbf{E}\left(e^{i\theta \frac{\sigma_f(t)}{\sigma_h(t)} \check{X}_t}\right) \mathbf{E}\left(e^{i\theta \frac{\sigma_g(t)}{\sigma_h(t)} \check{Y}_t}\right), \tag{3.1.11}$$

and therefore

$$\lim_{t\uparrow R}\mathbf{E}(e^{i\theta\check{Z}_t}) = e^{-\lambda^2\theta^2/2}e^{-\mu^2\theta^2/2} = e^{-\theta^2/2}.$$

Here we use that, for any $a > 0$ and any $t > a$, the characteristic function of the family $(\check{X}_t)$ is a continuous function, as function of two variables, see Section 6.1.4, and in particular Lemma 1.3.4.

We collect this result by the following Proposition.

**Proposition 3.1.24.** *Let $f, g \in \mathcal{K}$ power series with the same radius of convergence $R > 0$ and denote $h = fg$ the product power series. Assume that $f$ and $g$ are Gaussian and also that*

$$\lim_{t\uparrow R}\frac{\sigma_f^2(t)}{\sigma_f^2(t)+\sigma_g^2(t)} = \lambda^2, \qquad \textit{with } \lambda\in[0,1],$$

*then $h$ is a Gaussian power series.*

We also have the following result:

**Proposition 3.1.25.** *Let $f, g \in \mathcal{K}$ power series with the same radius of convergence $R > 0$ and denote $h = fg$ the product power series. Assume that $f$ is Gaussian and also that*

$$\lim_{t\uparrow R}\frac{\sigma_f^2(t)}{\sigma_f^2(t)+\sigma_g^2(t)} = 1,$$

*then $h$ is a Gaussian power series.*

*Proof.* With the same notations: for any $\theta\in\mathbb{R}$ we have

$$\left|\mathbf{E}\left(e^{i\theta\frac{\sigma_g(t)}{\sigma_h(t)}\check{Y}_t}\right)-1\right| \le |\theta|\frac{\sigma_g(t)}{\sigma_f(t)}\mathbf{E}(\check{Y}_t^2)^{1/2} = |\theta|\frac{\sigma_g(t)}{\sigma_h(t)},$$

here we use the inequality $|e^{i\theta}-1|\le|\theta|$, for any $\theta\in\mathbb{R}$, therefore, using that $\lim_{t\uparrow R}\sigma_g(t)/\sigma_h(t)=0$, we have that

$$\lim_{t\uparrow R}\mathbf{E}\left(e^{i\theta\frac{\sigma_g(t)}{\sigma_h(t)}\check{Y}_t}\right) = 1, \qquad \text{for any } \theta\in\mathbb{R}.$$

The result follows combining the Gaussianity of $f$ with equation (3.1.11). □

**B.2. Power series with different radius of convergence** Assume that $f, g \in \mathcal{K}$ are power series with finite radius of convergence $R > 0$ and $S > 0$, respectively, with $S > R$.

Denote by $h = fg$ the product of $f$ and $g$. This power series has radius of convergence $R > 0$, recall that $f, g \in \mathcal{K}$. Assume that $f$ is Gaussian and also that $\lim_{t\uparrow R}\sigma_f^2(t) = +\infty$, then

$$\lim_{t\uparrow R}\frac{\sigma_f^2(t)}{\sigma_f^2(t)+\sigma_g^2(t)} = 1,$$

and

$$\lim_{t\uparrow R} \frac{\sigma_g^2(t)}{\sigma_f^2(t)+\sigma_g^2(t)} = 0,$$

then, applying Lemma 1.3.4, we conclude that

$$\lim_{t\uparrow R} \mathbf{E}\left(e^{i\theta \frac{\sigma_g(t)}{\sigma_h(t)}\breve{Y}_t}\right) = 1, \qquad \text{for any } \theta \in \mathbb{R}\,.$$

We finally have, using independence, that

$$\lim_{t\uparrow R} \mathbf{E}(e^{i\theta \breve{Z}_t}) = \lim_{t\uparrow R} \mathbf{E}\left(e^{i\theta \frac{\sigma_f(t)}{\sigma_h(t)}\breve{X}_t}\right) \lim_{t\uparrow R} \mathbf{E}\left(e^{i\theta \frac{\sigma_g(t)}{\sigma_h(t)}\breve{Y}_t}\right) = e^{-\theta^2/2}, \qquad \text{for any } \theta \in \mathbb{R}\,.$$

Here we use, again, Lemma 1.3.4.

We collect the previous discussion by means of the following proposition.

**Proposition 3.1.26.** *Let $f, g$ be power series with radius of convergence $R > 0$ and $S > 0$, respectively, with $R > S$. Assume that $f$ is Gaussian and also that $\lim_{t\uparrow R} \sigma_f^2(t) = +\infty$, then $h = fg$ is Gaussian.*

### C. Powers of Gaussian power series

Let $f \in \mathcal{K}$ be a power series with radius of convergence $R > 0$. Fix an integer $N \geq 1$ and denote by $h(z) = f(z)^N$. Assume that $f$ is Gaussian, then we have that $f^N$ is also Gaussian.

Denote $(X_t)$ the Khinchin family associated to $f$ and $(Z_t)$ the Khinchin family associated to $h$, recall that

$$Z_t \stackrel{d}{=} X_t^{(1)} + \dots X_t^{(N)}, \qquad \text{for any } t \in [0, R).$$

Here $X_t^{(1)}, \dots, X_t^{(N)}$ are i.i.d copies of $X_t$.

The previous equality in distribution gives that

$$m_h(t) = N m_f(t), \qquad \text{and} \qquad \sigma_h^2(t) = N \sigma_f(t)^2, \qquad \text{for any } t \in [0, R),$$

moreover

$$\mathbf{E}(e^{i\theta \breve{Z}_t}) = \mathbf{E}(e^{i\theta \breve{X}_t/\sqrt{N}})^N, \qquad \text{for any } \theta \in \mathbb{R} \text{ and } t \in (0, R),$$

therefore

$$\lim_{t\uparrow R} \mathbf{E}(e^{i\theta \breve{Z}_t}) = \lim_{t\uparrow R} \mathbf{E}\left(e^{i\theta \breve{X}_t/\sqrt{N}}\right)^N = \left(e^{-\theta^2/2N}\right)^N = e^{-\theta^2/2}, \qquad \text{for any } \theta \in \mathbb{R}\,.$$

We collect this results by means of the following proposition.

**Proposition 3.1.27.** *Let $f \in \mathcal{K}$ be a power series with radius of convergence $R > 0$. Fix an integer $N \geq 1$ and denote by $h(z) = f(z)^N$ the $N$-th power of $f$. Assume that $f$ is Gaussian, then $h(z) = f(z)^N$ is Gaussian.*

### D. Shifts of Gaussian power series

Let $f \in \mathcal{K}$ a power series with radius of convergence $R > 0$. Let $C > 0$ be a constant, then, for any integer $M \geq 0$, the power series $h(z) = Cz^M f(z) \in \mathcal{K}_s$.

Denote $(X_t)$ the Khinchin family associated to $f$ and $(Y_t)$ the Khinchin family associated to $h$, recall that

$$Y_t \stackrel{d}{=} M + X_t, \quad \text{for any } t \in [0, R).$$

The previous equality in distribution gives that

$$m_h(t) = M + m_f(t), \quad \text{and} \quad \sigma_h^2(t) = \sigma_f^2(t), \quad \text{for any } t \in [0, R),$$

and also that

$$\breve{Y}_t \stackrel{d}{=} \breve{X}_t, \quad \text{for any } t \in (0, R).$$

We collect the previous discussion in the following Proposition.

**Proposition 3.1.28.** *Assume that $f \in \mathcal{K}$ is Gaussian, then, for any $C > 0$ and any integer $M \geq 0$, the power series $h(z) = Cz^M f(z)$ is Gaussian.*

### E. Composition of a Gaussian power series with a polynomial

Here we prove that under certain conditions the composition $P \circ f$ of a polynomial in $P \in \mathcal{K}$ with a Gaussian power series $f \in \mathcal{K}$ is Gaussian. First, we prove some auxiliary lemmas.

**Lemma 3.1.29.** *Let $P \in \mathcal{K}$ be a polynomial of degree $\deg(P) = N \geq 1$ and $f \in \mathcal{K}$ a power series with radius of convergence $R > 0$. Denote $h = P \circ f$ and assume that $\lim_{t \uparrow R} f(t) = +\infty$ and also that*

$$(\star) \quad \lim_{t \uparrow R} \frac{m_f^2(t)}{\sigma_f^2(t) f(t)} = 0,$$

*then*

$$\lim_{t \uparrow R} \frac{\sigma_h^2(t)}{\sigma_f^2(t)} = N, \quad \textit{and} \quad \lim_{t \uparrow R} \frac{m_h(t) - N m_f(t)}{\sigma_f(t)} = 0.$$

*Proof.* Using equation (1.3.28) we have that

$$(3.1.12) \qquad \sigma_h^2(t) = \sigma_P^2(f(t)) m_f^2(t) + m_P(f(t)) \sigma_f^2(t), \quad \text{for any } t \in [0, R),$$

then

$$\frac{\sigma_h^2(t)}{m_P(f(t)) \sigma_f^2(t)} = 1 + \frac{1}{m_P(f(t))} \frac{m_f^2(t)}{\sigma_f^2(t) f(t)} \sigma_P^2(f(t)) f(t).$$

Recall that for a polynomial $P \in \mathcal{K}$ with $\deg(P) = N \geq 1$ we have

$$\lim_{t\to+\infty} m_P(t) = N,$$

see equation (1.2.4) and the discussion in there. Since $\lim_{t\uparrow R} f(t) = +\infty$, we have that $\lim_{t\uparrow R} m_P(f(t)) = N$. We also have that there exists a constant $C > 0$ and certain integer $j \geq 1$ such that

$$(\dagger) \quad \sigma_P^2(t) \sim C\frac{1}{t^j}, \qquad \text{as } t \to +\infty,$$

see equation (3.1.2) for further details. The previous asymptotic formula gives that $\sigma_P^2(f(t))f(t)$ is bounded (and goes to 0 if $j > 1$).

Combining all these facts with $(\star)$ we conclude that

$$\lim_{t\uparrow R} \frac{\sigma_h^2(t)}{\sigma_f^2(t) m_P(f(t))} = 1,$$

and therefore

$$\lim_{t\uparrow R} \frac{\sigma_h^2(t)}{\sigma_f^2(t)} = N.$$

For the second claim observe that

$$\frac{m_h(t) - N m_f(t)}{\sigma_f(t)} = \frac{m_f(t)\sigma_P(f(t))}{\sigma_f(t)} \frac{(m_P(f(t)) - N)}{\sigma_P(f(t))}$$

here we use that $m_h(t) = m_P(f(t))m_f(t)$. Applying Lemma 3.1.3 we conclude that

$$\lim_{t\uparrow R} \frac{(m_P(f(t)) - N)}{\sigma_P(f(t))} = 0\,.$$

We also have, combining the asymptotic expression $(\dagger)$ with the hypothesis $\lim_{t\uparrow R} f(t) = +\infty$, that there exists a constant $C > 0$ such that

$$0 \leq \frac{m_f(t)\sigma_P(f(t))}{\sigma_f(t)} \leq C\frac{m_f(t)}{\sigma_f(t) f(t)^{1/2}}, \qquad \text{as } t \uparrow R,$$

and therefore using $(\star)$ we conclude that

$$\lim_{t\uparrow R} \frac{m_h(t) - N m_f(t)}{\sigma_f(t)} = 0.$$

□

**Lemma 3.1.30.** *Assume that $f \in \mathcal{K}$ is a transcendental entire function of finite order $\rho \geq 0$ such that $\lim_{t\to+\infty} \sigma_f^2(t) = +\infty$, then*

$$\lim_{t\uparrow R} \frac{\sigma_h^2(t)}{\sigma_f^2(t)} = N, \qquad \textit{and} \qquad \lim_{t\uparrow R} \frac{m_h(t) - N m_f(t)}{\sigma_f(t)} = 0.$$

*Proof.* Using that $f$ has finite order $\rho \geq 0$, we have that for any $\varepsilon > 0$, there are constants $C > 0$ and $T > 0$ such that

$$\frac{f'(t)}{f(t)} \leq C t^{\rho - 1 + \varepsilon}, \quad \text{for any } t \geq T\,. \tag{3.1.13}$$

See equation (2.5.2) and also Theorem 2.5.1 for further details.

Combining (3.1.2), that is using that there exist a constant $C > 0$ and $T > 0$ such that $\sigma_P^2(f(t)) \leq C/f(t)$, for any $t \geq T$, with the inequality (3.1.13) we find that (adjusting $T > 0$ if necessary) there are constants $C > 0$ and $M > 0$ such that

$$\sigma_P^2(f(t)) m_f(t)^2 \leq C \frac{t^2 f'(t)^2}{f(t)^3} \leq C \frac{t^M}{f(t)}, \quad \text{for any } t \geq T\,,$$

and therefore

$$(\dagger) \qquad \lim_{t \to +\infty} \sigma_P^2(f(t)) m_f(t)^2 = 0.$$

Using equation (3.1.12) we find that

$$\frac{\sigma_h^2(t)}{\sigma_f^2(t)} = \frac{\sigma_P^2(f(t)) m_f^2(t)}{\sigma_f^2(t)} + m_P(f(t)),$$

then using $(\dagger)$ and the hypothesis $\lim_{t\to\infty} \sigma_f^2(t) = +\infty$ we conclude that

$$\lim_{t \to +\infty} \frac{\sigma_h^2(t)}{\sigma_f^2(t)} = N\,.$$

□

The following lemma expresses the characteristic function of the normalized Khinchin family associated to the composition $h = P \circ f$, in terms of the characteristic function of the Khinchin family associated to $f$.

**Lemma 3.1.31.** *Let $P(z) = \sum_{j=0}^N a_j z^n \in \mathcal{K}$ be a polynomial of degree $\deg(P) = N \geq 1$ and $f \in \mathcal{K}$ a power series with radius of convergence $R > 0$. Denote $(X_t)$, $(Y_t)$ and $(Z_t)$ the Khinchin families of $f, P$ and $h = P \circ f$, respectively, then*

$$\mathbf{E}(e^{i\theta Z_t}) = \sum_{j=0}^N \mathbf{P}(Y_{f(t)} = j) \mathbf{E}(e^{i\theta X_t})^j\,,$$

*for any $t \in (0, R)$ and $\theta \in \mathbb{R}$, and*

$$\mathbf{E}(e^{i\theta \breve{Z}_t}) = \sum_{j=0}^N \mathbf{P}(Y_{f(t)} = j) \mathbf{E}(e^{i\theta X_t/\sigma_h(t)})^j e^{-i\theta m_h(t)/\sigma_h(t)}\,, \tag{3.1.14}$$

*for any $t \in (0, R)$ and $\theta \in \mathbb{R}$.*

*Proof.* For any $t \in [0, R)$ and $\theta \in \mathbb{R}$ we have

$$\mathbf{E}(e^{i\theta Z_t}) = \frac{P(f(te^{i\theta}))}{P(f(t))} = \sum_{j=0}^{N} \frac{a_j f(te^{i\theta})^j}{P(f(t))} = \sum_{j=0}^{N} \frac{a_j f(t)^j}{P(f(t))} \frac{f(te^{i\theta})^j}{f(t)^j}$$

therefore

$$(\dagger) \quad \mathbf{E}(e^{i\theta Z_t}) = \sum_{j=0}^{N} \mathbf{P}(Y_{f(t)} = j)\mathbf{E}(e^{i\theta X_t})^j.$$

This result is a particular case of Lemma 1.3.15. □

We finally state and prove the main result of this subsection.

**Theorem 3.1.32.** *Let $P(z) = \sum_{j=0}^{N} a_j z^j \in \mathcal{K}$ be a polynomial of degree $\deg(P) = N \geq 1$ and $f \in \mathcal{K}$ a Gaussian power series with radius of convergence $R > 0$. Assume that $\lim_{t\uparrow R} f(t) = +\infty$ and also that*

$$(\star) \quad \lim_{t\uparrow R} \frac{m_f^2(t)}{\sigma_f^2(t) f(t)} = 0,$$

*then $h = P \circ f$ is a Gaussian power series.*

See also [48, p. 84], there Hayman proves that for $f \in \mathcal{K}$ a Hayman function the composition of a polynomial $p(z) = b_0 + b_1 z + \ldots b_m z^m$, with real coefficients and $b_m > 0$, with $f$, i.e. the function $h = P \circ f$, is again a Hayman function, see also Section 3.3 for the definition of the Hayman class.

*Proof.* Denote $(X_t)$, $(Y_t)$ and $(Z_t)$ the Khinchin families of $f, P$ and $h = P \circ f$, respectively. Lemma 3.1.31, in particular equation (3.1.14), gives that for any $t \in (0, R)$ and any $\theta \in \mathbb{R}$ we have

$$\mathbf{E}(e^{i\theta \breve{Z}_t}) = \sum_{j=0}^{N} \mathbf{P}\left(Y_{f(t)} = j\right) \mathbf{E}(e^{i\theta X_t/\sigma_h(t)})^j e^{-i\theta \frac{m_h(t)}{\sigma_h(t)}}.$$

For any $j \geq 1$ we have

$$\left| \mathbf{E}(e^{i\theta X_t/\sigma_h(t)})^j e^{-i\theta \frac{m_h(t)}{\sigma_h(t)}} \right| \leq 1,$$

then

$$\left| \mathbf{P}\left(Y_{f(t)} = j\right) \mathbf{E}\left(e^{i\theta X_t/\sigma_h(t)}\right)^j e^{-i\theta \frac{m_h(t)}{\sigma_h(t)}} \right| \leq \mathbf{P}(Y_{f(t)} = j)\,.$$

Given that $\lim_{t\uparrow R} f(t) = +\infty$, we have $\lim_{t\uparrow R} \mathbf{P}(Y_t = j) = 0$, for any $0 \leq j \leq N-1$, see equation (3.1.1) and the discussion in there, and therefore

$$\lim_{t\uparrow R} \mathbf{P}\left(Y_{f(t)} = j\right) \mathbf{E}(e^{i\theta X_t/\sigma_h(t)})^j e^{-i\theta \frac{m_h(t)}{\sigma_h(t)}} = 0, \quad \text{for any integer } 0 \leq j < N.$$

Since $X_t = \sigma_h(t)\breve{X}_t + m_h(t)$, the previous reasoning gives that

$$\lim_{t\uparrow R} \mathbf{E}(e^{i\theta \breve{Z}_t}) = \lim_{t\uparrow R} \mathbf{E}(e^{i\theta X_t/\sigma_h(t)})^N e^{-i\theta \frac{m_h(t)}{\sigma_h(t)}}$$
$$= \lim_{t\uparrow R} \mathbf{E}\left(e^{i\theta \breve{X}_t \frac{\sigma_f(t)}{\sigma_h(t)}}\right)^N e^{-i\theta\left(\frac{m_h(t)-Nm_f(t)}{\sigma_h(t)}\right)}.$$

combining hypothesis $(\star)$ with Lemma 3.1.29 we conclude that

$$\lim_{t\uparrow R} \frac{m_h(t) - Nm_f(t)}{\sigma_h(t)} = 0\,.$$

We finally have, combining Lemma 1.3.4 with the previous discussion, and the Gaussianity of $f$, that

$$\lim_{t\uparrow R} \mathbf{E}(e^{i\theta \breve{Z}_t}) = \lim_{t\uparrow R} \mathbf{E}\left(e^{i\theta \breve{X}_t \frac{\sigma_f(t)}{\sigma_h(t)}}\right)^N = \left(e^{-\theta^2/(2N)}\right)^N = e^{-\theta^2/2}\,,$$

and therefore $h = P \circ f$ is a Gaussian power series. □

## F. Some applications

We apply some of the previous results to some particular cases. We give an alternative proof of the fact that the exponential of any polynomial in $\mathcal{K}$ is a Gaussian power series.

**F.1. Product of a polynomial and a Gaussian power series** Let $f \in \mathcal{K}$ be a Gaussian power series with radius of convergence $R > 0$ such that $\lim_{t\uparrow R} \sigma_f^2(t) = +\infty$ and $P \in \mathcal{K}$ a polynomial, then the product $h = Pf$ is a Gaussian power series.

In this case we have

$$\lim_{t\uparrow R} \frac{\sigma_P(t)}{\sigma_f(t)} = 0\,,$$

recall that for a polynomial $P \in \mathcal{K}$ we have $\lim_{t\to+\infty} \sigma_P^2(t) = 0$, therefore applying Propositions 3.1.25 and 3.1.24 we conclude that $h = Pf$ is a Gaussian power series.

**F.2. Exponential of a polynomial** Let $P$ be a polynomial in $\mathcal{K}$, then $f = e^P$ is a Gaussian power series.

Indeed: let $a, b > 0$ and $p > q > 1$ two integers. Consider the entire functions $f(z) = e^{az^p}$ and $g(z) = e^{bz^q}$, then we have $\sigma_f^2(t) = ap^2t^p$ and $\sigma_g^2(t) = bqt^q$ and this implies that

$$\lim_{t\to+\infty} \frac{\sigma_f^2(t)}{\sigma_f^2(t) + \sigma_g^2(t)} = 1,$$

recall that $p > q$. We know that $f$ and $g$ are Gaussian power series, this follows combining the Gaussinity of the exponential function with Proposition 3.1.23, therefore $fg$ is Gaussian, see Proposition 3.1.25.

In general for a polynomial $P(z) = \sum_{n=0}^{N} a_n z^n$ in $\mathcal{K}$, we have, applying this result, that $e^{P(z)} = e^{a_0} \prod_{n=1}^{N} e^{a_n z^n}$ is Gaussian.

We proved before, see Proposition 3.1.14, that the exponential of any polynomial $P$ such that $e^P \in \mathcal{K}$ is Gaussian. In Proposition 3.1.14, we only require $e^P \in \mathcal{K}$ rather than $P \in \mathcal{K}$, as we do in this case.

**F.3. Composition of an entire function of finite order with a polynomial** Let $f \in \mathcal{K}$ be an entire function of finite order $\rho \geq 0$ such that $\lim_{t\to+\infty} \sigma_f^2(t) = +\infty$, then we have

$$\lim_{t\to+\infty} \frac{m_f^2(t)}{\sigma_f^2(t) f(t)} = 0,$$

see Remark 3.1.30 for further details. Assume that $f$ is Gaussian, then, applying Theorem 3.1.32, we conclude that $P \circ f$ is a Gaussian power series.

## 3.2 Strongly Gaussian power series

Our interest in the Gaussian power series lies, fundamentally, in Hayman's identity (1.3.31). We want an asymptotic formula for the Taylor coefficients $a_n$ of a power series $f \in \mathcal{K}$ with Khinchin family $(X_t)_{t\in[0,R)}$ and $M_f = +\infty$. The identity (1.3.31) asks for analyzing the behavior of the integral

$$(\star) \qquad \int_{-\pi\sigma_f(t_n)}^{\pi\sigma_f(t_n)} \mathbf{E}(e^{i\theta \breve{X}_{t_n}}) d\theta,$$

as $n \to +\infty$. Recall that here $(t_n)_{n\geq 1} \subseteq (0, R)$ is a sequence such that $m_f(t_n) = n$, for any $n \geq 1$ (and therefore $t_n \to R$, as $n \to +\infty$).

For a Gaussian power series we have pointwise convergence of the integrand in $(\star)$ towards $e^{-\theta^2/2}$. Of course, this is not enough to conclude the convergence of the integral $(\star)$; we need some kind of additional bound to conclude the convergence of the integral.

We say that the power series $f \in \mathcal{K}$ with radius of convergence $R > 0$ and its Khinchin family $(X_t)$ are **strongly Gaussian,** or locally Gaussian, if the following two conditions hold:

$$(3.2.1) \qquad \lim_{t\uparrow R} \sigma_f^2(t) = +\infty \quad \text{and} \quad \lim_{t\uparrow R} \int_{-\pi\sigma_f(t)}^{\pi\sigma_f(t)} \left| \mathbf{E}(e^{i\theta \breve{X}_t}) - e^{-\theta^2/2} \right| d\theta = 0\,.$$

The strongly Gaussianity will allow us to pass to the limit in the integral term $(\star)$ and therefore obtain asymptotic formulas for the coefficients $a_n$, see Hayman's identity (1.3.31).

We will see later than the second condition in (3.2.1) does not follow from the Gaussianity of the power series $f$. See example 3.2.7 below. On the other hand we will also see, later, that strongly Gaussian power series are Gaussian.

The notion of strongly Gaussian power series, at least in this context, was introduced by Báez-Duarte, see [6].

Before continuing we test if some of the basic families are strongly Gaussian. Namely: $f(z) = e^z$, $f(z) = 1/(1+z)$ and $f(z) = 1 + z$.

### A. The Poisson family is strongly Gaussian

We already know that the Poisson family is Gaussian, in fact the exponential of any polynomial in $\mathcal{K}$ is Gaussian, see Section 3.1. The bound given by the next lemma will allow us to apply the dominated convergence theorem, and then obtain that the Poisson family is strongly Gaussian.

We collect some bounds by means of the following lemma: we will use these bounds repeatedly below these lines

**Lemma 3.2.1.** *We have*

*1)* $t\cos(\theta) - t \le -\frac{2}{\pi^2} t\theta^2$, *for any* $t \ge 0$ *and* $|\theta| \le \pi$.

*2)* $e^{t\cos(\theta)} - e^t \le -\frac{1}{2\pi^2} e^t\theta^2$, *for any* $t \ge 1$ *and* $|\theta| \le \pi$.

*Proof.* The inequality $\sin(\theta/2) \ge \theta/\pi$, for any $0 \le \theta \le \pi$ combined with identity $1 - \cos(\theta) = 2\sin^2(\theta/2)$ gives that

$$(\star) \quad \cos(\theta) - 1 \le -\frac{2}{\pi^2}\theta^2, \quad \text{for any } |\theta| \le \pi .$$

The bound $(\star)$ gives 1).

For the other inequality, using $(\star)$, and also that $t \ge 1$ and $e^{-y} - 1 \le \frac{-y}{4}$, for any $0 \le y \le 2$, we find that

$$e^{t\cos(\theta)} - e^t \le e^t\left(\exp\left(-2/\pi^2 t\theta^2\right) - 1\right) \le e^t\left(\exp\left(-2/\pi^2\theta^2\right) - 1\right) \le -\frac{1}{2\pi^2}e^t\theta^2 .$$

□

**Lemma 3.2.2.** *For* $f(z) = e^z$, *we have*

$$|\mathbf{E}(e^{i\theta \breve{X}_t})| \le e^{-2\theta^2/\pi^2}, \quad \textit{for any } \theta \in \mathbb{R} \textit{ such that } |\theta| < \pi\sigma_f(t) .$$

*Proof.* The characteristic function of $\breve{X}_t$ is given by

$$|\mathbf{E}(e^{i\theta \breve{X}_t})| = e^{t(\cos(\theta/\sqrt{t})-1)} ,$$

see equation (1.3.11). Recall that $\sigma_f(t) = \sqrt{t}$. The bound follows from Lemma 3.2.1. □

By virtue of Lebesgue dominated convergence theorem, here we use the Gaussianity of $f$, see Section 3.1, and also the bound given by the previous lemma, we conclude that $f(z) = e^z$ is a strongly Gaussian power series.

**B. The Pascal and the Bernoulli family are not strongly Gaussian**

The families given by $f(z) = 1 + z$, the Bernoulli family, and $f(z) = 1/(1 - z)$, the Pascal family, are not strongly Gaussian.

We know that these families are not Gaussian, see Chapter 3, and we will see, later on, that strongly Gaussian power series are Gaussian.

For the Bernoulli family, and in general for any polynomial in the class $\mathcal{K}$, we have that $\lim_{t\to+\infty} \sigma_f^2(t) = 0$, and therefore

$$\int_{-\pi\sigma_f(t)}^{\pi\sigma_f(t)} \left|\mathbf{E}(e^{i\theta \breve{X}_t}) - e^{-\theta^2/2}\right| d\theta \le 4\pi\sigma_f(t)\,,$$

tends to 0, as $t \to +\infty$, and

$$\int_{-\pi\sigma_f(t)}^{\pi\sigma_f(t)} e^{-\theta^2/2} d\theta \le 4\pi\sigma_f(t)\,,$$

also tends to 0, as $t \to +\infty$. In particular any polynomial in $\mathcal{K}$ is not strongly Gaussian, recall that polynomials in $\mathcal{K}$ are characterized by the property $\lim_{t\to+\infty} \sigma_f(t) = 0$.

For the Pascal family, that is, the Khinchin family associated to $f(z) = 1/(1 - z)$, we have $m_f(t) = t/(1-t)$ and $\sigma_f^2(t) = t/(1-t)^2$, both functions escaping to $+\infty$, as $t \uparrow 1$.

If we take $t_n = n/(n + 1)$, then we have $m_f(t_n) = n$, then

$$\begin{aligned}\int_{-\pi\sigma_f(t_n)}^{\pi\sigma_f(t_n)} \mathbf{E}(e^{i\theta \breve{X}_t}) d\theta &= \int_{-\pi\sigma_f(t_n)}^{\pi\sigma_f(t_n)} \frac{1-t_n}{1-t_n e^{i\theta/\sigma_f(t_n)}} e^{-i\theta m_f(t_n)/\sigma_f(t_n)} d\theta\\ &= \sigma_f(t_n) \int_{-\pi}^{\pi} \frac{1-t_n}{1-t_n e^{i\alpha}} e^{-i\alpha m_f(t_n)} d\alpha\\ &= \sqrt{t_n} \int_{-\pi}^{\pi} \frac{1-t_n}{1-t_n e^{i\alpha}} e^{-i\alpha n} d\alpha = \sqrt{t_n} 2\pi t_n^n\,.\end{aligned}$$

Here we use the identity (1.3.31). Recall that $t_n = n/(n+1)$, then $t_n^n = (1 - \frac{1}{n})^n$ and

$$\lim_{n\to+\infty} \int_{-\pi\sigma_f(t_n)}^{\pi\sigma_f(t_n)} \mathbf{E}(e^{i\theta \breve{X}_t}) d\theta = \frac{2\pi}{e}\,.$$

Hayman's identity (1.3.31) gives in this case, that is for $f(z) = 1/(1 - z)$, that

$$1 = a_n = \frac{f(t_n)}{2\pi t_n^n \sigma_f(t_n)} \int_{-\pi\sigma_f(t_n)}^{\pi\sigma_f(t_n)} \mathbf{E}(e^{i\theta \breve{X}_{t_n}}) d\theta, \qquad \text{for any } n \ge 1.$$

### 3.2.1 Hayman's local central limit Theorem

Hayman's local central limit theorem is a central result in this Theory. The proof in here is similar, at least in spirit, to the proof of the local central limit theorem for lattice random variables.

**Theorem 3.2.3.** *Let $f \in \mathcal{K}$ be a strongly Gaussian power series with radius of convergence $R > 0$, then*

$$\lim_{t\uparrow R} \sup_{n\in\mathbb{Z}} \left| \frac{a_n t^n}{f(t)} \sqrt{2\pi}\sigma_f(t) - e^{-(m(t)-n)^2/(2\sigma_f^2(t))} \right| = 0\,.$$

We take the supremum over $\mathbb{Z}$, here $a_n = 0$, for any integer $n < 0$. In other terms, if for each $t \in [0, R)$, we denote

$$\mathcal{V}_t = \left\{ \frac{n - m_f(t)}{\sigma_f(t)} : n \in \mathbb{Z} \right\},$$

then Theorem 3.2.3 claims that

$$\lim_{t\uparrow R} \sup_{x\in\mathcal{V}_t} \left| \mathbf{P}(\breve{X}_t = x)\sigma_f(t) - \frac{1}{\sqrt{2\pi}} e^{-x^2/2} \right| = 0.$$

*Proof of Theorem 3.2.3.* For each $t \in [0, R)$ we denote

$$S_t \triangleq \sup_{n\in\mathbb{Z}} \left| \frac{a_n t^n}{f(t)} \sqrt{2\pi}\sigma_f(t) - e^{-(m(t)-n)^2/(2\sigma_f^2(t))} \right| .$$

For any $n \in \mathbb{Z}$ and $t \in (0, R)$, Hayman's identity (1.3.31) gives that

$$\frac{a_n t^n}{f(t)} 2\pi\sigma_f(t) = \int_{-\pi\sigma_f(t)}^{\pi\sigma_f(t)} \mathbf{E}(e^{i\theta\breve{X}_t}) e^{i\theta(n-m_f(t))/\sigma_f(t)} d\theta$$

and

$$\int_{-\pi\sigma_f(t)}^{\pi\sigma_f(t)} \mathbf{E}(e^{i\theta\breve{X}_t}) e^{i\theta(n-m_f(t))/\sigma_f(t)} d\theta = \int_{-\pi\sigma_f(t)}^{\pi\sigma_f(t)} \left( \mathbf{E}(e^{i\theta\breve{X}_t}) - e^{-\theta^2/2} \right) e^{i\theta(n-m_f(t))/\sigma_f(t)} d\theta + \int_{-\pi\sigma_f(t)}^{\pi\sigma_f(t)} e^{-\theta^2/2} e^{i\theta(n-m_f(t))/\sigma_f(t)} d\theta\,.$$

Because $f$ is strongly Gaussian, the first integral at the right-hand side of the previous expression is $o(1)$ uniformly in $n$. For the second integral we have

$$\int_{-\pi\sigma_f(t)}^{\pi\sigma_f(t)} e^{-\theta^2/2} e^{i\theta(n-m_f(t))/\sigma_f(t)} d\theta = \int_{\mathbb{R}} e^{-\theta^2/2} e^{i\theta(n-m_f(t))/\sigma_f(t)} d\theta + o(1) = \sqrt{2\pi} e^{-(n-m_f(t))^2/(2\sigma_f^2(t))} + o(1).$$

Again the $o(1)$ is uniform in $n$. □

As a direct corollary of Theorem 3.2.3 we obtain that strongly Gaussian power series have $M_f = +\infty$.

**Corollary 3.2.4.** *Let $f \in \mathcal{K}$ be a strongly Gaussian power series, then $M_f = +\infty$.*

*Proof.* Restricting the supremum to integers $n \leq -1$ in Theorem 3.2.3, we have

$$\lim_{t\uparrow R} \sup_{n\leq -1} \left| e^{-(m_f(t)-n)^2/(2\sigma_f^2(t))} \right| = 0.$$

For each $t \in [0, R)$ the supremum is reached by taking $n = -1$, recall that $m_f(t) \geq 0$, then we have

$$\lim_{t\uparrow R} \exp\left( -\frac{(m_f(t)+1)^2}{\sigma_f^2(t)} \right) = 0\,,$$

therefore

$$(\dagger) \qquad \lim_{t\uparrow R} \frac{(m_f(t)+1)^2}{\sigma_f^2(t)} = +\infty.$$

Strongly Gaussian power series verify that $\lim_{t\uparrow R} \sigma_f^2(t) = +\infty$, then $(\dagger)$ gives that $M_f = \lim_{t\uparrow R} m_f(t) = +\infty\,.$ □

As a corollary of the previous result we have the following Theorem. Hayman enunciate the following Theorem for power series in the Hayman class, see Section 3.3, and not for strongly Gaussian power series as we do here, see [48].

**Corollary 3.2.5** (Hayman's asymptotic formulas)**.** *Let $f \in \mathcal{K}$ be a strongly Gaussian power series, then*

$$a_n \sim \frac{1}{\sqrt{2\pi}} \frac{f(t_n)}{\sigma_f(t_n) t_n^n}, \qquad \textit{as } n \to +\infty, \tag{3.2.2}$$

*here the sequence $(t_n)_{n\geq 0} \subseteq (0, R)$ is given by $m_f(t_n) = n\,.$*

Observe that strongly Gaussian power series have the property that the coefficients $a_n$ are positive after some value of $n$.

*Proof.* Let $f \in \mathcal{K}$ a strongly Gaussian power series, then Corollary 3.2.4 gives that $M_f = +\infty$. Let $(t_n)_{n\geq 0} \subseteq (0, R)$ be a sequence such that $m_f(t_n) = n$.

Denote

$$S_t \triangleq \sup_{n\in\mathbb{Z}} \left| \frac{a_n t^n}{f(t)} \sqrt{2\pi}\sigma_f(t) - e^{-(m(t)-n)^2/(2\sigma_f^2(t))} \right|.$$

Theorem 3.2.3 gives that $S_t \to 0$, as $t \uparrow R$. By definition we have

$$\left| \frac{a_n t_n^n}{f(t_n)} \sqrt{2\pi}\sigma_f(t_n) - 1 \right| \leq S_{t_n}, \qquad \text{for each } n \geq 1,$$

then

$$\lim_{n\to+\infty}\left|\frac{a_n t_n^n}{f(t_n)}\sqrt{2\pi}\sigma_f(t_n)-1\right|=0\,.$$

Alternatively, appealing to Hayman's identity (1.3.31) we have

$$a_n=\frac{f(t_n)}{2\pi\sigma_f(t_n)t_n^n}\int_{-\pi\sigma_f(t_n)}^{\pi\sigma_f(t_n)}\mathbf{E}(e^{i\theta\breve{X}_{t_n}})d\theta,\qquad \text{for each } n\geq 1,$$

and using that $f$ is strongly Gaussian we find that

$$\lim_{n\to\infty}\int_{-\pi\sigma_f(t_n)}^{\pi\sigma_f(t_n)}\mathbf{E}(e^{i\theta\breve{X}_{t_n}})d\theta=\int_R e^{-\theta^2/2}d\theta=\sqrt{2\pi}\,.$$

□

In the conditions of the previous corollary, and with the same argument, we obtain that if $\omega_n$ is a good approximation of $t_n$, meaning that

$$\lim_{n\to+\infty}\frac{m_f(\omega_n)-n}{\sigma_f(\omega_n)}=0\,,$$

then we obtain that

$$a_n\sim\frac{1}{\sqrt{2\pi}}\frac{f(\omega_n)}{\sigma_f(\omega_n)\omega_n^n},\qquad \text{as } n\to\infty.$$

We will come back to this variant later on.

Let's see some applications of Corollary 3.2.5. The following example give title to Hayman's paper: *A generalization of Stirling's formula*, see [48].

**Example 3.2.6** (Stirling's formula)**.** For $f(z)=e^z=\sum_{n=0}^{\infty}z^n/n!$ we have $R=+\infty$, $m_f(t)=t$ and $\sigma_f(t)=\sqrt{t}$. We have $M_f=+\infty$ and therefore $t_n=n$.

We already know that $f(z)=e^z$ is strongly Gaussian. Applying Hayman's asymptotic formula we find that

$$\frac{1}{n!}\sim\frac{1}{\sqrt{2\pi}}\frac{e^n}{\sqrt{n}n^n},\qquad \text{as } n\to\infty\,,$$

which is Stirling's formula. ⊡

**Example 3.2.7** (Power series in $\mathcal{K}$ which are Gaussian but not strongly Gaussian)**.** Consider the entire function $f(z)=e^{z^2}$, which is a power series in $\mathcal{K}$, that is, $R=+\infty$. The mean and the variance are given by $m_f(t)=t^2$ and $\sigma_f^2(t)=4t^2$. Proposition 3.1.23 gives that $f$ is Gaussian.

If we assume that $f$ is strongly Gaussian, then Hayman's formula will hold, but this is not possible: the coefficients of even index are all identically zero. ⊡

Let $f \in \mathcal{K}$ be a strongly Gaussian power series with radius of convergence $R > 0$. In order to apply Hayman's formula we need to solve the equation $m_f(t) = n$ (recall that $m_f$ is an increasing diffeomorphism) and then evaluate in that value $t_n$ both $f(t_n)$ and $\sigma_f(t_n)$.

If the function $f \in \mathcal{K}$ is intricate, the expression for $m_f(t)$ will also be intricate and solving $m_f(t) = n$ will be difficult. Also, any estimation of $t_n$ would also require to estimate the values of $f(t_n)$ and $\sigma_f(t_n)$, respectively.

The following theorem will be very relevant in future applications; it allows to circumvent the problem of finding the sequence $(t_n)_{n\geq 1}$ such that $m_f(t_n) = n$, which is some cases is almost impossible, see for instance equation 1.2.6.

**Theorem 3.2.8** (Baéz-Duarte asympotic formula)**.** *Let $f(z) = \sum_{n=0}^{\infty} a_n z^n \in \mathcal{K}$ a strongly Gaussian power series. Assume that $\tilde{m}_f(t)$ and $\tilde{\sigma}_f(t)$ are functions, defined on $(0, R)$, such that*

- $\tilde{m}_f(t)$ *grows to $+\infty$ continuously and monotonically, as $t \uparrow R$,*
- $\tilde{\sigma}_f(t) \sim \sigma_f(t)$, *as $t \uparrow R$,*
- $\displaystyle\lim_{t\uparrow R} \frac{\tilde{m}_f(t) - m_f(t)}{\sigma_f(t)} = 0,$

*then*

$$a_n \sim \frac{f(\tau_n)}{\sqrt{2\pi}\tilde{\sigma}_f(\tau_n)\tau_n^n}, \qquad \text{as } n \to \infty, \tag{3.2.3}$$

*where $\tau_n$ is given by $\tilde{m}_f(\tau_n) = n$, for each $n \geq 0$.*

Under the hypothesis of the previous lemma we have $\tilde{m}_f(t) \sim m_f(t)$, as $t \uparrow R$. Indeed: Corollary 3.2.10, below, gives that $\lim_{t\uparrow R} m_f(t)/\sigma_f(t) = +\infty$, then for any $M > 0$, there exists $T > 0$ such that

$$\frac{\tilde{m}_f(t) - m_f(t)}{\sigma_f(t)} = \left(\frac{\tilde{m}_f(t)}{m_f(t)} - 1\right)\frac{m_f(t)}{\sigma_f(t)} \geq M\left(\frac{\tilde{m}_f(t)}{m_f(t)} - 1\right), \qquad \text{for any } t \geq T.$$

and therefore $\tilde{m}_f(t) \sim m_f(t)$, as $t \uparrow R$.

*Proof of Theorem 3.2.8.* Using that $\tilde{m}$ is monotonic we conclude that $\tau_n \to \infty$, as $n \to \infty$. By definition of the sequence $(\tau_n)_{n\geq 0}$ we have that $\tilde{m}_f(\tau_n) = n$ and therefore

$$\lim_{n\to\infty} \frac{n - m_f(\tau_n)}{\sigma_f(\tau_n)} = 0.$$

From Theorem 3.2.3 we deduce that

$$a_n \sim \frac{f(\tau_n)}{\sqrt{2\pi}\sigma_f(\tau_n)\tau_n^n}, \qquad \text{as } n \to \infty,$$

and finally, given that $\sigma_f(t) \sim \tilde{\sigma}_f(t)$, as $t \uparrow R$, then we conclude

$$a_n \sim \frac{f(\tau_n)}{\sqrt{2\pi}\tilde{\sigma}_f(\tau_n)\tau_n^n}, \qquad \text{as } n \to \infty,$$

□

### 3.2.2 Strongly Gaussian power series are Gaussian

We prove here that strongly Gaussian power series are Gaussian.

**Theorem 3.2.9.** *Let $f \in \mathcal{K}$ be a strongly Gaussian power series with radius of convergence $R > 0$, then $f$ is Gaussian.*

*Proof.* Fix $b \in \mathbb{R}$ and $a < b$.

Denote

$$S_t \triangleq \sup_{n\in\mathbb{Z}} \left| \frac{a_n t^n}{f(t)} \sqrt{2\pi}\sigma_f(t) - e^{-(m(t)-n)^2/(2\sigma_f^2(t))} \right| , \qquad \text{for any } t \in [0, R) .$$

Theorem 3.2.3 gives that $\lim_{t\uparrow R} S_t = 0$.

Denote

$$\mathcal{D}_t = \{n \geq 0 : m_f(t) + a\sigma_f(t) < n < m_f(t) + b\sigma_f(t)\} ,$$

and

$$\mathcal{E}_t = \left\{ \frac{n - m_f(t)}{\sigma_f(t)} : n \in \mathcal{D}_t \right\} .$$

Observe that $\#\mathcal{E}_t \leq \sigma_f(t)|b - a| + 1$, then

$$\left| \sum_{n\in\mathcal{D}_t} \frac{a_n t^n}{f(t)} - \frac{1}{\sqrt{2\pi}\sigma_f(t)} \sum_{x\in\mathcal{E}_t} e^{-x^2/2} \right| \leq \frac{1}{\sqrt{2\pi}\sigma_f(t)} (\#\mathcal{E}_t) S_t, \qquad \text{for any } t \geq 0 .$$

Theorem 3.2.3 gives that $\lim_{t\uparrow R} S_t = 0$ and we also have that $\sigma_f(t)$ is bounded from below for large values of $t$, in fact $\lim_{t\uparrow R} \sigma_f(t) = +\infty$. The upper bound in the previous inequality tends to 0, as $t \uparrow R$.

Moreover

$$\sum_{n\in\mathcal{D}_t} \frac{a_n t^n}{f(t)} = \mathbf{P}(a < \breve{X}_t < b) .$$

Finally using that $\lim_{t\uparrow R} \sigma_f^2(t) = +\infty$, we have

$$\lim_{t\uparrow R} \frac{1}{\sigma_f(t)} \sum_{x\in\mathcal{E}_t} e^{-x^2/2} = \int_a^b e^{-x^2/2} dx,$$

therefore combining all the steps in previous analysis we find that

$$\lim_{t\uparrow R} \mathbf{P}(a < \breve{X}_t < b) = \frac{1}{\sqrt{2\pi}} \int_a^b e^{-x^2/2} dx .$$

For negative $a < 0$ we have, appealing to Chebyshev's inequality, that

$$\left|\mathbf{P}(\breve{X}_t \leq b) - \mathbf{P}(a < \breve{X}_t \leq b)\right| = \mathbf{P}(\breve{X}_t \leq a) \leq \mathbf{P}(|\breve{X}_t| \geq -a) \leq \frac{1}{a^2}\,.$$

From here we get that

$$\limsup_{t\uparrow R} \mathbf{P}(\breve{X}_t \leq b) \leq \frac{1}{\sqrt{2\pi}} \int_a^b e^{-x^2/2} dx + \frac{1}{a^2}$$

and therefore making $a \downarrow -\infty$, we find that

$$\limsup_{t\uparrow R} \mathbf{P}(\breve{X}_t \leq b) \leq \frac{1}{\sqrt{2\pi}} \int_{-\infty}^b e^{-x^2/2} dx\,.$$

By using an analogous argument we also find that

$$\limsup_{t\uparrow R} \mathbf{P}(\breve{X}_t \leq b) \geq \int_{-\infty}^b e^{-x^2/2} dx.$$

We have the same inequalities for the $\liminf$. □

### 3.2.3 Strongly Gaussian power series are Clans

As a consequence of Theorem 3.2.3 we also obtain that every strongly Gaussian function conforms a clan.

**Corollary 3.2.10.** *Let $f \in \mathcal{K}$ be strongly Gaussian, then $f$ conforms a clan.*

*Proof.* Restricting the supremum in the proof of Theorem 3.2.3 to $n = -1$, we get

$$\lim_{t\uparrow R} \exp\left(-\frac{(m_f(t)+1)^2}{2\sigma_f(t)^2}\right) = 0\,.$$

Since $\lim_{t\uparrow R} \sigma_f^2(t) = +\infty$, because $f$ is strongly Gaussian, we get that $\sigma_f(t) = o(m_f(t))$, as $t \uparrow R$. □

This corollary shows that a large collection of functions ranging from $f(z) = e^z$ to the generating functions of partitions $\prod_{j\geq 1}(1-z^j)^{-1}$, $|z| < 1$ are clans (see [17] and [18] for this and other interesting examples).

**Corollary 3.2.11.** *Let $f \in \mathcal{K}$ be a strongly Gaussian power series, with radius of convergence $R > 0$, then, for each $n \geq 0$, we have*

$$\lim_{t\uparrow R} \frac{f(t)}{t^n} = +\infty\,.$$

*In fact we have*

$$\lim_{t\uparrow R} \frac{f(t)}{\sigma_f(t) t^n} = +\infty\,.$$

*Proof.* Polynomials are not strongly Gaussian, recall that polynomials are characterized, among the functions in $\mathcal{K}$, by the property $\lim_{t\to+\infty}\sigma_f^2(t)=0$.

Fix $n\geq 0$ such that $a_n\neq 0$. Theorem 3.2.3 gives that

$$\lim_{t\uparrow R}\left|\frac{a_n t^n}{f(t)}\sqrt{2\pi}\sigma_f(t)-e^{-(m(t)-n)^2/(2\sigma_f^2(t))}\right|=0\,.$$

Corollary 3.2.10 gives that

$$\lim_{t\uparrow R}\frac{\sigma_f(t)}{m_f(t)}=0$$

therefore

$$\lim_{t\uparrow R}\frac{a_n t^n}{f(t)}\sqrt{2\pi}\sigma_f(t)=0,$$

using that $a_n\neq 0$ we conclude that

$$\lim_{t\uparrow R}\frac{t^n\sigma_f(t)}{f(t)}=0\,.$$

Because $f$ is strongly Gaussian there are an infinite amount of indices $n$ such that $a_n\neq 0$, and therefore

$$\lim_{t\uparrow R}\frac{t^n\sigma_f(t)}{f(t)}=0\,,\qquad\text{for any } n\geq 0.$$

Now using that $\lim_{t\uparrow R}\sigma_f(t)=+\infty$, we conclude that

$$\lim_{t\uparrow R}\frac{t^n}{f(t)}=0\,,\qquad\text{for any } n\geq 0.$$

□

There are power series in $\mathcal{K}$ that conform a clan but which are not strongly Gaussian. For instance, if $g\in\mathcal{K}$ is strongly Gaussian then $f(z)=g(z^N)$, with $N\in\mathbb{N}$, conforms a clan but it is not strongly Gaussian since $a_n\equiv 0$ when $n$ is not a multiple of $N$ and hence (3.2.2) does not hold.

### 3.2.4 Operations with strongly Gaussian power series

Here we prove that some operations with strongly Gaussian power series preserve the property of being strongly Gaussian.

**A. Subordination**

Let $Q \geq 1$ an integer. Denote $\mathcal{K}_Q$ the subfamily of power series in $g(z) = \sum_{n=0}^{\infty} a_n z^n \in \mathcal{K}$ such that $a_n = 0$ except for those index $n \geq 0$ which are multiples of $Q \geq 1$.

For such power series $g \in \mathcal{K}_Q$ we have that $Q$ divides $\gcd\{n \geq 1 : a_n \neq 0\}$, and reciprocally, if $Q$ divides $\gcd\{n \geq 1 : a_n \neq 0\}$, then $g \in \mathcal{K}_Q$. For any $g \in \mathcal{K}_Q$ there exists a power series $f(z) = \sum_{n=0}^{\infty} b_n z^n \in \mathcal{K}$ such that $g(z) = f(z^Q)$, with $\gcd\{n \geq 1 : b_n \neq 0\} = 1$.

Fix $g \in \mathcal{K}_Q$, then there exists $f \in \mathcal{K}$ such that $g(z) = f(z^Q)$. If $g$ has radius of convergence $R > 0$, then $f$ has radius of convergence $R^Q > 0$.

Denote $(X_t)_{t\in[0,R)}$ the Khinchin family of $g$ and $(Y_t)_{t\in[0,R^Q)}$ the Khinchin family of $f$. Recall that in this case

$$(\star) \quad X_t \stackrel{d}{=} Q \cdot Y_{t^Q}, \qquad \text{for any } t \in [0, R),$$

and therefore $m_g(t) = Q \cdot m_f(t^Q)$ and $\sigma_g(t) = Q \cdot \sigma_f(t^Q)$, for any $t \in [0, R)$.

We also have, using $(\star)$, that

$$\mathbf{E}(e^{i\theta \breve{X}_t}) = \mathbf{E}(e^{i\theta \breve{Y}_{t^Q}}), \qquad \text{for any } \theta \in \mathbb{R} \text{ and } t \in (0, R).$$

From the previous discussion: we say that a power series $g \in \mathcal{K}_Q$ is $Q$-strongly Gaussian if

$$\lim_{t \uparrow R} \sigma_g(t) = +\infty, \qquad \text{and} \qquad \lim_{t \uparrow R} \int_{-\pi\sigma_g(t)/Q}^{\pi\sigma_g(t)/Q} |\mathbf{E}(e^{i\theta \breve{Y}_{t^Q}}) - e^{-\theta^2/2}| d\theta = 0.$$

We have that $f$ is strongly Gaussian if and only if $g$ is $Q$-strongly Gaussian.

If $g$ is $Q$-strongly Gaussian, or equivalently, if $f$ is strongly Gaussian, the asymptotic formula for the coefficients $b_n$ of $f$, turns into an asymptotic formula for the coefficients $a_{Qn}$ of $g$. This asymptotic formula, written in terms of $g$, is

$$b_{nQ} \sim Q \frac{1}{\sqrt{2\pi}} \frac{g(s_{Qn})}{s_{Qn}^{nQ} \sigma_g(s_{Qn})}, \qquad \text{as } n \to \infty\,.$$

Here the sequence $(t_n)_{n\geq 0}$ is given by $m_f(t_n) = n$ and $s_{nQ} = t_n^{1/Q}$.

**B. Products of strongly Gaussian power series**

We have seen before, see Section 3.1, that under some simple conditions the product of two Gaussian power series in $\mathcal{K}$ is a Gaussian power series.

Now we study the same question, but this time for strongly Gaussian power series.

**B.1 Product of power series with the same radius of convergence** Assume that $f, g \in \mathcal{K}$ are power series with radius of convergence $R > 0$. The product $h = fg$ has radius of convergence $R > 0$, and is, of course, in $\mathcal{K}$.

Denote $(X_t)$, $(Y_t)$ and $(Z_t)$ the respective Khinchin families of $f, g$ and $h = fg$. Recall that

$$(\star) \quad Z_t \stackrel{d}{=} X_t \oplus Y_t\,, \qquad \text{for any } t \in [0, R)\,,$$

here $\oplus$ denotes the sum of independent copies. In particular, the previous inequality gives that $\sigma_h^2(t) = \sigma_f^2(t) + \sigma_g^2(t)$ and $m_h(t) = m_f(t) + m_g(t)$, for each $t \in [0, R)$.

Assume that $f$ and $g$ are strongly Gaussian and also that

$$\liminf_{t \uparrow R} \frac{\sigma_f^2(t)}{\sigma_f^2(t) + \sigma_g^2(t)} > 0, \qquad \text{and} \qquad \liminf_{t \uparrow R} \frac{\sigma_g^2(t)}{\sigma_f^2(t) + \sigma_g^2(t)} > 0\,.$$

These conditions are always verified if $f \equiv g$.

We always have $\lim_{t \uparrow R} \sigma_h^2(t) = +\infty$, this follows from the fact that $f$ and $g$ are strongly Gaussian power series.

For each $t \in [0, R)$ we write $\lambda(t) = \sigma_f(t)/\sigma_h(t)$ and $\mu(t) = \sigma_g(t)/\sigma_h(t)$, then we have

$$\liminf_{t \uparrow R} \lambda(t) > 0, \qquad \text{and} \qquad \liminf_{t \uparrow R} \mu(t) > 0\,.$$

Observe that $\lambda^2(t) + \mu^2(t) \equiv 1$, for any $t \in [0, R)$.

Using $(\star)$ we find that

$$\mathbf{E}(e^{i\theta \breve{Z}_t}) = \mathbf{E}(e^{i\theta\lambda(t)\breve{X}_t})\mathbf{E}(e^{i\theta\mu(t)\breve{Y}_t}), \qquad \text{for any } \theta \in \mathbb{R} \text{ and } t \in (0, R).$$

Now using the equality

$$\begin{aligned} \Big|\mathbf{E}(e^{i\theta\lambda(t)\breve{X}_t})\mathbf{E}(e^{i\theta\mu(t)\breve{Y}_t}) - e^{-\theta^2/2}\Big| = \Big|&\mathbf{E}(e^{i\theta\lambda(t)\breve{X}_t}) \left(\mathbf{E}(e^{i\theta\mu(t)\breve{Y}_t}) - e^{-\mu(t)^2\theta^2/2}\right) \\ &+ e^{-\mu(t)^2\theta^2/2} \left(\mathbf{E}(e^{i\theta\lambda(t)\breve{X}_t}) - e^{-\lambda(t)^2\theta^2/2}\right)\Big|, \end{aligned} \tag{3.2.4}$$

we conclude that

$$\Big|\mathbf{E}(e^{i\theta\lambda(t)\breve{X}_t})\mathbf{E}(e^{i\theta\mu(t)\breve{Y}_t}) - e^{-\theta^2/2}\Big| \;\le \Big|\mathbf{E}(e^{i\theta\lambda(t)\breve{X}_t}) - e^{-\lambda^2(t)\theta^2/2}\Big| + \Big|\mathbf{E}(e^{i\theta\mu(t)\breve{Y}_t}) - e^{-\mu^2(t)\theta^2/2}\Big|\,.$$

for any $\theta \in \mathbb{R}$ and any $t \in (0, R)$.

Using the previous inequality we find that

$$\begin{aligned} \int_{|\theta|<\pi\sigma_h(t)} \Big|\mathbf{E}(e^{i\theta\breve{Z}_t}) - e^{-\theta^2/2}\Big|\, d\theta &\le \int_{|\theta|\le\pi\sigma_h(t)} \Big|\mathbf{E}(e^{i\theta\lambda(t)\breve{X}_t}) - e^{-\lambda^2(t)\theta^2/2}\Big|\, d\theta \\ &+ \int_{|\theta|\le\pi\sigma_h(t)} \Big|\mathbf{E}(e^{i\theta\mu(t)\breve{Y}_t}) - e^{-\mu^2(t)\theta^2/2}\Big|\, d\theta. \\ &= \frac{1}{\lambda(t)} \int_{|\theta|\le\pi\sigma_f(t)} \Big|\mathbf{E}(e^{i\theta\breve{X}_t}) - e^{\theta^2/2}\Big|\, d\theta \\ &+ \frac{1}{\mu(t)} \int_{|\theta|\le\pi\sigma_g(t)} \Big|\mathbf{E}(e^{i\theta\breve{Y}_t}) - e^{\theta^2/2}\Big|\, d\theta \end{aligned}$$

We collect the previous discussion by means of the following theorem.

**Theorem 3.2.12.** *Let $f, g \in \mathcal{K}$ be power series with radius of convergence $R > 0$. Assume that $f, g$ are strongly Gaussian and also that*

$$\liminf_{t\uparrow R} \frac{\sigma_f^2(t)}{\sigma_f^2(t)+\sigma_g^2(t)} > 0, \qquad \textit{and} \qquad \liminf_{t\uparrow R} \frac{\sigma_g^2(t)}{\sigma_f^2(t)+\sigma_g^2(t)} > 0\,.$$

*then $h = fg$ is strongly Gaussian.*

**B.2 Product of power series with different radius of convergence** Now we study the case where $f \in \mathcal{K}$ has radius of convergence $R > 0$ and $g \in \mathcal{K}$ has radius of convergence $S > R$. Assume that $f$ is strongly Gaussian.

Denote $(X_t)$, $(Y_t)$ and $(Z_t)$ the respective Khinchin families of $f, g$ and $h = fg$. Given that $f$ is strongly Gaussian, and using that $\lim_{t\uparrow R} \sigma_f(t) = +\infty$ combined with the fact that the radius of convergence of $g$ is $S > R$ we have

$$\lim_{t\uparrow R} \lambda(t) = \lim_{t\uparrow R} \frac{\sigma_f(t)}{\sigma_h(t)} = 1$$

and therefore

$$\lim_{t\uparrow R} \mu(t) = \lim_{t\uparrow R} \frac{\sigma_g(t)}{\sigma_h(t)} = 0.$$

Now we write

$$\begin{aligned}\mathbf{E}(e^{i\theta \breve{Z}_t}) - e^{-\theta^2/2} &= \Big(\mathbf{E}(e^{i\theta\lambda(t)\breve{X}_t})\mathbf{E}(e^{i\theta\mu(t)\breve{Y}_t}) - e^{-\theta^2/2}\Big)\\ &= \mathbf{E}(e^{i\theta\mu(t)\breve{Y}_t})\Big(\mathbf{E}(e^{\theta\lambda(t)\breve{X}_t}) - e^{-\theta^2/2}\Big) + e^{-\theta^2/2}\Big(\mathbf{E}(e^{\theta\mu(t)\breve{Y}_t}) - 1\Big)\,,\end{aligned}$$

therefore

$$\begin{aligned}\int_{|\theta|\le\pi\sigma_h(t)} \Big|\mathbf{E}(e^{i\theta\breve{Z}_t}) - e^{-\theta^2/2}\Big|\,d\theta &\le \int_{|\theta|\le\pi\sigma_h(t)} e^{-\theta^2/2}\Big|\mathbf{E}(e^{i\theta\mu(t)\breve{Y}_t}) - 1\Big|\,d\theta\\ &+ \int_{|\theta|\le\pi\sigma_h(t)} \Big|e^{-\theta^2/2} - e^{-\lambda(t)^2\theta^2/2}\Big|\,d\theta\\ &+ \frac{1}{\lambda(t)}\int_{|\theta|\le\pi\sigma_f(t)} \Big|\mathbf{E}(e^{i\theta\breve{X}_t}) - e^{\theta^2/2}\Big|\,d\theta\end{aligned}$$

Lebesgue's dominated convergence Theorem gives that the first and second term at the right-hand side of the previous inequality tends to 0, as $t \uparrow R$. Here we use the following inequality: because $\lambda(t) \to 1$, as $t \uparrow R$, we have that there exists certain $0 < T < R$ such that $\lambda(t) > 1/2$, for $T < t < R$, then

$$|e^{-\theta^2/2} - e^{-\lambda(t)^2\theta^2/2}| \leq e^{-\theta^2/2} + e^{-\theta^2/8} \leq 2e^{-\theta^2/8}, \quad \text{for any } \theta \in \mathbb{R}.$$

Combining that $f$ is strongly Gaussian with the hypothesis $\lim_{t\uparrow R} \lambda(t) = 1$ we also conclude that the third term tends to 0, as $t \uparrow R$.

This argument is also valid for the case $S = R$, that is, we have the same result for power series $f, g \in \mathcal{K}$ with the same radius of convergence.

We collect the previous discussion by means of the following theorem.

**Theorem 3.2.13.** *Let $f, g \in \mathcal{K}$ be power series with radius of convergence $S \geq R > 0$, respectively. Assume that $f, g$ are strongly Gaussian and also that*

$$\lim_{t\uparrow R} \frac{\sigma_f^2(t)}{\sigma_f^2(t) + \sigma_g^2(t)} = 1, \quad \textit{and} \quad \lim_{t\uparrow R} \frac{\sigma_g^2(t)}{\sigma_f^2(t) + \sigma_g^2(t)} = 0\,.$$

*then $h = fg$ is strongly Gaussian.*

### C. Powers of strongly Gaussian power series

Assume that $f \in \mathcal{K}$ is a strongly Gaussian power series, with radius of convergence $R > 0$, and denote $h = f^N$, for integer $N \geq 1$, then $h = f^N$ is a strongly Gaussian power series.

In this case we have

$$m_h(t) = N \cdot m_f(t) \quad \text{and} \quad \sigma_h(t) = N \cdot \sigma_f(t), \quad \text{for any } t \in [0, R),$$

therefore

$$\lim_{t\uparrow R} \lambda(t) = \lim_{t\uparrow R} \frac{\sigma_f(t)}{\sigma_h(t)} = \frac{1}{N} > 0 \quad \text{and} \quad \lim_{t\uparrow R} \mu(t) = \lim_{t\uparrow R} \frac{\sigma_{f^{N-1}}(t)}{\sigma_h(t)} = \frac{N-1}{N} > 0,$$

then, by using an induction argument, which simply applies the discussion above for power series with the same radius of convergence, we conclude that $h = f^N$ is strongly Gaussian.

If we denote $(Z_t)$ and $(X_t)$ the Khinchin families of $f^N$ and $f$ respectively, then we have

$$Z_t \stackrel{d}{=} X_t^{(1)} \oplus X_t^{(2)} \oplus \cdots \oplus X_t^{(N)}, \quad \text{for any } t \in [0, R),$$

here $X_t^{(1)}, \ldots, X_t^{(N)}$ are i.i.d copies of $X_t$. This equality also gives that $\sigma_{f^N}(t) = \sqrt{N} \cdot \sigma_f(t)$ and $m_{f^N}(t) = N \cdot m_f(t)$. We also have that

$$\breve{Z}_t \stackrel{d}{=} \frac{1}{\sqrt{N}} \left( \breve{X}_t^{(1)} \oplus \cdots \oplus \breve{X}_t^{(N)} \right), \quad \text{for any } t \in (0, R)\,,$$

therefore

$$\mathbf{E}(e^{i\theta \breve{Z}_t}) = \mathbf{E}\left(e^{i\theta \breve{X}_t/\sqrt{N}}\right)^N, \quad \text{for any } \theta \in \mathbb{R} \text{ and } t \in (0, R).$$

The previous discussion gives that

$$\text{(3.2.5)}\quad \lim_{t\uparrow R}\int_{|\theta|\le\pi\sigma_{f^N}(t)}\left|\mathbf{E}(e^{i\theta\breve{Z}_t})-e^{-\theta^2/2}\right|d\theta=\lim_{t\uparrow R}\int_{|\theta|\le\pi\sigma_f(t)\sqrt{N}}\left|\mathbf{E}(e^{i\theta\breve{X}_t/\sqrt{N}})^N-e^{-\theta^2/2}\right|d\theta=0.$$

This observation will be useful later when proving that the composition of a strongly Gaussian power series with a polynomial is a strongly Gaussian power series.

## D. Shifts of strongly Gaussian power series

Let $f\in\mathcal{K}$ be a strongly Gaussian power series with radius of convergence $R>0$. Denote $C>0$ a constant and $M\ \geq 0$ a non-negative integer. We define $h(z)=Cz^Mf(z)$ which is a power series in $\mathcal{K}_s$.

Denote $(X_t)$ and $(Y_t)$ the Khinchin families of $h$ and $f$ respectively, then we have

$$X_t\stackrel{d}{=}M+Y_t,\qquad \text{for any } t\in[0,R)$$

this implies that $m_h(t)=M+m_f(t)$ and $\sigma_h^2(t)=\sigma_f^2(t)$, for any $t\in[0,R)$, and also that

$$\mathbf{E}(e^{i\theta\breve{X}_t})=\mathbf{E}(e^{i\theta\breve{Y}_t}),\qquad \text{for any } \theta\in\mathbb{R} \text{ and } t\in(0,R),$$

therefore $h$ is a strongly Gaussian power series.

## E. Composition of a strongly Gaussian power series with a polynomial

Here we prove that the composition $h=P\circ f$ of a polynomial of degree $N\geq 1$ which is in $\mathcal{K}$ with a strongly Gaussian power series $f\in\mathcal{K}$ is a strongly Gaussian power series.

See also [48, p. 84], there Hayman proves that the composition of a polynomial of degree $m$, with real coefficients and positive $m$-th coefficient, with a Hayman function is again a Hayman function, see also Section 3.3 for the definition of the Hayman class.

**Theorem 3.2.14.** *Let $f\in\mathcal{K}$ be a power series with radius of convergence $R>0$ and $P\in\mathcal{K}$ a polynomial of degree $N\geq 1$. Assume that $f$ is strongly Gaussian and also that*

$$(\star)\quad \lim_{t\uparrow R}\frac{m_f^2(t)}{\sigma_f^2(t)f(t)}=0,$$

*then $h=P\circ f$ is strongly Gaussian.*

Theorem 3.2.14 is stronger than Theorem 3.1.32. Recall that strongly Gaussian power series are Gaussian.

*Proof.* For $f\in\mathcal{K}$ be a strongly Gaussian power series with radius of convergence $R>0$ we always have $\lim_{t\uparrow R}f(t)=+\infty$, see Corollary 3.2.11.

Denote $(Z_t),(X_t)$ and $(Y_t)$ the Khinchin families of $h,f$ and $P$ respectively. Recall that condition $(\star)$ implies that $\lim_{t\uparrow R}\sigma_h^2(t)/\sigma_f^2(t)=N$, see Lemma 3.1.29.

The mean and variance of $h$ are given by

$$m_h(t) = m_P(f(t))m_f(t) \quad \text{for any } t \in [0, R),$$

and

$$\sigma_h^2(t) = \sigma_P^2(f(t))m_f^2(t) + m_P(f(t))\sigma_f^2(t), \quad \text{for any } t \in [0, R),$$

see equations (1.3.27) and (1.3.28) for further details. Combining the expression for $\sigma_h^2(t)$ with the fact that $f$ is strongly Gaussian and also with the inequality $\sigma_h^2(t) \geq m_P(f(t))\sigma_f^2(t)$, for any $t \in (0, R)$, we find that

$$\lim_{t\uparrow R} \sigma_h^2(t) = +\infty.$$

This is the variance condition for strongly Gaussian power series.

Proposition 3.1.31 gives that

$$\mathbf{E}(e^{i\theta \breve{Z}_t}) = \sum_{j=0}^{N} \mathbf{P}\left(Y_{f(t)} = j\right) \mathbf{E}(e^{i\theta X_t/\sigma_h(t)})^j e^{-i\theta \frac{m_h(t)}{\sigma_h(t)}},$$

for any $\theta \in \mathbb{R}$ and any $t \in (0, R)$, therefore

$$\left|\mathbf{E}(e^{i\theta \breve{Z}_t}) - \mathbf{E}(e^{i\theta X_t/\sigma_h(t)})^N e^{-i\theta m_h(t)/\sigma_h(t)}\right| \leq (1 - \mathbf{P}(Y_{f(t)} = N)) \left|\mathbf{E}(e^{i\theta X_t/\sigma_h(t)})\right| \tag{3.2.6}$$

inequality (3.2.6) gives that

$$\begin{aligned}\int_{|\theta|\leq \pi\sigma_h(t)} \left|\mathbf{E}(e^{i\theta \breve{Z}_t}) - \mathbf{E}\left(e^{i\theta \frac{X_t}{\sigma_h(t)}}\right)^N e^{-i\frac{\theta m_h(t)}{\sigma_h(t)}}\right| d\theta &\leq (1 - \mathbf{P}(Y_{f(t)} = N)) \int_{|\theta|\leq \pi\sigma_h(t)} \left|\mathbf{E}(e^{i\theta \check{X}_t \frac{\sigma_f(t)}{\sigma_h(t)}})\right| d\theta \\ &= (1 - \mathbf{P}(Y_{f(t)} = N)) \frac{\sigma_h(t)}{\sigma_f(t)} \int_{|\theta|\leq \pi\sigma_f(t)} \left|\mathbf{E}(e^{i\theta \check{X}_t})\right| d\theta\end{aligned}$$

Now using that $\lim_{t\uparrow R} f(t) = +\infty$ we have $\lim_{t\uparrow R} \mathbf{P}(Y_{f(t)} = N) = 1$, hypothesis $(\star)$ combined with Lemma 3.1.29 gives that $\lim_{t\uparrow R} \sigma_h(t)/\sigma_f(t) = 1/\sqrt{N}$. Using that $f$ is strongly Gaussian we obtain that

$$\int_{|\theta|\leq \pi\sigma_f(t)} \left|\mathbf{E}(e^{i\theta \check{X}_t})\right| d\theta = O(1), \quad \text{as } t \uparrow R.$$

The previous analysis gives that

$$\lim_{t\uparrow R} \int_{|\theta|\leq \pi\sigma_h(t)} \left|\mathbf{E}(e^{i\theta \breve{Z}_t}) - \mathbf{E}(e^{i\theta X_t/\sigma_h(t)})^N e^{-i\theta m_h(t)/\sigma_h(t)}\right| d\theta = 0\,.$$

Observe that for any $t \in (0, R)$ and any $\theta \in \mathbb{R}$ we have

$$\mathbf{E}\left(e^{i\theta X_t/\sigma_h(t)}\right)^N e^{-i\theta \frac{m_h(t)}{\sigma_h(t)}} = \mathbf{E}(e^{i\theta \check{X}_t \frac{\sigma_f(t)}{\sigma_h(t)}})e^{-i\theta \frac{(m_h(t) - N m_f(t))}{\sigma_h(t)}}$$

then

$$\int_{|\theta|\le\pi\sigma_h(t)}\left|\mathbf{E}\left(e^{i\theta\check{X}_t\frac{\sigma_f(t)}{\sigma_h(t)}}\right)^N e^{-i\theta\frac{(m_h(t)-Nm_f(t))}{\sigma_h(t)}}-e^{-\theta^2/2}\right|d\theta$$

$$=\frac{\sigma_h(t)}{\sigma_f(t)}\int_{|\phi|\le\pi\sigma_f(t)}\left|\mathbf{E}(e^{i\phi\check{X}_t})^N e^{-i\phi\frac{(m_h(t)-Nm_f(t))}{\sigma_f(t)}}-e^{-\frac{\sigma_h^2(t)}{\sigma_f^2(t)}\phi^2/2}\right|d\phi$$

where we change variables in the second integral: $\theta/\sigma_h(t)=\phi/\sigma_f(t)$. Recall that combining $(\star)$ and Lemma 3.1.29 we have $\lim_{t\uparrow R}\sigma_h(t)/\sigma_f(t)=1/\sqrt{N}$. Now we study the integral in the previous equality.

We have

$$I=\int_{|\phi|\le\pi\sigma_f(t)}\left|\mathbf{E}(e^{i\phi\check{X}_t})^N e^{-i\phi\frac{(m_h(t)-Nm_f(t))}{\sigma_f(t)}}-e^{-\frac{\sigma_h^2(t)}{\sigma_f^2(t)}\phi^2/2}\right|d\phi$$

$$=\int_{|\phi|\le\pi\sigma_f(t)}\left|\mathbf{E}(e^{i\phi\check{X}_t})^N-e^{-\frac{\sigma_h^2(t)}{\sigma_f^2(t)}\phi^2/2}e^{i\phi\frac{(m_h(t)-Nm_f(t))}{\sigma_f(t)}}\right|d\phi$$

$$\le\int_{|\phi|\le\pi\sigma_f(t)}\left|\mathbf{E}(e^{i\phi\check{X}_t})^N-e^{-\phi^2N/2}\right|d\phi+\int_{|\phi|\le\pi\sigma_f(t)}\left|e^{-\phi^2N/2}-e^{-\frac{\sigma_h^2(t)}{\sigma_f^2(t)}\phi^2/2}e^{i\phi\frac{(m_h(t)-Nm_f(t))}{\sigma_f(t)}}\right|d\phi$$

observe that

$$\int_{|\phi|\le\pi\sigma_f(t)}\left|\mathbf{E}(e^{i\phi\check{X}_t})^N-e^{-\phi^2N/2}\right|d\phi\le N\int_{\pi\sigma_f(t)}\left|\mathbf{E}(e^{i\phi\check{X}_t})-e^{-\phi^2/2}\right|d\phi$$

here we use the inequality $|z^N-w^N|\le N|z-w|$, for $|z|\le1$, $|w|\le1$, therefore, using that $f$ is strongly Gaussian we have

$$\lim_{t\uparrow R}\int_{|\phi|\le\pi\sigma_f(t)}\left|\mathbf{E}(e^{i\phi\check{X}_t})^N-e^{-\phi^2N/2}\right|d\phi=0,$$

this also follows from the fact that for $f\in\mathcal{K}$ a strongly Gaussian power series we have that $f^N$ is strongly Gaussian, for any integer $N\ \ge1$.

Now observe that

$$\lim_{t\uparrow R}\left(e^{-\phi^2N/2}-e^{-\frac{\sigma_h^2(t)}{\sigma_f^2(t)}\phi^2/2}e^{i\phi\frac{(m_h(t)-Nm_f(t))}{\sigma_f(t)}}\right)=0,$$

here we combine $(\star)$ with Lemma 3.1.29 to obtain that

$$\lim_{t\uparrow R}\frac{\sigma_h^2(t)}{\sigma_f^2(t)}=N,\quad\text{and}\quad\lim_{t\uparrow R}\frac{m_h(t)-Nm_f(t)}{\sigma_f(t)}=0.$$

Finally dominated convergence gives that

$$\lim_{t\uparrow R}\int_{|\phi|\le\pi\sigma_f(t)}\left|e^{-\phi^2N/2}-e^{-\frac{\sigma_h^2(t)}{\sigma_f^2(t)}\phi^2/2}\,e^{i\phi\frac{(m_h(t)-Nm_f(t))}{\sigma_f(t)}}\right|d\phi=0,$$

and therefore $h=P\circ f$ is strongly Gaussian. □

### D. Some applications

Here we apply some of the previous results to some particular cases.

**D.1. Product of a polynomial and a strongly Gaussian power series** Observe that the product of a polynomial $P\in\mathcal{K}$ and a strongly Gaussian power series $f\in\mathcal{K}$ gives a strongly Gaussian power series $h=Pf$.

If $f$ is an entire function, then we apply the discussion above for power series with the same radius of convergence. In this case we have

$$\lim_{t\to+\infty}\frac{\sigma_P(t)}{\sigma_h(t)}=0,\quad\text{and}\quad\lim_{t\uparrow R}\frac{\sigma_f(t)}{\sigma_h(t)}=1.$$

Recall that polynomials are characterized, among entire functions in $\mathcal{K}$, by $\lim_{t\to+\infty}\sigma_P^2(t)=0$, that is, $f\in\mathcal{K}$, an entire function, is a polynomial if and only if $\lim_{t\to+\infty}\sigma_f^2(t)=0$.

If $f\in\mathcal{K}$ has finite radius of convergence $R>0$, then we apply the discussion above for power series with different radius of convergence, to obtain that $h=Pf$ is a strongly Gaussian power series.

**D.2. Composition of a polynomial with an exponential** We already know that the Poisson family, that is the Khinchin family associated to $f(z)=e^z\in\mathcal{K}$, is strongly Gaussian. Recall that $m_f(t)=t$ and $\sigma_f^2(t)=t$. Fix $P\in\mathcal{K}$ a polynomial of degree $N\ge1$. In this case we have that

$$\lim_{t\to+\infty}\frac{m_f^2(t)}{\sigma_f^2(t)f(t)}=\lim_{t\to+\infty}\frac{t}{e^t}=0.$$

Applying Theorem 3.2.14 we conclude that $h(z)=P(e^z)$ is strongly Gaussian.

**D.3. Composition of a polynomial with an entire function of finite order $\rho\ge0$** Let $f\in\mathcal{K}$ be an entire function of finite order $\rho\ge0$ which is strongly Gaussian. We already know that in this case

$$\lim_{t\to+\infty}\frac{m_f^2(t)}{\sigma_f^2(t)f(t)}=0,$$

see Lemma 3.1.30 for further details. Applying Theorem 3.2.14 we conclude that $h=P\circ f$ is strongly Gaussian.

## 3.3 Hayman Class

We have seen that the functions $f \in \mathcal{K}$ which are strongly Gaussian have excellent asymptotic properties, as $t \uparrow R$; local central limit theorem and asymptotic formula for the coefficients. We have seen certain number of examples of strongly Gaussian power series.

Hayman's power series, that we introduce now, are power series satisfying certain, concrete, and verifiable properties and, as we will see later on, are strongly Gaussian power series. This result, that we prove in this Section, will allow us to verify that, for instance, the Bell family is strongly Gaussian, and obtain Moser-Wiman's asymptotic formula for the Bell numbers. Later on, we will also use Hayman's power series to obtain Hardy-Ramanujan asymptotic formula.

We say that a power series $f \in \mathcal{K}$ is in the Hayman class if

$$\text{(3.3.1)} \qquad \text{(variance condition)} \qquad \lim_{t\uparrow R} \sigma_f(t) = +\infty.$$

and if for certain function $h : [0, R) \to (0, \pi]$, which we refer as the cut (between a major arc and a minor arc), there hold

$$\text{(3.3.2)} \qquad \text{(major arc)} \qquad \lim_{t\uparrow R} \sup_{|\theta|\le h(t)\sigma_f(t)} \left|\mathbf{E}(e^{i\theta \breve{X}_t})e^{\theta^2/2} - 1\right| = 0\,.$$

$$\text{(3.3.3)} \qquad \text{(minor arc)} \qquad \lim_{t\uparrow R} \sigma_f(t) \sup_{\sigma_f(t)h(t)\le|\theta|\le\pi\sigma_f(t)} |\mathbf{E}(e^{i\theta \breve{X}_t})| = 0\,.$$

The cut $h$ of a function $f$ in the Hayman class is not uniquely determined.

The functions in the Hayman class are called admissible by Hayman [48]; they are also called Hayman-admissible, or just H-admissible.

**Remark 3.3.1** (Conditions to belong in the Hayman class)**.**

- Condition (3.3.2) gives that the characteristic function of $\breve{X}_t$ is approximated uniformly by $e^{-\theta^2/2}$ on the major arc $|\theta| \le h(t)\sigma_f(t)$.
- Condition (3.3.3) gives that the characteristic function of $\breve{X}_t$ is $o(1/\sigma_f(t))$, uniformly on the minor arc $\sigma_f(t)h(t) \le |\theta| \le \pi\sigma_f(t)$.

$\boxtimes$

The major arc is comparatively, and in general, smaller than the minor arc. The appellative *major* is used to highlight that its contribution is, relatively, more important.

As first illustration we prove that the exponential function $f(z) = e^z$ is in the Hayman class. This is the seminal, and primigenious, example of a function in the Hayman class.

**Example 3.3.2.** For the exponential $f(z) = e^z$, we have $m_f(t) = t$ and $\sigma_f^2(t) = t$. Let's see how to choose a cut function $h$ with values in $(0, \pi]$.

The major arc condition, see (3.3.2), requires that

$$\lim_{t\to\infty} \sup_{|\theta|\le h(t)} \left| e^{t(e^{i\theta}-1-i\theta+\theta^2/2)} - 1 \right| = 0.$$

We have the inequality

$$|e^{i\theta} - 1 - i\theta + \theta^2/2| \le \frac{|\theta|^3}{6}, \qquad \text{for any } \theta \in \mathbb{R},$$

now using that $|e^z - 1| \le e^{|z|}|z|$, for any $z \in \mathbb{C}$, we find that

$$\sup_{|\theta|\le h(t)} \left| e^{t(e^{i\theta}-1-i\theta+\theta^2/2)} - 1 \right| \le e^{th(t)^3} \frac{1}{6} th(t)^3,$$

and therefore we impose the condition $\lim_{t\to\infty} th(t)^3 = 0$ on $h$.

The minor arc condition, see (3.3.3), requires that

$$\lim_{t\to\infty} \sigma_f(t) \sup_{h(t)\le|\theta|\le\pi} e^{t(\cos(\theta)-1)} = 0,$$

that is,

$$\lim_{t\to\infty} \sqrt{t} e^{t(\cos(h(t))-1)} = 0.$$

Now using that $\cos(\theta) - 1 \le -\frac{2}{\pi^2}\theta^2$, for $|\theta| \le \pi$, we have that

$$e^{t(\cos(h(t))-1)} \le e^{-2/\pi^2 th(t)^2},$$

and therefore, we observe that for the minor arc condition to hold we need to impose on $h$ the condition

$$\lim_{t\to\infty} \sqrt{t} e^{-2/\pi^2 th(t)^2} = 0\,.$$

If we take $h(t) = \min\{\pi, t^{-\alpha}\}$, with $\alpha \in (1/3, 1/2)$, then both conditions, the major arc and the minor arc condition, hold simultaneously. ⊡

### 3.3.1 Size of the cut

The value $\theta = h(t)\sigma_f(t)$ set the cut between the major and the minor arc, there both conditions hold, and this impose certain restrictions on $h$.

**Lemma 3.3.3.** *Let $f \in \mathcal{K}$ a power series in the Hayman class, then*

$$\lim_{t\uparrow R} h(t)\sigma_f(t) = +\infty.$$

*In fact, for certain $t_0 \in (0, R)$ we have*

$$(\star) \quad h(t) \ge \frac{\sqrt{2\ln(\sigma_f(t))}}{\sigma_f(t)}, \qquad \textit{for any } t \in [t_0, R).$$

*Proof.* For $\theta = h(t)\sigma_f(t)$, condition (3.3.2), written in terms of $f$, gives that

$$\lim_{t\uparrow R} \frac{|f(te^{ih(t)})|}{f(t)} e^{h(t)^2\sigma_f^2(t)/2} = 1,$$

while condition (3.3.3), also written in terms of $f$, gives that

$$\lim_{t\uparrow R} \sigma_f(t) \frac{|f(te^{ih(t)})|}{f(t)} = 0.$$

Therefore we have

$$\lim_{t\uparrow R} \frac{e^{h(t)^2\sigma_f^2(t)/2}}{\sigma_f(t)} = \infty.$$

Given that $\lim_{t\uparrow R} \sigma_f(t) = +\infty$, then $\lim_{t\uparrow R} e^{h(t)^2\sigma_f^2(t)/2} = +\infty$.

For certain $t_0 \in (0, R)$, we have that $e^{h(t)^2\sigma_f^2(t)/2}/\sigma_f(t) \geq 1$, for any $t \geq t_0$, and therefore we deduce $(\star)$. □

On the other hand, we can assume that given $A > 1$ the cut function verifies that for certain $t_0 \in (0, R)$ we have

$$(\star\star) \quad h(t) \leq A\frac{\sqrt{2\ln(\sigma_f(t))}}{\sigma_f(t)}, \qquad \text{for any } t \in (t_0, R).$$

Indeed. Take $t_0 \in (0, R)$ such that

$$\text{(3.3.4)} \qquad \frac{|f(te^{i\theta})|}{f(t)} \leq 2e^{-\theta^2\sigma_f^2(t)/2}, \qquad \text{for any } t \in [t_0, R) \text{ and } |\theta| \leq h(t),$$

here we use condition (3.3.2).

Fix $t \in [t_0, R)$ and assume that $h(t) > A\sqrt{2\ln(\sigma_f(t))}/\sigma_f(t)$. We will see that in this case we can enlarge the minor arc (and reduce the major arc) putting the cut at $A\sqrt{2\ln(\sigma_f(t))}/\sigma_f(t)$. We need to check that if

$$A\frac{\sqrt{2\ln(\sigma_f(t))}}{\sigma_f(t)} \leq |\theta| \leq h(t),$$

then a similar uniform bound to that in the minor arc condition (3.3.3) holds. Observe that such $\theta$ are in the original major arc of $f$ and therefore (3.3.4) holds. Then, for those $\theta$, we have

$$\frac{|f(te^{i\theta})|}{f(t)} \leq 2e^{-\theta^2\sigma_f^2(t)/2} \leq 2e^{-A^2\ln(\sigma_f(t))} = \frac{2}{\sigma_f(t)^{A^2}}.$$

Denote $A(t) = A\sqrt{2\ln(\sigma_f(t))}/\sigma_f(t)$. Using the previous inequality we find that

$$\sigma_f(t) \sup_{A(t)\leq|\theta|\leq h(t)} \frac{|f(te^{i\theta})|}{f(t)} \leq \frac{2}{\sigma_f(t)^{A^2-1}}$$

We conclude then that we can assume that, for any $A > 1$, there exists $t_0 \in (0, R)$ such that

$$\frac{\sqrt{2\ln(\sigma_f(t))}}{\sigma_f(t)} \le h(t) \le A\frac{\sqrt{2\ln(\sigma_f(t))}}{\sigma_f(t)}, \qquad \text{for any } t \in [t_0, R).$$

In particular this implies that we can take a cut $h$ verifying that

$$h(t) \asymp \frac{\sqrt{\ln(\sigma_f(t))}}{\sigma_f(t)}, \qquad \text{as } t \uparrow R\,,$$

and therefore, in this case, we have that

$$h(t)\sigma_f(t) \asymp \sqrt{\ln(\sigma_f(t))}, \qquad \text{as } t \uparrow R\,.$$

### 3.3.2 Global criteria for the major arc condition for non-vanishing power series

Following the ideas in the proof of Theorem 3.1.4 we give a global criteria for non-vanishing power series $f \in \mathcal{K}$ to verify the major arc condition for some cut function $h$.

**Theorem 3.3.4.** *Let $f \in \mathcal{K}$ be a power series with radius of convergence $R > 0$ and having no zeros in $\mathbb{D}(0, R)$. For $h : [0, R) \to (0, \pi]$ a cut function suppose that*

$$\lim_{s\uparrow\ln(R)} \left(\sup_{|\theta|\le h(t)\sigma_f(t)} |F'''(s+i\theta)|\right) h(e^s)^3 = 0\,,$$

*then the major arc condition (3.3.2) holds.*

*Proof.* Denote $t = e^s$, for $s < \ln(R)$ and $|\theta| \le h(t)\sigma_f(t)$ we have

$$|F(s+i\theta) - F(s) - F'(s)i\theta + F''(s)\frac{\theta^2}{2}| \le \sup_{|\phi|\le h(t)\sigma_f(t)} |F'''(s+i\phi)|\frac{|\theta|^3}{6},$$

therefore

$$(\triangle) \qquad \left|\ln\left(\mathbf{E}(e^{i\theta\breve{X}_t})\right) + \frac{\theta^2}{2}\right| \le \frac{\sup_{|\phi|\le h(t)\sigma_f(t)} |F'''(s+i\phi)|}{\sigma_f(t)^3}\frac{|\theta|^3}{6}.$$

Using the numerical inequality $|e^z - 1| \le |z|e^{|z|}$, for any $z \in \mathbb{C}$, we find that

$$(3.3.5) \qquad |\mathbf{E}(e^{i\theta\breve{X}_t})e^{\theta^2/2} - 1| \le |\ln(\mathbf{E}(e^{i\theta\breve{X}_t})e^{\theta^2/2})| \exp(|\ln(\mathbf{E}(e^{i\theta\breve{X}_t})e^{\theta^2/2})|)\,.$$

If we denote $\Omega(t) = \left(\sup_{|\phi|\le h(t)\sigma_f(t)} |F'''(s+i\phi)|\right) h(t)^3$, inequality $(\triangle)$ gives that

$$\sup_{|\theta|\ \le h(t)\sigma_f(t)} |\mathbf{E}(e^{i\theta\breve{X}_t})e^{\theta^2/2} - 1| \le \Omega(t)e^{\Omega(t)},$$

therefore our hypothesis, that is, $\lim_{t\uparrow R}\Omega(t) = 0$, gives that the major arc condition (3.3.2) holds.

□

We also have the following corollary.

**Corollary 3.3.5.** *With the same conditions than in the previous theorem, assume that*

$$|F'''(s+i\theta)| \leq |F'''(s)|, \qquad \text{for any } \theta \in \mathbb{R} \text{ and } s < \ln(R)$$

*and also that*

$$\lim_{s \uparrow \ln(R)} |F'''(s)| h(e^s)^3 = 0,$$

*then the major arc condition (3.3.2) holds.*

An instance of the situation presented in Corollary 3.3.5 is given by the exponentials $f = e^g$, with $g \in \mathcal{K}_s$. This example will be relevant later on, as encodes some relevant combinatorial structures, see Chapter 4 for further details.

## A. Applications of the global criteria for the major arc condition for non-vanishing power series

Now we apply the previous criteria to some of our prominent examples.

**A.1. Partitions of integers.** For the generating function of the partitions of integers we have

$$P(z) = \prod_{j=1}^{\infty} \frac{1}{1-z^j}, \qquad \text{for any } |z| < 1.$$

In fact we know, see equation (3.1.8), that

$$F'''(-s) \sim \zeta(2)\Gamma(4)s^{-4}, \qquad \text{as } s \downarrow 0.$$

If we take $h(e^{-s}) = \min\{\pi, s^\alpha\}$, with $\alpha > 0$, then

$$F'''(-s)h(e^{-s})^3 \sim \zeta(2)\Gamma(4)\frac{h(e^{-s})^3}{s^4}, \qquad \text{as } s \downarrow 0.$$

this implies that $\alpha \in (4/3, +\infty)$ and therefore the major arc condition (3.3.2) with cut $h$ holds for $\alpha > 4/3$.

**A.2. Exponential function.** In this case we have $f(z) = e^z$, then $F(s) = e^s$, and therefore

$$|F'''(s+i\theta)|h(e^s)^3 = F'''(s)h(t)^k = e^s h(e^s)^3.$$

Denote $h(t) = \min\{\pi, t^{-\alpha}\}$, for $\alpha > 0$, and then $1 - 3\alpha < 0$, that is, $\alpha \in (1/3, +\infty)$. The previous reasoning gives that the major arc condition with cut $h$ holds for $\alpha > 1/3$.

**A.3. Exponential of a polynomial not necessarily in $\mathcal{K}$.** Let $p(z) = a_N z^N + \cdots + a_0$ be a polynomial of degree $N \geq 1$ such that $e^p \in \mathcal{K}$; observe that $p$ is not necessarily in $\mathcal{K}$. We want to prove that $f = e^p$ verifies the major arc condition for some cut function $h$.

Given that $p$ is a polynomial, there exists a constant $C > 0$ such that

$$|F'''(s + i\theta)| \sim Ce^{sN}, \quad \text{as } s \to \infty\,,$$

then, in order to apply Theorem 3.3.4, we want the cut function $h$ to verify that

$$\lim_{s\to\infty} e^{sN} h(e^s)^3 = 0.$$

If we take $h(t) = \min\{\pi, t^{-\alpha}\}$, with $\alpha > 0$, then $\alpha \in (N/3, +\infty)$ and therefore the major arc condition with cut $h$ holds for $\alpha > N/3$. Observe that for $p(z) = z$, a polynomial of degree $N = 1$, we recover the previous example: the exponential function $f(z) = e^z$.

**Remark 3.3.6** (Heuristic: *guessing* the upper bound for $\alpha$)**.** A way of guessing the upper bound for $\alpha$ in the previous examples is applying Lemma 3.3.3 to the candidate for cut function, that is, to impose the necessary condition for functions in the Hayman class:

$$\lim_{t\uparrow R} \sigma_f(t) h(t) = +\infty.$$

Observe that this is simply a necessary condition and doesn't imply that $f$ is in the Hayman class. If $f$ is in the Hayman class with cut given by $h$, depending on a parameter $\alpha$, then $\alpha$ should take certain specific values, this is the claim; Outside this interval of possible values for $\alpha$, the function $h$ is not a cut function for $f$.

In the previous examples we have taken $h$ depending on a parameter $\alpha$ (a family of cuts, recall that the cut function is not unique) and the possible values of $\alpha$ which are candidates for giving a cut $h$ such that $f$ is in the Hayman class are the following:

- Exponential function: we impose to the candidate for cut $h$ the necessary condition

  $$\lim_{t\uparrow R} \sigma_f(t) h(t) = +\infty.$$

  Observe that $\sigma_f(t) = t^{1/2}$, then $h(t)\sigma_f(t) = t^{1/2-\alpha}$, for $t \geq 1$, therefore if $f(z) = e^z$ is in the Hayman class with cut $h$, then $\alpha \in (1/3, 1/2)$, and this is the case, the reciprocal is also true; see Example 3.3.2 above.

- Generating function for the partitions of integers: again we impose to the candidate for cut function $h$ the necessary condition

  $$\lim_{t\uparrow R} \sigma_f(t) h(t) = +\infty.$$

  In this case we have

  $$\sigma_f(e^{-s}) \sim (\zeta(2)\Gamma(3))^{1/2} s^{-3/2}, \quad \text{as } s \downarrow 0,$$

  see equation (3.1.8) for further details, then $h(e^{-s})\sigma_f(e^{-s}) \sim \zeta(2)\Gamma(3)s^{\alpha-3/2}$ as $s \downarrow 0$. The previous discussion gives that if $P$ is in the Hayman class with cut $h$, then $\alpha \in (4/3, 3/2)$. We will see later on that the reciprocal is also true.

- Exponential of a polynomial not in $\mathcal{K}$: we impose that

$$\lim_{t\uparrow R} \sigma_f(t)h(t) = +\infty.$$

Notice that $a_N > 0$ is a positive real number, because $f \in \mathcal{K}$ and the variance $\sigma_f(t) > 0$, for any $t > 0$. Moreover, there exists a constant $C_N > 0$ such that

$$\sigma_f(t) \sim C_N t^{N/2}, \qquad \text{as } t \to +\infty,$$

therefore $\sigma_f(t)h(t) \sim C_N t^{N/2} h(t)$, as $t \to +\infty$. The previous discussion gives that if $e^P$ is in the Hayman class with cut $h$, then $\alpha \in (N/3, N/2)$. We will see later on that, after imposing certain conditions on the polynomial $P$, this is the case and the reciprocal is also true. Observe that for $P(z) = z$ we obtain $f(z) = e^z$ and therefore we recover the interval $(1/3, 1/2)$, see Example 3.3.2 above.

⊠

### 3.3.3 Hayman class and strongly Gaussian power series

We see here that functions in the Hayman class are strongly Gaussian.

**Theorem 3.3.7** (Hayman)**.** *If $f \in \mathcal{K}$ is in the Hayman class, then $f$ is strongly Gaussian.*

*Proof.* The power series $f$ is in the Hayman class, therefore the variance condition gives that $\lim_{t\uparrow R} \sigma_f(t) = +\infty$.

Denote

$$I_t = \int_{-\pi\sigma_f(t)}^{\pi\sigma_f(t)} |\mathbf{E}(e^{i\theta \breve{X}_t}) - e^{-\theta^2/2}| d\theta, \qquad \text{for any } t \in (0, R).$$

We divide the integral $I_t$ in two integrals, one in the major arc, that we denote by $J_t$, and one in the minor arc that we denote by $K_t$, that is:

$$J_t = \int_{|\theta| \le h(t)\sigma_f(t)} |\mathbf{E}(e^{i\theta \breve{X}_t}) - e^{-\theta^2/2}| d\theta\,,$$

$$K_t = \int_{h(t)\sigma_f(t)\ \le |\theta| \le \pi\sigma_f(t)} |\mathbf{E}(e^{i\theta \breve{X}_t}) - e^{-\theta^2/2}| d\theta\,.$$

For the integral $K_t$ we have:

$$K_t \le 2\pi\sigma_f(t) \sup_{h(t)\sigma_f(t) \le |\theta| \le \pi\sigma_f(t)} |\mathbf{E}(e^{i\theta \breve{X}_t})| + \int_{|\theta| \ge h(t)\sigma_f(t)} e^{-\theta^2/2} d\theta.$$

The first term at the right-hand side tends to 0, by virtue of the minor arc condition (3.3.3), while the second term tends to 0 because of Lemma 3.3.3.

For the integral $J_t$ we have:

$$J_t = \int_{|\theta|\le h(t)\sigma_f(t)} e^{-\theta^2/2}|\mathbf{E}(e^{i\theta\breve{X}_t})e^{\theta^2/2}-1|d\theta$$
$$\le \left(\sup_{|\theta|\le h(t)\sigma_f(t)} |\mathbf{E}(e^{i\theta\breve{X}_t})e^{\theta^2/2}-1|\right)\int_{\mathbb{R}} e^{-\theta^2/2}d\theta$$

using the major arc condition (3.3.2) the previous integral tends to 0, as $t\uparrow R$. □

### 3.3.4 Operations with Hayman functions

Here we prove that some operations between functions in the Hayman class preserve the property of being in the Hayman class.

#### A. Product of functions in the Hayman class

Now we analyze if the product of functions in the Hayman class is also in the Hayman class.

**A.1. Power series with the same radius of convergence** Let $f, g \in \mathcal{K}$ power series with radius of convergence $R > 0$ which are in the Hayman class with cuts $h_f$ and $h_g$, respectively. Denote $L = fg$ the product of $f$ and $g$, we prove here that $L$ is also a Hayman function.

Denote $(X_t)$, $(Y_t)$ and $(Z_t)$ the Khinchin families of $f, g$ and $L$, respectively. Recall that

$$Z_t \stackrel{d}{=} X_t \oplus Y_t, \qquad \text{for any } t\in[0,R). \tag{3.3.6}$$

The previous equality in distribution gives that for any $t\in[0,R)$ we have that

$$m_L(t) = m_f(t)+m_g(t) \quad \text{and} \quad \sigma_L^2(t) = \sigma_f^2(t)+\sigma_g^2(t),$$

and therefore, using that $f$ and $g$ are Hayman functions, we have that $\lim_{t\uparrow R}\sigma_L^2(t) = +\infty$ and the variance condition holds.

Define $h(t) = \min\{h_f(t), h_g(t)\}$, the candidate for cut function of the product $L = fg$ and denote

$$\varepsilon_f(t) = \sup_{|\theta|\le h_f(t)\sigma_f(t)} \left|\mathbf{E}(e^{i\theta\breve{X}_t})e^{\theta^2/2}-1\right|, \qquad \text{for any } t\in(0,R),$$

and, analogously for $g$, we denote $\varepsilon_g(t)$ the corresponding quantity for $g$. The major arc conditions for $f$ and $g$ can be summarized as $\lim_{t\uparrow R}(\varepsilon_f(t)+\varepsilon_g(t)) = 0$.

In addition we denote

$$\eta_f(t) = \sigma_f(t)\sup_{h_f(t)\le|\theta|\le\pi}\frac{|f(te^{i\theta})|}{f(t)}, \qquad \text{for any } t\in(0,R),$$

and, anagously for $g$, we denote $\eta_g(t)$ the corresponding quantity for $g$. The minor arc conditions of $f$ and $g$ can be summarized as $\lim_{t\uparrow R}(\eta_f(t)+\eta_g(t)) = 0$.

Lemma 3.3.3 gives that there exists $0 < T < R$ such that $\sigma_f(t)h_f(t) > 1$ and $\sigma_g(t)h_g(t) > 1$, simultaneously, for any $T < t < R$.

We study first the major arc for the product $L = fg$ with cut $h$: using the equality in distribution (3.3.6) we find that

$$\mathbf{E}(e^{i\theta \check{Z}_t}) = \mathbf{E}(e^{i\theta(\sigma_f(t)/\sigma_L(t))\check{X}_t})\mathbf{E}(e^{i\theta(\sigma_g(t)/\sigma_L(t))\check{Y}_t})$$

therefore

$$\mathbf{E}(e^{i\theta \check{Z}_t})e^{\theta^2/2} = \mathbf{E}(e^{i\theta(\sigma_f(t)/\sigma_L(t))\check{X}_t})e^{(\sigma_f(t)/\sigma_L(t))\theta^2/2}\mathbf{E}(e^{i\theta(\sigma_g(t)/\sigma_L(t))\check{Y}_t})e^{(\sigma_g(t)/\sigma_L(t))\theta^2/2}.$$

For $|\theta| \leq h(t)\sigma_L(t)$, that is, for $\theta$ in the potential major arc of $R$, we have that

$$\frac{\sigma_f(t)}{\sigma_L(t)}|\theta| \leq h(t)\sigma_f(t) \leq h_f(t)\sigma_f(t), \quad \text{and} \quad \frac{\sigma_g(t)}{\sigma_L(t)}|\theta| \leq h(t)\sigma_g(t) \leq h_g(t)\sigma_g(t),$$

and therefore

$$|\mathbf{E}(e^{i\theta \check{Z}_t})e^{\theta^2/2} - 1| \leq \varepsilon_f(t) + \varepsilon_g(t) + \varepsilon_f(t)\varepsilon_g(t),$$

which gives the major arc condition for $L = fg$.

Now we study the minor arc, that is, $\sigma_L(t)h(t) \leq |\theta| \leq \sigma_L(t)\pi$; we analyze the quantity

$$\sigma_L(t)|\mathbf{E}(e^{i\theta \check{Z}_t})| = \sigma_L(t)\frac{|f(te^{i\theta/\sigma_L(t)})|}{f(t)}\frac{|g(te^{i\theta/\sigma_L(t)})|}{g(t)}.$$

Assume that $\sigma_f(t) \geq \sigma_g(t)$; we argue analogously if $\sigma_f(t) \leq \sigma_g(t)$. We have $\sigma_f(t) \leq \sigma_L(t) \leq \sqrt{2}\sigma_f(t)$ and, moreover, $|g(te^{i\theta})| \leq g(t)$, we focus our attention in finding an upper bound for the quantity

$$A(t,\theta) = \sigma_f(t)\frac{|f(te^{i\theta})|}{f(t)}.$$

We distinguish two cases:

1. Assume first that $|\theta| \geq h_f(t)$, then we have $A(t,\theta) \leq \eta_f(t)$.
2. Assume now that $|\theta| \leq h_f(t)$, then $h(t) = h_g(t)$ and $h_g(t) \leq |\theta|$. This is the intermediate zone, intersection of the major arc of $f$ and the minor arc of $g$.

   Because $h(t) = h_g(t)$, the border between the major arc and the minor arc of $g$, we have, using the major arc condition, that

   $$(\star) \qquad \frac{|g(te^{ih(t)})|}{g(t)} = e^{-h(t)^2\sigma_g(t)^2/2}(1 + \tilde{\varepsilon}(t))$$

   where $\tilde{\varepsilon}(t) \leq \varepsilon_g(t)$.

We multiply by $\sigma_g(t)$ and apply the minor arc bound to obtain that

$$\sigma_g(t)e^{-h(t)^2\sigma_g^2(t)/2}(1+\tilde{\varepsilon}(t)) \leq \eta_g(t),$$

and therefore

$$\sigma_g(t)e^{-h(t)^2\sigma_g^2(t)/2} \leq \frac{\eta_g(t)}{(1+\tilde{\varepsilon}(t))}.$$

For $\gamma > 0$ fixed, the function $s \mapsto e^{-\gamma s^2/2}$ is decreasing for $s^2 > 1/\gamma$, then given that $\sigma_g(t)h_g(t) > 1$ and $\sigma_f(t) \geq \sigma_g(t)$, then we have

$$\sigma_f(t)e^{-h(t)^2\sigma_f^2(t)/2} \leq \frac{\eta_g(t)}{1+\tilde{\varepsilon}(t)}.$$

As $\theta$ is in the major arc of $f$ and also $|\theta| \geq h(t)$, we have that

$$\frac{|f(te^{i\theta})|}{f(t)} \leq e^{-\theta^2\sigma_f^2(t)/2}(1+\varepsilon_f(t)) \leq e^{-h(t)^2\sigma_f^2(t)/2}(1+\varepsilon_f(t)),$$

and then we conclude that

$$\sigma_f(t)\frac{|f(te^{i\theta})|}{f(t)} \leq \frac{\eta_g(t)}{1+\tilde{\varepsilon}(t)}(1+\varepsilon_f(t)),$$

which gives the minor arc condition.

We collect the previous discussion by means of the following theorem. See [48][Thm VII]

**Theorem 3.3.8** (Hayman). *Let $f, g \in \mathcal{K}$ power series with radius of convergence $R > 0$. Assume that $f$ and $g$ are in the Hayman class with cuts $h_f$ and $h_g$, respectively, then $L = fg$ belongs to the Hayman class with cut function $h = \min\{h_f, h_g\}$.*

**A.2. Power series with different radius of convergence** Now we consider the case where $f$ is a power series in $\mathcal{K}$ with radius of convergence $R > 0$ and $g$ is power series in $\mathcal{K}$ with radius of convergence $S \geq R > 0$. Assume that $f$ is in the Hayman class and also that

$$\Sigma_g = \sup_{t \in [0,R)} \sigma_g^2(t) < +\infty.$$

The power series $g$ is not necessarily in the Hayman class. This situation includes the case in which $S > R$ and also the case where $g$ is a polynomial.

We denote $(X_t)$, $(Y_t)$ and $(Z_t)$ the Khinchin families of $f, g$ and $L = fg$, respectively. In this case we have

$$Z_t \stackrel{d}{=} X_t \oplus Y_t, \qquad \text{for any } t \in [0, R),$$

then we have $\sigma_L^2(t) = \sigma_f^2(t)+\sigma_g^2(t)$, and jtherefore $\lim_{t\uparrow R} \sigma_L^2(t) = +\infty$; this is the variance condition. In particular this means that the product $L = fg$ has radius $R > 0$.

Observe that

$$\lim_{t\uparrow R}\frac{\sigma_g(t)}{\sigma_L(t)}=0,\qquad\text{and}\qquad \lim_{t\uparrow R}\frac{\sigma_f(t)}{\sigma_L(t)}=1.$$

Denote $h_f$ the cut function of $f$. We take as a candidate for cut function of the product $L=fg$ the cut function $h_f$. We have seen before, see Subsection 3.3.1, that we can take a cut function $h$ such that $\lim_{t\uparrow R} h(t)=0$.

For the minor arc, that is, for $h(t)\leq|\theta|\leq\pi$, we have

$$\sigma_L(t)\sup_{h(t)\leq|\theta|\leq\pi}\frac{|f(te^{i\theta})|}{f(t)}\frac{|g(te^{i\theta})|}{g(t)}\leq\frac{\sigma_L(t)}{\sigma_f(t)}\sigma_f(t)\sup_{h(t)\leq|\theta|\ \leq\pi}\frac{|f(te^{i\theta})|}{f(t)},$$

therefore, taking limits, as $t\uparrow R$, we obtain the minor arc condition for $L=fg$.

Now we study the major arc condition, that is, $|\theta|\leq h(t)\sigma_L(t)$. Denote $\lambda(t)=\sigma_f(t)/\sigma_L(t)$ and $\mu(t)=\sigma_g(t)/\sigma_L(t)$. In this case we have

$$\mathbf{E}(e^{i\theta\breve{Z}_t})e^{\theta^2/2}=\mathbf{E}(e^{i\theta\lambda(t)\breve{X}_t})e^{\lambda(t)\theta^2/2}\cdot\mathbf{E}(e^{i\theta\mu(t)\breve{Y}_t})e^{\mu(t)\theta^2/2}$$

We define the functions:

$$\begin{aligned}A(t,\theta)&=\mathbf{E}(e^{i\theta\lambda(t)\breve{X}_t}e^{\lambda(t)\theta^2/2}),\\ B(t,\theta)&=\mathbf{E}(e^{i\theta\mu(t)\breve{Y}_t}),\\ C(t,\theta)&=e^{\mu(t)\theta^2/2}.\end{aligned}$$

We prove now that that these functions converge uniformly to 1, when $|\theta|\leq h(t)\sigma_L(t)$, and of course, as $t\uparrow R$. We collect the previous discussion by means of the following Theorem.

If $|\theta|\leq h(t)\sigma_L(t)$, then $|\theta|\lambda(t)\leq h(t)\sigma_f(t)$, and the convergence of $A(t,\theta)$ follows from the major arc condition for $f$.

For $B(t,\theta)$ and $|\theta|\leq h(t)\sigma_L(t)$ we have, using that $|e^{i\phi}-1|\leq|\phi|$, for any $\phi\in\mathbb{R}$ and that $\mathbf{E}(|\breve{Y}_t|)\leq\mathbf{E}(|\breve{Y}_t|^2)^{1/2}$, that

$$|B(t\theta)-1|\leq\mathbf{E}(|e^{i\theta\mu(t)\breve{Y}_t}-1|)\leq\mu(t)|\theta|\leq h(t)\sigma_g(t)\leq\Sigma_g h(t),$$

and because $\lim_{t\uparrow R} h(t)=0$, we conclude that $B(t,\theta)$ converges uniformly to 1, in the required range, as $t\uparrow R$.

For $C(t,\theta)$ and $|\theta|\leq h(t)\sigma_L(t)$, we have that $\mu(t)|\theta|\leq h(t)\sigma_g(t)\leq\Sigma_g h(t)$, then

$$|C(t,\theta)-1|\leq e^{h(t)^2\Sigma_g^2/2}-1,$$

bound that tends to 0, as $t\uparrow R$.

We collect the previous discussion by means of the following theorem.

**Theorem 3.3.9.** *Let $f,g\in\mathcal{K}$ be a power series with radius of convergence $R>0$ and $S\geq R>0$, respectively. Assume that $f$ is in the Hayman class and also that*

$$\Sigma_g=\sup_{t\in[0,R)}\sigma_g^2(t)<+\infty.$$

*then $L=fg$ is in the Hayman class.*

### B. Powers of functions in the Hayman class

Using the previous criteria for the product of power series in the Hayman class we will obtain that the powers of a Hayman power series are also Hayman power series.

**Theorem 3.3.10.** *Let $f \in \mathcal{K}$ a power series with radius of convergence $R > 0$. Assume that $f$ is in the Hayman class with cut $h$, then $f^N$, with integer $N \geq 1$, is in the Hayman class with cut function $h$.*

The previous result follows applying Theorem 3.3.8.

### C. Some other operations

May be the reason for the relevance of the Hayman class in combinatorial questions is that it enjoys some closure and small perturbation properties.

For instance, see Hayman [48]:

a) If $f$ is in the Hayman class, then $e^f$ is in the Hayman class, see [48][Thm. VI.]

b) If $f$ is in the Hayman class and if $R$ is a polynomial in $\mathcal{K}$, then the product $Rf$ is in the Hayman class, see [48][Thm. VIII]

c) Let $B(z) = \sum_{n=0}^{N} b_n z^n$ be a non-constant polynomial with real coefficients such that for each $d > 1$ there exists $m$, not a multiple of $d$ such that $b_m \neq 0$, and such that if $m(d)$ is the largest such $m$, then $b_{m(d)} > 0$. If $e^B \in \mathcal{K}$, then $e^B$ is in the Hayman class. See [48][Thm. X]

## 3.3.5 Dominated criteria for strongly Gaussian power series

We see now an alternative criteria, a variant, of Theorem 3.3.7 which claims that Hayman functions are strongly Gaussian.

We still have a cut function, and its major and minor arc, but in contrast with the conditions to belong in the Hayman class in this case we assume that the power series $f$ is Gaussian. This relax the condition on the major arc for Hayman functions, which this time is only a uniform bound, instead of uniform convergence. See [6], of Baéz-Duarte, and [83], of Rosenbloom.

**Theorem 3.3.11** (Dominated criteria)**.** *Let $f \in \mathcal{K}$ be a Gaussian power series and assume that*

- *we have a cut function $h : [0, R) \to (0, \pi]$ and*
- *we have a function $\Phi$ in $[0, +\infty)$ positive and integrable such that:*

- **[Major H]**:

$$\frac{|f(te^{i\theta})|}{f(t)} \leq \Phi(|\theta|\sigma_f(t)), \quad \textit{for any } |\theta| \leq h(t) \textit{ and } t \in [0, R).$$

*In terms of characteristic functions, this condition is equivalent to*

(3.3.7) $$\left|\mathbf{E}(e^{i\theta\check{X}_t})\right| \le \Phi(|\theta|), \quad \textit{for any } |\theta| \le h(t)\sigma_f(t) \textit{ and } t \in (0,R).$$

- **[Minor H]**:

$$\lim_{t\uparrow R} \sigma_f(t) \sup_{h(t)\le|\theta|\le\pi} \frac{|f(te^{i\theta})|}{f(t)} = 0.$$

*In terms of characteristic functions, this condition is equivalent to*

(3.3.8) $$\lim_{t\uparrow R} \sigma_f(t) \sup_{h(t)\le|\theta|\le\pi} |\mathbf{E}(e^{i\theta\check{X}_t})| = 0.$$

- **[Variance Condition]**:

(3.3.9) $$\lim_{t\uparrow R} h(t)\sigma(t) = +\infty\,.$$

*then $f$ is strongly Gaussian.*

**Remark 3.3.12.** The previous result is also true if the Major H condition holds for large values of $t$, that is, if there exists $0 < T < R$ such that the Major H condition hold for $t > T$. ⊠

The Minor H condition is the same than the minor arc condition for the power series in the Hayman class. Hayman functions satisfy, automatically, the Major H condition. Combining this with Lemma 3.3.3 we have that Hayman functions verify the previous criteria.

For the Major H condition on the previous criteria, the function $\Phi(\theta)$ is usually given by $\Phi(\theta) = e^{-\delta\theta^2}$, with $\delta > 0$, or $\Phi(\theta) = e^{-\delta\theta^2+\omega|\theta|^p}$, with $\delta, \omega > 0$ and $p < 2$.

The proof of the previous criteria is analogous to that of Theorem 3.3.7.

*Proof.* Define

$$I_t = \int_{|\theta|\le\pi\sigma_f(t)} |\mathbf{E}(e^{i\theta\check{X}_t}) - e^{-\theta^2/2})|d\theta.$$

We need to prove that $\lim_{t\uparrow R} I_t = 0$. We divide the integral $I_t$ into two integrals, the first integral $J_t$ on the major arc, and the second integral $K_t$ on the minor arc.

In $J_t$ the integrand is bounded by $\Phi(|\theta|) + e^{-\theta^2/2}$, therefore using the Gaussinity of $f$, and applying dominated convergence, we conclude that $\lim_{t\uparrow R} J_t = 0$.

For the integral $K_t$ we have

$$K_t \le \sigma_f(t) \sup_{\sigma_f(t)h(t)\le|\theta|\le\pi\sigma_f(t)} |\mathbf{E}(e^{i\theta\check{X}_t})| + \int_{|\theta|\ge h(t)\sigma_f(t)} e^{-\theta^2/2}d\theta.$$

The first summand at the right-hand side of the previous inequality tends to 0, as $t \uparrow R$, using the hypothesis on the minor arc, while the second summand tends to 0 because of the the variance condition; recall that $\lim_{t\uparrow R} h(t)\sigma_f(t) = +\infty$. □

At this point we collect the following chain of implications:

$$\text{Hayman class} \Rightarrow \text{dominated criteria} \Rightarrow \text{strongly Gaussian}$$

We also have that strongly Gaussian power series are clans and Gaussian power series. This chain of implications give criteria to prove that certain power series is not strongly Gaussian or, even further, that is not in the Hayman class.

## 3.4 Uniformly Gaussian power series

Here we introduce the concept of uniformly Gaussian power series. The main difference with the previous concept of Gaussianity is that now we have two parameters which can escape to infinity.

We will use the notation $[n\to\infty \vee t\uparrow R]$ to mean that $1 \le n \to \infty$ or (inclusive) for some $0 < t_0 \le t \uparrow R$. The restriction $t \ge t_0 > 0$ is there to exclude the possibility of $t = 0$ and $n \to \infty$, simultaneously.

Let $f(z) = \sum_{n=0}^{\infty} b_n z^n$ be a power series in $\mathcal{K}$ with radius of convergence $R > 0$ and let $(X_t)_{t\in[0,R)}$ be its Khinchin family.

A power series $f \in \mathcal{K}$ and its Khinchin family $(X_t)_{t\in[0,R)}$ are called *uniformly Gaussian* if

$$\lim_{[n\to\infty \vee t\uparrow R]} \mathbf{E}\left(e^{\imath\theta \breve{X}_t/\sqrt{n}}\right)^n = e^{-\theta^2/2}, \qquad \text{for any } \theta \in \mathbb{R}.$$

Recall that by $[n\to\infty \vee t\uparrow R]$ we mean that $1 \le n \to \infty$ or (inclusive) $0 < t_0 \le t \uparrow R$.

We verify now that the exponential $f(z) = e^z$ is uniformly Gaussian. We have $\sigma_f(t) = \sqrt{t}$, for $t > 0$ and so $\sigma_f(t)\sqrt{n} = \sigma_f(nt)$, for $t > 0$ and $n \ge 1$. For the characteristic function of the normalized variable $\breve{X}_t$ we have that

$$\mathbf{E}\big(e^{\imath\theta\breve{X}_t}\big) = \exp\left(t\left(e^{\imath\theta/\sqrt{t}} - 1 - \imath\theta/\sqrt{t}\right)\right),$$

and thus that

$$\mathbf{E}\big(e^{\imath\theta\breve{X}_t/\sqrt{n}}\big)^n = \mathbf{E}\big(e^{\imath\theta\breve{X}_{nt}}\big),$$

therefore, using that $f(z) = e^z$ is Gaussian, we obtain that

$$\lim_{[n\to\infty \vee t\uparrow R]} \mathbf{E}\big(e^{\imath\theta\breve{X}_t/\sqrt{n}}\big)^n = \lim_{[n\to\infty \vee t\uparrow R]} \mathbf{E}\big(e^{\imath\theta\breve{X}_{nt}}\big) = e^{-\theta^2/2}, \qquad \text{for any } \theta \in \mathbb{R}.$$

## 3.5 Uniformly strongly Gaussian families

Let $f(z) = \sum_{n=0}^{\infty} b_n z^n$ be a power series in $\mathcal{K}$ with radius of convergence $R > 0$ and let $(X_t)_{t\in[0,R)}$ be its Khinchin family.

We have encountered two integral convergence results: the integral form of the Local Central Limit, see (†) below, and also Section 6.1 in Chapter 6, and the notion of strongly Gaussian power series, definition (3.2.1).

• Assume that $Q_f = \gcd\{n \geq 1 : b_n > 0\} = 1$. *For each fixed* $t \in (0, R)$, the normalized variable $\breve{X}_t$ is a lattice random variable with gauge function $1/\sigma_f(t)$, since $Q_f = 1$. Because of local central limit theorem for lattice variables then we have that

$$(\dagger) \qquad \lim_{n\to\infty} \int\limits_{|\theta|\leq \pi\sigma_f(t)\sqrt{n}} \left|\mathbf{E}\Big(e^{\imath\theta\breve{X}_t/\sqrt{n}}\Big)^n - e^{-\theta^2/2}\right| d\theta = 0\,.$$

• If the Khinchin family $(X_t)_{t\in[0,R)}$ is strongly Gaussian, then we have *for each fixed* $n \geq 1$, that

$$(\ddagger) \qquad \lim_{t\uparrow R} \int\limits_{|\theta|\leq \pi\sigma_f(t)\sqrt{n}} \left|\mathbf{E}\Big(e^{\imath\theta\breve{X}_t/\sqrt{n}}\Big)^n - e^{-\theta^2/2}\right| d\theta = 0\,.$$

This fact follows from the bound

$$\begin{aligned}\int\limits_{|\theta|\leq \pi\sigma_f(t)\sqrt{n}} \left|\mathbf{E}\Big(e^{\imath\theta\breve{X}_t/\sqrt{n}}\Big)^n - e^{-\theta^2/2}\right| d\theta &= \sqrt{n} \int\limits_{|\varphi|\leq \pi\sigma_f(t)} \left|\mathbf{E}\Big(e^{\imath\varphi\breve{X}_t}\Big)^n - e^{-\varphi^2 n/2}\right| d\varphi \\ &\leq n^{3/2} \int\limits_{|\varphi|\leq \pi\sigma_f(t)} \left|\mathbf{E}\Big(e^{\imath\varphi\breve{X}_t}\Big) - e^{-\varphi^2/2}\right| d\varphi\,,\end{aligned}$$

where, after the change of variables $\theta = \varphi\sqrt{n}$, we have used that for complex numbers $z, w$ such that $|z|, |w| \leq 1$ we have that $|z^n - w^n| \leq n|z - w|$.

Thus, for a Gaussian power series the integral

$$I_n(t) = \int\limits_{|\theta|\leq \pi\sigma_f(t)\sqrt{n}} \left|\mathbf{E}\Big(e^{\imath\theta\breve{X}_t/\sqrt{n}}\Big)^n - e^{-\theta^2/2}\right| d\theta$$

converges to 0 as $n \to \infty$ with $t$ fixed and as $t \uparrow R$ with $n$ fixed.

Power series $f$ in $\mathcal{K}$ with Khinchin family $(X_t)_{t\in[0,R)}$ are power series for which (†) and (‡) hold *simultaneously* in the sense that the involved integrals converges to 0, as $n \to \infty$ or $t \uparrow R$.

**Definition 2.** *A power series $f \in \mathcal{K}$ and its Khinchin family $(X_t)_{t\in[0,R)}$ are called uniformly strongly Gaussian if the following two conditions are satisfied:*

$$a) \quad \lim_{t\uparrow R} \sigma_f(t) = \infty \quad \textit{and} \quad b) \quad \lim_{[n\to\infty\,\vee\, t\uparrow R]} \int\limits_{|\theta|\leq \pi\sigma_f(t)\sqrt{n}} \left|\mathbf{E}\Big(e^{\imath\theta\breve{X}_t/\sqrt{n}}\Big)^n - e^{-\theta^2/2}\right| d\theta = 0\,.$$

By $[n\to\infty \vee t\uparrow R]$ we mean that $1 \leq n \to \infty$ or (inclusive) $0 \leq t_0 \leq t \uparrow R$. The restriction $t \geq t_0 > 0$ is there to exclude the possibility of $t = 0$ and $n \to \infty$, simultaneously.

By fixing $n = 1$ and letting $t \uparrow R$, we observe that uniformly strongly Gaussian power series are strongly Gaussian. In particular, if $f$ is a uniformly strongly Gaussian power series then

$$M_f = \lim_{t \uparrow R} m_f(t) = \infty\,.$$

Moreover, the coefficients $b_n$ of $f$ satisfy Hayman's asymptotic formula and in particular $Q_f = 1$.

Let us verify that the exponential $f(z) = e^z$ is uniformly strongly Gaussian. Again we have, see Section 3.4 above, $\sigma_f(t) = \sqrt{t}$, for $t > 0$ and so $\sigma_f(t)\sqrt{n} = \sigma_f(nt)$, for $t > 0$ and $n \geq 1$. For the characteristic function of the normalized variable $\breve{X}_t$ we have that

$$\mathbf{E}\big(e^{\imath\theta \breve{X}_t}\big) = \exp\left(t\left(e^{\imath\theta/\sqrt{t}} - 1 - \imath\theta/\sqrt{t}\right)\right),$$

and thus that

$$\mathbf{E}\big(e^{\imath\theta \breve{X}_t/\sqrt{n}}\big)^n = \mathbf{E}\big(e^{\imath\theta \breve{X}_{nt}}\big)\,.$$

Therefore,

$$\int\limits_{|\theta| \leq \pi\sigma_f(t)\sqrt{n}} \left|\mathbf{E}\Big(e^{\imath\theta \breve{X}_t/\sqrt{n}}\Big)^n - e^{-\theta^2/2}\right| d\theta = \int\limits_{|\theta| \leq \pi\sigma_f(nt)} \left|\mathbf{E}\Big(e^{\imath\theta \breve{X}_{nt}}\Big) - e^{-\theta^2/2}\right| d\theta\,.$$

Since $e^z$ is strongly Gaussian, this integral tends to 0 as $(nt) \to \infty$, and thus as $[n \to \infty \,\vee\, t \uparrow +\infty]$. Observe that in this case $R = +\infty$.

## 3.6 Uniformly Hayman power series

The notion of uniformly Hayman power series which we are about to introduce generalizes that of Hayman power series, see Subsection 3.3. We will verify shortly that uniformly Hayman power series are uniformly strongly Gaussian, much like power series in the Hayman class are strongly gaussian. In Theorem 4.3.22, see Chapter 4 provides us with an ample class of power series which are uniformly Hayman.

**Definition 3.** *Let $f(z)$ be a power series in $\mathcal{K}$ with radius of convergence $R > 0$ and let $(X_t)_{t \in [0,R)}$ be its Khinchin family.*

*We say that $f(z)$ and $(X_t)_{t \in [0,R)}$ are uniformly Hayman if for each $n \geq 1$ and $t \in (0, R)$ there exists $h(n,t) \in (0, \pi)$ (called cuts) such that the following requirements are satisfied:*

$$\sup_{|\theta| \leq h(n,t)\sigma_f(t)\sqrt{n}} \left|\mathbf{E}(e^{\imath\theta \breve{X}_t/\sqrt{n}})^n e^{\theta^2/2} - 1\right| \to 0, \quad \text{as } [n \to \infty \,\vee\, t \uparrow R]\,, \tag{3.6.1}$$

$$\sqrt{n}\sigma_f(t) \sup_{h(n,t)\sigma_f(t)\sqrt{n} \leq |\theta| \leq \pi\sigma_f(t)\sqrt{n}} \left|\mathbf{E}(e^{\imath\theta \breve{X}_t/\sqrt{n}})^n\right| \to 0, \quad \text{as } [n \to \infty \,\vee\, t \uparrow R]\,, \tag{3.6.2}$$

$$\lim_{t \to R} \sigma_f(t) = \infty\,. \tag{3.6.3}$$

Condition (3.6.2) may be written equivalently as

$$\sqrt{n}\sigma_f(t) \sup_{h(n,t)\le|\theta|\le\pi} \left|\mathbf{E}(e^{\imath\theta X_t})^n\right| \to 0, \quad \text{as } [n\to\infty \vee t\uparrow R]\,. \tag{3.6.4}$$

As announced:

**Theorem 3.6.1.** *Uniformly Hayman power series are uniformly strongly Gaussian.*

*Proof.* Denote $\theta(n,t) = h(n,t)\sigma_f(t)\sqrt{n}$. First we show, that

$$\theta(n,t) \to \infty\,, \quad \text{as } [n\to\infty \vee t\uparrow R]\,. \tag{3.6.5}$$

Abbreviate $\widehat{\theta} = \theta(n,t)$. By (3.6.1), we have that

$$\mathbf{E}(e^{\imath\widehat{\theta}\breve{X}_t/\sqrt{n}})e^{\widehat{\theta}^2/2} \to 1\,, \quad \text{as } [n\to\infty \vee t\uparrow R]\,,$$

while, from (3.6.2) we obtain that

$$\sqrt{n}\sigma_f(t)\mathbf{E}\big(e^{\imath\widehat{\theta}\breve{X}_t/\sqrt{n}}\big)^n \to 0\,, \quad \text{as } [n\to\infty \vee t\uparrow R]\,.$$

From these two limits we deduce that

$$\sqrt{n}\sigma_f(t)e^{-\widehat{\theta}^2/2} \to 0\,, \quad \text{as } [n\to\infty \vee t\uparrow R]\,,$$

and thus, since $\sqrt{n}\sigma_f(t) \to \infty$ as $[n\to\infty \vee t\uparrow R]$ we deduce that

$$\widehat{\theta} \to \infty\,, \quad \text{as } [n\to\infty \vee t\uparrow R]\,.$$

Denote by $A(n,t), B(n,t)$, respectively, the supremum in (3.6.1) and (3.6.2).

We bound

$$\begin{aligned}
&\int\limits_{|\theta|\le h(n,t)\sigma_f(t)\sqrt{n}} \left|\mathbf{E}(e^{\imath\theta\breve{X}_t/\sqrt{n}})^n - e^{-\theta^2/2}\right| d\theta \\
&= \int\limits_{|\theta|\le h(n,t)\sigma_f(t)\sqrt{n}} \left|\mathbf{E}(e^{\imath\theta\breve{X}_t/\sqrt{n}})^n e^{\theta^2/2} - 1\right| e^{-\theta^2/2}\, d\theta \\
&\le A(n,t)\sqrt{2\pi}\,,
\end{aligned}$$

and

$$\begin{aligned}
&\int\limits_{h(n,t)\sigma_f(t)\sqrt{n}\le|\theta|\le\pi\sigma_f(t)\sqrt{n}} \left|\mathbf{E}(e^{\imath\theta\breve{X}_t/\sqrt{n}}) - e^{-\theta^2/2}\right| d\theta \\
&\le 2\pi\sigma_f(t)\sqrt{n}B(n,t) + \int_{|\theta|\ge h(n,t)\sigma_f(t)\sqrt{n}} e^{-\theta^2/2}\, d\theta\,.
\end{aligned}$$

These two bounds and conditions (3.6.1) and (3.6.2) combined with (3.6.5) give the result. □

Chapter 4

# Set Constructions and Exponentials

## Contents

---

The set constructions of combinatorial classes fits nicely within the framework of Khinchin families. As a general reference for combinatorial classes and operations with them, we strongly suggest:

- Chapter II in Flajolet, P. and Sedgewick, R. Analytic combinatorics. See [32].
- Arriata, R., Barbour, A. D. and Tavaré, S. Logarithmic Combinatorial Structures: a Probabilistic Approach. See [4].

Most of the (ordinary or exponential) generating functions of sets of combinatorial classes are exponentials of power series with non-negative coefficients, the object of interest of this chapter.

We will verify by means of Theorems 4.3.15 and 4.3.17 that most of the set constructions give rise to generating functions which are in the Hayman class, so that there are asymptotic formulas for the coefficients given by Theorems 3.2.5 and 3.2.8. This chapter is mainly based on the papers:

- Maciá, V.J. et al. Khinchin families and Hayman class. *Comput. Methods Funct. Theory* **21** (2021), 851–904, see [17],
- Maciá, V.J. et al. Khinchin families, set constructions, partitions and exponentials, *Mediterr. J. Math.*, 21, 39 (2024), see [18].

## 4.1 Set construction of combinatorial classes

Exponentials of power series with positive coefficients are very relevant, in particular, in Combinatorics since they codify (most) generating functions of the set construction, which includes among them generating functions of partitions of many sorts. Here we discuss the set construction of combinatorial classes.

### 4.1.1 Labeled combinatorial classes and Sets

If $g(z) = \sum_{n=1}^{\infty}(b_n/n!)$ is the exponential generating function EGF of a labeled combinatorial class $\mathcal{G}$ (no object of size 0), then $f(z) = e^{g(z)} = \sum_{n=0}^{\infty}(a_n/n!)z^n$ is the EGF of the labeled class of sets formed with the objects of $\mathcal{G}$, termed assemblies in [4] .

Next, we consider first sets of the basics classes: sets, lists and cycles, and the sets of rooted trees and functions.

### A. Sets of Sets

Let $\mathcal{G}$ be the labeled class of non-empty sets: $b_n = 1$, for each $n \geq 1$, and $b_0 = 0$. Its EGF is $g(z) = e^z - 1$. Then

$$f(z) = \exp(e^z - 1) = \sum_{n=0}^{\infty} \frac{\mathcal{B}_n}{n!} z^n$$

is the EGF of sets of sets, or, equivalently, of partitions of sets, see [32, p. 107]. Here $\mathcal{B}_n$ is the $n$-th Bell number, which counts the number of partitions of the set $\{1, 2, \ldots, n\}$. In this case $R = \infty$ and the mean and variance functions are $m_f(t) = te^t$ and $\sigma_f^2(t) = t(t+1)e^t$, see also Subsection 1.1.2.

The characteristic function of the Khinchin family of $f$, that we denote by $(X_t)$, is given by

$$\mathbf{E}(e^{i\theta X_t}) = \exp(e^{te^{i\theta}} - e^t), \quad \text{for any } \theta \in \mathbb{R} \text{ and } t \geq 0,$$

and therefore

$$\mathbf{E}(e^{i\theta \breve{X}_t}) = \exp(e^{te^{i\theta e^{-t^2}/\sqrt{t(t+1)}}} - e^t - i\theta\sqrt{t/(t+1)e^{t/2}}) \tag{4.1.1}$$

**Remark 4.1.1** (Pointed sets)**.** The class $\mathcal{G}$ of pointed sets (i.e. sets with one marked element) has EGF $g(z) = ze^z$. The class of sets of pointed sets is isomorphic to the class of idempotent maps; its EGF $f$ is given by $f(z) = e^{ze^z}$. See [32, p.131]. ⊠

### B. Sets of Lists

For the labeled class $\mathcal{G}$ of (non-empty) lists, the function $g$ is just $g(z) = z/(1-z)$, and therefore $f(z) = e^{z/(1-z)}$ is the EGF of the sets of lists, the so called fragmented permutations. See [32, p.125]

### C. Sets of Cycles

The labeled class $\mathcal{G}$ of (non-empty) cycles has $g(z) = \ln(1/(1-z))$ as EGF. The function $f(z) = e^{g(z)} = 1/(1-z)$ is the EGF of the sets of cycles, or, in other terms of the permutations.

The length of the cycles could be restricted. Thus for integer $k \geq 1$,

$$\exp\Big(\sum_{n\leq k} z^n/n\Big) \quad \text{or} \quad \exp\Big(\sum_{n\geq k} z^n/n\Big)$$

are the EGF of permutations such that all the cycles in their cycle decomposition have length at most $k$, or at least $k$, respectively. We may also consider, for $k \geq 1$,

$$\exp\Big(\sum_{d\geq 1, d|k} z^d/d\Big)$$

which is the EGF of the permutations $\sigma$ such that $\sigma^k$ is the identity.

**D. Sets of Trees and Sets of Functions**

Here we study the generating functions of the sets of trees and the sets of functions.

**D.1. Sets of trees (forests)** The class $\mathcal{G}$ of rooted (labeled) trees has EGF $g(z) = \sum_{n=1}^{\infty}(n^{n-1}/n!)z^n$, see [32], Section II.5.1; this is Cayley's Theorem. The class of forests (sets) of rooted (labeled) trees has then EGF $f = e^g$. See also subsection 1.1.2-C.

Cayley's Theorem also shows that the EGF of the clas of unrooted (labeled) trees is $g(z) = \sum_{n=1}^{\infty}(n^{n-2}/n!)z^n$.

**D.2. Sets of functions** The class $\mathcal{G}$ of functions has EGF $g(z) = \sum_{n=1}^{\infty}(n^n/n!)z^n$. See [32], Section II.5.2. The class of sets of functions then has EGF $f = e^g$. See also Subsection 1.1.2-C.

### 4.1.2 Unlabeled combinatorial classes and Sets

We split the discussion into multisets and (proper) sets. See Section I.2.2 in [32] and [4] as general references for the constructions of unlabeled combinatorial classes.

**A. Multisets of Unlabeled Combinatorial Classes**

If $C(z) = \sum_{n=1}^{\infty} c_n z^n$ is the ordinary generating function OGF of a combinatorial (unlabeled) class $\mathcal{G}$ (no object of size 0), then

$$f(z) = \prod_{j=1}^{\infty} \frac{1}{(1-z^j)^{c_j}} = \exp\left(\sum_{n=1}^{\infty} C(z^n)/n\right),$$

is the OGF of the class of sets formed with the objects of $\mathcal{G}$, termed multisets in [4].

We may write $f(z) = e^{g(z)}$, where $g$ is the power series

$$g(z) = \sum_{j,k\geq 1} c_j \frac{z^{kj}}{k} = \sum_{m=1}^{\infty}\left(\sum_{j|m} jc_j\right)\frac{z^m}{m} \triangleq \sum_{m=1}^{\infty} b_m z^m.$$

Observe that the power series $g$ has non-negative coefficients:

$$b_n = \frac{1}{m}\sum_{j|m} jc_j, \qquad \text{for } m \geq 1.$$

An instance of this situation is given by the OGF $P$ of the partitions of integers. Recall that $P$ is given by the infinite product

$$P(z) = \prod_{j=1}^{\infty} \frac{1}{1-z^j} = \sum_{n=0}^{\infty} p(n)z^n, \qquad \text{for } z \in \mathbb{D},$$

the power series $P$ is in $\mathcal{K}$. Here $p(n)$ is the number of partitions of the integer $n \geq 1$, and $p(0) = 1$. The mean and variance functions of its Khinchin family $(X_t)$ are not so direct, see equation (1.2.7).

In this instance, the function $C$ is given by $C(z) = z/(1-z)$, as $c_n = 1$, for $n \geq 1$ (one object of size $n$, for each $n \geq 1$). And thus $g$ has coefficients $b_m = \sigma_1(m)/m$, where $\sigma_1(m)$ is the sum of the divisors of the integer $m$, that is,

$$\sigma_1(m) = \sum_{j|m} j.$$

For general $C(z) = \sum_{n=1}^{\infty} c_n z^n$ as above, the corresponding $f$ is the OGF of the colored partitions; $c_j$ different colors for part $j$.

### B. Sets of Unlabeled Combinatorial Classes

If again we write $C(z) = \sum_{n=1}^{\infty} c_n z^n$ for the ordinary generating function OGF of a combinatorial (unlabeled) class $\mathcal{G}$ (no objects of size 0), then

$$f(z) = \prod_{j=1}^{\infty} (1+z^j)^{c_j} = \exp\left(\sum_{n=1}^{\infty} (-1)^{n+1} C(z^n)/n\right)$$

is the OGF of the class of (proper) sets formed with the objects of $\mathcal{G}$, termed *selections* in [4].

We may write $f(z) = e^{g(z)}$, where $g$ is the power series

$$g(z) = \sum_{j,k\geq 1} c_j \frac{(-1)^{k+1} z^{kj}}{k} = \sum_{m=1}^{\infty} \left(\sum_{jk=m} j c_j (-1)^{k+1}\right) \frac{z^m}{m}.$$

In general, the power series $g$ could have negative coefficients, this is the case, for instance, for $f(z) = (1+z)^5(1+z^2)$.

But for the sequence of coefficients $c_j = j^{c-1}$, for $j \geq 1$ and $c$ is an integer $c \geq 1$, the coefficients of $g$ are non-negative; in fact, the $m$-th coefficient of $g$ is

$$\frac{1}{m}\sigma_c^{\text{odd}}\omega(m),$$

where

$$\sigma_c^{\text{odd}}(m) = \sum_{j|m,\ j \text{ odd}} j^c \quad \text{and} \quad \omega(m) = \frac{2^c-2}{2^c-1} 2^{\chi(m)c} + \frac{1}{2^c-1}, \quad \text{for any } m \geq 1,$$

here $\chi(m)$ is the highest integer exponent so that $2^{\chi(m)}|m$, and $\sigma_c^{\text{odd}}(m)$ is the sum of the $c$-th powers of the odd divisors of $m$. This is so because of the identity

$$\sum_{jk=m} j^c (-1)^{k+1} = \left(\sum_{j|m,\ j\text{odd}} j^c\right) \omega(m), \quad \text{for } m \geq 1. \tag{4.1.2}$$

To verify (4.1.2), observe first that as functions of $m$, both summations in (4.1.2) are multiplicative. Write $m$ as $m = 2^{\chi(m)}s$, with $s$ odd and observe that for $m = 2^r$, with $r \geq 1$, the summation on the left is $\omega(m)$ and the summation on the right is 1, while for $m = s$ odd, the two sums coincide and $\omega(s) = 1$.

We may also write

$$\sum_{jk=m} (-1)^{k+1} j^c = \sigma_c(m) - 2\sigma_c(m/2), \qquad \text{for } m \geq 1,$$

with the understanding that $\sigma_c(m/2) = 0$ if $m$ is odd. For the identity (4.1.2) and a variety of relations among a number of diverse sums of divisors, we refer to [37].

More generally, if the coefficient $c_j$ of $C$ is given by $c_j = R(j)$, for $j \geq 1$, where $R$ is a polynomial with non-negative integer coefficients, then the coefficients of $g$ are non-negative.

An instance of this situation is given by the OGF of partitions with distinct parts, given by

$$Q(z) = \prod_{j=1}^{\infty}(1 + z^j) = \sum_{n=0}^{\infty} q(n) z^n, \qquad \text{for any } z \in \mathbb{D}\,.$$

This power series is in $\mathcal{K}$. Here $q(n)$ is the number of partitions into distinct parts of the integer $n \geq 1$.

In this particular instance, the power series $C$ is given by $C(z) = z/(1 - z)$, as $c_n = 1$, for $n \geq 1$, and the power series $g$ is

$$g(z) = \sum_{m=1}^{\infty} \frac{\sigma_1(m) - 2\sigma_1(m/2)}{m} z^m.$$

## 4.2 Khinchin families of exponentials $f = e^g$, with $g \in \mathcal{P}$

We study now the Khinchin family associated to an exponential power series, that is, $f = e^g$, with $g$ a power series with non-negative coefficients and radius of convergence $R > 0$.

We define the class $\mathcal{P}$ as the class of non-constant power series $g$ with non-negative coefficients and having radius of convergence $R > 0$.

### 4.2.1 Mean and variance

For the mean and variance of a power series $f = e^g$, with $g \in \mathcal{P}$ we have the following lemma.

**Lemma 4.2.1.** *Let $f = e^g$ with $g \in \mathcal{P}$ a power series with radius of convergence $R > 0$, then*

$$m_f(t) = t g'(t) \qquad \text{and} \qquad \sigma_f^2(t) = t g'(t) + t^2 g''(t),$$

*for any $t \in (0, R)$.*

*Proof.* For any $g \in \mathcal{P}$ we have that $f = e^g \in \mathcal{K}$. Applying the general formula for the mean and the variance of a power series $f \in \mathcal{K}$, see equation (1.2.1), we have

$$m_f(t) = \frac{tf'(t)}{f(t)}, \qquad \text{and} \qquad \sigma_f^2(t) = tm_f'(t), \qquad \text{for any } t \in (0, R)\,,$$

then, for our particular case, i.e. $f = e^g$ with $g \in \mathcal{P}$, we have

$$m_f(t) = tg'(t), \qquad \text{for any } t \in (0, R),$$

and therefore

$$\sigma_f^2(t) = tg'(t) + t^2 g''(t) \qquad \text{for any } t \in (0, R).$$

□

**Range of the mean and the variance**

The mean and variance functions of $f = e^g$, with $g \in \mathcal{P}$, are very special. As we have seen before in Subsection 1.2.5 the limit $\lim_{t\uparrow R} \sigma_f^2(t)$ does not always exists for general $f \in \mathcal{K}$.

In this case, because of the expressions for the mean and the variance functions given by Lemma 4.2.1, both functions are monotone increasing functions (this was the case for the mean in the general case, but not for the variance) and therefore the limits

$$M_f = \lim_{t\uparrow R} m_f(t) \qquad \text{and} \qquad \Sigma_f = \lim_{t\uparrow R} \sigma_f^2(t)$$

always exist.

For the mean function we have that $\lim_{t\uparrow R} m_f(t) = +\infty$ if and only if $\lim_{t\uparrow R} tg'(t) = +\infty$. This is a particular case of Lemma 1.2.1.

For the variance we have the following lemma.

**Lemma 4.2.2** (Range of the variance)**.** *Let $f \in \mathcal{K}$ be such that $f = e^g$, where $g \in \mathcal{P}$ is a power series with radius of convergence $R > 0$. We distinguish two cases:*

1. *If $R < +\infty$, then*

$$\lim_{t\uparrow R} \sigma_f^2(t) = +\infty \qquad \text{if and only if} \qquad \lim_{t\uparrow R} g''(t) = +\infty.$$

2. *If $R = +\infty$, then*

$$\lim_{t\to+\infty} \sigma_f^2(t) = +\infty.$$

*Proof.* Lemma 4.2.1 gives that

$$\sigma_f^2(t) = tg'(t) + t^2 g''(t), \qquad \text{for any } t \in (0, R). \tag{4.2.1}$$

Assume first that $R < +\infty$ and also that $\lim_{t\uparrow R} g''(t) < +\infty$, comparing coefficients we find the inequality

$$tg'(t) \leq t^2 g''(t) + g'(0)t, \qquad \text{for any } t \in (0, R).$$

The previous inequality combined with the fact that $R < +\infty$ gives that $\lim_{t\uparrow R} g'(t) < +\infty$, and therefore we conclude that $\lim_{t\uparrow R} \sigma_f^2(t) < +\infty$. Assume now that $\lim_{t\uparrow R} g''(t) = +\infty$, then, equation (4.2.1) gives that $\lim_{t\uparrow R} \sigma_f^2(t) = +\infty$.

For entire power series $g \in \mathcal{P}$ we always have $\lim_{t\to\infty} tg'(t) = +\infty$, and therefore we conclude from (4.2.1) that $\lim_{t\to\infty} \sigma_f^2(t) = +\infty$. □

The previous lemma is a particular case of Corollary 4.2.5 below.

The mean and the variance functions of a exponential power series also verify the following relation.

**Lemma 4.2.3.** *Let $f = e^g$, with $g \in \mathcal{P}$ be an exponential power series with radius of convergence $R > 0$, then*

$$m_f(t) \leq \sigma_f^2(t), \qquad \textit{for any } t \in [0, R).$$

*Proof.* The expressions for the mean and the variance given in Lemma 4.2.1 give that

$$m_f(t) = tg'(t) \leq \sigma_f^2(t) = tg'(t) + t^2 g''(t), \qquad \text{for any } t \in [0, R).$$

Here we use that $g \in \mathcal{P}$ has non-negative coefficients.

□

### 4.2.2 Fulcrum

Here we discuss the fulcrum of an exponential power series. We adapt some of the results from Chapter 3, in particular Section 3.1, to the case of exponential power series $f = e^g$, with $g \in \mathcal{K}$.

#### A. Fulcrum of a general exponential power series

For any $f \in \mathcal{K}$, the holomorphic function $f$ does not vanish on the real interval $[0, R)$, and so, it does not vanish in a simply connected region containing that interval. We may consider $\ln(f)$, a holomorphic branch of the logarithm of $f$ which is real on $[0, R)$, and the function $F$, the fulcrum of $f$, see Chapter 1, in particular Subsection 1.3.4, is defined and holomorphic in a simply connected region containing $(-\infty, \ln(R))$. Recall that

$$F(z) = \ln(f(e^z)).$$

If $f$ does not vanish anywhere in the disk $\mathbb{D}(0,R)$, then the fulcrum $F$ of $f$ is defined on the whole half-plane $\{z \in \mathbb{C} : \Re(z) < \ln(R)\}$. In this case $f(z) = e^{g(z)}$, where $g$ is a holomorphic function in $\mathbb{D}(0,R)$ and $g(0) \in \mathbb{R}$, and therefore the fulcrum $F$ of $f$ is given by

$$F(z) = g(e^z), \quad \text{for } \Re(z) < \ln(R).$$

The mean and the variance may be expressed in terms of its fulcrum $F$ as

$$m_f(t) = F'(s) \quad \text{and} \quad \sigma_f^2(t) = F''(s), \quad \text{for } s < \ln(R) \text{ and } t = e^s.$$

These identities give an alternative proof of Lemma 4.2.1.

## B. Fulcrum of $e^g$ with $g \in \mathcal{P}$

Let $f$ be the power series $f = e^g$ with $g \in \mathcal{P}$ a power series with radius of convergence $R > 0$. This is a particular case of the case presented in A.

The interest in power series of the form $f = e^g$ with $g \in \mathcal{P}$ comes from the fact that these exponentials, as explained at the beginning of this chapter, appear in many combinatorial applications.

The fulcrum $F$ of power series $f = e^g$, with $g \in \mathcal{P}$, verifies that

$$|F^{(k)}(s + i\theta)| \le F^{(k)}(s), \qquad \text{for any } \theta \in \mathbb{R} \text{ and } s < \ln(R),$$

in fact

$$\sup_{\theta \in \mathbb{R}} |F^{(k)}(s + i\theta)| = F^{(k)}(s), \qquad \text{for any } s < \ln(R).$$

Our first lemma expresses the derivatives of the fulcrum $F$ as finite linear combination of terms of the form $e^{js} g^{(j)}(s)$, where $j \ge 1$.

**Lemma 4.2.4.** *Let $f \in \mathcal{K}$ be such that $f = e^g$, with $g \in \mathcal{P}$ a power series with non-negative coefficients and radius $R > 0$, then, for any $k \ge 1$, there are integer constants $c_1, \dots, c_{k-1} \ge 0$ such that*

$$F^{(k)}(s) = e^{ks} g^{(k)}(e^s) + \sum_{j=1}^{k-1} c_j (e^{sj} g^{(j)}(e^s)), \quad \textit{for any } s < \ln(R),$$

*Proof.* We proceed by induction on $n \ge 1$. For $n = 1$ we have

$$F'(s) = e^s g'(e^s), \quad \text{for any } s < \ln(R).$$

Assume that for some $n \ge 1$, there are integer constants $c_1, \dots, c_{n-1} \ge 0$ such that

$$F^{(n)}(s) = e^{ns} g^{(n)}(e^s) + \sum_{j=1}^{n-1} c_j (e^{sj} g^{(j)}(e^s)), \quad \text{for any } s < \ln(R),$$

then, taking derivatives, we find that

$$F^{(n+1)}(s) = e^{(n+1)s}g^{(n+1)}(e^s) + ne^{ns}g^{(n)}(e^s) + \sum_{j=1}^{n-1} c_j(je^{js}g^{(j)}(e^s) + e^{(j+1)s}g^{(j+1)}(e^s))$$

and then we have proved that there are integer constants $\tilde{c}_1, \ldots, \tilde{c}_n \geq 0$ such that

$$F^{(n+1)}(s) = e^{(n+1)s}g^{(n+1)}(e^s) + \sum_{j=1}^{n} \tilde{c}_j(e^{sj}g^{(j)}(e^s)), \quad \text{for any } s < \ln(R).$$

□

From the previous lemma we obtain asymptotic comparisons for the fulcrum $F$ and any of its derivatives, as $s \uparrow \ln(R)$.

**Corollary 4.2.5.** *Let $f \in \mathcal{K}$ be such that $f = e^g$, with $g \in \mathcal{P}$ a non-constant power series with non-negative coefficients and radius $R > 0$, then:*

- *If $g$ is not a polynomial, then for any $k \geq 1$, there exists a constants $C > 0$ and a polynomial $P_{k-1}$ of degree $k-1$ such that*

$$e^{ks}g^{(k)}(e^s) \leq F^{(k)}(s) \leq Ce^{ks}g^{(k)}(e^s) + P_{k-1}(e^s), \quad \textit{for any } s < \ln(R),$$

  *and therefore*

$$F^{(k)}(s) \asymp e^{ks}g^{(k)}(e^s), \quad \textit{as } s \uparrow \ln(R).$$

- *If $g(z) = a_N z^n + \cdots + a_0$ is a polynomial of degree $N \geq 1$, then, for any $k \geq 1$, we have*

$$F^{(k)}(s) \asymp a_N e^{sN}, \quad \textit{as } s \to +\infty.$$

*Proof.* Assume first that $g(z) = a_N z^N + \cdots + a_0 \in \mathcal{P}$ is a polynomial of degree $N \geq 1$, then $F(s) = g(e^s) = a_N e^{Ns} + \cdots + a_0$ and therefore, for any $k \geq 1$, we have

$$F^{(k)}(e^s) \asymp a_N e^{Ns}, \quad \text{as } s \to \infty.$$

Now we study the case where $g \in \mathcal{P}$ not a polynomial. Applying Lemma 4.2.4, for any $k \geq 1$, there are integer constants $c_1, \ldots, c_{n-1} \ \geq 0$ such that

$$(\star) \quad F^{(k)}(s) = e^{ks}g^{(k)}(e^s) + \sum_{j=1}^{k-1} c_j(e^{sj}g^{(j)}(e^s)), \quad \text{for any } s < \ln(R),$$

and therefore, because $g$ has non-negative coefficients, we have

$$e^{ks}g^{(k)}(e^s) \leq F^{(k)}(s), \quad \text{for any } s < \ln(R). \tag{4.2.2}$$

Comparing coefficients, for any $k \geq 0$, we have the inequality

$$t^k g^{(k)}(t) \leq t^{k+1} g^{(k+1)}(t) + t^k (g^{(k)}(0)/k!), \quad \text{for any } t \in (0, R), \tag{4.2.3}$$

and therefore, applying inequality (4.2.3), iteratively, to each of the terms in the sum at the right-hand side of $(\star)$, we find that there exists a constant $C > 0$ and a polynomial $P_{k-1}(z)$ of degree $k-1$ such that

$$F^{(k)}(s) \leq C e^{ks} g^{(k)}(e^s) + P_{k-1}(e^s), \quad \text{for any } s < \ln(R).$$

Now using that $g$ is not a polynomial, we conclude that

$$F^{(k)}(s) \asymp e^{ks} g^{(k)}(e^s), \quad \text{as } s \uparrow \ln(R).$$

□

The next corollary follows from a direct application of Corollary 4.2.5. This corollary will be relevant, later on, when giving criteria for an exponential power series to be Gaussian. In fact the following corollary gives, by means of Theorem 3.1.19, that the exponential of a polynomial in $\mathcal{P}$ is always a Gaussian power series.

**Corollary 4.2.6.** *Let $f \in \mathcal{K}$ be such that $f = e^g$, with $g \in \mathcal{P}$ a power series with radius of convergence $R > 0$. Assume that $g \in \mathcal{P}$ is not a polynomial, then, for any $k \geq 3$, we have*

$$\frac{\sup_{\theta \in \mathbb{R}} |F^{(k)}(s + i\theta)|}{F''(s)^{k/2}} = \frac{F^{(k)}(s)}{F''(s)^{k/2}} \asymp \frac{g^{(k)}(e^s)}{g''(s)^{k/2}}, \quad \text{as } s \uparrow \ln(R).$$

*In the particular case that $g \in \mathcal{P}$ is a polynomial, then*

$$\lim_{s \to \infty} \frac{\sup_{\theta \in \mathbb{R}} |F^{(k)}(s + i\theta)|}{F''(s)^{k/2}} = \lim_{s \to +\infty} \frac{F^{(k)}(s)}{F''(s)^{k/2}} = 0, \quad \text{for any } k \geq 3.$$

*Proof.* This result follows combining that $g \in \mathcal{P}$ has non-negative Taylor coefficients with Corollary 4.2.5. □

### C. The fulcrum of $f = e^g$ and the moments of $g$

Let $f = e^g$ a power series, with $g(z) = \sum_{n=0}^{\infty} b_n z^n \in \mathcal{K}_s$ a power series with radius of convergence $R > 0$. Denote $(X_t)$, $(Y_t)$ and $(Z_t)$ the Khinchin families of $f, g$ and $e^z$, respectively.

In this case the fulcrum $F$ of $f = e^g$ is given by

$$F(z) = g(e^z) = \sum_{n=0}^{\infty} b_n e^{nz}, \quad \text{for any } \Re(z) < \ln(R) \text{ and } \Im(z) \in \mathbb{R},$$

then, if we denote $t = e^s$, for $s < \ln(R)$, and for any integer $k \geq 1$ we have

$$F^{(k)}(s) = \sum_{n=0}^{\infty} n^k b_n e^{sn} = g(t)\mathbf{E}(Y_t^k) = \mathbf{E}(Z_{g(t)})\mathbf{E}(Y_t^k), \quad \text{for any } s < \ln(R).$$

and therefore, for any $k \geq 1$:

$$\frac{F^{(k)}(s)}{F''(s)^{k/2}} \geq 0, \quad \text{for any } s < \ln(R). \tag{4.2.4}$$

In fact we have

$$(\star) \quad \frac{F^{(k)}(s)}{F''(s)^{k/2}} = \frac{\mathbf{E}(Y_t^k)}{g(t)^{k/2-1}\mathbf{E}(Y_t^2)^{k/2}}, \quad \text{for } t = e^s \text{ and } s < \ln(R).$$

Inequality (4.2.4) is a particularity of power series of the form $f = e^g$, with $g \in \mathcal{K}_s$. For any $k \geq 3$, this is a manifestation of the fact that the variance, and all its derivatives, are non-negative and therefore the derivatives of the variance, and the variance itself, are monotone increasing functions.

The inequality (4.2.4) is not true for general power series $f \in \mathcal{K}$: for the polynomial $f(z) = 1 + z \in \mathcal{K}$ we have $F(s) = \ln(1 + e^s)$ and

$$F'''(s) = -\frac{e^s(e^s - 1)}{(e^s + 1)^3} \leq 0 \quad \text{for any } s \geq 0.$$

Let $g \in \mathcal{K}_s$ be a power series with radius of convergence $R > 0$ which is not a polynomial. Combining $(\star)$ with Corollary 4.2.5, for any $k \geq 1$, we have that

$$\frac{F^{(k)}(s)}{F''(s)^{k/2}} = \frac{\mathbf{E}(Y_t^k)}{g(t)^{k/2-1}\mathbf{E}(Y_t^2)^{k/2}} \asymp \frac{g^{(k)}(t)}{g''(t)^{k/2}}, \quad \text{as } t \uparrow R.$$

Let $g(z) = a_N z^N + \cdots + a_0 \in \mathcal{K}_s$ a polynomial of degree $N \geq 1$, then, for any $k \geq 1$, we have

$$\frac{F^{(k)}(s)}{F''(s)^{k/2}} = \frac{\mathbf{E}(Y_t^k)}{g(t)^{k/2-1}\mathbf{E}(Y_t^2)^{k/2}} \asymp e^{Ns(1-k/2)}, \quad \text{as } s \to \infty,$$

a direct application of Theorem 3.1.4 gives that the exponential of a polynomial $g \in \mathcal{K}_s$ is Gaussian.

### 4.2.3 Power series comparisons

In the proof of Theorems 4.3.15 and 4.3.17, we will resort to the following asymptotics for a couple of real power series. We have already stated Lemma 4.2.7 and Proposition 4.2.9 in Chapter 2, but we will write them again in here for easier reference. See Chapter 2, in particular Subsection 2.3, for proofs.

**Lemma 4.2.7.** *For $\beta > 0$,*

$$\sum_{n=1}^{\infty} n^{\beta-1} t^n \sim \Gamma(\beta)\frac{1}{(1-t)^\beta}, \quad \text{as } t \uparrow 1. \tag{4.2.5}$$

**Remark 4.2.8** (Moments of geometric random variables)**.** The power series comparison of Proposition 4.2.7 translated into a comparison of moments of the Khinchin family $(X_t)_{t\in[0,1)}$ of $1/(1-z)$, i.e, the geometric variables. Namely, the following: for $\beta > 0$,

$$\mathbf{E}(X_t^\beta) \sim \Gamma(\beta+1)\mathbf{E}(X_t)^\beta, \qquad \text{as } t \uparrow 1.$$

For $\beta \in (-1, 0)$, the same result holds for the variable $\tilde{X}_t$, the $X_t$ conditioned at being positive, given by $\mathbf{P}(\tilde{X}_t = k) = t^{k-1}(1-t)$, for $k \geq 1$ and $t \in (0,1)$. ⊠

**Proposition 4.2.9.** *For $\beta \geq 0$, we have that*

$$\sum_{n=0}^{\infty} n^\beta \frac{t^n}{n!} \sim t^\beta e^t, \qquad \textit{and} \qquad \sum_{n=1}^{\infty} \frac{1}{n^\beta}\frac{t^n}{n!} \sim \frac{e^t}{t^\beta}, \qquad \textit{as } t \to \infty.$$

For negative exponents, we refer to [55] and to Remark 4.2.10 below.

**Remark 4.2.10** (Moments of Poisson random variables)**.** Proposition 4.2.9 is actually a statement about positive moments of the Poisson variable $X_t$, and also about negative moments of the conditioned Poisson variables $\tilde{X}_t$, as $t$ tends to $\infty$. The conditional Poisson variables $\tilde{X}_t$, that is the variable $X_t$ conditioned at being positive, are given by

$$\mathbf{P}(\tilde{X} = k) = \frac{e^{-t}}{1-e^{-t}}\frac{t^k}{k!}, \qquad \text{for } k \geq 1 \text{ and } t > 0.$$

This moment estimation is the following: if $(X_t)_{t\in[0,+\infty]}$ is the Khinchin family of $e^z$, then for any $\beta > 0$,

$$\mathbf{E}(X_t^\beta) \sim t^\beta, \qquad \text{as } t \to \infty,$$

and for any $\beta < 0$,

$$\mathbf{E}(\tilde{X}_t^\beta) \sim t^\beta, \qquad \text{as } t \to \infty.$$

See also [55] for precise asymptotic *expansions* of the negative moments of Poisson variables. ⊠

## 4.3 Gaussian power series for $f = e^g$, with $g \in \mathcal{P}$

Now we study some of the classes of functions introduced in Chapter 3 in the case $f = e^g$, with $g \in \mathcal{P}$. We give criteria for these power series to be Gaussian and also to be in the Hayman class.

### 4.3.1 Gaussian power series

Power series $f = e^g \in \mathcal{K}$, with $g \in \mathcal{P}$ a power series with radius of convergence $R > 0$, are non-vanishing, then we can apply the criteria of Chapter 3, Subsection 3.1.1. Observe that in this case, because $g$ has non-negative coefficients, we have, for any $k \geq 0$, that

$$|F^{(k)}(s+i\theta)| \leq F^{(k)}(s), \qquad \text{for any } \theta \in \mathbb{R} \text{ and } s < \ln(R),$$

in fact

$$\sup_{\theta\in\mathbb{R}} |F^{(k)}(s+i\theta)| = F^{(k)}(s), \qquad \text{for any } s < \ln(R).$$

Combining the previous equality with Lemma 3.1.11 it is possible to translate the results of subsection 3.1.1 to this particular case.

### A. Global criteria to be Gaussian

Our first criteria, a direct translation of Theorem 3.1.4, for $f = e^g$, with $g \in \mathcal{P}$, is the following.

**Theorem 4.3.1.** *Let $f \in \mathcal{K}$ be such that $f = e^g$, where $g \in \mathcal{P}$ is a power series with radius of convergence $R > 0$. If $g$ is a polynomial or $g$ is not a polynomial and satisfies that*

$$(\star) \quad \lim_{t\uparrow R} \frac{g'''(t)}{g''(t)^{3/2}} = 0,$$

*then $f$ is Gaussian.*

*Proof.* We apply Theorem 3.1.4 to $f = e^g$, with $g \in \mathcal{P}$. Using that $g$ has non-negative coefficients we have that

$$\frac{\sup_{\theta\in\mathbb{R}} |F'''(s+i\theta)|}{F''(s)^{3/2}} = \frac{F'''(s)}{F''(s)^{3/2}}.$$

First assume that $g$ is a polynomial, then Corollary 4.2.6 gives that

$$\lim_{s\to\infty} \frac{\sup_{\theta\in\mathbb{R}} |F^{(k)}(s+i\theta)|}{F''(s)^{k/2}} = \lim_{s\to+\infty} \frac{F^{(k)}(s)}{F''(s)^{k/2}} = 0,$$

therefore Theorem 3.1.4 gives that $f = e^g$ is Gaussian.

Assume now that $g$ is not a polynomial, then, applying Corollary 4.2.6, we obtain that

$$\frac{F'''(s)}{F''(s)^{3/2}} \asymp \frac{e^{3s} g'''(e^s)}{(e^{2s} g''(e^s))^{3/2}} = \frac{g'''(e^s)}{g''(e^s)^{3/2}}, \qquad \text{as } s \uparrow \ln(R).$$

therefore combining hypothesis $(\star)$ with Theorem 3.1.4 we obtain that $f = e^g$ is Gaussian. □

We can generalize the previous result combining Theorem 3.1.6 with Corollary 4.2.6.

**Theorem 4.3.2.** *Let $f \in \mathcal{K}$ be such that $f = e^g$, where $g \in \mathcal{P}$ is a power series with radius of convergence $R > 0$. If $g$ is not and polynomial and satisfies that*

$$(\dagger) \quad \lim_{t\uparrow R} \frac{g^{(4)}(t)}{g''(t)^2} = 0,$$

*then*

$$f \text{ is Gaussian} \qquad \text{if and only if} \qquad \lim_{t\uparrow R} \frac{g'''(t)}{g''(t)^{3/2}} = 0.$$

*Proof.* We want to apply Theorem 3.1.6 to $f = e^g$, with $g \in \mathcal{P}$ not a polynomial. Corollary 4.2.6 gives that

$$\frac{F'''(s)}{F''(s)^{3/2}} \asymp \frac{g'''(s)}{g''(s)^{3/2}} \quad \text{and} \quad \frac{F^{(4)}(s)}{F''(s)^2} \asymp \frac{g^{(4)}(s)}{g''(s)^2}, \qquad \text{as } s \uparrow \ln(R), \tag{4.3.1}$$

therefore, using (†), we obtain that

$$\limsup_{s\uparrow\ln(R)} \frac{\sup_{\phi\in\mathbb{R}} |F^{(4)}(s+i\phi)|}{F''(s)^2} = \lim_{s\uparrow\ln(R)} \frac{F^{(4)}(s)}{F''(s)^2} = 0$$

Applying Theorem 3.1.6 we conclude that

$$f \text{ is Gaussian} \qquad \text{if and only if} \qquad \lim_{s\uparrow\ln(R)} \frac{F'''(s)}{F''(s)^{3/2}} = 0,$$

which is equivalent, by means of (4.3.1), to

$$f \text{ is Gaussian} \qquad \text{if and only if} \qquad \lim_{s\uparrow\ln(R)} \frac{g'''(s)}{g''(s)^{3/2}} = 0.$$

□

The following result is Corollary 3.1.8 applied to $f = e^g$, with $g \in \mathcal{P}$ not a polynomial.

**Corollary 4.3.3.** *Let $f \in \mathcal{K}$ be such that $f = e^g$, where $g \in \mathcal{P}$ is a power series with radius of convergence $R > 0$. If $g$ is not a polynomial and satisfies that*

$$\limsup_{t\uparrow R} \frac{g^{(4)}(t)}{g''(s)^2} < +\infty,$$

*then*

$$f \textit{ is Gaussian} \qquad \textit{implies that} \quad \lim_{t\uparrow R} \frac{g'''(t)}{g''(t)^{3/2}} = 0.$$

*Proof.* This result follows combining Corollary 4.2.6 with Corollary 3.1.8. □

**Remark 4.3.4** (Generalization of the previous Theorems)**.** Any of the previous Theorems can be generalized using the same ideas. Let $g \in \mathcal{P}$ be a power series with radius of convergence $R > 0$ which is not a polynomial. First assume that for some $k \geq 4$ we have

$$\limsup_{t\uparrow R} \frac{g^{(k)}(s)}{g''(s)^{k/2}} = 0,$$

then

$$f \text{ is Gaussian} \qquad \text{if and only if} \quad \lim_{t\uparrow R} \frac{g^{(j)}(s)}{g''(s)^{j/2}} = 0, \;\; \text{for any } 3 \leq j \leq k-1.$$

Assume now that for some $k \geq 4$ we have

$$\limsup_{t\uparrow R} \frac{g^{(k)}(s)}{g''(s)^{k/2}} < +\infty,$$

then

$$f \text{ is Gaussian} \qquad \text{implies that} \qquad \lim_{t\uparrow R} \frac{g^{(j)}(s)}{g''(s)^{j/2}} = 0, \ \text{for any } 3 \leq j \leq k-1.$$

⊠

### B. Moment criteria for Gaussianity

We translate Theorem 3.1.19 to the case $f = e^g$ with $g \in \mathcal{P}$ not a polynomial.

**Theorem 4.3.5.** *Let $f \in \mathcal{K}$ be such that $f = e^g$, where $g \in \mathcal{P}$ has radius of convergence $R > 0$. If $g$ is not a polynomial and satisfies*

$$\lim_{t\uparrow R} \frac{g^{(k)}(t)}{g''(t)^{k/2}} = 0, \qquad \textit{for any } k \geq 3,$$

*then $f = e^g$ is Gaussian.*

*Proof.* This result follows combining Corollary 4.2.6 with Theorem 3.1.19. □

**Remark 4.3.6.** Let $g \in \mathcal{P}$ be a power series with radius of convergence $R > 0$ which is not a polynomial. Assume that

$$\limsup_{t\uparrow R} \frac{g^{(k)}(t)}{g''(t)^{k/2}} < +\infty, \qquad \text{for any } k \geq 3,$$

then

$$f \text{ is Gaussian} \qquad \text{if and only if} \qquad \lim_{t\,\uparrow R} \frac{g^{(k)}(t)}{g''(t)^{k/2}} = 0, \qquad \text{for any } k \geq 3.$$

This result follows combining Theorem 4.3.5 with Remark 4.3.4. ⊠

### C. Applications

Now we apply the previous criteria for Gaussianity to a variety of combinatorial examples which have EGF given by an exponential power series, that is, generating functions of the form $f = e^g$, with $g \in \mathcal{P}$.

**C.1. EGF of the class of sets of (non-empty) sets** The exponential generating function of the class of sets of (non-empty) sets is given by $f(z) = \exp(e^z - 1)$, in this case the radius of convergence $R = \infty$ and $g(z) = (e^z - 1) \in \mathcal{P}$.

We have

$$\frac{g'''(t)}{g''(t)^{3/2}} = \frac{e^t}{e^{3t/2}} = e^{-t/2},$$

and therefore

$$\lim_{t\to\infty} \frac{g'''(t)}{g''(t)^{3/2}} = 0.$$

Applying Theorem 4.3.1 we conclude that the EGF of the sets of (non-empty) sets is Gaussian.

**C.2. EGF of the class of sets of (non-empty) lists** The exponential generating function of the class of sets of (non-empty) lists is given by

$$f(z) = \exp(z/(1-z)), \qquad \text{for any } |z| < 1.$$

In this case $g(z) = z/(1-z) \in \mathcal{P}$ and

$$g''(t) = \frac{2}{(1-t)^3} \quad \text{and} \quad g'''(t) = \frac{6}{(1-t)^4}, \quad \text{for any } t \in (0,1),$$

therefore

$$\lim_{t\uparrow 1} \frac{g'''(t)}{g''(t)^{3/2}} = 0.$$

Applying Theorem 4.3.1 we conclude that the EGF of the sets of (non-empty) lists is Gaussian.

**C.3. EGF of the class of sets of pointed sets** The EGF of the class of sets of pointed sets is given by $f(z) = e^{ze^z}$, then $g(z) = ze^z \in \mathcal{P}$ has radius of convergence $R = +\infty$.

We have

$$\lim_{t\to\infty} \frac{g'''(t)}{g''(t)^{3/2}} = \lim_{t\to\infty} \frac{e^t(t+3)}{(t+2)^{3/2}e^{3t/2}} = 0.$$

In fact, for any $k \geq 3$, we have

$$\lim_{t\to\infty} \frac{g^{(k)}(t)}{g''(t)^{k/2}} = \lim_{t\to\infty} \frac{e^t(t+k)}{(t+2)^{k/2}e^{kt/2}} = 0.$$

Then applying Theorem 4.3.1 or Theorem 4.3.5 we obtain that EGF of the class of sets of pointed sets is Gaussian.

**C.4. EGF of the class of sets of cycles** The EGF of the class of sets of cycles is given by

$$f(z) = \exp(\ln(1/(1-z))) = 1/(1-z), \quad \text{for } |z| < 1.$$

We already know that this power series is not Gaussian, see equation (3.1.4) in the examples at the beginning of Chapter 3. An alternative argument to prove that $f$ is not Gaussian is the following: the fulcrum $F$ of $f$ is given by

$$F(z) = \ln(1/(1-z)) = \sum_{n=1}^{\infty} \frac{z^n}{n}, \quad \text{for any } |z| < 1,$$

therefore, for any $k \geq 1$, we have

$$\sup_{\theta \in \mathbb{R}} |F^{(k)}(s + i\theta)| = F^{(k)}(s), \quad \text{for any } s < 0\,.$$

In fact, for any $s < 0$, we have

$$F''(s) = \frac{1}{(1-e^s)^2}, \quad F'''(s) = \frac{2}{(1-e^s)^3}, \quad \text{and} \quad F^{(4)}(s) = \frac{6}{(1-e^s)^4}\,.$$

Therefore

$$\lim_{s \uparrow 0} \frac{F^{(4)}(s)}{F''(s)^2} = 6 < +\infty \quad \text{and} \quad \lim_{s \uparrow 0} \frac{F'''(s)}{F''(s)^{3/2}} = 2\,.$$

Applying, for instance, Corollary 3.1.8 we conclude that $f$ is not Gaussian.

**C.5. Some parametric families**

- The power series $f(z) = \exp(1/(1-z)^\gamma)$, with $\gamma > 0$, are all Gaussian. Observe that $g(z) = 1/(1-z)^\gamma$ has radius of convergence $R = 1$ and

$$\frac{g'''(t)}{g''(t)^{3/2}} \asymp (1-t)^{\gamma+1}, \quad \text{as } t \uparrow 1,$$

  then, appealing to Theorem 4.3.1, we conclude that $f$ is Gaussian.

- Denote by $f = e^g$ with

$$g(z) = \sum_{n=1}^{\infty} n^\alpha z^n, \quad \text{for } |z| < 1.$$

  *The power series $f = e^g$ with $\alpha > -1$ are also Gaussian.* Indeed: observe that

$$t^3 g'''(t) = \sum_{n=1}^{\infty} n^{\alpha+3} t^n, \quad \text{and} \quad t^2 g''(t) = \sum_{n=1}^{\infty} n^{\alpha+2} t^n, \quad \text{for } t \in (0,1).$$

Appealing to Lemma 4.2.7 we find that

$$g'''(t) \asymp \frac{1}{(1-t)^{\alpha+4}} \quad \text{and} \quad g''(t) \asymp \frac{1}{(1-t)^{\alpha+3}}, \quad \text{as } t \uparrow 1,$$

then

$$\frac{g'''(t)}{g''(t)^{3/2}} \asymp (1-t)^{(\alpha+1)/2}, \quad \text{as } t \uparrow 1,$$

and therefore, for any $\alpha > -1$, we have

$$\lim_{t\uparrow 1} \frac{g'''(t)}{g''(t)^{3/2}} = 0\,.$$

By virtue of Theorem 4.3.1 we have that $f$ is Gaussian, for any $\alpha > -1$.

- Fix $\delta \geq 0$. Let $f$ be a power series given by

$$f(z) = \exp\left(\sum_{n=1}^{\infty} \frac{n^{n-\delta}}{n!} z^n\right), \quad \text{for any } |z| < 1/e\,.$$

Notice that for $\delta = 0$ we obtain the EGF of the class of sets of functions, for $\delta = 1$ the EGF of the class of forests of rooted trees and for $\delta = 2$ the EGF for the sets of trees.

In this case we have that

$$g(z) = \sum_{n=1}^{\infty} \frac{n^{n-\delta}}{n!} z^n, \quad \text{for any } |z| < 1/e, \tag{4.3.2}$$

and $f = e^g$.

*We prove now that $f$ is Gaussian if and only if* $0 \leq \delta < 1/2$.

Our fist lemma tell us the values of $\delta \geq 0$ for which $g$ is continuous on the closure of $\mathbb{D}(0, 1/e)$.

**Lemma 4.3.7.** *For any $\delta > 1/2$, the power series $g$ defined in (4.3.2) is continuous on the closure of* $\mathbb{D}(0, 1/e)$.

*Proof.* Using Stirling's formula we find that

$$\frac{n^{n-\delta}}{n!e^n} \asymp \frac{1}{n^{1/2+\delta}}, \quad \text{as } n \to \infty\,, \tag{4.3.3}$$

then, for any $\delta > 1/2$, uniform convergence of the partial sum in the closure of $\mathbb{D}(0, 1/e)$, gives that the power series $g$ is continuous on the closure of $\mathbb{D}(0, 1/e)$. □

For the variance $\sigma_f^2$ we have the following lemma.

**Lemma 4.3.8.** *For the variance of $f = e^g$ we have:*

- $\lim_{t\uparrow 1/e} \sigma_f^2(t) = +\infty$, *for any* $0 \le \delta \le 5/2$.
- $\lim_{t\uparrow 1/e} \sigma_f^2(t) = \sigma_f^2(1/e) < +\infty$, *for any* $\delta > 5/2$.

*Proof.* Appealing to Proposition 4.2.7, for any $0 \le \delta < 5/2$, we have

$$g''(t/e) \asymp \frac{1}{(1-t)^{5/2-\delta}}, \quad \text{as } t \uparrow 1, \tag{4.3.4}$$

then $\lim_{t\uparrow 1/e} \sigma_f^2(t) = +\infty$. For $\delta = 5/2$ we have

$$g''(t/e) \asymp \ln\left(\frac{1}{1-t}\right), \quad \text{as } t \uparrow 1,$$

then

$$\lim_{t\uparrow 1/e} \sigma_f^2(t) = \infty.$$

Finally we have, by virtue of equation (4.3.4), that

$$\lim_{t\uparrow 1/e} g''(t) < +\infty, \quad \text{for any } \delta > 5/2,$$

which is equivalent to $\lim_{t\uparrow 1/e} \sigma_f^2(t) < +\infty$, for any $\delta > 5/2$, see Lemma 4.2.2. □

**Remark 4.3.9.** For exponential power series, that is, power series of the form $f = e^g$, with $g \in \mathcal{P}$, we have the inequality $m_f(t) \le \sigma_f^2(t)$, see Lemma 4.2.3, then, Lemma 4.3.8 also gives that $\lim_{t\uparrow 1/e} m_f(t) = m_f(1/e) < +\infty$, for any $\delta > 5/2$.

Now we apply the previous analysis to the characteristic functions of $\breve{X}_t$:

- For $1/2 < \delta \le 5/2$ we have that $f$ is continuous on the closure of $\mathbb{D}(0, 1/e)$ and $\lim_{t\uparrow 1/e} \sigma_f^2(t) = +\infty$, then

  $$\lim_{t\uparrow 1/e} |\mathbf{E}(e^{i\theta \breve{X}_t})| = \lim_{t\uparrow 1/e} \frac{|f(te^{i\theta/\sigma_f(t)})|}{f(t)} = 1, \quad \text{for any } \theta \in \mathbb{R}.$$

  Therefore the power series $f$ is not Gaussian, for any $1/2 < \delta \le 5/2$.

  In fact, for any $0 \le \delta < 3/2$, we have, appealing again to Proposition 4.2.7, that

  $$m_f(t/e) = (t/e)g'(t/e) \asymp \frac{1}{(1-t)^{3/2-\delta}}, \quad \text{as } t \uparrow 1.$$

and therefore

$$\frac{m_f(t/e)}{\sigma_f(t/e)} \asymp (1-t)^{\delta/2-1/4}, \qquad \text{as } t \uparrow 1 \tag{4.3.5}$$

from the previous analysis we obtain that, for any $1/2 < \delta < 3/2$, we have

$$\lim_{t\uparrow 1/e} \mathbf{E}(e^{i\theta \breve{X}_t}) = \lim_{t\uparrow 1/e} \frac{f(te^{i\theta/\sigma_f(t)})}{f(t)} e^{-i\theta m_f(t)/\sigma_f(t)} = 1, \qquad \text{for any } \theta \in \mathbb{R},$$

and this implies that, for any $1/2 < \delta < 3/2$, the family $(\breve{X}_t)$ tends in distribution towards the constant 0, as $t \uparrow 1/e$.

Similarly, but using that $3/2 < \delta \leq 5/2$, we have that $\lim_{t\uparrow 1/e} m_f(t)$ is finite, see the comparison (4.3.3), while $\lim_{t\uparrow 1/e} \sigma_f^2(t) = +\infty$, see Lemma 4.3.8. For $\delta = 3/2$ we have that

$$m_f(t) \asymp \ln\left(\frac{1}{1-et}\right) \qquad \text{and} \qquad \sigma_f(t) \asymp \frac{1}{\sqrt{1-et}}, \qquad \text{as } t \uparrow 1/e.$$

Then, for any $3/2 \leq \delta \leq 5/2$ we have that

$$\lim_{t\uparrow 1/e} \frac{m_f(t)}{\sigma_f(t)} = 0,$$

and therefore the normalized family $(\breve{X}_t)$ tends in distribution towards the constant 0, as $t \uparrow 1/e$. Here we use again that $g$ is continuous on the closure of $\mathbb{D}(0, 1/e)$, for any $\delta > 1/2$.

– For $\delta > 5/2$ we have

$$\lim_{t\uparrow 1/e} |\mathbf{E}(e^{i\theta \breve{X}_t})| = \lim_{t\uparrow 1/e} \frac{|f(te^{i\theta/\sigma_f(t)})|}{f(t)} = \frac{|f((1/e)e^{i\theta/\sigma_f(1/e)})|}{f(1/e)}, \tag{4.3.6}$$

for any $\theta \in \mathbb{R}$. Here we use that $f$ is continuous on the closure of $\mathbb{D}(0, 1/e)$ and also that $\lim_{t\uparrow 1/e} \sigma_f^2(t) \triangleq \sigma_f^2(1/e) < +\infty$, see Lemma 4.3.7 and 4.3.8, respectively. The term at the right-hand side of the limit (4.3.6) is periodic, then $f$ is not Gaussian for any $\delta > 5/2$.

– Finally we study the case $0 \leq \delta \leq 1/2$. Appealing again to Proposition 4.2.7, for any $0 \leq \delta < 9/2$ we have

$$g^{(4)}(t/e) \asymp \frac{1}{(1-t)^{9/2-\delta}}, \qquad \text{as } t \uparrow 1\,,$$

therefore, for $0 \leq \delta < 5/2$, using equation (4.3.4), we conclude that

$$\frac{g^{(4)}(t)}{g''(t)^2} \asymp (1-et)^{1/2-\delta}, \qquad \text{as } t \uparrow 1/e.$$

For $0 \le \delta \le 1/2$, using an analogous argument, we find that

$$\frac{g'''(t)}{g''(t)^{3/2}} \asymp (1-et)^{1/4-\delta/2}, \qquad \text{as } t \uparrow 1/e.$$

and therefore, applying Theorem 4.3.2, we conclude that $f$ is Gaussian if and only if $0 \le \delta < 1/2$.

In summary: $f$ is Gaussian if and only if $0 \le \delta < 1/2$; this is the answer to Question 1 of [18]. For $1/2 < \delta \le 5/2$, the normalized family $\breve{X}_t$ tends in distribution towards the constant 0, as $t \uparrow 1/e$. For $\delta > 5/2$ the family of random variables $(\breve{X}_t)$ converges in distribution towards a non-constant random variable $X$ which is not a standard normal random variable.

Thus the EGF of the sets of functions ($\delta = 0$) is Gaussian, but the EGFs of the sets of trees ($\delta = 2$) and of the sets of rooted trees ($\delta = 1$) are not Gaussian.

### 4.3.2 Hayman class for exponentials $f = e^g$ with $g \in \mathcal{P}$

Along this section we consider again power series $f = e^g$, with $g \in \mathcal{P}$ a power series with radius of convergence $R > 0$. We are interested in conditions on $g$ which guarantee that $f = e^g$, a power series in $\mathcal{K}$, is in the Hayman class.

There are two results of Hayman, see [48], along this line:

(a) If $g$ is a power series in the Hayman class, then $e^g$ is in the Hayman class.

(b) If $B(z) = b_N z^N + \cdots + b_0$ is a non-constant polynomial with non-negative coefficients and such that $Q_B = \gcd\{n \ge 1 : b_n > 0\} = 1$, then $e^B$ is in the Hayman class.

**Remark 4.3.10.** In [72], it is shown that if $g$ is in the Hayman class, then $f = e^g$ satisfies the stronger (than Hayman) conditions of Harris and Schoenfeld, [47], which allow to obtain full asymptotic expansions for the coefficients and not just asymptotic formulas. ⊠

**Remark 4.3.11.** In Theorem X of [48], Hayman proves the following stronger version of (b): Let $B(z) = \sum_{j=0}^{N} b_n z^n$ be a non-constant polynomial with real coefficients such that for each $d > 1$ there exists $m$, not a multiple of $d$, such that $b_m \neq 0$ and such that if $m(d)$ is the largest such $m$, then $b_{m(d)} > 0$. If $e^B$ is in $\mathcal{K}$, then $e^B$ is in the Hayman class. For a direct proof of this result see Proposition 5.1 of [17]. ⊠

Theorem 4.3.12, stated below, which is [17, Theorem 4.1], exhibits conditions on the function $g$ which ensure that $f$ is in the Hayman class.

Further, we describe in Theorem 4.3.13 a 'practical' approach to verify that a power series satisfies the conditions of Theorem 4.3.12.

Finally, Theorems 4.3.15 and 4.3.17 give easily verifiable conditions on *the coefficients* $b_n$ *of* $g$ that imply that $f = e^g$ is in the Hayman class. These theorems may be directly applied to the generating functions of sets constructions (and other exponentials) to check that they belong to the Hayman class.

Let $(X_t)$ be the Khinchin family of $f = e^g$, with $g \in \mathcal{P}$ a power series with radius of convergence $R > 0$. The variance condition required for $f$ to be in the Hayman class translates readily in the following condition in terms of $g$:

$$\lim_{t \uparrow R} \big(t g'(t) + t^2 g''(t)\big) = +\infty\,.$$

Observe that for $g \in \mathcal{P}$, not a polynomial, the variance condition implies that $\lim_{t\uparrow R} g''(t) = +\infty$ and therefore $\lim_{t\uparrow R} g'''(t) = +\infty$, see Lemma 4.2.2.

### Criteria for $f = e^g$, with $g \in \mathcal{P}$, to be in the Hayman class

The following theorem is [17, Theorem 4.1].

**Theorem 4.3.12.** *Denote $f = e^g$, with $g \in \mathcal{P}$ a power series with radius of convergence $R > 0$ which is not a polynomial.*

*If the variance condition*

$$\lim_{t\uparrow R} \sigma_f^2(t) = \lim_{t\uparrow R} \big(t g'(t) + t^2 g''(t)\big) = +\infty\,. \tag{4.3.7}$$

*is satisfied and there is a cut function $h(t)$ so that*

$$\lim_{t\uparrow R} t^3 g'''(t)\, h(t)^3 = 0\,. \tag{4.3.8}$$

*and that*

$$\lim_{t\uparrow R} \sigma_f(t) \exp\Big(\sup_{h(t)\le|\theta|\le\pi} \Re g(te^{i\theta}) - g(t)\Big) = 0\,. \tag{4.3.9}$$

*hold, then $f = e^g$ is in the Hayman class.*

Condition (4.3.23) of Theorem 4.3.12, gives that the cut function $h$ fulfills condition (3.3.2) on the major arc. In terms of $h$ and $g$, condition (3.3.2) on the minor arc follows from condition 4.3.9 of Theorem 4.3.12.

*Proof.* The variance condition (4.3.22) follows combining the variance condition (3.3.1) with Lemma 4.2.1. The major arc condition (4.3.23) follows combining Theorem 3.3.4 with Corollary 4.2.5. Finally for the minor arc condition we have

$$|\mathbf{E}(e^{i\theta \breve{X}_t})| = \frac{|f(te^{i\theta/\sigma_f(t)})|}{f(t)} = \exp\big(\Re(g(te^{i\theta/\sigma_f(t)})) - g(t)\big)$$

and therefore

$$\begin{aligned}\sup_{h(t)\sigma_f(t)\le\ |\theta|\ \le\pi\sigma_f(t)} |\mathbf{E}(e^{i\theta\breve{X}_t})| &= \sup_{h(t)\le|\theta|\le\pi} \exp\big(\Re(g(te^{i\theta})) - g(t)\big)\\ &= \exp\Big(\sup_{h(t)\le|\theta|\le\pi} \big(\Re(g(te^{i\theta})) - g(t)\big)\Big)\end{aligned}$$

the previous equality gives that, in this case, that is, for $f = e^g$, the minor arc condition (3.3.3) is equivalent to condition (4.3.9). □

In practice, condition (4.3.9) on the cut $h(t)$ is verified by checking that there are positive functions $U, V$ defined in $(t_0, R)$, for some $t_0 \in (0, R)$, where $U$ takes values in $(0, \pi]$ and $V$ in $(0, \infty)$ and are such that

$$\sup_{|\varphi|\geq\omega} \big(\Re g(te^{\imath\varphi}) - g(t)\big) \leq -V(t)\,\omega^2\,, \quad \text{for } \omega \leq U(t) \text{ and } t \in (t_0, R)\,, \tag{4.3.10}$$

and thus that the cut $h$ is such that

$$h(t) \leq U(t), \text{ for } t \in (t_0, 1) \quad \text{and} \quad \lim_{t\uparrow R} \sigma(t)\, e^{-V(t)h(t)^2} = 0\,. \tag{4.3.11}$$

Thus we have the following corollary of Theorem 4.3.12, which we state as a theorem.

**Theorem 4.3.13.** *Let $g \in \mathcal{P}$ a power series with radius of convergence $R > 0$ which is not a polynomial, and let $h : [0, R) \to (0, \pi)$ be a cut function.*

*Assume that the variance condition (4.3.22) and the condition on the major arc (4.3.23) of Theorem 4.3.12 are satisfied.*

*Assume further that there are positive functions $U, V$ defined in $(t_0, R)$, for some $t_0 \in (0, R)$, where $U$ takes values in $(0, \pi]$ and $V$ in $(0, \infty)$, and such that*

$$\sup_{|\theta|\geq\omega} \big(\Re(g(te^{i\theta} - g(t)))\big) \leq -V(t)\omega^2, \quad \textit{for } \omega \leq U(t) \textit{ and } t \in (t_0, R), \tag{4.3.12}$$

*and*

$$h(t) \leq U(t), \quad \textit{for } t \in (t_0, R) \quad \textit{and} \quad \lim_{t\uparrow R} \sigma_f(t) e^{-V(t)h(t)^2} = 0. \tag{4.3.13}$$

*Then $f = e^g$ is in the Hayman class.*

*Proof.* The pair of conditions (4.3.12) and (4.3.13) together imply condition (4.3.9). This is so since the first half of (4.3.13) allows us to take $\omega = h(t)$ in (4.3.12) to obtain

$$\exp\left(\sup_{|\theta|\geq h(t)} \big(\Re(g(te^{i\theta})) - g(t)\big)\right) \leq \exp(-V(t)h(t)^2), \quad \text{for } t \in (t_0, R),$$

and to then to apply the second half 4.3.13. □

Notice further that if for a function $V$ defined in $(t_0, R)$, for some $t_0 \in (0, R)$, and taking values in $(0, +\infty)$, we have that

$$\Re(g(te^{i\theta})) - g(t) \leq -V(t)\theta^2, \quad \text{for } |\theta| \leq \pi \text{ and } t \in (t_0, R), \tag{4.3.14}$$

then we may set $U \equiv \pi$, and the cut $h$ is required to satisfy

$$\lim_{t\uparrow R} \sigma_f(t) e^{-V(t)h(t)^2} = 0. \tag{4.3.15}$$

Thus these two conditions (4.3.14) and 4.3.15 imply the conditions 4.3.12 and 4.3.13 of Theorem 4.3.13

### 4.3.3 Coefficient criteria for $g \in \mathcal{P}$ for $f = e^g$ to be in the Hayman class

Next we will exhibit easily verifiable requirements on the coefficients of $g$ to obtain functions $V$ and $U$ so that the conditions 4.3.12 and 4.3.13 hold, or simply, a functions $V$ so that the conditions (4.3.14) and (4.3.15) are satisfied.

The announced requirements differ if the power series $g$ is entire or has finite radius of convergence; we split the discussion accordingly.

#### A. Case $g$ entire, $R = \infty$

Here $g \in \mathcal{P}$ is a transcendental entire function.

**Lemma 4.3.14.** *If $g \in \mathcal{P}$, a transcendental entire function, satisfies, for some $\beta > 0$ and $B > 0$ that*

$$\Re(g(te^{i\theta})) - g(t) \le B(e^{\beta t \cos(\theta)} - e^{\beta t}), \qquad \text{for } t > 0 \text{ and } |\theta| \le \pi, \tag{4.3.16}$$

*then condition (4.3.14) is satisfied with $V(t) = Ce^{\beta t}$, for some constant $C > 0$ depending only on $B$ and $\beta$, and condition (4.3.15) requires that $\lim_{t\to\infty} \sigma_f(t) \exp\left(-Ce^{\beta t} h(t)^2\right) = 0$.*

*Proof.* For $t > 0$ and $\beta > 0$, the convexity of the exponential function $x \mapsto e^{\beta t x}$ in the interval $[-1, 1]$ gives

$$e^{\beta t} - e^{\beta t \cos(\theta)} \ge (1 - \cos(\theta)) \frac{e^{\beta t} - e^{-\beta t}}{2}, \qquad \text{for } \theta \in [-\pi, \pi].$$

Thus, we have, for $t \ge 1$, that

$$\begin{aligned} e^{\beta t \cos(\theta)} - e^{\beta t} \le (\cos(\theta) - 1) \frac{e^{\beta t} - e^{-\beta t}}{2} &\le e^{\beta t} \frac{1 - e^{-2\beta}}{2} (\cos(\theta) - 1) \\ &\le -\frac{1 - e^{-2\beta}}{\pi^2} e^{\beta t} \theta^2. \end{aligned}$$

In the second inequality, we have used that $t \ge 1$; the inequality $1 - \cos(\theta) \ge 2x^2/\pi^2$, for any $x \in [-\pi, \pi]$, yields the last step. ☐

**Theorem 4.3.15.** *Let $g(z) = \sum_{n=0}^{\infty} b_n z^n \in \mathcal{P}$ be a power series, not a polynomial, verifying*

$$B \frac{\beta^n}{n!} \le b_n \le L \frac{\lambda^n}{n!}, \qquad \text{for any } n \ge 1, \tag{4.3.17}$$

*for some constants $B, L > 0$ and $\beta, \lambda \ge 0$ such that $2\lambda < 3\beta$. Then $g$ satisfies the requirements of Theorem 4.3.13 and, therefore, $f = e^g$ is in the Hayman class.*

*Proof.* Condition (4.3.17) gives that the power series $g$ has radius of convergence $R = \infty$. Besides, condition (4.3.17) and Proposition 4.2.9 show that the variance $\sigma_f^2(t)$ satisfies that $t^2 e^{\beta t} = O(\sigma_f^2(t))$ and $\sigma_f^2(t) = O(t^2 e^{\lambda t})$. Appealing again to Proposition 4.2.9 we also find that $t^3 g'''(t) = O(t^3 e^{\lambda t})$, as $t \to \infty$.

The variance condition (4.3.22) is thus obviously satisfied. The major arc condition (4.3.9) holds with the cut $h(t) = e^{-\alpha t}$, where $\alpha \in (\lambda/3, \beta/2)$, since $3\alpha > \lambda$.

Also, for $t > 0$ and $|\theta| \leq \pi$, we have that

$$\begin{aligned}\Re(g(te^{i\theta})) - g(t) &= \sum_{n=1}^{\infty} b_n t^n (\cos(n\theta) - 1) \leq B \sum_{n=1}^{\infty} \frac{(\beta t)^n}{n!}(\cos(n\theta) - 1) \\ &= B(\Re(e^{\beta t e^{i\theta}}) - e^{\beta t}) \leq B(|e^{\beta t e^{i\theta}}| - e^{\beta t}) = B(e^{\beta t \cos(\theta)} - e^{\beta t})\end{aligned}$$

and, consequently, condition (4.3.16) of Lemma 4.3.14 is satisfied.

Condition (4.3.14) is satisfied with the choice $V(t) = Ce^{t\beta}$, $t > 1$, and (4.3.15) hold because $2\alpha < \beta$. As mentioned above, these two conditions imply condition (4.3.12) and (4.3.13) of Theorem 4.3.13, and thus we conclude that $f$ is in the Hayman class. □

### B. Case $g$ not entire, $R < \infty$

Along this section, we will repeatedly appeal to Proposition 4.2.7.

**Lemma 4.3.16.** *If a power series $g$ with non-negative coefficients and radius of convergence $R = 1$ satisfies, for some $\beta > 0$ and $B > 0$, that*

(4.3.18) $$\Re(g(te^{i\theta})) - g(t) \leq B\left(\frac{1}{|1 - te^{i\theta}|^{\beta}} - \frac{1}{(1-t)^{\beta}}\right), \quad \text{for } t \in (0,1) \text{ and } |\theta| \leq \pi,$$

*then*

$$\sup_{|\theta| \geq \omega} \Big(\Re(g(te^{i\theta})) - g(t)\Big) \leq -C\frac{1}{(1-t)^{2+\beta}}\omega^2$$

*for $t \in (1/2, 1)$ and $0 \leq \omega \leq D(1-t)$, where $C > 0$ and $D > 0$ depend only on $\beta$ and $B$.*

*And, in particular, if we set $V(t) = C/(1-t)^{2+\beta}$, and $U(t) = D(1-t)$, for $t \in (1/2, 1)$, then condition (4.3.12) is satisfied, and condition (4.3.13) requires that*

$$\lim_{t \uparrow 1} \sigma_f(t) \exp(-Ch(t)^2/(1-t)^{2+\beta}) = 0\,.$$

*Proof.* Let $\beta > 0$. For $t \in (0,1)$ and $|\theta| \leq \pi$, we have that

$$\Re(g(te^{i\theta})) - g(t) \leq \frac{B}{(1-t)^{\beta}}\left(\left(\left|\frac{1-t}{1-te^{i\theta}}\right|^2\right)^{\beta/2} - 1\right).$$

Thus, for $t \in (1/2, 1)$ and $0 \leq \omega \leq (1-t)$ we have that

$$\sup_{|\theta| \geq \omega} \big(\Re(g(te^{i\theta})) - g(t)\big) \leq \frac{B}{(1-t)^{\beta}}\left((1 - C\frac{1}{(1-t)^2}\omega^2)^{\beta/2} - 1\right),$$

and, for $\omega < D_\beta(1-t)$, for $D_\beta > 0$ appropriate small, and some constant $C_\beta > 0$, we have that

$$\sup_{|\theta| \geq \omega} \Big(\Re(g(te^{i\theta})) - g(t)\Big) \leq -\frac{B}{(1-t)^{\beta}} C_\beta \frac{1}{(1-t)^2}\omega^2\,.$$

□

**Theorem 4.3.17.** *Let $g(z) = \sum_{n=0}^{\infty} b_n z^n$ satisfy*

$$B\frac{n^\beta}{R^n} \le b_n \le L\frac{n^\lambda}{R^n}\,, \qquad \text{for } n \ge 1, \tag{4.3.19}$$

*for some constants $B, L > 0$, finite radius $R > 0$ and $\beta, \lambda > -1$ such that $2\lambda < 3\beta + 1$. Then $g$ satisfies the requirements of Theorem 4.3.13 and, therefore, $f = e^g$ is in the Hayman class.*

*Proof.* Inequality (4.3.19) gives that the power series $g$ has radius of convergence $R > 0$. By considering $g(Rz)$, we may assume that $R = 1$. Appealing to Proposition 4.2.7 we see that the variance of $f = e^g$ satisfies that

$$\frac{1}{(1-t)^{3+\beta}} = O(\sigma_f(t)^2), \qquad \text{and} \qquad \sigma_f^2(t) = O\left(\frac{1}{(1-t)^{3+\lambda}}\right), \qquad \text{as } t \uparrow 1\,,$$

and also that

$$t^3 g'''(t) = O(1/(1-t)^{4+\lambda})\,, \qquad \text{as } t \uparrow 1.$$

The variance condition (4.3.22) is obviously satisfied. We propose the cut $h(t) = (1-t)^\alpha$, where $\alpha \in (\lambda/3 + 4/3, \beta/2 + 3/2)$.

The major arc condition (4.3.9) holds since $3\alpha > \lambda + 4$

For $t \in (0,1)$ and $|\theta| \le \pi$, we have that

$$\begin{aligned}
\Re(g(te^{i\theta})) - g(t) &= \sum_{n=1}^{\infty} b_n t^n (\cos(n\theta) - 1) \le B \sum_{n=1}^{\infty} n^\beta t^n (\cos(n\theta) - 1)\\
&\le B_\beta \sum_{n=1}^{\infty} \frac{\Gamma(n+\beta+1)}{\Gamma(\beta+1) n!} t^n (\cos(n\theta) - 1)\\
&= B_\beta \left( \Re\left(\frac{1}{(1-te^{i\theta})^{\beta+1}}\right) - \left(\frac{1}{1-t}\right)^{\beta+1}\right)\\
&\le B_\beta \left( \left|\frac{1}{1-te^{i\theta}}\right|^{\beta+1} - \left(\frac{1}{1-t}\right)^{\beta+1}\right)
\end{aligned}$$

where we have appealed to the comparison (4.2.5) withing Proposition 4.2.7. Consequently, condition (4.3.18) of Lemma 4.3.16 does hold.

If we set $V(t) = C/(1-t)^{\beta+3}$, for $t \in (1/2, 1)$, (4.3.12) is satisfied and (4.3.13) requires that

$$\lim_{t \uparrow 1} \sigma_f(t) \exp\left(-Ch(t)^2/(1-t)^{\beta+3}\right) = 0,$$

which does hold since $2\alpha < \beta + 3$.

Thus, the conclusion follows from Theorem 4.3.13. □

### 4.3.4 Applications of Theorems 4.3.15 and 4.3.17

In this subsection we apply the previous criteria to a variety of exponential power series which are EGFs of set combinatorial classes

#### A. Set of Labeled Classes

- The EGF of the Bell is numbers, enumerating sets of (non-empty) sets, that is partitions, is given by $f(z) = e^{e^z-1}$. Here $g(z) = e^z - 1$ and $b_n = 1/n!$, for $n \geq 1$. This power series satisfies the hypothesis of Theorem 4.3.15, and therefore $f$ is in the Hayman class.

- The EGF of the sets of pointed sets is given by $f(z) = e^{ze^z}$. Here $g(z) = ze^z$, with $b_n = n/n!$, for $n \geq 1$. This power series satisfies the hypothesis of Theorem 4.3.15, and therefore $f$ is in the Hayman class.

  We collect the previous discussion, for later use, by means of the following theorem.

  **Theorem 4.3.18.** *The EGF of the Bell numbers is in the Hayman class.*

- The power series $f(z) = \exp(z/(1-z))$, EGF of sets of (non-empty) lists, has $g(z) = z/(1-z)$ with $b_n = 1$, for any $n \geq 1$. The coefficients of the function $g$ satisfy the hypothesis of Theorem 4.3.17 and, thus, in particular, $f = e^g$ is in the Hayman class.

  In general, using a similar argument, we obtain that $f(z) = \exp(z/(1-z)^\gamma)$, for $\gamma > 0$, is in the Hayman class. But $f(z) = \exp(\ln(1/(1-z))) = 1/(1-z)$, the EGF of the sets of cycles is not even Gaussian, as shown in Subsection 4.3.1.

- The EGF of the sets of functions is $f(z) = e^g$, with $g(z) = \sum_{n=1}^{\infty}(n^n/n!)z^n$, for $|z| < 1/e$. In this case, the coefficients $b_n = n^n/n!$ of $g$ satisfy (4.3.19) with $R = 1/e$ and $\beta = \lambda = -1/2$, and therefore $f$ is in the Hayman class.

  In general the exponential of $\sum_{n=1}^{\infty}(n^{n-\alpha}/n!)z^n$ is in the Hayman class if $0 \leq \alpha < 1/2$. As we have seen before this power series is Gaussian if and only if $0 \leq \alpha < 1/2$, which is, of course, consistent with the previous result. Recall that power series in the Hayman class are strongly Gaussian, and therefore, Gaussian.

#### B. Sets of Unlabeled Classes

- For the OGF of the partitions of integers $P$, we have that $P = e^g$, with

$$\sum_{n=1}^{\infty} \frac{\sigma_1(n)}{n} z^n$$

  where $\sigma_1(n)$ denotes the sum of the divisors of the positive integer $n$. The coefficients $b_n = \sigma_1(n)/n$ satisfy in this case the inequality

$$1 \leq b_n \leq D_\varepsilon n^\varepsilon,$$

for each $\varepsilon > 0$ and some constant $D_\varepsilon > 0$. See Theorem 322 in [45]. (Actually we may bound from above with $\ln(\ln(n))$, see Theorem 323 in [45]). Therefore the function $g$ satisfies the hypothesis of Theorem 4.3.17 and this means in particular that $P = e^g$ is in the Hayman class.

The same argument gives that the infinite product $\prod_{j=1}^{\infty} 1/(1-z^j)^{c_j}$, where the $c_j$ are integers satisfying $1 \leq c_j \leq c$, for some constant $c > 1$, is in the Hayman class.

We collect the previous discussion by means of the following theorem.

**Theorem 4.3.19.** *The ordinary generating function of the partitions of integers $P$ is in the Hayman class.*

- For the OGF of partitions into distinct parts $Q$, given by

$$Q(z) \overset{(1)}{=} \prod_{j=1}^{\infty}(1+z^j) \overset{(2)}{=} \prod_{j=1}^{\infty} \frac{1}{1-z^{2j+1}}, \qquad \text{for any } |z| < 1,$$

we have that $Q = e^g$, where

$$g(z) \overset{(1)}{=} \sum_{k,j\geq 1} \frac{(-1)^{k+1}}{k} z^{kj} \overset{(2)}{=} \sum_{k\geq 1; j\geq 0} \frac{1}{k} z^{k(2j+1)} = \sum_{n=1}^{\infty} \frac{\sigma_1^{\text{odd}}(n)}{n} z^n,$$

where $\sigma_1^{\text{odd}}$ registers the sum of the odd divisors of $n$.

The coefficient $b_n = \sigma_1^{\text{odd}}(n)/n$ satisfy the inequality

$$\frac{1}{n} \leq b_n \leq D_\varepsilon n^\varepsilon \tag{4.3.20}$$

for each $\varepsilon > 0$ and some constant $D_\varepsilon > 0$. The inequality on the left holds simply because $1|n$, while the inequality on the right-hand side holds because of the inequality $\sigma_1^{\text{odd}} \leq \sigma_1(n)$ and Theorem 322 in [45].

These bounds, though, are not within the reach of Theorem 4.3.17, and for the function $Q$ we will verify directly the hypothesis of Theorem 4.3.13.

Using (4.3.20), Proposition 4.2.7 and formula

$$\sigma_f^2(t) = tg'(t) + t^2 g''(t) = \sum_{n=1}^{\infty} n^2 b_n t^n,$$

we have for the variance function $\sigma_Q^2(t)$ that

$$\frac{1}{(1-t)^2} = O(\sigma_Q^2(t)) \quad \text{and} \quad \sigma_Q^2(t) = O\left(\frac{1}{(1-t)^{3+\varepsilon}}\right), \qquad \text{as } t \uparrow 1,$$

and for the function $g$ also that $t^3 g'''(t) = O(1/(1-t)^{4+\varepsilon})$, as $t \uparrow 1$.

The variance condition (4.3.22) is obviously satisfied. We propose a cut $h(t) = (1-t)^\alpha$, where $\alpha \in (4/3, 3/2)$. The major arc condition (4.3.23) holds since $3\alpha > 4+\varepsilon$, for appropriate small $\varepsilon > 0$.

For $t \in (0,1)$ and $|\theta| \le \pi$, denote $z = te^{i\theta}$, then we have

$$\left|\frac{1+te^{i\theta}}{1+t}\right| = 1 + \frac{2t(\cos(\theta)-1)}{(1+t)^2} \le 1 + \frac{t}{2}(\cos(\theta)-1) \le e^{t(\cos(\theta)-1)/2} = e^{(\Re(z)-|z|)/2},$$

and so

$$\frac{|Q(z)|}{Q(|z|)} \le \exp\left(\frac{1}{4}\left(\Re\left(\frac{z}{1-z}\right) - \frac{|z|}{1-|z|}\right)\right) \le \exp\left(\frac{1}{4}\left(\Re\left(\frac{1}{1-z}\right) - \frac{1}{1-|z|}\right)\right)$$

and, consequently,

$$\Re(g(z)) - g(|z|) = \ln\left(\frac{|Q(z)|}{Q(|z|)}\right) \le \frac{1}{4}\left(\Re\left(\frac{1}{1-z}\right) - \frac{1}{1-|z|}\right) \le \left|\frac{1}{1-z}\right| - \frac{1}{1-|z|}.$$

Lemma 4.3.16 gives that $V(t) = C/(1-t)^3$ and $U(t) = D(1-t)$. Now, condition (4.3.13) is satisfied since $\alpha < 3/2$.

Thus the function $g$ satisfies the conditions of Theorem 4.3.13, and, in particular $Q = e^g$ is in the Hayman class.

We collect the previous discussion by means of the following theorem.

**Theorem 4.3.20.** *The ordinary generating function of the partitions of integers into distinct parts $Q$ is in the Hayman class.*

- For the OGF of the plane partitions, see [32, p. 580], we have

$$M(z) = \prod_{j=1}^{\infty} \frac{1}{(1-z^j)^j},$$

we have $M = e^g$, with $g$ given by

$$g(z) = \sum_{n=1}^{\infty} \frac{\sigma_2(n)}{n} z^n,$$

where $\sigma_2(n)$ denotes the sum of the squares of the divisors of the integer $n \ge 1$. For each $\varepsilon > 0$, there is a constant $C_\varepsilon > 0$, such that

$$n \le \frac{\sigma_2(n)}{n} \le C_\varepsilon n^{1+\varepsilon}, \qquad \text{for each } n \ge 1$$

This follows since $\sigma_2(n) \le n\sigma_1(n)$ and $\sigma_1(n) \le C_\varepsilon n^{1+\varepsilon}$. Thus we see that $g$ satisfies the condition on Theorem 4.3.17, and, in particular, we obtain that $M = e^g$ is in the Hayman class.

We collect the previous discussion by means of the following theorem.

**Theorem 4.3.21.** *The ordinary generating function of the plane partitions $M$ is in the Hayman class.*

Likewise, and more generally, we see that for integer $c \geq 0$, the OGF of the colored partitions,

$$\prod_{j=1}^{\infty} \frac{1}{(1-z^j)^c},$$

where each part $j$ appears in $j^c$ different colors, is in the Hayman class. Observe that $n^c \leq \sigma_{c+1}(n)/n \leq C_\varepsilon n^{c+\varepsilon}$, for $n \geq 1$.

### 4.3.5 Uniformly Hayman exponentials

Let $g(z) = \sum_{n=0}^{\infty} b_n z^n$, with $b_n \geq 0$, for $n \geq 0$, and radius of convergence $R > 0$. Let $f \in \mathcal{K}$ be given by $f = e^g$.

One of the main results of [18], see [17, Theorem 4.1] and [18, Theorem F] gives conditions on the powers series $g$ that guarantees that $f = e^g$ is in the Hayman class, and thus strongly Gaussian, and, therefore, amenable to the Hayman asymptotic formula, see Corollary (3.2.5).

It turns out that these same conditions on $g$ are enough for $f$ being uniformly Hayman; this is the content of Theorem 4.3.22.

Denote

$$\omega_g(t) \triangleq \frac{1}{6}(b_1 t + 8 b_2 t^2 + \frac{9}{2} t^3 g'''(t)), \quad \text{for } t \in (0, R)\,. \tag{4.3.21}$$

This is the upper bound in Corollary 4.2.5, here we give the constant $C > 0$ and the polynomial $P_2(z)$ explicitly. Observe that this bound is also valid for $g \in \mathcal{P}$ a polynomial.

**Theorem 4.3.22.** *Let $g$ be a non-constant power series with radius of convergence $R > 0$ and non-negative coefficients.*

*Assume that the variance condition is satisfied*

$$\lim_{t \uparrow R} \big(t g'(t) + t^2 g''(t)\big) = +\infty\,. \tag{4.3.22}$$

*Assume further that there is a cut function $h(t)$ satisfying*

$$\lim_{t \uparrow R} \omega_g(t) h(t)^3 = 0\,. \tag{4.3.23}$$

*and that there are positive functions $U, V$ defined in $(t_0, R)$, for some $t_0 \in (0, R)$, where $U$ takes values in $(0, \pi]$ and $V$ in $(0, \infty)$ and are such that*

$$\sup_{|\varphi| \geq \omega} \big(\Re g(t e^{\imath\varphi}) - g(t)\big) \leq -V(t)\,\omega^2\,, \quad \textit{for } \omega \leq U(t) \textit{ and } t \in (t_0, R)\,, \tag{4.3.24}$$

*and thus that the cut $h$ is such that*

$$h(t) \leq U(t), \textit{ for } t \in (t_0, 1) \quad \textit{and} \quad \lim_{t \uparrow R} \sigma_f(t)\, e^{-V(t) h(t)^2} = 0\,. \tag{4.3.25}$$

*then the function $f = e^g$ is uniformly Hayman.*

*Proof.* Condition (4.3.22) is directly condition (3.6.3).

We shall verify the conditions on the cuts of the definition of uniformly Hayman power series (3.6.1) and (3.6.4) with cuts $h(n,t)$ given by

$$h(n,t) = h(t)\, n^{-\beta}$$

where the parameter $\beta$ satisfies $1/3 < \beta < 1/2$.

From the discussion in Subsection 3.3.2 we have that

$$\left| \ln \mathbf{E}(e^{\imath\theta \breve{X}_t}) + \frac{\theta^2}{2} \right| \le \omega_g(t) \frac{|\theta|^3}{\sigma^3(t)}\,, \quad \text{for } t \in (0,R) \text{ and } \theta \in \mathbb{R}\,,$$

and, thus, that

$$\left| n \ln \mathbf{E}(e^{\imath\theta \breve{X}_t/\sqrt{n}}) + \frac{\theta^2}{2} \right| \le \omega_g(t) \frac{|\theta|^3}{\sigma^3(t)\sqrt{n}}\,, \quad \text{for } t \in (0,R) \text{ and } \theta \in \mathbb{R}\,.$$

For $|\theta| \le h(n,t)\sigma_f(t)\sqrt{n}$ we deduce that

$$\left| n \ln \mathbf{E}(e^{\imath\theta \breve{X}_t/\sqrt{n}}) + \frac{\theta^2}{2} \right| \le \omega_g(t)\, h(t)^3\, n^{1-3\beta}$$

Hypothesis (4.3.23) on $h(t)$ and the fact that $1 - 3\beta < 0$ gives us that

$$\lim_{[n\to\infty \vee t\uparrow R]} \omega_g(t)\, h(t)^3\, n^{1-3\beta} = 0\,,$$

and, thus, that condition (3.6.1) is satisfied.

Since $h(n,t) \le h(t) \le U(t)$, condition (4.3.24) gives us that

$$n \sup_{h(n,t)\le|\theta|\le\pi} \big(\Re g(te^{\imath\theta}) - g(t)\big) \le -V(t)h(t)^2 n^{1-2\beta}\,,$$

and so that

$$\sup_{h(n,t)\le|\theta|\le\pi} \left|\mathbf{E}(e^{\imath\theta X_t})^n\right| \le \exp\Big( - V(t)h(t)^2 n^{1-2\beta}\Big)\,.$$

Since $\lim_{t\uparrow R} V(t)h(t)^2 = \infty$, we have, for a certain $t_0 \in (0,R)$ that $V(t)h(t)^2 \ge 1$, for $t \in (t_0, R)$, and thus that

$$V(t)h(t)^2 n^{1-2\beta} \ge V(t)h(t)^2 + n^{1-2\beta}\,.$$

We deduce that

$$\sqrt{n}\sigma_f(t) \sup_{h(n,t)\le|\theta|\le\pi} \left|\mathbf{E}(e^{\imath\theta X_t})^n\right| \le \sqrt{n}e^{-n^{1-2\beta}}\, \sigma_f(t) e^{-V(t)h(t)^2}\,,$$

and, because of hypothesis (4.3.25) and since $\beta < 1/2$, that condition (3.6.4) of the definition of uniformly Hayman is satisfied. □

In [18] a large number of exponentials $f = e^g$ where $g$ is a power series with non-negative coefficients which satisfy the conditions of Theorem 4.3.22 are exhibited. For instance, the EGF of the Bell numbers, or the generating functions $P$ of partitions or $Q$ of partitions into distinct parts, are actually uniformly Hayman, and also are uniformly Hayman related examples like the EGF of sets of pointed sets, the EGF of sets of functions or the OGF of plane partitions or of some colored partitions. See also Chapter 5.

**Exponential of polynomials**

Let $g$ be a polynomial with non-negative coefficients $g(z) = \sum_{n=0}^N b_n z^n$ and of degree $N$, so that $b_N > 0$.

Assume that $Q_g = \gcd\{1 \le n \le N : b_n > 0\} = 1$. Then $f = e^g$ is in the Hayman. This is a particular case of a result of Hayman [48, Theorem X]. See [17, Proposition 5.1] for a simpler proof of this particular case.

We are going to show next that $f = e^g$ *is actually uniformly Hayman* with an argument similar to the one used to show in [17, Proposition 5.1] that $f = e^g$ is in the Hayman class.

Observe first that

$$\sigma_f^2(t) = tg'(t) + t^2 g''(t) \sim N^2 b_N t^N, \quad \text{as } t \to \infty\,.$$

Thus, the variance condition (3.6.3) of being uniformly Hayman is satisfied.

We have

$$\omega_g(t) = \frac{1}{6}\big(b_1 t + 8 b_2 t^2 + \frac{9}{2} t^3 g'''(t)\big) = O(t^N)\,, \quad \text{as } t \to \infty\,.$$

For cuts we propose $h(n,t) = h(t) n^{-\beta} = t^{-N\alpha} n^{-\beta}$, with $\alpha, \beta$ in the interval $(1/3, 1/2)$. For concreteness, we take $\alpha = \beta = 5/12$.

From the proof of Theorem 4.3.22 we have that

$$\begin{aligned}\Big| n \ln \mathbf{E}(e^{\imath\theta \breve{X}_t/\sqrt{n}}) + \frac{\theta^2}{2}\Big| &\le \omega_g(t) t^{-5N/4} n^{-1/4}\\ &= O\big(t^{-N/4} n^{-1/4}\big)\,, \quad \text{for } |\theta| \le h(n,t)\sigma_f(t)\sqrt{n}\,,\end{aligned}$$

and thus, we see that condition (3.6.1) is satisfied.

Now, the proof that $f = e^g$ is in the Hayman class of [17] gives $\eta \in (0,\pi)$ and $t_0 > 0$, depending on $g$, so that

$$\sup_{|\theta| > \omega} \big(\Re g(te^{\imath\theta}) - g(t)\big) \le -C_g \min\big\{t, t^N \omega^2\big\}, \quad \text{for } \omega \ge \eta \text{ and } t > t_0\,,$$

for some constant $C_g$ depending on $g$.

Thus for some $t_1 > t_0$, so that $h(t) < \eta$, for $t > t_1$, we have that

$$\begin{aligned} n \sup_{h(n,t) \le |\theta| \le \pi} \big(\Re g(te^{\imath\theta}) - g(t)\big) &\le -C \min\big\{tn, t^N t^{-5N/6} n^{1-5/6}\big\}\\ &= -C \min\big\{tn, t^{N/6} n^{1/6}\big\}\,.\end{aligned}$$

With $\delta = \min\{1, N/6\}$, we have, for $t \geq 1$, that

$$\min\{tn, t^{N/6}n^{1/6}\} \geq t^{\delta}n^{1/6} \geq (1/2)(t^{\delta} + n^{1/6}),$$

and thus, since $|\mathbf{E}(e^{\imath\theta X_t})| = e^{\Re g(te^{\imath\theta})-g(t)}$, that

$$\begin{aligned}\sqrt{n}\sigma_f(t) \sup_{h(n,t)\leq|\theta|\leq\pi} \big|\mathbf{E}(e^{\imath\theta X_t})^n\big| &= O\Big(\sqrt{n}t^{N/2}\exp\big(-C\min\{tn, t^{N/6}n^{1/6}\}\big)\Big) \\ &= O\Big(\sqrt{n}\exp\big(-(C/2)n^{1/6}\big)\Big)\,O\Big(t^{N/2}\exp\big(-(C/2)t^{\delta}\big)\Big),\end{aligned}$$

and, therefore, condition (3.6.2) is satisfied.

# Chapter 5

# Partitions of integers and partitions of sets

## Contents

In this chapter we discuss ordinary generating functions (OGFs) and exponential generating functions (EGFs) of different sort of partitions: partitions of sets and partitions of integers. For the partitions of integers we study the OGF of partitions of integers $P$, the OGF of the partitions into distinct parts $Q$, the OGF of the partitions with parts in an arithmetic progression $P_{a,b}$ and also the OGF of the plane and some colored partitions. Later on we study the EGFs of the Bell numbers, that is, partitions of non empty sets.

We want to apply the machinery developed in previous chapters to these OGFs. We will prove that the mentioned OGFs are in the Hayman class, and therefore that these power series are strongly Gaussian. Later on, using Baéz-Duarte asymptotic formula, see Theorem 3.2.8, we will provide with an asymptotic formula for the coefficients of each of the mentioned OGFs.

These asymptotic formulae of coefficients $a_n$ of partitions of integers turn out to be always of the form

$$a_n \sim \alpha \frac{1}{n^\beta} e^{\gamma n^\delta}, \qquad \text{as } n \to \infty,$$

for appropriate $\alpha, \beta, \gamma, \delta > 0$.

This chapter is mainly based on the papers:

- Maciá, V.J. et al. Khinchin families and Hayman class. *Comput. Methods Funct. Theory* **21** (2021), 851–904, see [17] ,
- Maciá, V.J. et al. Khinchin families, set constructions, partitions and exponentials. *Mediterr. J. Math.*, 21, 39 (2024), see [18].

## 5.1 Partitions of integers

In this section we give asymptotic formulas for the coefficients of the OGFs of different sort of partitions of integers. In particular we give, by means of Baéz-Duarte asymptotic formula, see Theorem 3.2.8, asymptotic formulas for the partitions, partitions into distinct parts, partitions with parts in an arithmetic progression, plane partitions and colored partitions.

### 5.1.1 Euler summation

We will use repeatedly along this chapter Euler summation's formula, we refer to [26] for the specifics of this summation formula.

For $f \in C^k[0, N]$, with $k \geq 1$, the Euler summation formula of order $k$ reads:

$$\sum_{j=0}^{N} f(j) = \int_0^N f(t)dt + \frac{1}{2} f(0) + \frac{1}{2} f(N) + \sum_{j=1}^{k-1} (-1)^{j+1} \frac{B_{j+1}}{(j+1)!} (f^{(j)}(N) - f^{(j)}(0)) + (-1)^{k+1} \int_0^N f^{(k)}(t) \frac{B_k(\{t\})}{k!} dt \tag{5.1.1}$$

Here $\{t\}$ denotes the fractional part of $t$, $B_j$, for $j \geq 0$, denotes the $j-th$ Bernoulli number (with $B_1 = -1/2$) while $B_j(x)$, for $j \geq 0$, stands for the Bernoulli polynomials. See, for instance, [26].

### 5.1.2 Approximation of Means and Variances

We collect here approximations of the mean and the variance functions. We will prove that these asymptotic approximations verify the conditions of Baéz-Duarte Theorem 3.2.8.

The precise values of the integrals in the following lemma will be invoked frequently later on.

**Lemma 5.1.1.** *We have*

*a)* $\int_0^\infty s^u \ln\left(\frac{1}{1-e^{-s}}\right) = \zeta(u+2)\Gamma(u+1)$, *for* $u > -1$,

*b)* $\int_0^\infty s^u \frac{e^{-s}}{1-e^{-s}} ds = \zeta(u+1)\Gamma(u+1)$, *for* $u > 0$,

*c)* $\int_0^\infty s^u \frac{e^{-s}}{(1-e^{-s})^2} ds = \zeta(u)\Gamma(u+1)$, *for* $u > 1$.

### A. OGF of the partitions of integers $P$

Consider first the OGF of the partitions of integers

$$P(z) = \prod_{j=1}^{\infty} \frac{1}{1-z^j}, \qquad \text{for any } |z| < 1,$$

with Khinchin family $(X_t)_{t\in[0,1)}$ and mean and variance functions

$$m_P(t) = \sum_{j=1}^{\infty} \frac{jt^j}{1-t^j} \quad \text{and} \quad \sigma_P^2(t) = \sum_{j=1}^{\infty} \frac{j^2 t^j}{(1-t^j)^2}, \qquad \text{for any } t \in (0,1).$$

We want an asymptotic formula for the mean and the variance of $P$, as $s \uparrow 1$. For convenience, we consider $m_P(e^{-s})$, for $s > 0$, and then as $s \downarrow 0$. We have

$$m_P(e^{-s}) = \frac{1}{s} \sum_{j=1}^{\infty} (js) \frac{e^{-js}}{1-e^{-js}}.$$

Let $\phi$ be the function

$$\phi(x) = \frac{xe^{-x}}{1-e^{-x}}, \qquad \text{for any } x > 0.$$

If we set $\phi(0) = 0$, then $\phi \in C^\infty[0, +\infty)$ and, by virtue of Lemma 5.1.1, we obtain that

$$\int_0^\infty \phi(s) ds = \zeta(2).$$

In addition we have that $\lim_{x\to+\infty} \phi(x) = 0$ and $\int_0^\infty |\phi'(s)| ds < +\infty$.

Fixing $k = 1$ in Euler summation formula (5.1.1) we obtain that

$$\begin{aligned} s m_P(e^{-s}) = \sum_{j=1}^{\infty} \phi(js) &= \int_0^\infty \phi(st) dt + s \int_0^\infty \phi'(st) B_1(\{t\}) dt \\ &= \frac{1}{s} \int_0^\infty \phi(x) dx + \int_0^\infty \phi'(x) B_1(\{x/s\}) dx\,, \end{aligned}$$

then there exists a constant $C > 0$ such that

$$\left| s m_P(e^{-s}) - \frac{1}{s}\int_0^\infty \phi(x)dx \right| \le C \int_0^\infty |\phi'(x)|dt < +\infty,$$

and therefore

$$s m_P(e^{-s}) - \frac{\zeta(2)}{s} = O(1)\,, \qquad \text{as } s \downarrow 0.$$

In fact we have $\lim_{s\downarrow 0} s^2 m_P(e^{-s}) = \zeta(2)$, and so

$$m_P(e^{-s}) \sim \frac{\zeta(2)}{s^2}, \qquad \text{as } s \downarrow 0.$$

For the variance function $\sigma_P^2(t)$ we obtain, analogously, also using Lemma 5.1.1, that

$$\sigma_P^2(e^{-s}) \sim \frac{\Gamma(3)\zeta(2)}{s^3}, \qquad \text{as } s \downarrow 0.$$

We collect the previous discussion by means of the following lemma.

**Lemma 5.1.2.** *Define $\tilde{m}_P(e^{-s}) = \pi^2/(6s^2)$ and $\tilde{\sigma}_P^2(e^{-s}) = \pi^2/(3s^3)$, for $s > 0$, then*

*a) $m_P(e^{-s}) \sim m_P(e^{-s})$, as $s \downarrow 0$,*

*b) $\sigma_P(e^{-s}) \sim \tilde{\sigma}_P(e^{-s})$, as $s \downarrow 0$,*

*c) $\frac{m_P(e^{-s}) - \tilde{m}_P(e^{-s})}{\tilde{\sigma}_P(e^{-s})} = O(\sqrt{s}) = o(1)$, as $s \downarrow 0$.*

### B. OGF of the partitions into distinct parts $Q$

Now we study the mean and variance of the OGF of the partitions into distinct parts $Q$. Recall that

$$Q(z) = \prod_{j=1}^{\infty}(1 + z^j) = \prod_{j=1}^{\infty} \frac{1}{1 - z^{2j+1}}, \qquad \text{for any } |z| < 1.$$

Applying Lemma 5.1.1 combined with the identity

$$\ln(1 + x) = \ln(1 - x^2) - \ln(1 - x), \qquad \text{for } x \in (0, 1),$$

and its first two derivatives, we find that

$$m_Q(e^{-s}) = \sum_{j=1}^{\infty} j \frac{e^{-js}}{1 + e^{-js}} \sim \frac{\zeta(2)}{2s^2}, \text{ as } s \downarrow 0$$

$$\sigma_Q^2(e^{-s}) = \sum_{j=1}^{\infty} j^2 \frac{e^{-sj}}{(1 + e^{-js})^2} \sim \frac{\zeta(2)}{s^3}, \text{ as } s \downarrow 0$$

and again, using Euler's summation formula (5.1.1) in a similar fashion, we obtain the following lemma.

**Lemma 5.1.3.** *Define* $\tilde{m}_Q(e^{-s}) = \zeta(2)/(2s^2)$ *and* $\tilde{\sigma}_Q^2(e^{-s}) = \zeta(2)/s^3$, *for* $s > 0$, *then*

*a)* $m_Q(e^{-s}) \sim m_Q(e^{-s})$, *as* $s \downarrow 0$,

*b)* $\sigma_Q(e^{-s}) \sim \tilde{\sigma}_Q(e^{-s})$, *as* $s \downarrow 0$,

*c)* $\frac{m_Q(e^{-s})-\tilde{m}_Q(e^{-s})}{\tilde{\sigma}_Q(e^{-s})} = O(\sqrt{s}) = o(1)$, *as* $s \downarrow 0$.

### C. OGF of the partitions with parts in arithmetic progression $P_{a,b}$

Now we study the mean and variance of the OGF of the partitions in the arithmetic progression $\{aj + b : j \geq 0\}$, for integers $a, b \geq 1$ that is, the OGF $P_{a,b}$. Recall that

$$P_{a,b}(z) = \prod_{j=1}^{\infty} \frac{1}{1 - z^{aj+b}}, \qquad \text{for any } |z| < 1.$$

For the mean and the variance we have, analogously,

$$m_{P_{a,b}}(e^{-s}) = \sum_{j=0}^{\infty} (aj + b) \frac{e^{-(aj+b)s}}{1 - e^{-(aj+b)s}} \sim \frac{\zeta(2)}{as^2}, \qquad \text{as } s \downarrow 0, \tag{5.1.2}$$

$$\sigma^2_{P_{a,b}}(e^{-s}) = \sum_{j=0}^{\infty} (aj + b)^2 \frac{e^{-(aj+b)s}}{(1 - e^{-(aj+b)s})^2} \sim \frac{2\zeta(2)}{as^3}, \qquad \text{as } s \downarrow 0, \tag{5.1.3}$$

again, applying Euler's summation formula (5.1.1), we find the following lemma.

**Lemma 5.1.4.** *Define* $\tilde{m}_{P_{a,b}}(e^{-s}) = \zeta(2)/(as^2)$ *and* $\tilde{\sigma}^2_{P_{a,b}}(e^{-s}) = (2\zeta(2))/(as^3)$, *for* $s > 0$. *Then*

*a)* $m_{P_{a,b}}(e^{-s}) \sim \tilde{m}_{P_{a,b}}(e^{-s})$, *as* $s \downarrow 0$,

*b)* $\sigma_{P_{a,b}}(e^{-s}) \sim \tilde{\sigma}_{P_{a,b}}(e^{-s})$, *as* $s \downarrow 0$,

*c)* $\frac{m_{P_{a,b}}(e^{-s})-\tilde{m}_{P_{a,b}}(e^{-s})}{\tilde{\sigma}_{P_{a,b}}(e^{-s})} = O(\sqrt{s}) = o(1)$, *as* $s \downarrow 0$.

### D. Generalized OGF $W_a^b$ (including colored and plane partitions)

For integers $a \geq 1$ and $b \geq 0$ consider the infinite product

$$W_a^b(z) = \prod_{j=1}^{\infty} \left( \frac{1}{1 - z^{j^a}} \right)^{j^b}, \qquad \text{for any } |z| < 1.$$

We study the mean and variance associated to $W_{a,b}$. Consider the real function

$$\phi(x) = x^b \frac{x^a e^{-x^a}}{1 - e^{-x^a}}, \qquad \text{for } x \geq 0,$$

with $\phi(0) = 0$. The derivatives of order less than $b$ of this function vanish at $x = 0$ and $x = +\infty$, and

$$\int_0^\infty \phi(x)dx = \frac{1}{a}\zeta\left(1 + \frac{b+1}{a}\right)\Gamma\left(1 + \frac{b+1}{a}\right).$$

Arguing as in the previous cases, and using Euler's summation formula (5.1.1) of order $b+1$, we get that

$$m_{W_a^b}(e^{-s}) = \frac{1}{a}\zeta\left(1 + \frac{b+1}{a}\right)\Gamma\left(1 + \frac{b+1}{a}\right)\frac{1}{s^{1+(b+1)/a}} + O(1/s), \qquad \text{as } s \downarrow 0.$$

Analogously,

$$\sigma^2_{W_a^b}(e^{-s}) \sim \frac{1}{a}\zeta\left(1 + \frac{b+1}{a}\right)\Gamma\left(2 + \frac{b+1}{a}\right)\frac{1}{s^{2+(b+1)/a}}, \qquad \text{as } s \downarrow 0.$$

For any $s > 0$, denote

$$\tilde{m}_{W_a^b}(e^{-s}) = \frac{1}{a}\zeta\left(1 + \frac{b+1}{a}\right)\Gamma\left(1 + \frac{b+1}{a}\right)\frac{1}{s^{1+(b+1)/a}},$$

and

$$\tilde{\sigma}^2_{W_a^b}(e^{-s}) = \frac{1}{a}\zeta\left(1 + \frac{b+1}{a}\right)\Gamma\left(2 + \frac{b+1}{a}\right)\frac{1}{s^{2+(b+1)/a}}.$$

We collect the previous discussion by means of the following lemma.

**Lemma 5.1.5.** *For the functions $\tilde{m}_{W_a^b}(e^{-s})$ and $\tilde{\sigma}^2_{W_a^b}(e^{-s})$, we have*

*a) $m_{W_a^b}(e^{-s}) \sim \tilde{m}_{W_a^b}(e^{-s})$, as $s \downarrow 0$,*

*b) $\sigma_{W_a^b}(e^{-s}) \sim \tilde{\sigma}_{W_a^b}(e^{-s})$, as $s \downarrow 0$,*

*c) $\frac{m_{W_a^b}(e^{-s}) - \tilde{m}_{W_a^b}(e^{-s})}{\tilde{\sigma}_{W_a^b}(e^{-s})} = O(s^{(b+1)/(2a)}) = o(1)$, as $s \downarrow 0$.*

### 5.1.3 Asymptotic for the OGFs of partitions on (0,1)

In this subsection we give asymptotic formulas for the OGFs of different kind of partitions, evaluated at $z = e^{-s}$, and as $s \downarrow 0$.

**A. OGF of the partitions of integers $P$**

We start with the OGF of partitions $P$. We provide a proof since its ingredients are to be used below to handle other OGFs of partitions.

Fix $s > 0$, we want to apply Euler's summation formula (5.1.1) of order 2 to the function $f(x) = -\ln(1 - e^{-sx})$. Observe that $\lim_{x\to\infty} f(x) = \lim_{x\to\infty} f'(x) = 0$. We write

$$\ln(P(e^{-s})) = \sum_{j=1}^{\infty} f(j) = \int_1^{\infty} f(x)dx + \frac{1}{2}f(1) - \frac{1}{12}f'(1) - \int_1^{\infty} f''(x)\frac{B_2(\{x\})}{2!}dx$$

where $B_2(y)$ denote the second Bernoulli polynomial (the sum starts at $j = 1$, then the formula is modified to fit with these indices).

Now, applying Lemma 5.1.1, we find that

$$\begin{aligned}\int_1^{\infty} \ln\left(\frac{1}{1-e^{-sx}}\right) dx &= \frac{1}{s}\int_s^{\infty} \ln\left(\frac{1}{1-e^{-x}}\right) dx \\ &= \frac{1}{s}\int_0^{\infty} \ln\left(\frac{1}{1-e^{-x}}\right) dx - \frac{1}{s}\int_0^{s} \ln\left(\frac{1}{1-e^{-x}}\right) dx \\ &= \frac{\zeta(2)}{s} + \ln(s) - 1 + O(s), \qquad \text{as } s \downarrow 0.\end{aligned}$$

Also,

$$\frac{1}{2}f(1) - \frac{1}{12}f'(1) = -\frac{1}{2}\ln(s) + \frac{1}{12} + O(s), \qquad \text{as } s \downarrow 0.$$

and

$$\int_1^{\infty} f''(x)\frac{B_2(\{x\})}{2!}dx = \int_1^{\infty} \frac{(sx)^2 e^{sx}}{(e^{sx}-1)^2}\frac{B_2(\{x\})}{2!}\frac{1}{x^2}dx.$$

Since the function $y \mapsto (y^2 e^y)/(e^y - 1)^2$ is bounded in $[0, +\infty)$ and tends to 1 as $y \downarrow 0$, dominated convergence gives that

$$\int_1^{\infty} f''(x)\frac{B_2(\{x\})}{2!}dx = \int_1^{\infty} \frac{B_2(\{x\})}{2!}\frac{1}{x^2}dx + o(1), \text{ as } s \downarrow 0.$$

Thus,

$$\ln(P(e^{-s})) = \frac{\zeta(2)}{s} + \frac{1}{2}\ln(s) - 1 + \frac{1}{12} - \int_1^{\infty} \frac{B_2(\{x\})}{2!}\frac{1}{x^2}dx + o(1), \qquad \text{as } s \downarrow 0.$$

We may identify the constant term in the previous expression by using Euler's summation of order 2 between 1 and $N$ applied to the function $\ln(x)$ to obtain this way Stirling's approximation in the following precise (and standard) form

$$\begin{aligned}\ln(N!) &= N\ln(N) - N + 1 + \frac{1}{2}\ln(N) + \frac{1}{12}\left(\frac{1}{N} - 1\right) + \int_1^N \frac{B_2(\{x\})}{2!}\frac{1}{x^2}dx \\ &= \ln\left(\frac{N^N}{\sqrt{N}e^N}\right) + 1 - \frac{1}{12} + \frac{1}{12}\frac{1}{N} + \int_1^N \frac{B_2(\{x\})}{2!}\frac{1}{x^2}dx\end{aligned}$$

then

$$1 - \frac{1}{12} + \int_1^\infty \frac{B_2(\{x\})}{2!}\frac{1}{x^2}dx = \ln(\sqrt{2\pi}),$$

and

$$\ln(P(e^{-s})) = \frac{\zeta(2)}{s} + \frac{1}{2}\ln(s) - \ln(\sqrt{2\pi}) + o(1), \qquad \text{as } s \downarrow 0.$$

Finally, for the ordinary generating function of the partitions of integers $P$ we have the well-known asymptotic formula in the interval $(0,1)$.

**Lemma 5.1.6.** *For the ordinary generating function of the partitions of integers $P$ we have*

$$\ln(P(e^{-s})) = \frac{\zeta(2)}{s} + \frac{1}{2}\ln(s) - \ln(\sqrt{2\pi}) + o(1), \qquad \text{as } s \downarrow 0\,,$$

*and thus*

$$P(e^{-s}) \sim \frac{1}{\sqrt{2\pi}}\sqrt{s}e^{\pi^2/(6s)}, \qquad \text{as } s \downarrow 0\,.$$

See Hardy-Ramanujan [44], Section 3.2 (the formula for $g(x)$ there should have $(1-x)^{3/2}$ instead of $\sqrt{1-x}$). See also, de Bruijn [26], Ex. 3 of Chapter 3, or, even, Báez-Duarte [6], formula (2.21).

### B. OGF of the partitions into distinct parts $Q$

Consider now the ordinary generating function of the partitions into distinct parts $Q$. Recall that

$$Q(z) = \prod_{j=1}^\infty (1 + z^j) = \prod_{j=1}^\infty \frac{1}{1 - z^{2j+1}}, \qquad \text{for any } |z| < 1.$$

From the infinite product for $Q$ we deduce that

$$Q(z)P(z^2) = P(z), \qquad \text{for any } |z| < 1,$$

then

$$Q(e^{-s}) = P(e^{-s})/P(e^{-2s}), \qquad \text{as } s \downarrow 0.$$

Using Lemma 5.1.6 we obtain that

$$Q(e^{-s}) \sim \frac{1}{\sqrt{2}}e^{\pi^2/(12s)} \qquad \text{as } s \downarrow 0.$$

We collect the previous discussion by means of the following lemma.

**Lemma 5.1.7.** *For the ordinary generating function of the partitions of integers into distinct parts $Q$ we have*

$$Q(e^{-s}) \sim \frac{1}{\sqrt{2}} e^{\pi^2/(12s)}, \qquad \textit{as } s \downarrow 0.$$

### C. OGF of the partitions with parts in arithmetic progression $P_{a,b}$

For the ordinary generating function of the partitions with parts in the arithmetic progression $\{aj+b : j \geq 0\}$, with $a, b \geq 1$, we have

$$P_{a,b}(z) = \prod_{j=0}^{\infty} \frac{1}{1 - z^{aj+b}}, \qquad \text{for any } |z| < 1.$$

If $\gcd(a,b) > 1$, then $P_{a,b}$ is not strongly Gaussian (although is Gaussian and $d$-strongly Gaussian, see Subsection 3.2.4) since we can write $P_{a,b}(z) = G(z^d)$, for some holomorphic function in $\mathbb{D}$, and thus the coefficients of $P_{a,b}$ satisfy no asymptotic formula.

In fact, if $\gcd(a,b) = d > 1$, then we have

$$P_{a,b}(z) = P_{a',b'}(z^d)\,, \quad \text{for } z \in \mathbb{D}\,,$$

where $a' = a/d$ and $b' = b/d$. Observe that $\gcd(a', b') = 1$.

Now we prove that $P_{a,b}$, with $a, b \geq 1$ and $\gcd(a,b) = 1$ is in the Hayman class. In this case we cannot use the coefficient criteria, Theorem 4.3.17, for exponential power series.

We can write $P_{a,b} = \exp(g)$ where

$$g(z) = \sum_{k\geq 1; j \geq 0} \frac{1}{k} z^{k(aj+b)} = \sum_{k=1}^{\infty} \frac{1}{k} \frac{z^{bk}}{1 - z^{ak}} = \sum_{k=1}^{\infty} \frac{1}{k} U(z^k)\,, \quad \text{for } z \in \mathbb{D}\,,$$

where $U$ is the rational function

$$U(z) = \frac{z^b}{1 - z^a}\,,$$

holomorphic in $\mathbb{D}$. We want to apply Theorem 4.3.12 to $P_{a,b}$, in order to apply this result we need some preliminary lemmas.

**Lemma 5.1.8.** *Let $\omega \in (0, \pi)$, and $t \in (0,1)$. Then*

$$\sup_{\omega \leq |\theta| \leq \pi} \Re\Big(\frac{1-t}{1-te^{\imath\theta}}\Big) = \Re\Big(\frac{1-t}{1-te^{\imath\omega}}\Big) \quad \textit{and} \quad \sup_{\omega \leq |\theta| \leq \pi} \Big|\frac{1-t}{1-te^{\imath\theta}}\Big| = \Big|\frac{1-t}{1-te^{\imath\omega}}\Big|\,.$$

**Lemma 5.1.9.** *Let $h(t)$ be defined in the interval $(0,1)$ with values in $(0,\pi)$ and such that $\lim_{t\uparrow 1} h(t)/(1-t) = 0$. Let $k$ be an integer $k \geq 1$. Then we have*

$$\lim_{t\uparrow 1} \frac{(1-t)^2}{h(t)^2} \Big( \Big|\frac{1-t}{1-te^{\imath h(t)}}\Big|^k - 1\Big) = -\frac{k}{2}\,.$$

Observe that the Taylor coefficients of $U(z) = z^b/(1-z^a)$ around $z = 0$ are non negative.

We search now for a cut function $h(t) = (1-t)^\alpha$ with appropriate $\alpha > 0$. Since $g'''(t) = O(1/(1-t)^4)$, as $t \uparrow 1$, the condition (4.3.23) of Theorem 4.3.12 for the major arc is satisfied if $\alpha > 4/3$.

Next, the minor arc. This is where we use the condition that $\gcd(a, b) = 1$. We have

$$\Re g(z) - g(|z|) = \sum_{k=1}^{\infty} \frac{1}{k}\big(\Re(U(z^k)) - U(|z|^k) \le \Re(U(z)) - U(|z|)\big)\,, \quad \text{for } z \in \mathbb{D}\,.$$

The inequality holds since all the summands in the series above are non positive.

The function $U(z)$ has simple poles at the $a$-th roots of unity: $\gamma_j = e^{2\pi \imath j/a}$ with $0 \le j < a$. The residue of $U$ at $\gamma_j$ is $-\gamma_j^{b+1}/a$.

Let us denote by $\mathcal{R}_j$ the region $\mathcal{R}_j = \big\{z = re^{\imath\phi} : 1/2 \le r < 1, |\phi - 2\pi j/a| \le \pi/a\big\}$ and let $\widetilde{\mathcal{R}} = \bigcup_{j=1}^{a-1} \mathcal{R}_j$; observe that we do not include $\mathcal{R}_0$ in $\widetilde{\mathcal{R}}$.

For a certain constant $M > 0$ we have that

$$(\flat) \quad \left|U(z) + \frac{\gamma_j^{b+1}}{a(z-\gamma_j)}\right| \le M\,, \quad \text{for } z \in \mathcal{R}_j \text{ and } 0 \le j < a\,.$$

We start analyzing $\Re U(z)$ for $z \in \widetilde{\mathcal{R}}$.

**Lemma 5.1.10.** *For $0 < j < a$, we have that*

$$\limsup_{t\uparrow 1} \sup_{\substack{z\in\mathcal{R}_j;\\ t\le|z|\le 1}} \Re\Big(U(z)\frac{1-|z|^a}{|z|^b}\Big) \le \frac{1}{2}(1 + \cos(2\pi j\, b/a))\,.$$

*Besides, for certain $\delta \in (0,1)$ and $t_0 \in (1/2, 1)$ we have that*

$$\Re U(z)\frac{1-|z|^a}{|z|^b} \le \delta\,, \quad \textit{for every } z \in \widetilde{\mathcal{R}} \textit{ with } |z| > t_0\,.$$

*Proof.* Fix $0 < j < a$ and $z \in \mathcal{R}_j$. Let $w \in \mathcal{R}_0$ such that $z = \gamma_j w$. By appealing to $(\flat)$ we have that

$$(\flat\flat) \quad \Re U(z)\frac{1-|z|^a}{|z|^b} \le \frac{1-|z|^a}{a(1-|z|)|z|^b}\,\Re\Big(\gamma_j^b\frac{1-|w|}{1-w}\Big) + M\frac{1-|z|^a}{|z|^b}\,.$$

Now, for every $w \in \mathbb{D}$ one has that

$$\frac{1-|w|}{1-w} \in \mathbb{D}\Big(\frac{1}{2}, \frac{1}{2}\Big) \cup \{1\}\,,$$

and, in particular,

$$\Re\big(e^{\imath\eta}\frac{1-|w|}{1-w}\big) \le \frac{1}{2}\big(1 + \cos(\eta)\big)\,, \quad \text{for all } \eta \in \mathbb{R} \text{ and } w \in \mathbb{D}\,.$$

Therefore, we obtain from (♭♭) that

$$(\flat\flat\flat)\qquad \Re\Big(U(z)\frac{1-|z|^a}{|z|^b}\Big) \le \frac{1-|z|^a}{a(1-|z|)|z|^b}\,\delta_j + M\frac{1-|z|^a}{|z|^b}\,.$$

with

$$\delta_j = \frac{1}{2}\big(1+\cos(2\pi\imath jb/a)\big)\,.$$

This bound (♭♭♭) gives the lim sup of the statement. The bound now follows since $\gcd(a,b)=1$ and $0<j<a$, implies that $\cos(2\pi\imath jb/a)<1$. □

As a consequence of Lemma 5.1.10 we have

$$(\natural)\quad \Re U(z) - U(t) \le (\delta-1)U(t) \le -\frac{1-\delta}{a}\,\frac{1}{1-t}\,,\qquad \text{for } z\in\widetilde{\mathcal{R}} \text{ with } |z|=t\in(t_0,1)\,.$$

It remains to analyze $\Re U(z)$ in the region $\mathcal{R}_0$. For $z\in\mathcal{R}_0$ we have that

$$\Big|U(z)-\frac{1}{a(1-z)}\Big| \le M\,,$$

and, therefore, for $z\in\mathcal{R}_0$ with $|z|=t$, we have that

$$\Re U(z) - U(t) \le M + \frac{1}{a}\Big(\Re\frac{1}{1-z}-\frac{1}{1-t}\Big) + \frac{1}{a(1-t)} - \frac{t^b}{(1-t^a)}\,.$$

Now, $\dfrac{1}{a(1-t)} - \dfrac{t^b}{(1-t^a)} \le \dfrac{b}{a}$, for each $t\in(0,1)$, and thus from Lemmas 5.1.8 and 5.1.9 and arguing as in the case above of $P$ and $Q$ and increasing the $t_0$ above if necessary we deduce that for certain $\delta>0$,

$$(\natural\natural)\qquad \sup_{\substack{z=te^{\imath\theta}\in\mathcal{R}_0;\\ h(t)\le|\theta|}} \Big(\Re U(z)-U(t)\Big) \le M + \frac{b}{a} - \delta\frac{h(t)^2}{(1-t)^3}\,,\qquad \text{for } t\in(t_0,1)\,.$$

The minor arc for $\{|z|=t\}$ comprises $\widetilde{\mathcal{R}}\cap\{|z|=t\}$ (where we have the bound (♮)) and the two subarcs in $\mathcal{R}_0\cap\{|z|=t\}$ of those $z=te^{\imath\theta}$ where $|\theta|\ge h(t)$.

From (♮) and (♮♮) and the moderate growth of $\sigma$, see (5.1.2), we conclude that $\alpha<3/2$ suffices to guarantee the minor arc condition of Theorem 4.3.12, and consequently that, under the assumption that $\gcd(a,b)=1$, the ogf $P_{a,b}$ is strongly Gaussian and that its coefficients satisfy Hayman's formula.

We collect the previous discussion by means of the following theorem.

**Theorem 5.1.11.** *Let $a,b\ge 1$ with $\gcd(a,b)=1$, then the ordinary generating function of the partitions in arithmetic progression $P_{a,b}$ is in the Hayman class.*

By a similar argument to that used for $P$ we obtain the asymptotic formula

**Theorem 5.1.12.** *For the ogf $P_{a,b}$ of partitions with parts in the arithmetic progression $\{aj+b; j \geq 0\}$ with integers $a, b \geq 1$ that*

$$P_{a,b}(e^{-s}) \sim \frac{1}{\sqrt{2\pi}}\,\Gamma\Big(\frac{b}{a}\Big)(as)^{b/a-1/2}\exp\Big(\frac{\zeta(2)}{as}\Big)\,, \quad \text{as } s \downarrow 0\,. \tag{5.1.4}$$

Recall that $P_{1,1} \equiv P$ and $P_{2,1} \equiv Q$, and then the previous theorem contains the asymptotic formulas for $P$ and $Q$.

### D. Generalized OGF $W_a^b$ (including colored and plane partitions)

Now we turn to the ogf $W_a^b$, with integers $a \geq 1$ and $b \geq 0$.

Fix $s > 0$. We will apply Euler summation of order $b+1$ to the function

$$f(x) = x^b \ln \frac{1}{1 - e^{-sx^a}}\,, \quad \text{for } x > 0\,.$$

Note that

$$\ln W_a^b(e^{-s}) = \sum_{j=1}^{\infty} f(j)\,.$$

Observe that

$$\begin{aligned}
\int_1^\infty f(x)dx &= \frac{1}{a}\frac{1}{s^{(b+1)/a}}\int_s^\infty y^{\frac{b+1}{a}} \ln\frac{1}{1-e^{-y}}\,\frac{dy}{y} = \\
&= \frac{1}{a}\frac{1}{s^{(b+1)/a}}\int_0^\infty y^{\frac{b+1}{a}} \ln\frac{1}{1-e^{-y}}\,\frac{dy}{y} - \frac{1}{a}\frac{1}{s^{(b+1)/a}}\int_0^s y^{\frac{b+1}{a}} \ln\frac{1}{1-e^{-y}}\,\frac{dy}{y} \\
&= \frac{1}{a}\zeta\big(1+\frac{b+1}{a}\big)\,\Gamma\Big(\frac{b+1}{a}\Big)\frac{1}{s^{(b+1)/a}} - \frac{1}{b+1}\ln\frac{1}{s} - \frac{a}{(b+1)^2} + O(s)\,.
\end{aligned}$$

For the function $f$ at $\infty$ we have that $\lim_{x\to\infty} f^{(j)}(x) = 0$, for $j \geq 0$. We need the values of $f$ and its derivatives up to order $b+1$ at $x = 1$. We have

$$f(1) = \ln\frac{1}{1-e^{-s}} = \ln\frac{1}{s} + O(s)\,.$$

Write $f(x) = \dfrac{1}{s^{b/a}} g(s^{1/a}x)$, where

$$g(y) = y^b \ln\frac{1}{1-e^{-y^a}} = \underbrace{y^{b+a} + y^b \ln\frac{y^a}{e^{y^a}-1}}_{=\omega(y)} - a\,\underbrace{(y^b \ln y)}_{=\eta(y)}\,.$$

Observe that $f^{(j)}(1) = \dfrac{1}{s^{(b-j)/a}} g^{(j)}(s^{1/a})$, for $j \geq 1$.

The function $\omega$ is holomorphic near 0, and its Taylor expansion starts with $(1/2)y^{b+a}$, while

$$\eta^{(j)}(y) = y^{b-j}\frac{b!}{(b-j)!}\Big(\ln y + \sum_{i=1}^{j}\frac{1}{b+1-i}\Big), \quad \text{for } 1 \le j \le b \text{ and } y > 0\,.$$

Also, $\eta^{(b+1)}(y) = \dfrac{b!}{y}$ and $\eta^{(b+2)}(y) = -\dfrac{b!}{y^2}$, for $y > 0$. And, besides,

$$f^{(b+1)}(1) = s^{1/a}g^{(b+1)}(s^{1/a}) = s^{1/a}\omega^{(b+1)}(s^{1/a}) - as^{1/a}\eta^{(b+1)}(s^{1/a}) = -b!\,a + O(s)\,.$$

With all this, we have that

$$\begin{aligned}\frac{1}{2}(f(1)+f(\infty))+&\sum_{j=1}^{b+1}(-1)^{j+1}\frac{B_{j+1}}{(j+1)!}(f^{(j)}(\infty)-f^{(j)}(1))\\ &=\frac{1}{2}\ln\frac{1}{s}+\Big(\sum_{j=1}^{b}(-1)^{j+1}\frac{B_{j+1}}{(j+1)!}\frac{b!}{(b-j)!}\Big)\ln s\\ &+a\sum_{j=1}^{b}(-1)^{j+1}\frac{B_{j+1}}{(j+1)!}\frac{b!}{(b-j)!}\sum_{i=1}^{j}\frac{1}{b+1-i}+a(-1)^{b+2}\frac{B_{b+2}}{(b+2)!}b!+O(s)\,.\end{aligned}$$

Applying Euler summation of order $b+1$ to $x^b$, with $b \ge 1$, just in the interval $[0,1]$ and taking into account that its $b$-th derivative is a constant ($b!$, actually), we deduce that

$$1 = \frac{1}{b+1}+\frac{1}{2}+\sum_{j=1}^{b-1}(-1)^{j+1}\frac{B_{j+1}}{(j+1)!}\frac{b!}{(b-j)!}\,.$$

Thus,

$$\sum_{j=1}^{b}(-1)^{j+1}\frac{B_{j+1}}{(j+1)!}\frac{b!}{(b-j)!} = (-1)^{b+1}\frac{B_{b+1}}{b+1}+\frac{1}{2}-\frac{1}{b+1}\,,$$

for $b \ge 1$, and for $b = 0$ also.

So, we may simplify and write

$$\begin{aligned}&\frac{1}{2}(f(1)+f(\infty))+\sum_{j=1}^{b+1}(-1)^{j+1}\frac{B_{j+1}}{(j+1)!}(f^{(j)}(\infty)-f^{(j)}(1))\\ =&\Big(-\frac{1}{b+1}+(-1)^{b+1}\frac{B_{b+1}}{b+1}\Big)\ln s\\ &+a\sum_{j=1}^{b}(-1)^{j+1}\frac{B_{j+1}}{(j+1)!}\frac{b!}{(b-j)!}\sum_{i=1}^{j}\frac{1}{b+1-i}+a(-1)^{b+2}\frac{B_{b+2}}{(b+2)!}b!+O(s)\,.\end{aligned}$$

Finally,

$$(-1)^{b+3}\int_1^\infty f^{(b+2)}(x)\frac{B_{b+2}(\{x\})}{(b+2)!}dx = (-1)^{b+3}\int_1^\infty (s^{1/a}x)^2 g^{(b+2)}(s^{1/a}x)\frac{B_{b+2}(\{x\})}{(b+2)!}\frac{dx}{x^2}\,.$$

The function $y \mapsto y^2 g^{(b+2)}(y)$ is bounded in $(0,\infty)$ and as $y \downarrow 0$ tends to $b!\,a$. We deduce from dominated convergence that

$$(-1)^{b+3}\int_1^\infty f^{(b+2)}(x)\frac{B_{b+2}(\{x\})}{(b+2)!}dx = (-1)^{b+3} b!\,a\int_1^\infty \frac{B_{b+2}(\{x\})}{(b+2)!}\frac{dx}{x^2} + o(1)\,.$$

With all this we have

$$\begin{aligned}\ln W_a^b(e^{-s}) &= \frac{1}{a}\zeta\big(1+\frac{b+1}{a}\big)\,\Gamma\big(\frac{b+1}{a}\big)\frac{1}{s^{(b+1)/a}}\\ &+ (-1)^{b+1}\frac{B_{b+1}}{b+1}\ln s\\ &- \frac{a}{(b+1)^2} + a\sum_{j=1}^{b}(-1)^{b+1}\frac{B_{j+1}}{(j+1)!}\frac{b!}{(b-j)!}\sum_{i=1}^{j}\frac{1}{b+1-i} + a(-1)^{b+2}\frac{B_{b+2}}{(b+2)!}b!\\ &+ (-1)^{b+3}b!\,a\int_1^\infty \frac{B_{b+2}(\{x\})}{(b+2)!}\frac{dx}{x^2} + o(1)\end{aligned}$$

The constant term of the above expression may be compacted by appealing to the so called (generalized) Glaisher-Kinkelin constants appearing (in fact, defined by) the asymptotic formula of (generalized) hyperfactorials.

If we apply Euler summation of order $b+2$ to the function $x^b \ln x$ between 1 and $N$, we obtain that

$$(\ddagger)\quad\begin{aligned}\sum_{n=1}^{N} n^b \ln n &= \frac{1}{b+1}N^{b+1}\ln N - \frac{1}{(b+1)^2}N^{b+1} + \frac{1}{2}N^b \ln N\\ &+ \sum_{j=1}^{b}(-1)^{j+1}\frac{B_{j+1}}{(j+1)!}\frac{b!}{(b-j)!}N^{b-j}\ln N + \sum_{j=1}^{b-1}(-1)^{j+1}\frac{B_{j+1}}{(j+1)!}\frac{b!}{(b-j)!}N^{b-j}\sum_{i=1}^{j}\frac{1}{b-i+1}\\ &+ (-1)^{b+1}\frac{B_{b+1}}{b+1}H_b + \frac{1}{(b+1)^2} - \sum_{j=1}^{b}(-1)^{j+1}\frac{B_{j+1}}{(j+1)!}\frac{b!}{(b-j)!}\sum_{i=1}^{j}\frac{1}{b-i+1}\\ &- (-1)^{b+2}\frac{B_{b+2}}{(b+2)!}b! - (-1)^{b+3}\int_1^\infty \frac{b!}{(b+2)!}B_{b+2}(\{x\})\frac{dx}{x^2} + O(1/N)\,.\end{aligned}$$

Here, $H_b$ denotes the $b$-th harmonic number.

This expression ($\ddagger$) is the analogue of Stirling's formula but now for the hyperfactorials $\prod_{j=1}^N j^{j^b}$ of order $b$; the factorials being the case $b=0$. The constant term (the sum of the terms not depending on $N$) in the above expression is $\ln GK_b$. These $GK_b, b \geq 1$ are the generalized Glaisher-Kinkelin constants introduced by Bendersky, [8], and identified in the following closed form

$$(\star)\quad \ln GK_b = \frac{(-1)^{b+1}B_{b+1}}{b+1}H_b - \zeta'(-b)\,,\quad \text{for } b \geq 1$$

by Choudhury, [20] and Adamchik, [3]. See also, Wang, [99]. The proper (original) Glaisher-Kinkelin constant is $GK_1$. By taking $H_0 = 0$, formula ($\star$) is also valid for $b=0$. The constant $GK_0$ is actually $\sqrt{2\pi}$.

We collect the previous discussion by means of the following theorem.

**Theorem 5.1.13.** *For the ordinary generating function $W_a^b$ with integers $a \geq 1$ and $b \geq 0$ we have*

$$\begin{aligned}\ln W_a^b(e^{-s}) = \frac{1}{a}\zeta\Big(1+\frac{b+1}{a}\Big)\,\Gamma\Big(\frac{b+1}{a}\Big)\frac{1}{s^{(b+1)/a}} - \zeta(-b)\ln s \\ + a\,\zeta'(-b) + o(1)\,, \quad \textit{as } s \downarrow 0\,,\end{aligned} \tag{5.1.5}$$

*where for compactness we have used that $\zeta(-b) = (-1)^b \dfrac{B_{b+1}}{b+1}$, for $b \geq 0$, and therefore*

$$W_a^b(e^{-s}) \sim \frac{e^{a\,\zeta'(-b)}}{s^{\zeta(-b)}}\exp\left(\frac{1}{a}\zeta\Big(1+\frac{b+1}{a}\Big)\,\Gamma\Big(\frac{b+1}{a}\Big)\frac{1}{s^{(b+1)/a}}\right), \qquad \textit{as } s \downarrow 0. \tag{5.1.6}$$

### 5.1.4 Asymptotic for the partitions of integers

In this section we give asymptotic formulas for different families of partitions of integers. We already know, see Chapter 4, in particular the applications of Theorem 4.3.17 in there, that some of the ordinary generating functions of the partitions of integers (with different restrictions) are in the Hayman class, and therefore that these OGFs are strongly Gaussian power series.

Now we are in position to apply Baéz-Duarte asymptotic formula, see Theorem 3.2.8: notice that we previously found approximations for the mean and variance functions of each of the generating functions for the partitions of integers and also asymptotic formulas for the generating functions in the interval $(0, 1)$. We also proved that some of the generating functions above are in the Hayman class. These are all the ingredients we need to find an asymptotic formula for each of these restricted partitions of integers, as $n \to \infty$, see Theorem 3.2.8.

#### A. OGF of the partitions of integers $P$: Theorem of Hardy-Ramanujan

We already know that the generating function of the partitions $P$ is in the Hayman class, see Theorem 4.3.19 and discussion in there, and therefore $P$ is strongly Gaussian.

Lemma 5.1.2 gives approximations $\tilde{m}_P$ and $\tilde{\sigma}_P$ of $m_P$ and $\sigma_P$ verifying the conditions of Baéz-Duarte asymptotic formula, see Theorem 3.2.8. Recall that

$$\tilde{m}_P(e^{-s}) = \frac{\pi^2}{6s^2} \quad \text{and} \quad \tilde{\sigma}_P(e^{-s}) = \frac{\pi}{\sqrt{3}s^{3/2}}$$

solving $m_P(e^{-s_n}) = n$, where $\tau_n = e^{-s_n}$, we find that $s_n = \pi/\sqrt{6n}$ and therefore $\tau_n = e^{-\pi/\sqrt{6n}}$. From the previous discussion we find that

$$\tilde{\sigma}_P(\tau_n) = \frac{(6n)^{3/4}}{\sqrt{3\pi}} \quad \text{and} \quad \tau_n^n = e^{-\pi\sqrt{n}/\sqrt{6}}.$$

Finally Lemma 5.1.6 gives an asymptotic formula for $P(\tau_n)$:

$$P(\tau_n) \sim \frac{1}{\sqrt{2\pi}} \sqrt{s_n} e^{\pi^2/(6s_n)} = \frac{1}{2^{3/4}(3n)^{1/4}} e^{\pi\sqrt{n}/\sqrt{6}}, \qquad \text{as } n \to \infty,$$

and therefore, applying Theorem 3.2.8, we obtain the following theorem.

**Theorem 5.1.14** (Hardy-Ramanujan)**.**

$$p(n) \sim \frac{P(\tau_n)}{\sqrt{2\pi}\tilde{\sigma}_P(\tau_n)\tau_n^n} \sim \frac{1}{4\sqrt{3}} \frac{1}{n} e^{\pi\sqrt{2/3}\sqrt{n}}, \qquad \textit{as } n \to \infty,$$

this is Hardy-Ramanujan partition Theorem, see [44].

### B. OGF of the partitions into distinct parts $Q$

Now we obtain, by means of Baéz-Duarte asymptotic formula, see Theorem 3.2.8, an asymptotic formula for the partitions of integers into distinct parts $q(n)$.

The generating function $Q$ is in the Hayman class, see Theorem 4.3.20, and this implies that $Q$ is strongly Gaussian. Lemma 5.1.3 gives approximations $\tilde{m}_Q$ and $\tilde{\sigma}_Q$ for the mean and variance $m_Q$ and $\sigma_Q$. This lemma gives that these approximations verify the hypothesis of Baéz-Duarte asymptotic formula, see Theorem 3.2.8.

Recall that $\tilde{m}_Q(e^{-s}) = \pi^2/(12s^2)$ and $\tilde{\sigma}_Q(e^{-s}) = \pi/(\sqrt{6}s^{3/2})$, then $s_n =$ and $\tau_n =$.

Using Lemma 5.1.7 we obtain that

$$Q(e^{-s}) \sim \frac{1}{\sqrt{2}} e^{\pi^2/(12s)}, \qquad \text{as } s \downarrow 0\,,$$

and therefore we have the following theorem

**Theorem 5.1.15.** *For the partitions of integers into distinct parts $q(n)$ we have*

$$q(n) \sim \frac{Q(\tau_n)}{\sqrt{2\pi}\tilde{\sigma}_Q(\tau_n)\tau_n^n} = \frac{1}{4 \cdot 3^{1/4}} \frac{1}{n^{3/4}} e^{\pi\sqrt{1/3}\sqrt{n}}, \qquad \textit{as } n \to \infty\,.$$

### C. OGF of the partitions in arithmetic progression $P_{a,b}$: Theorem of Ingham

We want to apply Baéz-Duarte asymptotic formula, see Theorem 3.2.8, to the ordinary generating function $P_{a,b}$, with $a, b \geq 1$ and $\gcd(a, b) = 1$.

We proved in Theorem 5.1.11 that $P_{a,b}$, with $a, b \geq 1$ and $\gcd(a, b) = 1$, is in the Hayman class, and therefore that this is a strongly Gaussian power series.

Lemma 5.1.4 gives approximations $\tilde{m}_{P_{a,b}}$ and $\tilde{\sigma}_{P_{a,b}}$ for the mean and variance $m_{P_{a,b}}$ and $\sigma_{P_{a,b}}$. In fact we have

$$\tilde{m}_{P_{a,b}}(e^{-s}) = \zeta(2)/(as^2) \quad \text{and} \quad \tilde{\sigma}^2_{P_{a,b}}(e^{-s}) = (2\zeta(2))/(as^3), \qquad \text{for any } s > 0.$$

These approximations verify the hypothesis of Baéz-Duarte asymptotic formula, see Theorem 3.2.8.

Solving $\tilde{m}_{P_{a,b}}(e^{-s_n}) = \tilde{m}_{P_{a,b}}(\tau_n) = n$ we find that $s_n = (\zeta(2)/(an))^{1/2}$ and therefore $\tau_n = e^{-s_n} = e^{-(\zeta(2)}/(an))^{1/2}$. Observe that

$$\tilde{\sigma}(\tau_n) = \tilde{\sigma}(e^{-s_n}) = \sqrt{\frac{2a}{\zeta(2)}} n^{3/4} = \frac{\sqrt{12a}}{\pi} n^{3/4}.$$

Combining the sequence $\tau_n$, the sequence $\tilde{\sigma}(\tau_n)$, the asymptotic formula given by Theorem 5.1.12 and Baéz-Duarte asymptotic formula (3.2.3), we obtain the following theorem.

**Theorem 5.1.16** (Ingham)**.** *For integers $a, b \geq 1$ with $\gcd(a,b) = 1$, the number of partitions $p_{a,b}(n)$ of integer $n \geq 1$ with parts drawn for the arithmetic sequence $\{aj + b : j \geq 0\}$ satisfies that*

$$p_{a,b}(n) \sim \left[\frac{1}{\sqrt{2}\,2\pi}\Gamma\Big(\frac{b}{a}\Big) a^{b/(2a)-1/2}\Big(\frac{\pi^2}{6}\Big)^{b/(2a)}\right] \frac{1}{n^{1/2+b/(2a)}} \exp\Big(\pi\sqrt{\frac{2}{3a}}\sqrt{n}\Big), \tag{5.1.7}$$

*as $n \to \infty$.*

## D. Plane and Colored Partitions: Wright's Theorem

Again, we apply Baéz-Duarte asymptotic formula to the generating function $W_1^b$ where $b \geq 0$, the OGF of the colored partitions with coloring sequence $(j^b)_{j\geq 1}$.

First we discuss the plane partitions, that is, the case $b = 1$.

**D.1. Plane partitions** The MacMahon function

$$W_1^1(z) = M(z) = \prod_{j=1}^{\infty} \frac{1}{(1-z^j)^j} = \sum_{n=0}^{\infty} M_n z^n, \quad \text{for any } |z| < 1,$$

codifies the plane partitions: $M_n$ is the number of plane partitions of the integer $n$, for $n \geq 1$, with $M_0 = 1$.

We already know that $M$ is in the Hayman class, see Theorem 4.3.21, and therefore that $M$ is strongly Gaussian. We also have approximations $\tilde{m}_M$ and $\tilde{\sigma}_M$ for the mean $m_M$ and the variance $\sigma_M$, see Lemma 5.1.5. These approximations verify the hypothesis of Baéz-Duarte Theorem 3.2.8, see Lemma 5.1.5.

In fact we have

$$\tilde{m}_M(e^{-s}) = \zeta(3)\,\Gamma(3)\,\frac{1}{s^3}, \quad \text{and} \quad \tilde{\sigma}_M^2(e^{-s}) = \zeta(3)\,\Gamma(4)\,\frac{1}{s^4}, \quad \text{for any } s > 0,$$

and therefore $s_n = (\zeta(3)\Gamma(3))^{1/3}/n^{1/3}$ and $\tau_n = e^{-s_n}$ are so that $\tilde{m}_M(\tau_n) = n$ and $\tilde{\sigma}_M(\tau_n) = (\zeta(3)\Gamma(3))^{-1/6} n^{2/3}$.

We also have at our disposal, see equation (5.1.5), a suitable asymptotic formula for $W_1^1$ on the interval $(0,1)$:

$$W_1^1(e^{-s}) \sim e^{\zeta'(-1)}\, s^{1/12}\, e^{\zeta(3)/s^2} = \frac{e^{1/12}}{GK_1}\, s^{1/12}\, e^{\zeta(3)/s^2}\,, \quad \text{as } s \downarrow 0\,. \tag{5.1.8}$$

Now we are in position to apply Theorem 3.2.8, this way we obtain the following theorem.

**Theorem 5.1.17** (Wright)**.** *For number of plane partitions of an integer $n$ we have the asymptotic formula*

$$M_n \sim \left[\frac{e^{\zeta'(-1)}}{\sqrt{12\pi}}\, \zeta(3)^{7/36} 2^{25/36}\right] \frac{1}{n^{25/36}}\, \exp\left(3(\zeta(3)/4)^{1/3} n^{2/3}\right), \quad \textit{as } n \to \infty\,, \tag{5.1.9}$$

which is Wright's asymptotic formula for plane partitions, [104].

**D.2. Colored partitions with coloring sequence** $(j^b)_{j\geq 1}$ The OGF $W_1^b$, with integer $b \geq 0$, is given by the infinite product

$$W_1^b(z) = \prod_{j=1}^{\infty} \frac{1}{(1-z^j)^{j^b}}, \qquad \text{for any } |z| < 1.$$

First we prove that the OGF $W_1^b$ is in the Hayman class. Let $g$ be the power series

$$g(z) = -\sum_{j=1}^{\infty} j^b \ln(1-z^j) = \sum_{n=1}^{\infty} \frac{\sigma_{b+1}(n)}{n} z^n, \qquad \text{for any } |z| < 1\,,$$

where $\sigma_{b+1}(n) = \sum_{jk=n} j^{b+1}$ denotes the sum of the $(b+1)$-th powers of the divisors of $n$. For this power series $g$ we have $W_1^b = e^g$. We have already encountered $\sigma_1(n)$ in Subsection 4.3.4-B.

We want to apply Theorem 4.3.17 to $g$. Fix $\varepsilon > 0$ such that $\varepsilon < (b+1)/2$, for instance $\varepsilon = 1/4$, then there exists a constant $D_\varepsilon > 0$ such that

$$n^b \leq \frac{\sigma_{b+1}(n)}{n} \leq D_\varepsilon n^{b+\varepsilon}, \qquad \text{for any } n \geq 1\,,$$

the inequality at the right-hand side follows since $\sigma_{b+1}(n) \leq n^b \sigma_1(n)$, finally applying Theorem 322 in [45] we find that for any $\varepsilon > 0$, there exists a constant $D_\varepsilon > 0$ such that $\sigma_1(n) \leq D_\varepsilon n^{1+\varepsilon}$, for any $n \geq 1$. See also 4.3.4-B. The inequality on the left holds because $n$ always divides $n$.

The previous coefficient inequality gives that $W_1^b$ is in the Hayman class, see Theorem 4.3.17, and therefore $W_1^b$ is strongly Gaussian.

Lemma 5.1.5, and the discussion in there, give approximations $\tilde{m}_{W_1^b}$ and $\tilde{\sigma}_{W_1^b}$ for the mean $m_{W_1^b}$ and the variance $\sigma_{W_1^b}$. These approximations verify the hypothesis of Baéz-Duarte Theorem. In fact we have

$$\tilde{m}_{W_1^b}(e^{-s}) = \zeta\,(b+2)\,\Gamma\,(b+2)\,\frac{1}{s^{b+2}}\,,$$

and

$$\tilde{\sigma}^2_{W_a^b}(e^{-s}) = \zeta\,(b+2)\,\Gamma\,(b+3)\,\frac{1}{s^{b+3}}.$$

We also have an asymptotic formula for $W_1^b(e^{-s})$, as $s \downarrow 0$, see Theorem 5.1.13. Combining all these pieces we find the following theorem.

**Theorem 5.1.18.** *For the colored partitions with colors in the sequence $(j^b)_{j\geq 1}$ and integer $b \geq 0$ we have*

$$\text{COEFF}_{[n]}(W_1^b(z)) \sim \alpha_b \frac{1}{n^{\beta_b}} e^{\gamma_b\, n^{\frac{b+1}{b+2}}}\ , \quad \textit{as } n \to \infty\ , \tag{5.1.10}$$

*where*

$$\alpha_b = \frac{1}{\sqrt{2\pi}} e^{\zeta'(-b)} \frac{1}{\sqrt{b+2}} \Big[\Gamma(b+2)\zeta(b+2)\Big]^{\frac{-2\zeta(-b)+1}{2(b+2)}}$$

$$\beta_b = \frac{-2\zeta(-b)+b+3}{2(b+2)}$$

$$\gamma_b = \frac{b+2}{b+1}\big(\Gamma(b+2)\zeta(b+2)\big)^{1/(b+2)}.$$

If in the previous theorem we fix $b=0$ (Hardy-Ramanujan Theorem), we obtain the asymptotic for the partitions of integers and if we fix $b=1$ we obtain the asymptotic formula for the plane partitions (Wright's Theorem).

For asymptotic formulae (and further asymptotic expansions) in the cases $W_2^0$ and $W_a^0$, of partitions into squares and $a$-th powers, obtained via the circle method we refer to Vaughan [97], Gafni [35] and the primordial Wright [106].

## 5.2 Partitions of sets

Recall that the EGF of the partitions of (non-empty) sets, which is the EGF of the Bell numbers $B_n$, is given by $B(z) = e^{e^z-1} = \sum_{n=0}^{\infty}(\mathcal{B}_n/n!)z^n$.

We already know that $B$ is in the Hayman class, see Theorem 4.3.18 and the discussion in there, in particular $B$ is a strongly Gaussian power series. The mean and variance functions of $B$ are given by

$$m_B(t) = te^t \quad \text{and} \quad \sigma_B(t) = (t+1)e^t.$$

Let $(t_n)$ be the sequence given by $m_B(t_n) = t_n e^{t_n} = n$, then $t_n = W(n)$, where $W$ is the Lambert function. Observe that $\sigma_B(t_n) = \sqrt{(W(n)+1)e^{W(n)}}$ and $B(t_n) = \exp(e^{W(n)}-1)$.

Applying Baéz-Duarte asymptotic formula, see Theorem 3.2.8, we find the following theorem.

**Theorem 5.2.1** (Moser-Wyman). *For the Bell numbers $B_n$ we have*

$$\frac{\mathcal{B}_n}{n!} \sim \frac{1}{\sqrt{2\pi}} \frac{e^{e^{W(n)}-1}}{\sqrt{W(n)(W(n)+1)e^{W(n)}}W(n)^n}, \quad \text{as } n \to \infty,$$

which is the asymptotic formula for Bell numbers of Moser-Wyman [68].

In contrast with the cases studied above, partitions of integers with different restrictions, this asymptotic formula follows in a very direct way. First we can solve the equation $m_f(t_n) = n$ in a very direct way by using the Lambert function. Then, the function $f$ itself, the mean and variance functions have simple expressions which give, in a very direct way, the previous asymptotic formula.

Chapter 6

# Large Powers asymptotics and Lagrangian distributions

## Contents

---

In this chapter we introduce a unified approach to deal with *the asymptotic behavior of the coefficients of large powers of power series with non-negative coefficients.*

The framework of this approach arises from the use of Local Central Limit Theorems for individual lattice random variables (Theorem 6.1.5) and also for continuous families of lattice random variables (Theorem 6.1.8) as a tool within the theory of Khinchin families of random variables associated with power series.

Large powers arise while considering sums of independent copies of a random variables: if $\psi(z)$ is the probability generating function of a random variable $Y$ taking values in $\{0, 1, \ldots\}$ and $Y_1, \ldots, Y_n$ are independent copies of $Y$ then

$$\mathbf{P}\Big(\sum_{j=1}^{n} Y_j = k\Big) = \text{COEFF}_{[k]}(\psi(z)^n)\,,$$

and if $(X_t)_{t\in(0,R)}$ is the Khinchin family associated with some power series $f$, and $X_t^{(1)}, \ldots, X_t^{(n)}$ are $n$ independent copies of $X_t$ then

$$\begin{aligned}\mathbf{P}\Big(\sum_{j=1}^{n} X_t^{(j)} = k\Big) &= \text{COEFF}_{[k]}\bigg(\Big(\frac{f(tz)}{f(t)}\Big)^n\bigg)\\ &= \text{COEFF}_{[k]}(f^n(z))\,\frac{t^k}{f^n(t)}\,, \quad \text{for any } k \geq 0 \text{ and } t \in (0, R)\,.\end{aligned}$$

Large powers also appear in Combinatorics. For instance, if $\psi(z)$ is the generating function of an unlabelled combinatorial class $\mathcal{C}$, then $\text{COEFF}_{[k]}(\psi(z))$ is the number of the objects in the class with weight $k$, while $\text{COEFF}_{[k]}(\psi(z)^n)$ counts the number of lists of length $n$ of objects of $\mathcal{C}$ with total weight $k$. For general background on (Analytic) Combinatorics we refer to the comprehensive treatise [32].

Along this chapter, we reserve $k$ to signify *index of the coefficient* and $n$ to denote *the power to which a power series $f$ is raised.* The power $n$ tends to infinity. Different results appear depending on how $k$ behaves with $n$, or actually, on how $k/n$ behaves as $n \to \infty$. As we shall see the Local Central Limit theorem for continuous families is used to handle the large exponent $n$ while the Khinchin families, particularly, of Gaussian (or strongly Gaussian) power series are used to handle the coefficient index $k$. This chapter is mainly based on the preprint:

- Maciá, V.J. et al. Large Powers asymptotics, Khinchin families and Lagrangian distributions, ArXiv:2201.11746, see [31].

## 6.1 Local Central Limit Theorem and continuous families

To obtain asymptotic results for the coefficients of large powers of power series with non-negative coefficients we will write down formulae (to be presented in Section 6.2.2) for the coefficients of the large powers in terms of the associated Khinchin families and Hayman's identities.

To handle the Hayman's identities of a Khinchin family we shall use the Local Central Limit Theorem for lattice variables in some integral form. This is Theorem 6.1.5, as applied to a single lattice variable, and Theorem 6.1.8, which is a version of the Local Central Limit Theorem for continuous families of lattice variables. We provide a detailed proof of the latter result which requires some preliminary material which we will discuss in Section 6.1.4.

### 6.1.1 Approximation and bound of characteristic functions

For a random variable $Y$ and integer $n \geq 1$, we have for $\mathbf{E}\big(e^{\imath\theta Y/\sqrt{n}}\big)^n$ the following well known approximation, Theorem 6.1.1, and bound, Theorem 6.1.2, involving the third moment of $Y$.

**Theorem 6.1.1.** *There exists a constant $C > 0$, such that for any random variable $Y$ verifying that $\mathbf{E}(Y) = 0$, $\mathbf{E}(Y^2) = 1$ and $\mathbf{E}(|Y|^3) < \infty$, it holds that*

$$\left|\mathbf{E}\big(e^{\imath\theta Y/\sqrt{n}}\big)^n - e^{-\theta^2/2}\right| \leq C\frac{\mathbf{E}(|Y|^3)}{\sqrt{n}}\,, \quad \textit{if } |\theta| \leq \frac{\sqrt{n}}{\mathbf{E}(|Y|^3)} \textit{ and } n \geq 1\,.$$

See, for instance, [42, Lemma 6.2, Chapter 7]. In particular,

$$\lim_{n\to\infty} \mathbf{E}\big(e^{\imath\theta Y/\sqrt{n}}\big)^n = e^{-\theta^2/2}\,, \quad \text{for each } \theta \in \mathbb{R}\,, \tag{6.1.1}$$

which is the Central Limit Theorem for sums of independent identically distributed copies of the variable $Y$ with the extra assumption of finite third absolute moment.

**Theorem 6.1.2.** *For any random variable $Y$ such that $\mathbf{E}(Y) = 0$, $\mathbf{E}(Y^2) = 1$ and $\mathbf{E}(|Y|^3) < \infty$, we have that*

$$\left|\mathbf{E}\big(e^{\imath\theta Y/\sqrt{n}}\big)^n\right| \leq e^{-\theta^2/3}, \quad \textit{for any } \theta \textit{ such that } |\theta| \leq \frac{\sqrt{n}}{4\mathbf{E}(|Y|^3)} \textit{ and } n \geq 1\,.$$

See, for instance, step 1 in the proof of [42, Lemma 6.2, Chapter 7].

### 6.1.2 Lattice distributed variables

A random variable $Z$ is said to be *lattice distributed* or a *lattice random variable* if $Z$ is nonconstant and for a certain $a \in \mathbb{R}$ and $h > 0$ one has that $\mathbf{P}(Z \in a + h\mathbb{Z}) = 1$. The largest $h$ such that $\mathbf{P}\big(Z \in a+h\mathbb{Z}\big) = 1$ is called the *gauge* of the lattice variable $Z$. For gauge $h$, the variable $(Z-a)/h$ takes integer values and has gauge 1.

Let $f(z) = \sum_{n=0}^\infty a_n z^n$ be a power series in $\mathcal{K}$. Let $\mathcal{I}_f = \{n \geq 1 : a_n \neq 0\}$. Recall the notation $Q_f = \gcd\{\mathcal{I}_f\}$ from Section 1.3.6.

If $Q_f = 1$, then each $X_t$ of the Khinchin family of $f$ is a lattice variable of gauge 1 and each $\breve{X}_t$ is a lattice variable of gauge $1/\sigma_f(t)$. The condition $Q_f = 1$ is equivalent to the requirement that for each $d \geq 1$ there is $n \in \mathcal{I}_f$ such that $d \nmid n$.

In general, if $Q_f \geq 1$, then each $X_t$ of the Khinchin family of $f$ is a lattice variable of gauge $Q_f$ and each $\breve{X}_t$ is a lattice variable of gauge $Q_f/\sigma_f(t)$.

**Lemma 6.1.3.** *If $Z$ is a lattice random variable with gauge $h$, then*

$$|\mathbf{E}(e^{\imath\theta Z})| < 1\,, \quad \textit{for } 0 < \theta < 2\pi/h\,.$$

*In particular, for each $0 < \delta < \pi/h$, there exists $\omega = \omega_\delta < 1$ such that*

$$|\mathbf{E}(e^{\imath\theta Z})| \leq \omega\,, \quad \textit{for } |\theta - \pi/h| \leq \delta\,.$$

See [30, Section 3.5]. The second part of the statement follows simply from continuity of the characteristic function.

The next result is a standard inversion lemma whereas the mass function of a variable $U$ is expressed in terms of its characteristic function.

**Lemma 6.1.4.** *Let $U$ be a lattice random variable with gauge $h$, so that for certain $a \in \mathbb{R}$ we have that $\mathbf{P}(U \in a + h\mathbb{Z}) = 1$. Then*

$$\mathbf{P}(U = a + hk) = \frac{h}{2\pi}\int_{-\pi/h}^{\pi/h} \mathbf{E}(e^{\imath\theta U})e^{-\imath\theta(a+kh)}\,d\theta\,, \quad \textit{for each } k \in \mathbb{Z}\,.$$

In other terms, for each attainable value $x$ of the variable $U$ we have that

$$\mathbf{P}(U = x) = \frac{h}{2\pi}\int_{-\pi/h}^{\pi/h} \mathbf{E}(e^{\imath\theta U})e^{-\imath\theta x} d\theta\,.$$

*Proof.* We may assume that $a = 0$ and $h = 1$. For each $k \in \mathbb{Z}$, let $p_k = \mathbf{P}(U = k)$. Then, we have that $\sum_{k\in\mathbb{Z}} p_k = 1$ and

$$\mathbf{E}(e^{\imath\theta U}) = \sum_{k\in\mathbb{Z}} p_k e^{\imath k\theta}\,, \quad \text{for each } \theta \in \mathbb{R}\,,$$

and therefore, that $p_k = \dfrac{1}{2\pi}\displaystyle\int_{-\pi}^{\pi} \mathbf{E}(e^{\imath\theta U})e^{-\imath k\theta} d\theta\,, \quad \text{for each } k \in \mathbb{Z}\,.$ □

### 6.1.3 Local Central Limit Theorem for lattice variables: integral form

**Theorem 6.1.5.** *Let $Z$ be a lattice random variable with gauge $h$ and such that $\mathbf{E}(Z) = 0$ and $\mathbf{E}(Z^2) = 1$. Then*

$$\lim_{n\to\infty} \int\limits_{|\theta|\leq\pi\sqrt{n}/h} \left|\mathbf{E}(e^{\imath\theta Z/\sqrt{n}})^n - e^{-\theta^2/2}\right| d\theta = 0\,.$$

From the convergence of this integral, the usual statement of the Local Central Limit Theorem follows, see Section 6.1.3. In fact, Theorem 6.1.5 is what is actually shown (as an intermediate step) in the usual proofs of the Local Central Limit Theorem, see, for instance, [38, Section 43] and the proof of [30, Theorem 3.5.3].

We give a proof below with the extra assumption that $\mathbf{E}(|Y|^3) < \infty$ anticipating and as preparation for the corresponding result, Theorem 6.1.8, for continuous families of lattice variables.

*Proof.* The Central Limit Theorem, as in equation (6.1.1), gives that, as $n \to \infty$, the integrand converges towards 0 for each $\theta \in \mathbb{R}$.

Denote $\tau \triangleq 1/(2\mathbf{E}(|Z|^3)$. Theorem 6.1.2 provides us with the bound

$$|\mathbf{E}(e^{\imath\theta Z/\sqrt{n}})^n| \le e^{-\theta^2/6}\,, \quad \text{si } |\theta| \le \tau\sqrt{n}\,.$$

Thus, dominated convergence implies that

$$(\flat) \quad \lim_{n\to\infty} \int_{|\theta|\le\tau\sqrt{n}} \left|\mathbf{E}(e^{\imath\theta Z/\sqrt{n}})^n - e^{-\theta^2/2}\right| d\theta = 0\,.$$

Lemma 6.1.3 gives, in particular, a bound $\omega \in (0,1)$ such that

$$|\mathbf{E}(e^{\imath\theta Z})| \le \omega\,, \quad \text{si } \tau \le |\theta| \le \pi/h\,,$$

and so that

$$|\mathbf{E}(e^{\imath\theta Z/\sqrt{n}})^n| \le \omega^n\,, \quad \text{si } \tau\sqrt{n} \le |\theta| \le \pi\sqrt{n}/h\,.$$

The bound

$$\int_{\tau\sqrt{n}\le|\theta|\le\pi\sqrt{n}/h} \left|\mathbf{E}(e^{\imath\theta Z/\sqrt{n}})^n - e^{-\theta^2/2}\right| d\theta \le \int_{\tau\sqrt{n}\le|\theta|\le\pi\sqrt{n}/h} \left(\omega^n + e^{-\theta^2/2}\right) d\theta\,,$$

tends to 0 as $n \to \infty$. This together with $(\flat)$ gives the result. □

### Local Central Limit Theorem

Let $Y$ be a random variable such that $\mathbf{E}(Y) = 0$ and $\mathbf{E}(Y^2) = 1$. Let $Y_1, Y_2, \ldots$ be a sequence of independent random variables each one of them distributed like $Y$.

Assume that $Y$ is a lattice variable with gauge $h$ and that $\mathbf{P}(Y \in a + h\mathbb{Z}) = 1$, for some fixed $a \in \mathbb{R}$.

For each $n \ge 1$, we denote $S_n = (Y_1 + \cdots + Y_n)/\sqrt{n}$. The variable $S_n$ is likewise a lattice variable, with gauge $h/\sqrt{n}$.

Let $\mathcal{V}_n$ denote the set of values that the variable $S_n$ can take, i.e.

$$\mathcal{V}_n = \{a\sqrt{n} + kh/\sqrt{n} : k \in \mathbb{Z}\}.$$

Applying Lemma 6.1.4 to the variable $S_n$ and using that $\mathbf{E}(e^{\imath\theta S_n}) = \mathbf{E}(e^{\imath\theta Y/\sqrt{n}})^n$, we have that

$$\mathbf{P}(S_n = x) = \frac{h}{2\pi\sqrt{n}} \int_{-\pi\sqrt{n}/h}^{\pi\sqrt{n}/h} \mathbf{E}(e^{\imath\theta Y/\sqrt{n}})^n e^{-\imath\theta x}\, d\theta\,, \quad \text{for every } x \in \mathcal{V}_n\,.$$

Since

$$\frac{1}{\sqrt{2\pi}} e^{-x^2/2} = \frac{1}{2\pi} \int_{\mathbb{R}} e^{-\theta^2/2} e^{-\imath\theta x}\, d\theta\,, \quad \text{for every } x \in \mathbb{R}\,,$$

we deduce that

$$\sup_{x\in\mathcal{V}_n} \left| \frac{\sqrt{n}}{h}\mathbf{P}(S_n = x) - \frac{1}{\sqrt{2\pi}} e^{-x^2/2} \right| \le \frac{1}{2\pi} \int\limits_{|\theta|\le\pi\sqrt{n}/h} \left| \mathbf{E}(e^{\imath\theta Y/\sqrt{n}})^n - e^{-\theta^2/2} \right| d\theta$$
$$+ \frac{1}{2\pi} \int\limits_{|\theta|>\pi\sqrt{n}/h} e^{-\theta^2/2}\, d\theta\,.$$

Theorem 6.1.5 allows us to conclude from the bound above that

**Corollary 6.1.6** (Local Central Limit Theorem)**.** *With the notations above,*

$$\lim_{n\to\infty} \sup_{x\in\mathcal{V}_n} \left| \frac{\sqrt{n}}{h}\mathbf{P}(S_n = x) - \frac{1}{\sqrt{2\pi}} e^{-x^2/2} \right| = 0\,.$$

### 6.1.4 Continuous families of lattice variables

In this subsection we introduce some results for continuous families of random variables that will be used in the context of large powers asymptotics, see Chapter 1, in particular Subsection 6.1.4, for the definition, and properties, of continuous family.

Recall that the (un-normalized) Khinchin family $(X_t)$ of any $f$ in $\mathcal{K}$ with radius of convergence $R$ is continuous and bounded in any interval $[0, b]$ with $0 < b < R$. If $M_f < +\infty$ (and $R < +\infty$) the family $(X_t)_{t\in[0,R]}$ extended to the closed interval $[0, R]$ is continuous and bounded; observe that $\mathbf{E}(|X_t|) = \mathbf{E}(X_t) \le M_f$, for $t \in [0, R]$. See Chapter 1, in particular Subsection 6.1.4.

Let $(Z_s)_{s\in[a,b]}$ be a continuous family of lattice variables. For each $s \in [a, b]$, let $h(s)$ denote the gauge of $Z_s$. The gauge function $h$ is upper semicontinuous:

$$\limsup_{u\to s} h(u) \le h(s)\,, \quad \text{for each } s \in [a, b]\,.$$

In general, the gauge of a continuous family is not a continuous function. For instance, for $s \in [0, 1/2]$, let the variable $Z_s$ take the values $0, 1/2$ and $1$ with respective probabilities $1/2 - s, 2s, 1/2 - s$. The family $(Z_s)$ is continuous, and $h(s) = 1/2$, for $s \neq 0$, but $h(0) = 1$.

The next lemma is a counterpart of Lemma 6.1.3 for continuous and bounded families; it follows from Lemma 6.1.3 by continuity and compactness.

**Lemma 6.1.7.** *Let $(Z_s)_{s\in[a,b]}$ be a continuous and bounded family of lattice random variables, with continuous gauge function $h(s)$.*

*Let $H > 0$, be such that $h(s) \le H$, for each $s \in [a, b]$.*

*Then,*

$$|\mathbf{E}(e^{\imath\theta Z_s})| < 1, \quad \textit{if } s \in [a,b] \textit{ and } 0 < |\theta| < 2\pi/h(s)\,.$$

*In particular, for each $0 < \tau < \pi/H$, there exists $\omega_\tau < 1$ such that*

$$|\mathbf{E}(e^{\imath\theta Z_s})| \le \omega_\tau, \quad \textit{if } s \in [a,b] \textit{ and } \tau \le |\theta| \le \pi/h(s)\,.$$

The next theorem is a uniform analogue of Theorem 6.1.5 for continuous families of lattice variables:

**Theorem 6.1.8.** *Let $(Z_s)_{s\in[a,b]}$ be a continuous family of lattice variables with gauge function $h(s)$ continuous in $[a,b]$ and such that for each $s \in [a,b]$ we have that $\mathbf{E}(Z_s) = 0$, $\mathbf{E}(Z_s^2) = 1$ and $\mathbf{E}(|Z_s|^3) \le \Gamma$, for a certain constant $\Gamma > 0$.*

*Then,*

$$\lim_{n\to\infty} \sup_{s\in[a,b]} \int_{-\pi\sqrt{n}/h(s)}^{\pi\sqrt{n}/h(s)} \left|\mathbf{E}(e^{\imath\theta Z_s/\sqrt{n}})^n - e^{-\theta^2/2}\right| d\theta = 0\,.$$

**Remark 6.1.9.** The hypothesis $\sup_{s\in[a,b]} \mathbf{E}(|Z_s|^3) < \infty$ may be replaced by the assumption that $\sup_{s\in[a,b]} \mathbf{E}(|Z_s|^{2+\delta}) < \infty$, for some $\delta > 0$, or even by assuming that the family $(|Z_s|^2)_{s\in[a,b]}$ is uniformly integrable. ⊠

*Proof.* First, we pay attention to the range of integration $|\theta| \le \sqrt{n}/(4\Gamma)$.

Theorem 6.1.2 gives us that

$$\left|\mathbf{E}\big(e^{\imath\theta Z_s/\sqrt{n}}\big)^n\right| \le e^{-\theta^2/3}, \quad \text{if } s \in [a,b] \text{ and } |\theta| \le \frac{\sqrt{n}}{4\Gamma}\,,$$

while Theorem 6.1.1 gives us that

$$\left|\mathbf{E}\big(e^{\imath\theta Z_s/\sqrt{n}}\big)^n - e^{-\theta^2/2}\right| \le C\frac{\Gamma}{\sqrt{n}}\,, \quad \text{if } s \in [a,b] \text{ and } |\theta| \le \frac{\sqrt{n}}{\Gamma}\,.$$

Therefore, for any sequence $(s_n)_{n\ge1}$ extracted from $[a,b]$ we deduce, from dominated convergence, that

$$\lim_{n\to\infty} \int_{-\sqrt{n}/(4\Gamma)}^{\sqrt{n}/(4\Gamma)} \left|\mathbf{E}(e^{\imath\theta Z_{s_n}/\sqrt{n}})^n - e^{-\theta^2/2}\right| d\theta = 0\,,$$

and, therefore, that

$$(\star) \qquad \lim_{n\to\infty} \sup_{s\in[a,b]} \int_{-\sqrt{n}/(4\Gamma)}^{\sqrt{n}/(4\Gamma)} \left|\mathbf{E}(e^{\imath\theta Z_s/\sqrt{n}})^n - e^{-\theta^2/2}\right| d\theta = 0\,.$$

Fix $0 < J < H$ such that $J \le h(s) \le H$, for each $s \in [a,b]$. Take $\tau > 0$ such that $\tau < 1/(4\Gamma)$ and $\tau < \pi/H$.

Since, for every $s \in [a,b]$, we have that $\mathbf{E}(|Z_s|) \leq \mathbf{E}(|Z_s|^2)^{1/2} = 1$, the family $(Z_s)_{s\in[a,b]}$ is bounded. Thus, Lemma 6.1.7 shows that there exists $\omega_\tau < 1$ such that

$$|\mathbf{E}(e^{\imath\theta Z_s})| \leq \omega_\tau, \quad \text{if } s \in [a,b] \text{ and } \tau \leq |\theta| \leq \pi/h(s)\,,$$

and so such that

$$|\mathbf{E}(e^{\imath\theta Z_s/\sqrt{n}})^n| \leq \omega_\tau^n, \quad \text{if } s \in [a,b] \text{ and } \tau\sqrt{n} \leq |\theta| \leq \pi\sqrt{n}/h(s)\,.$$

Therefore, for each $s \in [a,b]$ we have that

$$\int_{\tau\sqrt{n}\leq|\theta|\leq\pi\sqrt{n}/h(s)} \left|\mathbf{E}(e^{\imath\theta Z_s/\sqrt{n}})^n - e^{-\theta^2/2}\right| d\theta \leq \omega_\tau^n \pi\sqrt{n}/J + \int_{|\theta|\geq\tau\sqrt{n}} e^{-\theta^2/2} d\theta\,.$$

Consequently,

$$(\star\star) \quad \lim_{n\to\infty} \sup_{s\in[a,b]} \int_{\tau\sqrt{n}\leq|\theta|\leq\pi\sqrt{n}/h(s)} \left|\mathbf{E}(e^{\imath\theta Z_s/\sqrt{n}})^n - e^{-\theta^2/2}\right| d\theta = 0\,.$$

Combining $(\star)$ and $(\star\star)$ we obtain the result since $\tau < 1/(4\Gamma)$. □

## 6.2 Coefficients of large powers

We now turn our attention to the study of the asymptotic behavior of coefficients of large powers of functions $\psi$ in the class $\mathcal{K}$: the behavior of the coefficient $\text{COEFF}_{[k]}(\psi^n(z))$, as $n \to \infty$.

The objective of this section is to establish appropriate expressions for the Hayman's identities for the coefficients $\text{COEFF}_{[k]}(\psi^n(z))$ of (large) powers in terms of the Khinchin family of $\psi$: Lemma 6.2.2 and Lemma 6.2.3.

In forthcoming sections we will apply these lemmas to study the behavior of the coefficient $\text{COEFF}_{[k]}(\psi^n(z))$ under a variety of assumptions upon the joint behavior of the index $k$ and of the power $n$. At a later stage, in Section 6.6, we will consider the asymptotic behavior of the coefficients of $h(z)\psi^n(z)$ where $h(z) \in \mathcal{K}$.

**Remark 6.2.1.** The results which follow about large powers actually cover the general situation of $\psi$ with non-negative coefficients, and not just $\psi \in \mathcal{K}$. To see this we may reason as follows.

If $\psi$ has non-negative coefficients and it is not in $\mathcal{K}$, then one of the two following possibilities occurs: 1) $\psi$ is a constant or a monomial like $\psi(z) = bz^m$, for some $b \neq 0$ and integer $m \geq 1$, which are trivial situations as coefficients and large powers are concerned, 2) $\psi(0) = 0$ and $\psi$ has at least two non-negative coefficients, in which case, $\psi \in \mathcal{K}_s$ and for some integer $l \geq 1$ we have that $\varphi(z) \triangleq \psi(z)/z^l \in \mathcal{K}$, and thus $\text{COEFF}_{[k]}(\psi^n) = \text{COEFF}_{[k-nl]}(\varphi^n)$, for $k \geq nl$. ⊠

Henceforth, $\psi(z) = \sum_{j=0}^\infty b_j z^j$ will denote a power series in $\mathcal{K}$ with radius of convergence $R > 0$, and from now on we let $(Y_t)_{t\in[0,R)}$ denote its Khinchin family.

We reserve $k$ to denote index of coefficient and $n$ to denote power of $\psi$. In all asymptotic discussions below the power $n$ tends to $\infty$ (large powers), while the index $k$ could tend to $\infty$ in such a way that either $k \asymp n$ (including the case when $k/n$ tends to a finite nonzero limit) or $k = o(n)$ (including the possibility that $k$ could remain fixed) or $n = o(k)$ (which would require, as we shall see, that the power series $\psi$ is uniformly strongly Gaussian). These three cases are discussed in Sections 6.3, 6.4 and 6.5, respectively.

### 6.2.1 Auxiliary function $\phi$

Recall the notation of Section 1.3.6:

$$Q_f = \gcd\{n \geq 1 : a_n \neq 0\} = \lim_{N\to\infty} \gcd\{1 \leq n \leq N : a_n \neq 0\}$$

for any power series $f(z) = \sum_{n=0}^{\infty} a_n z^n$ in $\mathcal{K}$.

If $Q_\psi > 1$, then $\psi(z) = \phi(z^{Q_\psi})$ for a certain auxiliary power series $\phi \in \mathcal{K}$ which has radius of convergence $R^{Q_\psi}$. Observe that $Q_\phi = 1$ and that

$$\textsc{coeff}_{[kQ_\psi]}(\psi^n(z)) = \textsc{coeff}_{[k]}(\phi^n(z^{Q_\psi})),$$

while for $q$ not a multiple of $Q_\psi$, we have that

$$\textsc{coeff}_{[q]}(\psi^n(z)) = 0\,.$$

Denote by $(Z_t)_{t\in(0,R^Q)}$ the Khinchin family associated with this auxiliary function $\phi$, then, see Section 1.3.6,

$$Y_t \stackrel{d}{=} Q_\psi \cdot Z_{t^{Q_\psi}}\,, \qquad \text{for any } t \in (0, R).$$

The mean and variance functions of $\psi$ and $\phi$ are related by

$$m_\psi(t) = Q_\psi \cdot m_\phi(t^{Q_\psi}) \quad \text{and} \quad \sigma_\psi^2(t) = Q_\psi^2 \cdot \sigma_\phi^2(t^{Q_\psi})\,.$$

For each $t \in (0, R)$, the variable $Y_t$ is a lattice random variable with gauge $Q_\psi$. Likewise, the normalized variable $\breve{Y}_t$ is a lattice random variable with gauge $h(t) = Q_\psi/\sigma_\psi(t)$.

In subsequent analysis, we shall obtain asymptotic formulae for $\textsc{coeff}_{[k]}(\psi(z)^n)$ first under the assumption that $Q_\psi = 1$, and then in the general case when $Q_\psi \geq 1$, by considering the auxiliary power series $\phi$, with $Q_\phi = 1$, and translating the results back to $\psi$.

### 6.2.2 Hayman's identities and large powers

We may express $\textsc{coeff}_{[k]}(\psi(z)^n)$ in terms of the characteristic function of the normalized variables $(\breve{Y}_t)_{t\in[0,R)}$ by Hayman's identity as follows.

**Lemma 6.2.2.** *With the notations above,*

$$\textsc{coeff}_{[k]}(\psi(z)^n) = \frac{1}{2\pi}\frac{\psi^n(t)}{t^k}\frac{1}{\sqrt{n}}\frac{1}{\sigma_\psi(t)} \int\limits_{|\theta|\leq \pi\sigma_\psi(t)\sqrt{n}} \mathbf{E}\big(e^{\imath\theta \breve{Y}_t/\sqrt{n}}\big)^n e^{\imath(nm_\psi(t)-k)\theta/(\sigma_\psi(t)\sqrt{n})} d\theta\,, \tag{6.2.1}$$

*for any index $k \geq 1$, power $n \geq 1$ and for all $t \in (0, R)$.*

*Proof.* For $t \in (0, R)$, Cauchy's integral formula gives that

$$\begin{aligned}\textsc{coeff}_{[k]}(\psi(z)^n) &= \frac{1}{2\pi\imath}\int_{|z|=t}\frac{\psi(z)^n}{z^{k+1}}dz = \frac{1}{2\pi}\frac{\psi(t)^n}{t^k}\int_{|\theta|<\pi}\frac{\psi(te^{\imath\theta})^n}{\psi(t)^n}e^{-\imath\theta k}d\theta\\ &= \frac{1}{2\pi}\frac{\psi(t)^n}{t^k}\int_{|\theta|<\pi}\mathbf{E}(e^{\imath\theta Y_t})^n e^{-\imath\theta k}d\theta\\ &= \frac{1}{2\pi}\frac{\psi^n(t)}{t^k}\frac{1}{\sqrt{n}}\frac{1}{\sigma_\psi(t)}\int_{|\theta|<\pi\sigma_\psi(t)\sqrt{n}}\mathbf{E}(e^{\imath\theta\breve{Y}_t/\sqrt{n}})^n e^{\imath(nm_\psi(t)-k)\theta/(\sigma_\psi(t)\sqrt{n})}d\theta.\end{aligned}$$

□

Observe that the integral expression of (6.2.1) is greatly simplified if the radius $t$ is such that ($\star$) $n\, m_\psi(t) = k$. Whenever possible we will select and use that value of $t$ such that ($\star$) holds exactly or at least approximately.

In the next sections, we shall use the formula (6.2.1) to study the asymptotic behavior of $\textsc{coeff}_{[k]}(\psi^n(z))$, as $n \to \infty$, while $k \asymp n$, $k/n \to 0$ or $k/n \to \infty$.

A minor variation of the proof of Lemma 6.2.2 above gives that if $H(z)$ is a function holomorphic in $\mathbb{D}(0, R)$, then, with the notations above, we have that

$$\begin{aligned}&\textsc{coeff}_{[k]}(H(z)\psi(z)^n)\\ (6.2.2)\qquad &= \frac{1}{2\pi}\frac{\psi^n(t)}{t^k}\frac{1}{\sqrt{n}}\frac{1}{\sigma_\psi(t)}\int\limits_{|\theta|\le\pi\sigma_\psi(t)\sqrt{n}} H\big(te^{\imath\theta/(\sigma_\psi(t)\sqrt{n})}\big)\mathbf{E}\big(e^{\imath\theta\breve{Y}_t/\sqrt{n}}\big)^n e^{\imath(nm_\psi(t)-k)\theta/(\sigma_\psi(t)\sqrt{n})}d\theta\,,\end{aligned}$$

for any index $k \ge 1$, power $n \ge 1$ and $t \in (0, R)$.

If $M_\psi < \infty$ (and $R < \infty$), then $\psi$ extends to be continuous in $\mathrm{cl}(\mathbb{D}(0, R))$ and the Khinchin family $(Y_t)_{t\in[0,R)}$ extends to the closed interval by adding $Y_R$ given by

$$\mathbf{P}(Y_R = n) = b_n R^n/\psi(R)\,, \quad \text{for } n \ge 0\,.$$

See Section 1.2.4 and Section 1.3.7 and the notations therein. In this case, we can take $t = R$ to obtain the following formula for the coefficients of $\psi^n$.

**Lemma 6.2.3.** *With the notations above, if $R < \infty$ and $M_\psi = m_\psi(R) < \infty$ and $\sigma^2_\psi(R) < \infty$, then*

$$\begin{aligned}&\textsc{coeff}_{[k]}(\psi(z)^n)\\ (6.2.3)\qquad &= \frac{1}{2\pi}\frac{\psi^n(R)}{R^k}\frac{1}{\sqrt{n}}\frac{1}{\sigma_\psi(R)}\int\limits_{|\theta|\le\pi\sigma_\psi(R)\sqrt{n}}\mathbf{E}\big(e^{\imath\theta\breve{Y}_R/\sqrt{n}}\big)^n e^{\imath(nm_\psi(R)-k)\theta/(\sigma_\psi(R)\sqrt{n})}d\theta\,,\end{aligned}$$

*for any index $k \ge 1$ and power $n \ge 1$.*

## 6.3 Index $k$ and power $n$ are comparable

In this section we discuss the asymptotic behavior of $\text{COEFF}_{[k]}(\psi^n(z))$ when the index $k$ and the power $n$ are such that $k \asymp n$, as $n \to \infty$.

The plan is to analyze the integral term in the Hayman's identities of Section 6.2.2 when $k \asymp n$, as $n \to \infty$, by invoking the Local Central Limits Theorems of Section 6.1.

Along this section we assume that, for certain $A, B$, fixed, $0 < A < B$, the index $k$ and the power $n$ are such that $A \le k/n \le B$.

To deal with this case, we assume additionally that $M_\psi > B$. This being the case, we denote $a = m_\psi^{-1}(A)$ and $b = m_\psi^{-1}(B)$.

Assume first that $Q_\psi = 1$; a restriction to be lifted shortly in Section 6.3.1. The family of random variables $\breve{Y}_t$, where $t \in [a, b]$, is a continuous family of lattice random variables with gauge function $h(t) = 1/\sigma_\psi(t)$; see Section 6.1.4. In particular, $\max_{t\in[a,b]} \mathbf{E}(|\breve{Y}_t|^3) < +\infty$.

For each $n \ge 1$, define $\tau_n$ by $m_\psi(\tau_n) = k/n$. This choice is possible because $M_\psi > B$. For each $n \ge 1$, we have $\tau_n \in [a, b]$.

Taking $t = \tau_n$ in Hayman's identity (6.2.1), the integral term in there simplifies to

$$I_n \triangleq \int_{|\theta|\le \pi\sigma_\psi(\tau_n)\sqrt{n}} \mathbf{E}\big(e^{\imath\theta \breve{Y}_{\tau_n}/\sqrt{n}}\big)^n d\theta\,.$$

Since $\sigma_\psi(\tau_n) = 1/h(\tau_n)$, Theorem 6.1.8 gives that

$$\lim_{n\to\infty} \int_{|\theta|\le \pi\sigma_\psi(\tau_n)\sqrt{n}} \Big|\mathbf{E}\big(e^{\imath\theta \breve{Y}_{\tau_n}/\sqrt{n}}\big)^n - e^{-\theta^2/2}\Big| d\theta = 0\,. \tag{6.3.1}$$

In particular, since $\min_{t\in[a,b]} \sigma_\psi(t) > 0$, we conclude that $\lim_{n\to\infty} I_n = \sqrt{2\pi}$, and, thus, that

$$\text{COEFF}_{[k]}(\psi(z)^n) \sim \frac{1}{\sqrt{2\pi}} \frac{\psi^n(\tau_n)}{\tau_n^k} \frac{1}{\sqrt{n}} \frac{1}{\sigma_\psi(\tau_n)}\,, \quad \text{as } n \to \infty.$$

### 6.3.1 General $Q_\psi$

We now lift the assumption that $Q_\psi = 1$. To simplify notation, write $Q = Q_\psi \ge 1$. Consider the auxiliary function $\phi$, so that $\psi(z) = \phi(z^Q)$.

Recall that $m_\psi(t) = Qm_\phi(t^Q)$, let $A' = A/Q$ and $B' = B/Q$ and observe that $M_\phi = M_\psi/Q > B/Q = B'$.

Define $\tau_n'$ by $m_\phi(\tau_n') = k'/n$, for $k'$ such that $A' < k'/n < B'$.

Since $Q_\phi = 1$, we have that

$$(\natural) \quad \text{COEFF}_{[k']}(\phi^n(z)) \sim \frac{1}{\sqrt{2\pi}} \frac{\phi(\tau_n')^n}{(\tau_n')^{k'}} \frac{1}{\sqrt{n}} \frac{1}{\sigma_\phi(\tau_n')}\,, \quad \text{as } n \to \infty\,.$$

Let $\tau_n = (\tau_n')^{1/Q}$. Take $k = k'Q$, observe that $A < \dfrac{k}{n} < B$ and that

$$m_\psi(\tau_n) = Q\, m_\phi(\tau_n') = \frac{k'Q}{n} = \frac{k}{n}\,,$$

and

$$\sigma_\psi(\tau_n) = Q\, \sigma_\phi(\tau_n')\,.$$

Since $\text{COEFF}_{[k]}(\psi^n) = \text{COEFF}_{[k']}(\phi^n)$, translating ($\natural$) into terms of $\psi$, we get that

$$\text{COEFF}_k(\psi^n(z)) \sim \frac{Q_\psi}{\sqrt{2\pi}} \frac{\psi(\tau_n)^n}{\tau_n^k} \frac{1}{\sqrt{n}} \frac{1}{\sigma_\psi(\tau_n)}\,,$$
$$\text{as } n \to \infty \text{ while } A < k/n < B \text{ and } k \text{ is a multiple of } Q\,.$$

**Theorem 6.3.1.** *For $0 < A < B < M_\psi$, we have that*

$$\text{COEFF}_k(\psi^n(z)) \sim \frac{Q_\psi}{\sqrt{2\pi}} \frac{\psi(\tau_n)^n}{\tau_n^k} \frac{1}{\sqrt{n}} \frac{1}{\sigma_\psi(\tau_n)},$$

*as $n \to \infty$, while $A \le k/n \le B$ and $k$ is a multiple of $Q_\psi$, and where $\tau_n$ is given uniquely by $m_\psi(\tau_n) = k/n$.*

### 6.3.2 With further information on $k/n$

We assume now that for some $L$ such that $0 < L < +\infty$ and $L \le M_\psi$, we have that $k/n \to L$, as $n \to \infty$. With that information, we may sharpen the asymptotic formula for $\text{COEFF}_{[k]}(\psi(z)^n)$ of Theorem 6.3.1.

Assume that $k/n \to L$, where $0 < L < +\infty$ and $L \le M_\psi$, and, in fact, in such a way that for some $\omega \in \mathbb{R}$

$$\frac{k}{n} = L + \omega \frac{1}{\sqrt{n}} + o\Big(\frac{1}{\sqrt{n}}\Big)\,, \quad \text{as } n \to \infty\,.$$

or, equivalently, that

$$\lim_{n\to\infty} \frac{nL - k}{\sqrt{n}} = -\omega\,. \tag{6.3.2}$$

• *Suppose that $L < M_\psi$.* Choose $\tau \in (0, R)$ such that $m_\psi(\tau) = L$. Observe that now we have a fixed value of $\tau$, and no $\tau_n$ varying with $n$.

Assume first that $Q_\psi = 1$. If we choose $t = \tau$ in formula (6.2.1) and invoke the Local Central Limit Theorem 6.1.5 we readily find that

$$\text{COEFF}_{[k]}(\psi(z)^n) \sim \frac{1}{2\pi}\frac{\psi^n(\tau)}{\tau^k}\frac{1}{\sqrt{n}}\frac{1}{\sigma_\psi(\tau)}\int_{\mathbb{R}} e^{\imath\omega\theta/\sigma_\psi(\tau)}e^{-\theta^2/2}\,d\theta$$
$$= \frac{e^{-\omega^2/(2\sigma_f^2(\tau))}}{\sqrt{2\pi}\,\sigma_f(\tau)}\frac{1}{\sqrt{n}}\frac{\psi^n(\tau)}{\tau^k}\,, \text{ as } n \to \infty.$$

An argument much like the one in Section 6.3.1 allows us to deduce the general case $Q_\psi \geq 1$ from the case $Q_\psi = 1$.

**Theorem 6.3.2.** *If $L < M_\psi$ and*

$$\lim_{n\to\infty}\frac{nL-k}{\sqrt{n}} = -\omega\,,$$

*then*

$$\text{COEFF}_{[k]}(\psi(z)^n) \sim Q_\psi\frac{e^{-\omega^2/(2\sigma_\psi^2(\tau))}}{\sqrt{2\pi}\,\sigma_\psi(\tau)}\frac{1}{\sqrt{n}}\frac{\psi^n(\tau)}{\tau^k}\,, \textit{ as } n \to \infty.$$

*where $\tau$ is given by $m_\psi(\tau) = L$ and $k$ is a multiple of $Q_\psi$.*

The case where $k = n-1$ and thus $L = 1$ (and $\omega = 0$) shall be used in the discussion of the Otter-Meir-Moon Theorem, Theorem 6.7.1, on an asymptotic formula for the coefficients of the solutions of Lagrange's equation.

- *Suppose $L = M_\psi$ and $R < +\infty$.* Observe that $M_\psi < \infty$.

Now there is no $\tau \in (0, R)$ such that $m_\psi(\tau) = L$. But $\psi$ extends to be continuous in $\text{cl}(\mathbb{D}(0, R))$ and the family $(Y_t)_{t\in[0,R)}$ may be completed with a variable $Y_R$ which has $m_\psi(R) = M_\psi = L$, see Section 1.3.3. If further $\sigma_\psi^2(R) = \mathbf{V}(Y_R) < \infty$, which occurs if $\sum_{n=0}^\infty n^2 b_n R^n < \infty$, then we may use formula (6.2.3) instead of formula (6.2.1). By the same argument above, we then have

**Theorem 6.3.3.** *If $M_\psi = L$ and $R < +\infty$ and $\sum_{n=0}^\infty n^2 b_n R^n < \infty$. If*

$$\lim_{n\to\infty}\frac{nL-k}{\sqrt{n}} = -\omega\,,$$

*then*

$$\text{COEFF}_{[k]}(\psi(z)^n) \sim Q_\psi\frac{e^{-\omega^2/(2\sigma_\psi^2(R))}}{\sqrt{2\pi}\,\sigma_\psi(R)}\frac{1}{\sqrt{n}}\frac{\psi^n(R)}{R^k}\,,$$

*where $k$ is a multiple of $Q_\psi$.*

**Remark 6.3.4.** There remains the case where $L = M_f$ and $R = +\infty$. In this case, since $M_f < +\infty$, we have that $f$ is a polynomial of degree $\deg(f) = L = M_f$ and thus, in particular, we have that $L$ is an integer. In the extreme case when $k = Ln$, we have that $\text{COEFF}_{[k]}(\psi^n) = b_L^n$, where $b_L$ is the $L$-th coefficient $\psi(z)$. ⊠

**A. Binomial coefficients.**

As an illustration, we apply next Theorem 6.3.2 to binomial coefficients. We take $\psi(z) = 1+z$, which belongs to $\mathcal{K}$. In this case, we have

$$m_\psi(t) = \frac{t}{1+t}, \quad \text{and} \quad \sigma^2_\psi(t) = \frac{t}{(1+t)^2}, \quad \text{for } t \in (0,+\infty)\,.$$

In particular, $M_\psi = 1$.

Observe that

$$\text{COEFF}_{[k]}(\psi(z)^n) = \binom{n}{k}, \quad \text{for } n \geq 1 \text{ and } k \geq 1\,.$$

• Let $p \in (0,1)$. For each $n \geq 1/p$, let $k = \lfloor pn \rfloor$. We have that

$$0 \leq p - \frac{k}{n} \leq \frac{1}{n}\,.$$

Thus we may apply Theorem 6.3.2 with $L = p$ and $\omega = 0$, to deduce that

$$\binom{n}{\lfloor pn \rfloor} \sim \frac{1}{\sqrt{2\pi np(1-p)}} \frac{1}{(1-p)^{n-k}p^k}, \quad \text{as } n \to \infty\,.$$

If $p = 1/2$, we have that

$$\binom{n}{\lfloor n/2 \rfloor} \sim \frac{2^{n+1}}{\sqrt{2\pi}\sqrt{n}}, \quad \text{as } n \to \infty\,.$$

• Let $p \in (0,1)$ and $\lambda \in \mathbb{R}$. For $n \geq N$, let $k = \lfloor pn + \lambda\sqrt{n} \rfloor$, where $N$ is chosen so that $pn + \lambda\sqrt{n} \geq 1$, for $n \geq N$. Then

$$\frac{k}{n} = p + \frac{\lambda}{\sqrt{n}} + O\Big(\frac{1}{n}\Big), \quad \text{as } n \to \infty\,,$$

and Theorem 6.3.2 with $L = p$ and $\omega = \lambda$, gives us that

$$\binom{n}{\lfloor pn + \lambda\sqrt{n} \rfloor} \sim \frac{1}{\sqrt{2\pi np(1-p)}} \frac{1}{(1-p)^{n-k}p^k} e^{-\lambda^2/(2p(1-p))}, \quad \text{as } n \to \infty\,.$$

For $p = 1/2$ and $\lambda \in \mathbb{R}$ fixed, we have that

$$\binom{n}{\lfloor n/2 + \lambda\sqrt{n} \rfloor} \sim \frac{2^{n+1}}{\sqrt{2\pi}\sqrt{n}} e^{-2\lambda^2}, \quad \text{as } n \to \infty\,.$$

**B. Back to Local Central Limit Theorem**

Let us assume that $\psi$ is the probability generating function of a certain random variable $X$ (with values in $\{0, 1, 2 \ldots\}$) which has mean $\mu$ and standard deviation $s$. Assume further that $\psi$ has radius of convergence $R > 1$ and that $Q_\psi = 1$.

We have $M_\psi = \lim_{t\uparrow R} m_\psi(t) > m_\psi(1) = \mu$. For $\tau = 1 \in (0, R)$, we have $m_\psi(\tau) = \mu$ and $\sigma_\psi(\tau) = s$.

Let $Y$ be the random variable $Y = (X - \mu)/s$. This variable $Y$ is a lattice random variable with gauge $h = 1/s$, since $Q_\psi = 1$. The variable $Y$ has $\mathbf{E}(Y) = 0$ and $\mathbf{E}(Y^2) = 1$.

For each $n \geq 1$, let $S_n$ denote the random variable

$$S_n = \frac{Y_1 + \cdots + Y_n}{\sqrt{n}} = \frac{X_1 + \cdots + X_n}{s\sqrt{n}} - \Big(\frac{\mu}{s}\Big)\sqrt{n}\,,$$

where $X_1, X_2, \ldots$ are independent copies of $X$ and $Y_j = (X_j - \mu)/s$, for $j \geq 1$.

For $k \in \mathbb{Z}$, denote

$$v(k) = k\frac{1}{s\sqrt{n}} - \Big(\frac{\mu}{s}\Big)\sqrt{n}\,.$$

Let $\mathcal{V}_n$ denote the collection of values that $S_n$ can take:

$$\mathcal{V}_n = \Big\{v(k) : k \in \mathbb{Z}\Big\}\,.$$

For $x \in \mathcal{V}_n$, let $k_x = s\sqrt{n}x + \mu n$, so that $v(k_x) = x$. We have

$$\frac{k_x}{n} = \mu + \frac{s}{\sqrt{n}}x\,.$$

Moreover, with the notations of the hypothesis of Theorem 6.3.2, we have that $L = \mu$ and $\omega = sx$. Now we use that with $\tau = 1$, we have that $m_\psi(1) = \mu(= L)$, $\psi(1) = 1$, $\sigma_\psi(1) = s$ and $\omega^2/\sigma^2_\psi(1) = x^2$, and appealing to Theorem 6.3.2 we obtain that

$$\text{COEFF}_{[k]}(\psi(z)^n) \sim \frac{1}{\sqrt{2\pi}}\frac{1}{s\sqrt{n}}e^{-x^2/2}\,, \quad \text{as } n \to \infty\,.$$

Since

$$\mathbf{P}(S_n = x) = \mathbf{P}(X_1 + \cdots + X_n = k_x) = \text{COEFF}_{[k_x]}(\psi(z)^n)\,,$$

we deduce that

$$\lim_{n\to\infty} s\sqrt{n}\,\mathbf{P}(S_n = x) = \frac{1}{\sqrt{2\pi}}e^{-x^2/2}\,, \quad \text{as } n \to \infty\,,$$

which is "consistent" with the Local Central Limit Theorem as stated in Corollary 6.1.6.

## 6.4 Index $k$ is little 'o' of $n$

In this section we deal with the large powers asymptotics when $k/n \to 0$ as $k, n \to \infty$.

For the discussion of this case and *throughout this section we assume that* $\psi'(0) \neq 0$. This hypothesis implies that $Q_\psi = 1$.

We will use formula (6.2.1) with an appropriate choice of $t$. To precise the limit value of the integral term as $n \to \infty$, we shall appeal to the following lemma, akin to Lemma 2 in [65].

**Lemma 6.4.1.** *Suppose that $\psi \in \mathcal{K}$ and $\psi'(0) \neq 0$, then for each $\theta_0 \in (0,\pi)$, there exists $c>0$ and $r>0$ such that*

$$|\mathbf{E}(e^{\imath\theta\breve{Y}_t})| = \frac{|\psi(te^{\imath\theta})|}{\psi(t)} \leq e^{-ct}\,, \quad \textit{if } t \leq r \textit{ and } \theta_0 \leq |\theta| \leq \pi\,.$$

*Proof.* Without loss of generality we assume that $b_0 = 1$. We have

$$|\psi(te^{\imath\theta})|^2 = 1 + 2b_1 t\cos\theta + O(t^2)\,, \quad \text{as} \quad t\downarrow 0,$$

and

$$|\psi(t)|^2 = 1 + 2b_1 t + O(t^2)\,, \quad \text{as} \quad t\downarrow 0.$$

Therefore,

$$|\mathbf{E}(e^{\imath\theta\breve{Y}_t})|^2 = \frac{|\psi(te^{\imath\theta})|^2}{\psi(t)^2} = 1 + 2b_1 t\,(\cos\theta - 1) + O(t^2)\,, \quad \text{as} \quad t\downarrow 0.$$

For $\theta_0 \leq |\theta| \leq \pi$, we have $\cos\theta \leq 1-\delta$, for certain $\delta>0$, which may depend upon $\theta_0$. Since, by hypothesis, $b_1>0$, we then have, for certain $r\in(0,R)$, that

$$|\mathbf{E}(e^{\imath\theta\breve{Y}_t})|^2 = \frac{|\psi(te^{\imath\theta})|^2}{\psi(t)^2} \leq 1 - 2b_1 t\delta + O(t^2)\,, \quad \text{for all } t \leq r\,,$$

Therefore, for $t \leq r$ we have that

$$|\mathbf{E}(e^{\imath\theta\breve{Y}_t})|^2 = \frac{|\psi(te^{\imath\theta})|^2}{\psi(t)^2} \leq 1 - b_1 t\delta \leq e^{-b_1\delta t}\,,$$

as claimed. □

For $n$ large enough, we have that $k/n < M_\psi$ and thus we may define uniquely $\tau_n \in (0,R)$ such that $m_\psi(\tau_n) = k/n$. Observe that $\tau_n \to 0$, as $n\to\infty$.

We will use here the results from Section 2.2 in Chapter 2 about the moments of the Khinchin family $(Y_t)_{t\in[0,R)}$ when $t$ is close to 0. Since $b_1>0$, we have that

$$(\natural) \quad \lim_{t\downarrow 0} t^{1/2}\,\mathbf{E}(|\breve{Y}_t|^3) = \sqrt{\frac{b_0}{b_1}}\,.$$

As $m_\psi(t) = \dfrac{b_1}{b_0}t + O(t^2)$, we have that

$$\tau_n \sim \frac{b_0}{b_1}\frac{k}{n}\,, \quad \text{as } n\to\infty\,.$$

Besides,

$$\sigma_\psi^2(\tau_n) \sim \frac{k}{n}\,, \quad \text{as } n\to\infty\,.$$

For $t = \tau_n$ the expression ($\natural$) above tells us that

$$(\flat) \qquad \frac{\mathbf{E}(|\breve{Y}_{\tau_n}|^3)}{\sqrt{n}} \sim \sqrt{\frac{1}{k}}, \quad \text{as } k \to \infty.$$

If we choose $t = \tau_n$ in formula (6.2.1), the integral term simplifies to

$$I_n \triangleq \int_{|\theta| \le \pi \sigma_\psi(\tau_n)\sqrt{n}} \mathbf{E}\big(e^{\imath\theta \breve{Y}_{\tau_n}/\sqrt{n}}\big)^n d\theta .$$

Let us see that the integral $I_n$ converges to $\sqrt{2\pi}$ as $n \to \infty$.

Theorem 6.1.1 and ($\flat$) give that

$$\lim_{n\to\infty} \mathbf{E}\big(e^{\imath \breve{Y}_{\tau_n}\theta/\sqrt{n}}\big)^n = e^{-\theta^2/2}, \quad \text{for all } \theta \in \mathbb{R}.$$

Theorem 6.1.2 and ($\flat$) give that there exists $N \ge 1$ such that if $n \ge N$, then

$$\big|\mathbf{E}\big(e^{\imath \breve{Y}_{\tau_n}\theta/\sqrt{n}}\big)^n\big| \le e^{-\theta^2/3}, \quad \text{for all } |\theta| \le \sqrt{k}/5 .$$

By dominated convergence we have

$$\lim_{n\to\infty} \int_{|\theta|\le\sqrt{k}/5} \big|\mathbf{E}\big(e^{\imath \breve{Y}_{\tau_n}\theta/\sqrt{n}}\big)^n - e^{-\theta^2/2}\big|\, d\theta = 0 .$$

We apply lemma 6.4.1 with $\theta_0 = 1/10$ to ascertain that there exists $N \ge 1$ and $C > 0$ such that for all $n \ge N$ we have that

$$|\mathbf{E}(e^{\imath\theta Y_{\tau_n}})| = \frac{|\psi(\tau_n e^{\imath\theta})|}{\psi(\tau_n)} \le e^{-c\tau_n}, \quad \text{for } \tfrac{1}{10} \le |\theta| \le \pi ,$$

and, therefore, that

$$|\mathbf{E}(e^{\imath\theta Y_{\tau_n}/(\sigma_\psi(\tau_n)\sqrt{n})})|^n = \Big|\frac{\psi(\tau_n e^{\imath\theta/(\sigma_\psi(\tau_n)\sqrt{n})})}{\psi(\tau_n)}\Big|^n \le e^{-cn\tau_n},$$

for $\frac{1}{10}\sigma_\psi(\tau_n)\sqrt{n} \le |\theta| \le \pi\sigma_\psi(\tau_n)\sqrt{n}$.

Since for $n$ large enough we have that $\frac{1}{10}\sigma_\psi(\tau_n)\sqrt{n} \le \sqrt{k}/5$, we may bound

$$\int_{\sqrt{k}/5 \le |\theta| \le \pi\sigma_\psi(\tau_n)\sqrt{n}} \big|\mathbf{E}(e^{\imath\theta \breve{Y}_{\tau_n}/\sqrt{n}})^n - e^{-\theta^2/2}\big| d\theta \le e^{-cn\tau_n}(8\pi/3)\sqrt{k} + \int_{|\theta|\ge\sqrt{k}/3} e^{-\theta^2/2}\, d\theta ,$$

for $n$ large enough.

As $n \to \infty$, the first summand of this bound converges to 0, since $n\tau_n \sim (a_0/a_1)k$ and $k \to \infty$, while the second summand converges to 0 since $k \to \infty$. We conclude that

$$\lim_{n\to\infty} \int\limits_{\sqrt{k}\le|\theta|\le\pi\sigma_\psi(\tau_n)\sqrt{n}} \left|\mathbf{E}(e^{\imath\theta\check{Y}_{\tau_n}/\sqrt{n}})^n - e^{-\theta^2/2}\right| d\theta = 0\,,$$

and, therefore, that

$$\lim_{n\to\infty} \int\limits_{|\theta|\le\pi\sigma_\psi(\tau_n)\sqrt{n}} \left|\mathbf{E}(e^{\imath\theta\check{Y}_{\tau_n}/\sqrt{n}})^n - e^{-\theta^2/2}\right| d\theta = 0\,. \tag{6.4.1}$$

In particular

$$\lim_{n\to\infty} I_n = \lim_{n\to\infty} \int\limits_{|\theta|\le\pi\sigma_\psi(\tau_n)\sqrt{n}} \mathbf{E}(e^{\imath\theta\check{Y}_{\tau_n}/\sqrt{n}})^n \, d\theta = \sqrt{2\pi}\,.$$

In sum,

**Theorem 6.4.2.** *If $\psi'(0) \neq 0$, then*

$$\begin{aligned}\text{COEFF}_{[k]}(\psi(z)^n) &\sim \frac{1}{\sqrt{2\pi}}\frac{\psi^n(\tau_n)}{\tau_n^k}\frac{1}{\sqrt{n}}\frac{1}{\sigma_\psi(\tau_n)}\\ &\sim \frac{1}{\sqrt{2\pi}}\frac{\psi^n(\tau_n)}{\tau_n^k}\frac{1}{\sqrt{k}}\,, \quad \textit{as } n \to \infty\,,\end{aligned}$$

*where $k/n \to 0$ and $\tau_n$ is given by $m_\psi(\tau_n) = k/n$ (for $n$ large enough).*

### 6.4.1 With further information on $k/n$

If further information on how fast $k/n$ tends to 0 is available, then we may express the asymptotic formula (6.4.2) directly in terms of $k$ and $n$ and not on $\tau_n$.

We maintain the assumption that $\psi'(0) > 0$.

The function $m_\psi$ is holomorphic and injective near 0. Its inverse $m_\psi^{-1}$ is actually the solution of Lagrange's equation with data $z/m_\psi(z) = \psi(z)/\psi'(z)$. Besides, $\ln\psi$ is holomorphic near 0, since as $\psi \in \mathcal{K}$, we have $\psi(0) > 0$. Thus Lagrange's inversion formula gives the expansion

$$\ln\psi\big(m_\psi^{-1}\big)(z) = \ln b_0 + \sum_{j=1}^{\infty} B_j z^j,$$

where

$$B_j = \frac{1}{j}\text{COEFF}_{[j-1]}\left(\left(\frac{\psi'}{\psi}\right)\left(\frac{\psi}{\psi'}\right)^j\right) = \frac{1}{j}\text{COEFF}_{[j-1]}\left(\left(\frac{\psi}{\psi'}\right)^{j-1}\right), \quad \text{for } j \ge 1\,. \tag{6.4.2}$$

For each $j \ge 1$, the coefficient $B_j$ depends only on $b_0, b_1, \ldots, b_j$.

Likewise,

$$\ln \frac{m_\psi^{-1}(z)}{z} = \ln\Big(\frac{\psi}{\psi'}\Big)\big(m_\psi^{-1}\big)(z) = \ln \frac{b_0}{b_1} + \sum_{j=1}^{\infty} C_j z^j,$$

where

$$C_j = \frac{1}{j}\text{COEFF}_{[j]}\Big(\frac{\psi}{\psi'}\Big)^j = \frac{j+1}{j} B_{j+1}\,, \quad \text{for } j \ge 1\,. \tag{6.4.3}$$

Now, with the notations of the general discussion of the case $k/n \to 0$, we have that

$$\ln \psi(\tau_n)^n = n \ln \psi(\tau_n) = n \ln \psi\Big(m_\psi^{-1}\Big(\frac{k}{n}\Big)\Big) = n \ln b_0 + \sum_{j=1}^{\infty} B_j (k^j/n^{j-1})\,,$$

and also that,

$$\ln \tau_n^k = k \ln \tau_n = k \ln m_\psi^{-1}\Big(\frac{k}{n}\Big) = k \ln \frac{b_0}{b_1} + k \ln \frac{k}{n} + \sum_{j=1}^{\infty} C_j \frac{k^{j+1}}{n^j}.$$

Using formula (6.4.3), and that $B_1 = 1$, we may write

$$\ln \psi(\tau_n)^n - \ln \tau_n^k = (n-k)\ln b_0 + k \ln b_1 + k \ln \frac{n}{k} + k - \sum_{j=2}^{\infty} \frac{B_j}{j-1} \frac{k^j}{n^{j-1}}. \tag{6.4.4}$$

Let $\lambda = \limsup_{n\to\infty} \ln k/\ln n$ and let $J$ the smallest integer such that $J > \lambda/(1-\lambda)$; if $\lambda = 1$, we set $J = +\infty$. Thus

$$\sum_{j=J+1}^{\infty} \frac{B_j}{j-1}\frac{k^j}{n^{j-1}} = nO\Big(\frac{k}{n}\Big)^{J+1} = o(1)\,, \quad \text{as } n \to \infty\,,$$

and then

$$\frac{\psi(\tau_n)^n}{\tau_n^k} \sim b_0^{n-k} b_1^k \frac{n^k e^k}{k^k} \exp\Big(-\sum_{j=2}^{J} \frac{B_j}{j-1}\frac{k^j}{n^{j-1}}\Big)\,, \quad \text{as } n \to \infty\,.$$

Therefore

**Theorem 6.4.3.** *If $\psi'(0) \neq 0$, and $\lambda$ and $J$ as above, then*

$$\text{COEFF}_{[k]}(\psi(z)^n) \sim \frac{1}{\sqrt{2\pi}} b_0^{n-k} b_1^k \frac{n^k e^k}{k^k \sqrt{k}} \exp\Big(-\sum_{j=2}^{J} \frac{B_j}{j-1}\frac{k^j}{n^{j-1}}\Big)\,, \quad \textit{as } n \to \infty\,.$$

If $\lambda < 1/2$, then $J = 1$, and then

$$\text{COEFF}_{[k]}(\psi(z)^n) \sim b_0^{n-k} b_1^k \frac{n^k}{k!} \sim \text{COEFF}_{[k]}((b_0 + b_1 z)^n)\,, \quad \text{as } n \to \infty\,. \tag{6.4.5}$$

For the particular case where, $k = \lfloor \sqrt{n} \rfloor$, (with $\lambda = 1/2$, and $J = 2$), we have that

$$\text{COEFF}_{[k]}(\psi(z)^n) \sim b_0^{n-k} b_1^k \frac{n^k}{k!} e^{-B_2}\,, \quad \text{as } n \to \infty\,.$$

where $B_2 = (1/2) - (b_2 b_0)/b_1^2$.

### 6.4.2 Particular case: $k$ fixed, while $n \to \infty$

In this particular case, where $k$ remains fixed, while $n \to \infty$, the analysis of the asymptotic behavior of $\text{COEFF}_{[k]}(\psi(z)^n)$ just requires the multinomial theorem.

Observe that

$$\text{COEFF}_{[k]}(\psi(z)^n) = \text{COEFF}_{[k]}\Big(\sum_{j=0}^{k} b_j z^j\Big)^n .$$

The multinomial theorem gives then that

$$\text{COEFF}_{[k]}(\psi(z)^n) = \sum \frac{n!}{j_0! \dots j_k!} b_0^{j_0} \dots b_k^{j_k} ,$$

where the sum extends to all $k+1$-tuples $(j_0, \dots, j_k)$ of non-negative integers such that

$$\begin{aligned} j_0 + j_1 + \cdots + j_k &= n, \\ j_1 + 2j_2 + \cdots + kj_k &= k. \end{aligned}$$

For such a tuple, we have $n - k \le j_0 \le n$ and $j_i \le k$, for $1 \le i \le k$. We write $j_0 = n - l$, with $0 \le l \le k$ and classify by $l$:

$$\text{COEFF}_{[k]}(\psi(z)^n) = \sum_{l=0}^{k} \binom{n}{l} b_0^{n-l} \underbrace{\sum \binom{j_1 + \cdots + j_k}{j_1, \dots, j_k} b_1^{j_1} \dots b_k^{j_k}}_{\triangleq C_l} .$$

For each integer $l, 0 \le l \le k$, the interior sum $C_l$ extends to all $k$-tuples with $j_1 + \cdots + j_k = l$ and $j_1 + 2j_2 + \cdots + kj_k = k$. Observe that $C_l$ does not depend upon $n$. Therefore

$$\text{COEFF}_{[k]}(\psi(z)^n) = \sum_{l=0}^{k} \binom{n}{l} b_0^{n-l} C_l,$$

is a polynomial in $n$ whose degree $\gamma$ is given by the largest $l \le k$ such that $C_l \ne 0$, and

$$\text{COEFF}_{[k]}(\psi(z)^n) \sim \frac{1}{\gamma!} b_0^{n-\gamma} C_\gamma \, n^\gamma , \quad \text{as } n \to \infty .$$

- If $b_1 \ne 0$, then $C_k \ne 0$, and $\gamma = k$. In fact, the sum $C_k$ contains just one summand corresponding to $j_1 = k$ and $j_2 = \ldots = j_k = 0$, and thus $C_k = b_1^k$. We have in this case that

$$\text{COEFF}_{[k]}(\psi(z)^n) \sim \frac{1}{k!} b_0^{n-k} b_1^k \, n^k \sim \binom{n}{k} b_0^{n-k} b_1^k , \quad \text{as } n \to \infty .$$

In other terms and asymptotically speaking, the $k$-th coefficient of $\psi(z)^n$ behaves as the $k$-th coefficient of $(b_0 + b_1 z)^n$, as $n \to \infty$. Compare with the formula (6.4.5), where $k \to \infty$ with $n$, but slowly.

- If $b_1 = 0$, then the nonzero summands of $C_l$ must have $j_1 = 0$ and they must satisfy $j_2 + \cdots + j_k = l$ and $2j_2 + \cdots + kj_k = k$. Thus, if $C_l$ is not zero, then $2l \le k$, and so, $\gamma \le k/2$.

Assume further that $b_2 \neq 0$. We distinguish now between $k$ even and $k$ odd.

If $k = 2q$, then $C_q$ has only one nonzero summand, the one corresponding to $j_2 = q$. Thus $C_q = b_2^q$, $\gamma = q = k/2$, and

$$\text{COEFF}_{[k]}(\psi(z)^n) \sim \frac{1}{(k/2)!} b_0^{n-k/2} b_2^{k/2}\, n^{k/2}, \quad \text{as } n \to \infty.$$

If $k = 2q+1$, then $\gamma \le q$. The only nonzero summand of $C_q$ has $j_2 = q-1, j_3 = 1$. Thus $C_q = q b_2^{q-1} b_3$.

If $b_3 \neq 0$, then $C_q \neq 0$, $\gamma = q$ and

$$\text{COEFF}_{[k]}(\psi(z)^n) \sim \frac{1}{((k-3)/2)!}\, b_0^{n-(k-1)/2} b_2^{(k-3)/2} b_3\, n^{(k-1)/2}, \quad \text{as } n \to \infty.$$

If $b_3 = 0$, then $C_q = 0$ and actually the degree $\gamma$ is at most $k/3$.

- In general, if $b_1 = \ldots = b_{m-1} = 0$, and $b_m \neq 0$, then for each $h, 0 \le h \le m-1$, the coefficients of index $k \equiv h \mod (m)$, satisfy an asymptotic formula

$$\text{COEFF}_{[k]}(\psi(z)^n) \sim \Omega_{k,h} n^{\gamma_{k,h}}, \quad \text{as } n \to \infty,$$

where $\gamma_{k,h} \le (k-h)/m$ and $\Omega_{k,h} > 0$.

For $h = 0$, i.e., for those $k$ which are multiples of $m$, we have $\gamma_{k,0} = k/m$ and, actually, that

$$\text{COEFF}_{[k]}(\psi(z)^n) \sim \frac{1}{(k/m)!} b_0^{n-k/m} b_m^{k/m}\, n^{k/m}, \quad \text{as } n \to \infty.$$

In particular, if $Q = Q_\psi > 1$ and $k$ is a multiple of $Q$, then

$$\text{COEFF}_{[k]}(\psi(z)^n) \sim \frac{1}{(k/Q)!} b_0^{n-k/Q} b_Q^{k/Q}\, n^{k/Q}, \quad \text{as } n \to \infty.$$

Observe the $b_Q > 0$. If $Q$ is not a divisor of $k$, then, of course, $\text{COEFF}_{[k]}(\psi(z)^n) = 0$. for every $n \ge 1$.

## 6.5 Power $n$ is little 'o' of index $k$

To obtain a large powers asymptotic formula in this case where $k/n \to \infty$, as $n \to \infty$, we require more from $\psi$, namely, we will assume that $\psi$ is uniformly strongly Gaussian. In an a previous chapter of this document, we describe at some length these so called uniformly strongly Gaussian power series.

Since uniformly strongly Gaussian power series are strongly Gaussian, Hayman's asymptotic formula (3.2.2) is valid for $\psi$ and thus $Q_\psi = 1$.

The assumption of uniform strongly Gaussianity on $\psi$ is to be compared with the assumptions of Gardy in [36, Theorems 5 and 6], see also [36, Section 6].

The exponential generating function of set partitions, $e^{e^z-1}$, and the ordinary generating function of partitions of integers, $\prod_{j=1}^{\infty} 1/(1-z^j)$ are examples of uniformly strongly Gaussian power series.

Since $\psi$ is uniformly strongly Gaussian, we have that $M_\psi = +\infty$, and thus we may define (uniquely) $\tau_n \in (0, R)$ as being such that $m_\psi(\tau_n) = k/n$. Observe that $\tau_n \uparrow R$, as $n \to \infty$.

Now, formula (6.2.1) of Lemma 6.2.2, with $t = \tau_n$, gives us that

$$\text{COEFF}_{[k]}(\psi(z)^n) = \frac{1}{2\pi} \frac{\psi^n(\tau_n)}{\tau_n^k} \frac{1}{\sqrt{n}} \frac{1}{\sigma_\psi(\tau_n)} \int\limits_{|\theta| \le \pi\sigma_\psi(\tau_n)\sqrt{n}} \mathbf{E}\big(e^{\imath\theta \breve{Y}_{\tau_n}/\sqrt{n}}\big)^n d\theta\,.$$

Since $\tau_n \to R$, uniform Gaussianity of $\psi$ implies that

$$\lim_{n\to\infty} \int\limits_{|\theta| \le \pi\sigma_\psi(\tau_n)\sqrt{n}} \mathbf{E}\big(e^{\imath\theta \breve{Y}_{\tau_n}/\sqrt{n}}\big)^n d\theta = \sqrt{2\pi}\,,$$

and we conclude that

**Theorem 6.5.1.** *If $\psi$ is uniformly strongly Gaussian, then*

$$\text{COEFF}_{[k]}(\psi(z)^n) \sim \frac{1}{\sqrt{2\pi}} \frac{\psi^n(\tau_n)}{\tau_n^k} \frac{1}{\sqrt{n}} \frac{1}{\sigma_\psi(\tau_n)}\,, \quad \text{as } n \to \infty\,,$$

*where $k/n \to \infty$ as $n \to \infty$ and $\tau_n$ is such that $m_\psi(\tau_n) = k/n$.*

## 6.6 Coefficients of $h(z)\psi(z)^n$

Next we consider asymptotics of the coefficients of $h(z)\psi(z)^n$, as $n \to \infty$, where both $h$ and $\psi$ are in $\mathcal{K}$. This is to be used in Section 6.8 which deals with Lagrangian distributions and Galton-Watson processes: the power series $h$ would codify the initial offspring distribution while the power series $\psi$ would codify the probability distributions subsequent offspring distributions.

We assume here from the start that $Q_\psi = 1$.

Let $h(z) = \sum_{j=0}^{\infty} c_j z^j$ be a power series in $\mathcal{K}$. We assume throughout that the radius of convergence of $h$ is at least $R_\psi$, the radius of convergence of $\psi$.

We treat jointly the cases in which $k \asymp n$ and $k = o(n)$, with $k \to \infty$, by adjusting the arguments in the previous sections (whose notations we will use liberally) in which $h \equiv 1$.

Assume the hypothesis of Theorem 6.3.1, when $k \asymp n$, and of Theorem 6.4.2, when $k = o(n)$ and $k \to \infty$.

The $k$-th coefficient of $h(z)\psi(z)^n$ is given, see (6.2.2), by

$$\frac{1}{2\pi} \frac{\psi^n(\tau_n)}{\tau_n^k} \frac{1}{\sqrt{n}} \frac{1}{\sigma_\psi(\tau_n)} \int\limits_{|\theta| \le \pi\sigma_\psi(\tau_n)\sqrt{n}} h(\tau_n e^{\imath\theta/(\sigma_\psi(\tau_n)\sqrt{n})}) \mathbf{E}\big(e^{\imath\theta \breve{Y}_{\tau_n}/\sqrt{n}}\big)^n d\theta\,,$$

where we are taking $\tau_n$ such that $m_\psi(\tau_n) = k/n$.

Denote by $I_n$ the integral on the right hand-side of the previous expression.

For $0 \le s \le r < R$ and $\phi \in \mathbb{R}$ we have

$$\left|h(se^{\imath\phi}) - h(s)\right| \le \max_{|z|\le r} |h'(z)|\, s\, |\phi|\,.$$

For $k \asymp n$, we have $\tau_n \le m_\psi^{-1}(B) < R$, and for $k = o(n)$ we have $\tau_n \to 0$, so, in both cases, there exists a constant $K > 0$ such that

$$(\nabla) \qquad \left|h(\tau_n e^{\imath\theta/(\sigma_\psi(\tau_n)\sqrt{n})}) - h(\tau_n)\right| \le K \frac{|\theta|}{\sigma_\psi(\tau_n)\sqrt{n}}\,.$$

Denote by $\mathcal{D}_n$ the interval $\{\theta \in \mathbb{R} : |\theta| \le \pi\sigma_\psi(\tau_n)\sqrt{n}\}$.

Using the bound $(\nabla)$ and writing

$$\mathbf{E}\big(e^{\imath\theta\breve{Y}_t/\sqrt{n}}\big)^n = \Big(\mathbf{E}\big(e^{\imath\theta\breve{Y}_t/\sqrt{n}}\big)^n - e^{-\theta^2/2}\Big) + e^{-\theta^2/2}\,,$$

we conclude that

$$\begin{aligned} I_n &= h(\tau_n)\int_{\mathcal{D}_n} e^{-\theta^2/2}d\theta + h(\tau_n)\int_{\mathcal{D}_n}\Big(\mathbf{E}\big(e^{\imath\theta\breve{Y}_{\tau_n}/\sqrt{n}}\big)^n - e^{-\theta^2/2}\Big)d\theta \\ &+ O\Big(\int_{\mathcal{D}_n}\Big|\mathbf{E}\big(e^{\imath\theta\breve{Y}_{\tau_n}/\sqrt{n}}\big)^n - e^{-\theta^2/2}\Big|\frac{|\theta|}{\sigma_\psi(\tau_n)\sqrt{n}}d\theta\Big) + O\Big(\int_{\mathcal{D}_n} e^{-\theta^2/2}\frac{|\theta|}{\sigma_\psi(\tau_n)\sqrt{n}}d\theta\Big)\,. \end{aligned}$$

For $\theta \in \mathcal{D}_n$, we have that $\dfrac{|\theta|}{\sigma_\psi(\tau_n)\sqrt{n}} \le \pi$, and thus we see that the third term in the sum above tends to 0, by virtue of (6.3.1) (if $k \asymp n$) or (6.4.1) (if $k = o(n)$).

The second term in the sum tends to 0, since $h(\tau_n)$ is bounded and the integral converges to 0, by (6.3.1) and (6.4.1), respectively.

The fourth term in the sum tends to 0 since $\sigma_\psi(\tau_n)\sqrt{n} \asymp \sqrt{k}$ and $\int_{\mathbb{R}} e^{-\theta^2/2}|\theta|d\theta = 2$.

The first term in the sum is $h(\tau_n)\sqrt{2\pi} + o(1)$, since $h(\tau_n)$ is bounded, and so

$$I_n = h(\tau_n)\sqrt{2\pi} + o(1)\,, \qquad \text{as } n \to \infty\,,$$

but, using that $h(\tau_n)$ is bounded from below by $h(0)$ and recalling that $h \in \mathcal{K}$, we conclude that

$$I_n = h(\tau_n)\sqrt{2\pi}(1 + o(1))\,, \qquad \text{as } n \to \infty\,.$$

We summarize the discussion above as follows.

**Theorem 6.6.1.** *With the notations above,*

- *For the case $k \asymp n$, one has that*

$$\text{COEFF}_{[k]}(h(z)\psi(z)^n) \sim \frac{1}{\sqrt{2\pi}}\frac{h(\tau_n)\psi^n(\tau_n)}{\tau_n^k}\frac{1}{\sqrt{n}}\frac{1}{\sigma_\psi(\tau_n)}\,, \qquad \text{as } n \to \infty\,.$$

- *For the case $k = o(n)$ and $k \to \infty$, one has that*

$$\text{COEFF}_{[k]}(h(z)\psi(z)^n) \sim \frac{h(0)}{\sqrt{2\pi}} \frac{\psi^n(\tau_n)}{\tau_n^k} \frac{1}{\sqrt{k}}, \quad \textit{as } n \to \infty .$$

Observe that, in the case $k = o(n)$ and $k \to \infty$, we have that $h(\tau_n) \to h(0)$ as $n \to \infty$.

If more information is available on how $k/n$ tends to 0, as discussed in Section 6.4.1, we could apply the argument and the conclusions obtained there to refine the above asymptotic formula.

Suppose that $k/n \to L > 0$, as $n \to \infty$ and that equation (6.3.2) is satisfied. Assume that $L < M_\psi$. Then for $\tau = m_\psi^{-1}(L)$, we deduce, using the Hayman's formula (6.2.2) and arguing as in the case $h \equiv 1$, that

$$\text{COEFF}_{[k]}(h(z)\psi(z)^n) \sim \frac{1}{\sqrt{2\pi}} \frac{h(\tau)\psi^n(\tau)}{\tau^k} \frac{1}{\sqrt{n}} \frac{1}{\sigma_\psi(\tau)}, \quad \text{as } n \to \infty . \tag{6.6.1}$$

Recall that the power series $h$ has radius of convergence $R_\psi$ and that $\tau \in (0, R_\psi)$.

To conclude, we consider the case in which $k$ is fixed and $n \to \infty$. Denote $Q = Q_\psi = \{n \geq 1 : b_n > 0\}$. Then

$$\begin{aligned}\text{COEFF}_{[k]}(h(z)\psi(z)^n) &= \sum_{j=0}^{k} c_j \text{COEFF}_{[k]}(z^j\psi(z)^n) \\ &= \sum_{j=0}^{k} c_j \text{COEFF}_{[k-j]}(\psi(z)^n).\end{aligned}$$

The only $(k-j)$-th coefficients of $\psi(z)^n$ which are nonzero are the coefficients with indices $j \equiv k$, mod $Q$. Appealing to the case $h \equiv 1$, we have that

$$\begin{aligned}\text{COEFF}_{[k]}(h(z)\psi(z)^n) &= \sum_{\substack{0 \leq j \leq k \\ j \equiv k \text{ ,mod} Q}} c_j \text{COEFF}_{[k-j]}(\psi(z)^n) \\ &= \sum_{\substack{0 \leq j \leq k \\ j \equiv k \text{ ,mod} Q}} c_j \frac{1}{((k-j)/Q)!} b_0^{n-(k-j)/Q} \, b_Q^{(k-j)/Q} \, n^{(k-j)/Q} \, (1 + o(1)) .\end{aligned}$$

Denoting $j_0 = \min\{0 \leq j \leq k : j \equiv k, \mod Q \text{ and } c_j \neq 0\}$, we have that

$$\text{COEFF}_{[k]}(h(z)\psi(z)^n) \sim c_{j_0} \frac{1}{((k-j_0)/Q)!} b_0^{n-(k-j_0)/Q} \, b_Q^{(k-j_0)/Q} \, n^{(k-j_0)/Q}, \quad \text{as } n \to \infty .$$

In particular, if $h(0) = c_0 \neq 0$ and if $Q$ is a divisor of $k$, then

$$\text{COEFF}_{[k]}(h(z)\psi(z)^n) \sim c_0 \frac{1}{(k/Q)!} \, b_0^{n-k/Q} \, b_Q^{k/Q} \, n^{k/Q}, \quad \text{as } n \to \infty .$$

Observe that in this case, $j_0 = 0$.

## 6.7 Coefficients of solutions of Lagrange's equation

Next we shall apply the asymptotic results about large powers of Section 6.2 to obtain asymptotic formulae for the coefficients of solutions of Lagrange's equation when the data $\psi$ of the equation belongs to $\mathcal{K}$.

### 6.7.1 Lagrange's equation

We start with a function $\psi(z)$ which is holomorphic in a neighborhood of $z = 0$, has radius of convergence $R > 0$ and is such that $\psi(0) \neq 0$. Consider *Lagrange's equation*:

$$(\dagger) \qquad g(w) = w\psi(g(w))\,.$$

The *solution* $g(w)$ of Lagrange's equation with *data* $\psi(z)$ is a holomorphic function $g(w) = \sum_{n=1}^{\infty} A_n w^n$ which satisfies ($\dagger$) in a neighborhood of $w = 0$.

Assume that the data $\psi(z)$ of Lagrange's equation has power series expansion

$$\psi(z) = \sum_{n=0}^{\infty} b_n z^n, \quad \text{for all } z \in \mathbb{D}(0, R).$$

The holomorphic function $g(w)$, solution of Lagrange's equation with data $\psi(z)$, is unique. In fact, the coefficients $A_n$ of the Taylor expansion of $g(w)$ around $w = 0$ are given by *Lagrange's inversion formula*:

$$A_n = \frac{1}{n}\text{COEFF}_{[n-1]}(\psi(z)^n)\,, \quad \text{for all } n \geq 1\,.$$

For $n = 0$, we have $A_0 = g(0) = 0$.

This formula is exact for each coefficient $A_n$ of $g$ such that $n \geq 1$.

In a more general setting, consider the coefficients of $H(g(z))$, where $g$ is the solution of Lagrange's equation with data $\psi$ and $H$ is a holomorphic function with non-negative coefficients around $z = 0$, then the *extended Lagrange's inversion formula* gives

$$(6.7.1) \qquad \text{COEFF}_{[n]}(H(g(z))) = \frac{1}{n}\text{COEFF}_{[n-1]}(H'(z)\psi(z)^n), \quad \text{for all } n \geq 1\,,$$

and $\text{COEFF}_{[0]}(H(g(z))) = H(0)$.

A particular instance of this formula, that will be of interest later on, is given by the choice $H(z) = z^q$ for an integer $q \geq 1$. In this case we have

$$(6.7.2) \qquad \text{COEFF}_{[n]}(g(z)^q) = \frac{q}{n}\text{COEFF}_{[n-q]}(\psi(z)^n), \quad \text{for all } n \geq 1\,,$$

and $\text{COEFF}_{[0]}(g(z)^q) = H(0) = 0$.

If $Q_\psi > 1$, then the only nonzero coefficients $b_n$ are those where the index $n$ is a multiple of $Q_\psi$.Therefore, for the solution $g(w)$ of Lagrange's equation with data $\psi$, the only nonzero coefficients $A_n$ are those where the index $n$ verifies that $n - 1$ is a multiple of $Q_\psi$.

If $\psi(0) = 0$, then the only solution of Lagrange's equation is $g \equiv 0$, but, recall, nonetheless that $\psi \in \mathcal{K}$ implies that $\psi(0) > 0$.

### 6.7.2 The Otter-Meir-Moon Theorem and some extensions

The Otter-Meir-Moon Theorem, Theorem 6.7.1 below, comes from [75, Theorem 4] and [64, Theorem 3.1]. It gives an asymptotic formula for the coefficients $A_n$ of $g$, when $n-1$ is a multiple of $Q_\psi$ and $n \to \infty$, using minimal (but crucial) information about the function $\psi$.

We distinguish three cases, according as $M_\psi$ is $> 1, = 1$ or $< 1$.

♦ $M_\psi > 1$. This is the original assumption of both Otter and Meir-Moon; it means that there is a unique $\tau \in (0, R)$ such that $m_\psi(\tau) = 1$.

**Theorem 6.7.1** (Otter-Meir-Moon theorem)**.** *Let $\psi(z)$ be a power series in $\mathcal{K}$ with radius of convergence $R > 0$ and let $\psi(z) = \sum_{n=0}^{\infty} b_n z^n$ be its power series expansion.*

*Assume that $M_\psi > 1$ and let $\tau \in (0, R)$ be given by $m_\psi(\tau) = 1$.*

*Then the coefficients $A_n$ of the solution $g(w)$ of Lagrange's equation verify that*

- *if $n \not\equiv 1 \mod Q_\psi$, then $A_n = 0$,*
- *for the indices $n$ such that $n \equiv 1 \mod Q_\psi$ we have the asymptotic formula*

(6.7.3)
$$A_n \sim \frac{Q_\psi}{\sqrt{2\pi}} \frac{\tau}{\sigma_\psi(\tau)} \frac{1}{n^{3/2}} \Big(\frac{\psi(\tau)}{\tau}\Big)^n, \quad \text{as } n \to \infty\,.$$

Theorem 6.7.1 of Otter and Meir-Moon follows readily from Theorem 6.3.2 with $L = 1$ and $\omega = 0$, since $k = n - 1$.

Given that

$$\sigma_\psi^2(t) = m_\psi(t)(1 - m_\psi(t)) + t^2 \frac{\psi''(t)}{\psi(t)}, \quad \text{for all } t \in [0, R)\,,$$

we have

(6.7.4)
$$\frac{\psi''(\tau)}{\psi(\tau)} = \frac{\sigma_\psi^2(\tau)}{\tau^2}\,.$$

We thus may rewrite the conclusion of Theorem 6.7.1 (in Laplace's method or saddle point approximation style) as

$$A_n \sim Q_\psi \sqrt{\frac{\psi(\tau)}{2\pi\psi''(\tau)}} \frac{1}{n^{3/2}} \Big(\frac{\psi(\tau)}{\tau}\Big)^n, \quad \text{when } n{-}1 \text{ is a multiple of } Q_\psi \text{ and } n \to \infty\,.$$

♦ $M_\psi = 1$ (and $R < +\infty$). For the case $M_\psi = 1$, there is an Otter-Meir-Moon like theorem under a certain non-degeneracy condition on $\psi$. Recall the notations of Section 1.2.4.

**Theorem 6.7.2.** *Let $\psi(z)$ be a power series in $\mathcal{K}$ with radius of convergence $R > 0$ and let $\psi(z) = \sum_{n=0}^{\infty} b_n z^n$ be its power series expansion.*

*Assume that $M_\psi = 1$ and that $R < \infty$, and besides that $\sum_{n=0}^{\infty} n^2 b_n R^n < +\infty$ or equivalently $\psi''(R) < +\infty$.*

*Then the coefficients $A_n$ of the solution $g(w)$ of Lagrange's equation verify that*

- *if* $n \not\equiv 1 \mod Q_\psi$, *then* $A_n = 0$,
- *for the indices* $n$ *such that* $n \equiv 1 \mod Q_\psi$ *we have the asymptotic formula*

$$A_n \sim \frac{Q_\psi}{\sqrt{2\pi}} \frac{R}{\sigma_\psi(R)} \frac{1}{n^{3/2}} \Big(\frac{\psi(R)}{R}\Big)^n, \quad \textit{as } n \to \infty\,. \tag{6.7.5}$$

The power series $\psi$ extends to be continuous in $\mathrm{cl}(\mathbb{D}(0,R))$ and the Khinchin family $(Y_t)_{t\in[0,R)}$ extends (continuously in distribution) to include a variable $Y_R$ with mean 1 and variance $\sigma^2_\psi(R) = \lim_{t\uparrow R} \sigma^2_\psi(t)$ which is finite because by hypothesis $\sum_{n=0}^\infty n^2 b_n R^n < +\infty$. See Section 1.2.4.

Theorem 6.7.2 follows readily from Theorem 6.3.3 with $L = 1$ and $\omega = 0$, since $k = n - 1$.

This limit case of the Otter-Meir-Moon Theorem appears as Remark 3.7 in [29, Chapter 3]; the proof suggested there, of a different nature than the one above, comes from the Appendix of [53]. See also [67, Section 2.3.1]. In there you may find the same asymptotic result under the assumption that $\psi$ is a probability generating function satisfying $\psi^{(4)}(R) = \lim_{t\uparrow R} \psi^{(4)}(t) < \infty$.

Observe that, in particular, it follows from Theorem 6.7.1 that when $M_\psi > 1$ and with $m_\psi(\tau) = 1$ the radius of convergence of the solution $g$ is $\tau/\psi(\tau)$. When $M_\psi = 1$ (and $R < +\infty$), the radius of convergence of $g$ is $R/\psi(R)$.

If $M_\psi = 1$ and $R = +\infty$, then $\psi$ is a polynomial of degree 1. Thus $\psi(z) = a + bz$, with both $a, b > 0$. In this case the solution $g(z)$ of Lagrange equation with data $\psi$ is simply

$$g(z) = \frac{az}{1 - bz} = \sum_{n=1}^\infty (ab^{n-1}) z^n\,.$$

and thus

$$A_n = ab^{n-1}\,, \quad \text{para } n \geq 1\,.$$

♦ $M_\psi < 1$. In this case, we necessarily have that $R < \infty$, since $R = \infty$ and $M_\psi < 1$ would imply that $\psi$ is a constant, which is not allowed. See Lemma 1.2.1.

We assume as above that $\sum_{n=0}^\infty n^2 b_n R^n < \infty$. As in the preceding case $M_\psi = 1$ above, the power series $\psi$ is continuous in $\mathrm{cl}(\mathbb{D}(0,R))$ and the Khinchin family extends to include continuously a random variable $Y_R$ with mean $M_\psi < 1$ and variance $\lim_{t\uparrow R} \sigma^2_\psi(t) = \sigma^2_\psi(R) < \infty$. Assume for simplicity that $Q_\psi = 1$

In this case, we do not obtain a proper asymptotic formula, just only that

$$\lim_{n\to\infty} A_n \frac{R^{n-1}}{\psi^n(R)} n^{3/2} = 0\,. \tag{6.7.6}$$

For $A_n$ we have that

$$A_n = \frac{1}{n} \text{COEFF}_{[n-1]}(\psi^n) = \frac{1}{2\pi} \frac{\psi^n(R)}{R^{n-1}} \frac{1}{n^{3/2}} \frac{1}{\sigma_\psi(R)} I_n$$

where

$$I_n \triangleq \int_{|\theta|\le \pi\sigma_\psi(R)\sqrt{n}} \mathbf{E}\big(e^{\imath\theta \check{Y}_R/\sqrt{n}}\big)^n e^{\imath\theta\alpha(n)} d\theta\,, \qquad \text{for } n\ge 1\,,$$

and

$$\alpha(n) = \frac{nM_\psi - (n-1)}{\sigma_\psi(R)\sqrt{n}}\,, \qquad \text{for } n \ge 1\,.$$

We shall now verify that $I_n$ tends towards 0 as $n\to\infty$. We split $I_n$ as

$$\begin{aligned} I_n = & \int_{|\theta|\le \pi\sigma_\psi(R)\sqrt{n}} \Big(\mathbf{E}\big(e^{\imath\theta \check{Y}_R/\sqrt{n}}\big)^n - e^{-\theta^2/2}\Big) e^{\imath\theta\alpha(n)} d\theta \\ & + \int_{\mathbb{R}} e^{-\theta^2/2} e^{\imath\theta\alpha(n)} d\theta \\ & - \int_{|\theta|> \pi\sigma_\psi(R)\sqrt{n}} e^{-\theta^2/2} e^{\imath\theta\alpha(n)} d\theta\,. \end{aligned}$$

The third summand obviously tends to 0, while the first summand tends to 0 on account of Theorem 6.1.5, the integral form of the Local Central Limit Theorem. The second summand may be written as

$$\int_{\mathbb{R}} e^{-\theta^2/2} e^{\imath\theta\alpha(n)} d\theta = e^{-\alpha(n)^2/2}$$

to observe that since $M_\psi < 1$, we have that $\lim_{n\to\infty}\alpha(n) = -\infty$, and thus that this second summand too converges to 0, as $n\to\infty$.

Minami in [67, Section 2.3.2] obtains a faster rate of convergence (higher power of $n$) in (6.7.6) under the stronger assumption that $\lim_{t\uparrow R}\psi^{(k)}(t) < \infty$ for some $k\ge 3$.

### 6.7.3 Coefficients of powers of solutions of Lagrange's equation

We assume $Q_\psi = 1$. We are interested now, see Meir-Moon [65], in asymptotic formulas for the coefficients of powers of the solution of Lagrange's equation.

For $q\ge 1$ and $n\ge 1$, the $n$-th coefficient of $g(w)^q$, which we will denote by $B_{n,q}$, is given by the exact formula

$$B_{n,q} = \frac{q}{n}\mathrm{COEFF}_{[n-q]}(\psi(z)^n)\,.$$

- For $q\ge 1$ fixed, Theorem 6.3.2 gives, under the assumption $M_\psi > 1$ and with $\tau$ given by $m_\psi(\tau) = 1$, that

$$B_{n,q} \sim \frac{q}{\sqrt{2\pi}}\frac{\tau^q}{\sigma_\psi(\tau)}\frac{1}{n^{3/2}}\Big(\frac{\psi(\tau)}{\tau}\Big)^n\,, \qquad \text{as } n\to\infty\,.$$

- If $q$ and $n$ are related by $q = \alpha n + \beta\sqrt{n} + o(\sqrt{n})$ as $n\to\infty$, where $\alpha\in[0,1)$ and $\beta\in\mathbb{R}$, then for $\tau_\alpha$ given by $m_\psi(\tau_\alpha) = 1-\alpha$, we obtain, analogously to (6.6.1), that

$$B_{n,q} \sim e^{-\beta^2/(2\sigma_\psi^2(\tau_\alpha))}\frac{q}{\sqrt{2\pi}}\frac{\tau_\alpha^q}{\sigma_\psi(\tau_\alpha)}\frac{1}{n^{3/2}}\Big(\frac{\psi(\tau_\alpha)}{\tau_\alpha}\Big)^n\,, \qquad \text{as } n\to\infty\,.$$

Observe that

$$\sigma_\psi(\tau_\alpha) = \alpha(1-\alpha) + \tau_\alpha^2 \frac{\psi''(\tau_\alpha)}{\psi(\tau_\alpha)}\,.$$

**Coefficients of functions of solutions of Lagrange's equation**

In a more general setting, let $H$ be a power series with non-negative coefficients.

Assume further that the radius of convergence of $H$ is at least the radius of convergence $R$ of $\psi$ and also that $M_\psi > 1$. Let $\tau$ be such that $m_\psi(\tau) = 1$.

It follows from formulas (6.7.1) and (6.6.1) that the coefficients of $H(g(z))$, where $g(z)$ is the solution of Lagrange's equation with data $\psi$, satisfy:

$$\text{COEFF}_{[n]}(H(g(z))) \sim \frac{1}{\sqrt{2\pi}} \frac{H'(\tau)\tau}{\sigma_\psi(\tau)} \frac{1}{n^{3/2}} \Big(\frac{\psi(\tau)}{\tau}\Big)^n, \quad \text{as } n \to \infty\,.$$

## 6.8 Probability generating functions

Let $\psi$ be the probability generating function of a random variable $X$. We assume for convenience that $\psi(0) \neq 0$ and that $\psi'(0) \neq 0$, and further that $Q_\psi = 1$, that $M_\psi = +\infty$ and that the radius of convergence of $\psi$ is $R > 1$. This last assumption implies, in particular, that the variable $X$ has finite moments of all orders.

Let $F$ be the fulcrum of $\psi$, see Section 1.3.4, given by $F(z) = \ln\psi(e^z)$, in a region containing $[0, R_\psi)$. Recall that in terms of the fulcrum $F$ we have that $m_\psi(t) = F'(s)$ and $\sigma_\psi^2(t) = F''(s)$, where $s$ and $t$ are related by $e^s = t$.

The function $e^{F(z)} = \psi(e^z) = \sum_{n=0}^\infty b_n e^{nz}$, which is holomorphic in a neighborhood of $z = 0$, is the *moment generating function* of $X$.

If $X_1, X_2, \ldots$ are independent copies of $X$, then

$$\text{COEFF}_{[k]}(\psi(z)^n) = \mathbf{P}\left(\frac{X_1 + \cdots + X_n}{n} = \frac{k}{n}\right).$$

Let $m_\psi(\tau_n) = k/n$ and let $s_n$ be given by $e^{s_n} = \tau_n$. Observe that $\psi(\tau_n) = e^{F(s_n)}$ and that $\sigma_\psi^2(\tau_n) = F''(s_n)$.

In terms of the fulcrum $F$ (or the ln of the moment generating function), we have that if $k \asymp n$, see Theorem 6.3.1 or $k/n \to 0$, see Theorem 6.4.2, that

$$\begin{aligned}\mathbf{P}\left(\frac{X_1 + \cdots + X_n}{n} = \frac{k}{n}\right) &\sim \frac{1}{\sqrt{2\pi n}\,\sigma_\psi(\tau_n)} \frac{\psi(\tau_n)^n}{\tau_n^k} \\ &= \frac{1}{\sqrt{2\pi n\, F''(s_n)}} e^{n[F(s_n) - s_n(k/n)]}, \quad \text{as } n \to \infty\,.\end{aligned}$$

If $\psi$ is uniformly strongly Gaussian, the asymptotic formula above holds also if $k/n \to \infty$.

### 6.8.1 Lagrangian distributions (and Galton-Watson processes)

We now turn our attention to Lagrangian probability distributions. We refer to [89] for a neat presentation of the basic theory of these probability distributions. See also [22], for a comprehensive treatment, and [76].

We start with a probability generating function $\phi$ of a variable $Y$ taking values in $\{0, 1, \ldots\}$. We assume that $\phi$ is in $\mathcal{K}$.

This variable $Y$ generates a cascade or Galton-Watson random process starting (initial stage) with a single individual. The variable $Y$ gives the random number of immediate descendants, the offsprings of each individual in every generation.

The random number of individuals in generation $n$ is denoted by $G_n$; the generation 0 is the initial stage. We have $G_0 \equiv 1$ and

$$(\star) \quad G_{n+1} = \sum_{j=1}^{G_n} Y_{n,j}\,, \quad \text{for } n \geq 0\,,$$

where the $Y_{n,j}$ are independent copies of $Y$.

Let us denote $\phi'(1) \triangleq \lim_{t\uparrow 1} \phi'(t) = \mathbf{E}(Y)$. Observe that $\phi'(1) = m_\phi(1)$.

The total progeny $Z$ of the single individual of the generation 0 (including itself) is the random variable $Z = \sum_{n=0}^{\infty} G_n$. This variable $Z$ could take the value $\infty$, but it is a proper random variable (i.e., $\mathbf{P}(Z < \infty) = 1$) if and only if $\psi'(1) \leq 1$.

Assume thus that $\psi'(1) \leq 1$. The probability generating function $g$ of the total progeny $Z$ is actually the solution of Lagrange's equation with data $\phi$.

Furthermore, let $f$ be a power series with non-negative coefficients which is the probability generating function of a random variable $X$ taking values in $\{0, 1, \ldots\}$.

We enhance the process by allowing the size of the initial stage (generation 0) to be randomly chosen following the distribution of $X$. Thus $G_0 = X$ and the evolution is determined by the recurrence $(\star)$. The total progeny $Z$ including the individuals of the initial stage has probability generating function $f \circ g$.

This size of progeny $Z$ is said to *follow a Lagrangian distribution* $\mathcal{L}(\phi, f)$; $\phi, f$ are called the generators of $\mathcal{L}(\phi, f)$.

If $f(z) \equiv z$, then $X \equiv 1$, and there is (deterministically) a single individual in the initial generation.

We have, see (6.7.1), that

$$\mathbf{P}(Z = n) = \text{COEFF}_{[n]}(f(g(z))) = \frac{1}{n}\text{COEFF}_{[n-1]}(f'(z)\phi^n(z))\,, \quad \text{for } n \geq 1\,,$$

and $\mathbf{P}(Z = 0) = f(0)$.

### 6.8.2 Lagrangian distributions and Khinchin families

We will have two ingredients: a Khinchin family of offspring probability distributions and a Khinchin family of probability distribution for the initial distribution.

(A) Let $\psi(z)=\sum_{n=0}^{\infty} b_n z^n$ be a power series in $\mathcal{K}$ with radius of convergence $R$. We assume from the outset that $Q_\psi=1$. We let $(Y_t)_{t\in[0,R)}$

Assume that either $M_\psi>1$ and then we let $\tau\in(0,R)$ be such that $m_\psi(\tau)=1$ or $M_\psi=1$ and $R<\infty$ and $\lim_{t\uparrow R}\sigma_\psi(t)=\sigma_\psi(R)<\infty$ and then we let $\tau=R$.

For each $t\in(0,\tau]$, we let $\psi_t$ denote the power series $\psi_t(z)=\psi(tz)/\psi(t)$. This $\psi_t$ is the probability generating function of $Y_t$, and besides, since $t\le\tau$, we have that

$$(\dagger)\quad \psi_t'(1)=m_\psi(t)\le 1\,.$$

For each $t\in(0,\tau]$, we let $g_t$ be the solution of Lagrange's equation with data $\psi_t$:

$$g_t(z)=z\psi_t(g_t(z))\,. \tag{6.8.1}$$

Because of $(\dagger)$ and the discussion above in Section 6.8.1, $g_t$ is the probability generating function of the random distribution of the total progeny of a single individual with offspring distribution $Y_t$.

We let $g$ denote the solution of Lagrange's equation with data $\psi$, then we may write each $g_t$ in terms of $g$ as follows

**Lemma 6.8.1.** *With the notations above*

$$g_t(z)=\frac{1}{t}g\big((t/\psi(t))z\big)\,,\quad \textit{for } t\le\tau\,.$$

*Proof.* This is a consequence of the uniqueness of solution of Lagrange's equation. Let $\widetilde{g}_t(z)=\frac{1}{t}g\big((t/\psi(t))z\big)$. Now,

$$\begin{aligned}
z\psi_t(\widetilde{g}_t(z)) &= \frac{z}{\psi(t)}\psi(t\,\widetilde{g}_t(z))=\frac{z}{\psi(t)}\psi(g\big((t/\psi(t))z\big)\\
&=\frac{1}{t}\frac{zt}{\psi(t)}\psi(g\big((t/\psi(t))z\big)=\frac{1}{t}g\big((t/\psi(t))z\big)\\
&=\widetilde{g}_t(z)\,.
\end{aligned}$$

Uniqueness gives that $\widetilde{g}_t(z)=g_t(z)$, as claimed. □

**Remark 6.8.2.** The $g_t$ are *not* the probability generating functions of a Khinchin family, but if we change parameters and substitute $t$ by $g(u)$ and let

$$\widetilde{g}_u(z)=g_{g(u)}(z)$$

we have that $\widetilde{g}_u(z)=g(uz)/g(u)$. The power series $g$ is not in $\mathcal{K}$ but $g(z)/z$ is in $\mathcal{K}$, since $g'(0)=\psi(0)>0$. If $(W_u)$ is the Khinchin family of $g(z)/z$ then $\widetilde{g}_u$ is the probability generating functions of $W_u+1$. ⊠

(B) Next, we let $f(z)=\sum_{n=0}^{\infty}a_n z^n$ be a non-constant power series with non-negative coefficients and radius of convergence $S>0$. We do not require $f$ to be in $\mathcal{K}$; in fact $f(z)=z^m$, for integer $m\geq 1$, is particularly relevant.

For $s\in(0,S)$ we denote by $f_s$ the power series

$$f_s(z)=\frac{f(sz)}{f(s)}$$

which is the probability generating function of a random variable $X_s$, say.

**Remark 6.8.3.** If the power series $f$ were not in $\mathcal{K}$, then since $f$ is not constant, there should exist at least one index $m\geq 1$ so that $a_m\neq 0$. Thus only two cases could occur. In the first case, $f(z)=a_m z^m$, for some integer $m\geq 1$. And in the second; there is $\phi\in\mathcal{K}$ and an integer $l\geq 0$ so that $f(z)=z^l\phi(z)$.

In the first case $(X_s)$ is not a Khinchin family, but $X_s\equiv m$, for all $s\in(0,S)$, and therefore there are exactly $m\geq 1$ nodes in the first generation; this is the deterministic case $f_s(z)=z^m$.

In the second case $(X_s)_{s\in[0,S)}$ is a shifted Khinchin family. If $f(0)>0$, then $(X_s)_{s\in[0,S)}$ is a proper Khinchin family. ⊠

Now, for $s\in(0,S)$ and $t\in(0,\tau]$, the composition $f_s(g_t(z))$ is the probability generating function of the total progeny $Z_{s,t}$ of a Galton-Watson process with initial distribution $f_s$ and offspring distribution $g_t$. The progeny $Z_{s,t}$ has Lagrange distribution $\mathcal{L}(\psi_t,f_s)$.

For radius $u>0$ such that $us<S$ and $ut\leq\tau$ and appealing to (6.7.1) we may write

$$\begin{aligned}\mathbf{P}(Z_{s,t}=n)&=\text{COEFF}_{[n]}\big(f_s(g_t(z))\big)=\frac{1}{n}\text{COEFF}_{[n-1]}\big(f_s'(z)\psi_t^n(z)\big)\\&=\frac{s}{f(s)}\frac{\psi(tu)^n}{\psi(t)^n}\frac{1}{u^{n-1}}\frac{1}{n}\frac{1}{2\pi}\int\limits_{|\theta|\leq\pi}f'(sue^{\imath\theta})\frac{\psi(tue^{\imath\theta})^n}{\psi(tu)^n}e^{-\imath(n-1)\theta}d\theta\,.\end{aligned}$$

If the parameter $s$ is further restricted to $s\tau<tS$, we may take $u=\tau/t$ in the expression above and write

$$\begin{aligned}\mathbf{P}(Z_{s,t}=n)&=\frac{1}{2\pi}\frac{s}{f(s)}\frac{\psi(\tau)^n}{\psi(t)^n}\left(\frac{t}{\tau}\right)^{n-1}\frac{1}{n}\int\limits_{|\theta|\leq\pi}f'\left(\frac{s\tau}{t}e^{\imath\theta}\right)\mathbf{E}(e^{\imath Y_\tau\theta})^n e^{-\imath(n-1)\theta}d\theta\\&=\frac{1}{2\pi}\frac{s}{f(s)}\frac{\psi(\tau)^n}{\psi(t)^n}\left(\frac{t}{\tau}\right)^{n-1}\frac{1}{n^{3/2}}\frac{1}{\sigma_\psi(\tau)}\\&\int\limits_{|\theta|\leq\pi\sigma_\psi(\tau)\sqrt{n}}f'\left(\frac{s\tau}{t}e^{\imath\theta/(\sigma_\psi(\tau)\sqrt{n})}\right)\mathbf{E}(e^{\imath\breve{Y}_\tau\theta/\sqrt{n}})^n e^{\imath\theta/(\sigma_\psi(\tau)\sqrt{n})}d\theta\end{aligned}$$

For $\theta$ fixed, we have that

$$\lim_{n\to\infty}f'\left(\frac{s\tau}{t}e^{\imath\theta/(\sigma_\psi(\tau)\sqrt{n})}\right)e^{\imath\theta/(\sigma_\psi(\tau)\sqrt{n})}=f'\left(\frac{s\tau}{t}\right)$$

and

$$\left|f'\left(\frac{s\tau}{t}e^{\imath\theta/(\sigma_\psi(\tau)\sqrt{n})}\right)e^{\imath\theta/(\sigma_\psi(\tau)\sqrt{n})}\right|\leq f'\left(\frac{s\tau}{t}\right)\,.$$

And thus taking into account the integral form of the Local Central Limit Theorem, Theorem 6.1.5, we deduce that

$$\lim_{n\to\infty} \int_{|\theta|\le\pi\sigma_\psi(\tau)\sqrt{n}} f'\left(\frac{s\tau}{t}e^{\imath\theta/(\sigma_\psi(\tau)\sqrt{n})}\right)\mathbf{E}(e^{\imath\breve{Y}_\tau\theta/\sqrt{n}})^n e^{\imath\theta/(\sigma_\psi(\tau)\sqrt{n})}d\theta = \sqrt{2\pi}f'\left(\frac{s\tau}{t}\right)$$

and, consequently, that

$$\mathbf{P}(Z_{s,t}=n)\sim\frac{1}{\sqrt{2\pi}}\frac{s}{f(s)}\frac{\psi(\tau)^n}{\psi(t)^n}\left(\frac{t}{\tau}\right)^{n-1}\frac{1}{n^{3/2}}\frac{1}{\sigma_\psi(\tau)}f'\left(\frac{s\tau}{t}\right), \quad \text{as } n\to\infty, \tag{6.8.2}$$

as long as $s\tau < tS$, which amounts to no restriction if $S=+\infty$.

(C) As an illustration, consider the case where $\psi(z)=e^z$ and $f(z)=z^j$, for some integer $j\ge 1$.

In this case, $R=S=\infty$, $m_\psi(t)=t$ and $\sigma^2_\psi(t)=t$. Also $M_\psi=\infty$ and $\tau=1$.

For $t\le 1=\tau$, we have that $\psi_t(z)=e^{t(z-1)}$ and for $s<\infty$, we have that $f_s(z)=z^j$. Observe that for any $s$, $f_s$ is the probability generating function of the constant $j$.

For $0<t\le 1$ and $0<s<\infty$, the variable $Z_{s,t}$ is the total progeny of a Galton-Watson process, where the initial generation consists of exactly $j$ individuals and the offspring of each individual is given by a Poisson variable of parameter $t$. This distribution, $\mathcal{L}(e^{t(z-1)}, z^j)$, is the Borel-Tanner distribution with parameters $t$ and $j$. The case $j=1$ is the Borel distribution. See [22], [79] and [89], and also the original sources [15], [43] and [92].

Using (6.8.2), we deduce that

$$\mathbf{P}(Z_{s,t}=n)\sim\frac{j}{\sqrt{2\pi}}\frac{1}{n^{3/2}}t^{n-j}e^{n(1-t)}, \quad \text{as } n\to\infty. \tag{6.8.3}$$

In fact, for the Borel-Tanner distribution with parameters $t$ and $j$ we have the exact formula

$$\mathbf{P}(Z_{s,t}=n)=\frac{j}{n}\frac{e^{-tn}(tn)^{n-j}}{(n-j)!}, \quad \text{for } n\ge j. \tag{6.8.4}$$

The asymptotic formula (6.8.3) follows then from Stirling's formula.

(D) Consider now the case $\psi(z)=e^z$ and $f(z)=e^z$, so that $R=S=\infty$, $m_\psi(t)=t$ and $\sigma^2_\psi(t)=t$ and, also, $M_\psi=\infty$ and $\tau=1$.

For $t\le 1=\tau$, we have that $\psi_t(z)=e^{t(z-1)}$ and for $s<\infty$, we have that $f_s(z)=e^{s(z-1)}$.

For $0<t\le 1$ and $0<s<\infty$, the variable $Z_{s,t}$ is the total progeny of a Galton-Watson process, where the size of the initial generation is drawn from a Poisson distribution of parameter $s$ individuals and the offspring of each individual is given by a Poisson variable of parameter $t$. This distribution, in Lagrangian distribution notation, is $\mathcal{L}(e^{t(z-1)}, e^{s(z-1)})$.

From (6.8.2), we have that

$$\mathbf{P}(Z_{s,t}=n)\sim\frac{1}{\sqrt{2\pi}}e^{s/t-s}st^{n-1}e^{n(1-t)}\frac{1}{n^{3/2}}, \quad \text{as } n\to\infty. \tag{6.8.5}$$

Conditioning on the size of the initial generation and using (6.8.4), we deduce for the distribution $\mathcal{L}(e^{t(z-1)}, e^{s(z-1)})$ that

$$\mathbf{P}(Z_{s,t} = n) = \frac{1}{n!} e^{-tn-s}(tn+s)^{n-1} s\,, \quad \text{para } n \geq 1\,.$$

See [89]. The asymptotic formula (6.8.5) follows then from Stirling's formula.

(E) *Back to Galton-Watson.* Assume that $\psi$ is a probability generating function and that $f(z) = z$. We assume that $\psi'(1) \leq 1$, so that the total progeny $Z$ of the Galton-Watson process starting with a single individual and with offspring distribution $\psi$ is a (proper) random variable.

Assume that $M_\psi > 1$ or $M_\psi = 1$ with $R < \infty$ and $\sigma_\psi(R) < \infty$. In the first case we take $\tau \in (0, R)$ such that $m_\psi(\tau) = 1$. Since $\psi'(1) \leq 1$, we have that $m_\psi(1) = 1$, and since $m_\psi$ is increasing we see that $1 \leq \tau < R$. In the second case we take $\tau = R$; observe that $R \geq 1$.

With $t = 1$ (so that $\psi_1(z) \equiv \psi(z)$) and $s = 1$ (an immaterial choice since $f(z) = z$), we have

$$\mathbf{P}(Z = n) \sim \frac{1}{\sqrt{2\pi}} \frac{\tau}{\sigma_\psi(\tau)} \Big(\frac{\psi(\tau)}{\tau}\Big)^n \frac{1}{n^{3/2}}\,, \quad \text{as } n \to \infty\,.$$

This is Theorem 6.7.1 applied to the solution of Lagrange's equation with data $\psi$. Recall that $\sigma_\psi^2(\tau) = \tau^2 \psi''(\tau)/\psi(\tau)$, see formula (6.7.4).

### Limit cases

By limit cases we mean $t = \tau$ and $s \to S$, with $S < \infty$. See [70] for related results.

For $t = \tau$ and $s < S$ we have the exact formula

$$(\flat) \quad \begin{aligned} \mathbf{P}(Z_{s,\tau} = n) &= \frac{1}{2\pi} \frac{s}{f(s)} \frac{1}{n^{3/2}} \frac{1}{\sigma_\psi(\tau)} \\ &\int\limits_{|\theta| \leq \pi\sigma_\psi(\tau)\sqrt{n}} f'\big(s e^{\imath\theta/(\sigma_\psi(\tau)\sqrt{n})}\big) \mathbf{E}(e^{\imath \check{Y}_\tau \theta/\sqrt{n}})^n e^{\imath\theta/(\sigma_\psi(\tau)\sqrt{n})} d\theta \end{aligned}$$

and the asymptotic formula.

$$\mathbf{P}(Z_{s,\tau} = n) \sim \frac{1}{\sqrt{2\pi}} m_f(s) \frac{1}{\sigma_\psi(\tau)} \frac{1}{n^{3/2}}\,, \quad \text{as } n \to \infty\,,$$

where with a slight abuse of notation we have written $sf'(s)/f(s) = m_f(s)$, see Remark 6.8.3.

To consider $s \to S$, we first rewrite ($\flat$) as

$$(\flat\flat) \quad \begin{aligned} \mathbf{P}(Z_{s,\tau} = n) &= \frac{1}{2\pi} \frac{1}{n^{3/2}} m_f(s) \frac{1}{\sigma_\psi(\tau)} \\ &\int\limits_{|\theta| \leq \pi\sigma_\psi(\tau)\sqrt{n}} \frac{f'\big(s e^{\imath\theta/(\sigma_\psi(\tau)\sqrt{n})}\big)}{f'(s)} \mathbf{E}(e^{\imath \check{Y}_\tau \theta/\sqrt{n}})^n e^{\imath\theta/(\sigma_\psi(\tau)\sqrt{n})} d\theta \end{aligned}$$

If $\sum_{n=0}^{\infty} na_n S^n < \infty$, then $f'$ extends continuously to $\mathrm{cl}(\mathbb{D}(0,S))$, and by appealing to Theorem 6.1.5, we readily see that

$$\lim_{\substack{s\uparrow S;\\ n\to\infty}} \frac{n^{3/2}}{m_f(s)} \mathbf{P}(Z_{\tau,s} = n) = \frac{1}{\sigma_\psi(\tau)\sqrt{2\pi}} \,.$$

More generally, if $s_n < S$ and $n$ are such that

$$(\sharp) \qquad \lim_{n\to\infty} \frac{f'(s_n e^{\imath\phi/\sqrt{n}})}{f'(s_n)} = 1 \,, \quad \text{for each } \phi \in \mathbb{R} \,,$$

then

$$\lim_{n\to\infty} \frac{n^{3/2}}{m_f(s_n)} \mathbf{P}(Z_{\tau,s_n} = n) = \frac{1}{\sigma_\psi(\tau)\sqrt{2\pi}} \,.$$

For instance for $f(z) = 1/(1-z)$, condition $(\sharp)$ is satisfied if $s_n$ and $n$ are related so that $\lim_{n\to\infty}(1-s_n)\sqrt{n} = \infty$.